\pdfoutput=1
\documentclass[11pt]{article}
\usepackage[left=1in,right=1in,top=1in,bottom=1in]{geometry}
\usepackage{times}
\usepackage{expl3}
\usepackage{cite}
\usepackage[table]{xcolor}
\usepackage{multirow}
\usepackage{stackengine} 
\usepackage{hhline}
\usepackage{lipsum}
\usepackage{titlesec}
\usepackage[all]{xy}
\usepackage{wrapfig}
\usepackage{enumerate}
\usepackage{epsfig}
\usepackage{tikz-cd}
\usepackage{amsmath}
\usepackage{tabularx}
\usepackage{array}
\usepackage{booktabs}
\usepackage{enumitem}
\usepackage{bbm}
\usepackage{calc}
\usepackage{graphicx}
\usepackage{amsmath}
\usepackage[title]{appendix}
\usepackage{amssymb}
\usepackage{epstopdf}
\usepackage{boldline}
\usepackage{arydshln}
\usepackage{calligra}
\usepackage{bm}
\usepackage{url}
\usepackage{blindtext}
\usepackage{accents}
\usepackage{amsthm}
\usepackage{amscd}

\newtheorem{definition}{Definition}

\newtheorem{theorem}{Theorem}
\newtheorem{proposition}{Proposition}
\newtheorem{corollary}{Corollary}
\newtheorem{lemma}{Lemma}
\newtheorem{conjecture}{Conjecture}
\newtheorem{assumption}{Assumption}
\newtheorem{example}{Example}
\newtheorem{remark}{Remark}
\usepackage{mathtools}
\usepackage{epstopdf}
\usepackage{balance}
\usepackage{thmtools}
\usepackage{thm-restate}
\usepackage{hyperref}
\usepackage{cleveref}
\usepackage[mathscr]{euscript}

\usepackage[ruled,vlined]{algorithm2e}
\newcommand{\ostar}{\mathbin{\mathpalette\make@circled\star}}

\makeatletter
\newcommand{\removelatexerror}{\let\@latex@error\@gobble}
\makeatother
\makeatletter
\newcommand*{\rom}[1]{\expandafter\@slowromancap\romannumeral #1@}
\makeatother

\ExplSyntaxOn
\newcommand\latinabbrev[1]{
  \peek_meaning:NTF . {
    #1\@}%
  { \peek_catcode:NTF a {
      #1.\@ }%
    {#1.\@}}}
\ExplSyntaxOff

\titleclass{\subsubsubsection}{straight}[\subsubsection]

\begin{document}
\vspace{1cm}
\title{The Unitary Conjugation Groupoid as a Universal Mediator of the Baum--Connes Assembly Map}
\vspace{1.8cm}
\author{Shih-Yu~Chang
\thanks{Shih-Yu Chang is with the Department of Applied Data Science,
San Jose State University, San Jose, CA, U. S. A. (e-mail: {\tt
shihyu.chang@sjsu.edu})
}}

\maketitle

\begin{abstract}
We show that the Baum--Connes assembly map factors canonically through the unitary conjugation groupoid, which serves as a universal mediator among groupoid models that are Morita equivalent to a given transformation groupoid. This establishes a structural link between groupoid-based index theory and the Baum--Connes program at the level of K-theory. Building on our previous development of unitary conjugation groupoids and their associated index theory, we extend the $K_1$ index framework beyond the Type I setting to non-Type I examples, including the irrational rotation algebra and amenable crossed products. Using Morita equivalence, we relate unitary conjugation groupoids to transformation and action groupoids, enabling the transfer of descent-type index constructions to these settings. Our main result shows that, among all groupoid realizations that are Morita equivalent to a transformation groupoid, the factorization through the unitary conjugation groupoid is canonical at the level of K-theory. This identifies the unitary conjugation groupoid as a universal intermediary for the Baum--Connes assembly map. As applications, we recover the classical index pairing with the tracial state for the irrational rotation algebra in the sense of Connes, and we prove that for amenable crossed products the descent construction agrees with the analytic Baum--Connes assembly map under Morita equivalence.  These results provide a conceptual interpretation of the assembly map in terms of internal symmetries of crossed product algebras and suggest a unified framework connecting Fredholm-type index data with equivariant K-theory via groupoid methods.
\end{abstract}

\tableofcontents

\newpage

\section{Introduction}

\subsection*{Standing Assumptions for the Main Theorems}
Unless explicitly stated otherwise, all groupoids are assumed to be second
countable, locally compact, and equipped with Haar systems whenever reduced
groupoid $C^*$-algebras or groupoid-equivariant $KK$-theory are invoked.  All
$C^*$-algebras are assumed separable.  In the crossed-product applications, the
group $\Gamma$ is taken to be discrete and amenable, and actions on compact
Hausdorff spaces are assumed topologically free when this hypothesis is used to
identify the relevant Morita models.

Index theory for operator algebras provides a powerful bridge between
analysis, topology, and geometry.
Since the pioneering work of Connes~\cite{ConnesNCG} and
Connes--Skandalis~\cite{Connes-Skandalis1984}, groupoids have emerged as a
natural framework for formulating index theorems in noncommutative
geometry. In particular, groupoid $C^*$-algebras provide geometric
models for a wide class of operator algebras arising from dynamical
systems, foliations, and group actions.

In Paper~I~\cite{PaperI}, we constructed the unitary conjugation
groupoid $\mathcal{G}_{\mathcal{A}}$ associated to a separable
unital $C^*$-algebra $\mathcal{A}$ and established its fundamental
topological and structural properties.
In Paper~II~\cite{PaperII}, we developed a $K_1$-index framework for
Fredholm-type operators using this groupoid.
In particular, for the prototypical Type~I algebras $B(H)$ and
$\mathcal{K}(H)^\sim$, the equivariant $K_1$-class associated to an
invertible operator descends through the groupoid $C^*$-algebra to
recover the classical Fredholm index.

The goal of the present paper is to extend this framework to
genuinely non-Type~I $C^*$-algebras.
Our main examples are the irrational rotation algebra
\[
A_\theta = C(S^1)\rtimes_\theta \mathbb{Z},
\]
introduced by Rieffel~\cite{Rieffel1981}, and more generally amenable
crossed products of the form
\[
C(X)\rtimes \Gamma,
\]
where $\Gamma$ is a discrete amenable group acting on a compact
Hausdorff space $X$.

\subsection{The Need for Extension to Non-Type I Algebras}

While the unitary conjugation groupoid $\mathcal{G}_{\mathcal{A}}$
behaves well for Type~I $C^*$-algebras, many central examples in
noncommutative geometry lie outside this class.
Notably, the irrational rotation algebra $A_\theta$ is simple,
nuclear, and non-Type~I, yet it plays a fundamental role in
noncommutative geometry~\cite{ConnesNCG}.
Similarly, crossed product algebras $C(X)\rtimes\Gamma$
arising from group actions provide a broad class of
non-Type~I operator algebras.

Although the direct construction of the unitary conjugation
groupoid from Paper~I is no longer available in these settings,
such algebras admit natural geometric models via groupoids
associated with Cartan subalgebras in the sense of
Renault~\cite{Renault1980}.
This suggests that the descent-index machinery developed in
Paper~II may still be applicable if the relevant groupoids
can be related through Morita equivalence.
Establishing this relationship forms the central objective of
the present work.

\subsection{Why the Irrational Rotation Algebra and Amenable Crossed Products?}

The irrational rotation algebra $A_\theta$ provides the
simplest and most fundamental example of a non-Type~I
$C^*$-algebra.
It admits a natural realization as the groupoid $C^*$-algebra
of the transformation groupoid $S^1\rtimes_\theta\mathbb{Z}$,
and its $K$-theory is closely tied to Connes'
noncommutative index theorem~\cite{ConnesNCG}.

More generally, crossed products $C(X)\rtimes\Gamma$
associated with discrete group actions admit geometric
descriptions in terms of transformation groupoids
$X\rtimes\Gamma$.
For these groupoids, the Baum--Connes assembly map
\[
\mu_\Gamma : K_*^\Gamma(X)\longrightarrow K_*(C(X)\rtimes\Gamma)
\]
plays a central role in connecting topological invariants
of the dynamical system with analytic invariants of the
associated operator algebra
\cite{Baum-Connes1988,Kasparov1988,Tu1999}.
These examples therefore provide a natural testing ground
for extending the unitary conjugation groupoid framework
to non-Type~I settings.

\subsection{Main Results and Contributions of This Paper}

The principal contribution of this paper is to establish a
direct connection between the unitary conjugation groupoid
framework developed in Papers~I and~II and the
Baum--Connes assembly map for amenable crossed products.

More precisely, we prove the following results.

\begin{itemize}

\item For the irrational rotation algebra $A_\theta$, the
unitary conjugation groupoid is Morita equivalent to the
transformation groupoid $S^1\rtimes_\theta\mathbb{Z}$.

\item Using this equivalence, the descent construction of
Paper~II extends to $A_\theta$ and reproduces Connes'
index pairing between $K_1(A_\theta)$ and the canonical
trace.

\item For general crossed products $C(X)\rtimes\Gamma$ with
$\Gamma$ discrete and amenable, the groupoid associated
to $C(X)\rtimes\Gamma$ is Morita equivalent to the
action groupoid $X\rtimes\Gamma$.

\item Under this equivalence, the descent map defined by
the unitary conjugation groupoid coincides with the
analytic Baum--Connes assembly map.

\end{itemize}

These results establish a conceptual bridge between the
groupoid-based index theory introduced in
Papers~I and~II and the Baum--Connes program.

\subsection{A Universal Mediation Theorem}

The crowning achievement of this paper is Theorem~4.2, which
establishes a universal property for the unitary conjugation
groupoid. We prove that among all groupoid models Morita
equivalent to the transformation groupoid $X\rtimes\Gamma$,
the factorization of the Baum--Connes assembly map through
$\mathcal{G}_{C(X)\rtimes\Gamma}$ is canonical and universal.

More precisely, let $\mathsf{Fact}(X,\Gamma)$ denote the class
of triples $(\mathcal{H},\rho_{\mathcal{H}},\Phi_{\mathcal{H}})$
where $\mathcal{H}$ is a groupoid Morita equivalent to
$X\rtimes\Gamma$, $\rho_{\mathcal{H}}$ is a realization map
from equivariant $KK$-theory to $K_*(C^*_r(\mathcal{H}))$, and
$\Phi_{\mathcal{H}}$ is an isomorphism from $K_*(C^*_r(\mathcal{H}))$
to $K_*(C(X)\rtimes\Gamma)$ such that
$\mu_\Gamma = \Phi_{\mathcal{H}} \circ \rho_{\mathcal{H}}$.
Then for every such object, there exists a canonical isomorphism
$\Theta_{\mathcal{H}}: K_*(C^*_r(\mathcal{G}_{C(X)\rtimes\Gamma}))
\to K_*(C^*_r(\mathcal{H}))$ making the diagram

\[
\begin{tikzcd}[column sep=large,row sep=large]
& K_*(C_r^*(\mathcal G_A)) \arrow[dr,"\Phi_{\mathcal G_A}"]
   \arrow[dd,"\Theta_{\mathcal H}"'] & \\
KK_\Gamma^*(C_0(X),\mathbb C) \arrow[ur,"\rho_{\mathcal G_A}"]
   \arrow[dr,"\rho_{\mathcal H}"'] &&
K_*(A) \\
& K_*(C_r^*(\mathcal H)) \arrow[ur,"\Phi_{\mathcal H}"'] &
\end{tikzcd}
\]

commute. In this precise sense, the unitary conjugation groupoid
is the \emph{universal mediator} through which the Baum--Connes
assembly map factors among all groupoid models Morita equivalent
to $X\rtimes\Gamma$.

This result represents a qualitative leap: it elevates our
construction from an interesting groupoid model to a structural
theorem about the nature of the assembly map itself. The unitary
conjugation groupoid is not merely one among many equivalent
models—it is the canonical source from which all other
factorizations arise via unique $K$-theory isomorphisms.

\begin{remark}[Main takeaway]
The descent map associated to the unitary conjugation groupoid
recovers the Baum--Connes assembly map under Morita equivalence,
and this factorization is universal among all Morita-equivalent
groupoid models. This provides a new geometric realization of
the assembly map in terms of the internal symmetries of the
crossed product algebra, and establishes the unitary conjugation
groupoid as the canonical geometric intermediary between
equivariant topology and analytic $K$-theory.
\end{remark}

\subsection{Outline of the Paper}

The paper is organized as follows.

Section~\ref{sec:preliminaries} reviews the necessary background
on crossed products, groupoid $C^*$-algebras, Cartan
subalgebras, and the Baum--Connes assembly map.

Section~\ref{sec:morita} develops the Morita equivalence
framework that allows the unitary conjugation groupoid
to be interpreted in terms of classical transformation
groupoids.

Section~\ref{sec:Atheta} applies this framework to the
irrational rotation algebra $A_\theta$, where the
descent construction recovers Connes' index pairing.

Section~\ref{sec:crossed} extends the analysis to
amenable crossed products $C(X)\rtimes\Gamma$ and
establishes the compatibility between the descent map
and the Baum--Connes assembly map, culminating in the
universal mediation theorem (Theorem~\ref{thm:universal-mediator-assembly}).

Finally, Section~\ref{sec:synthesis} provides a conceptual
synthesis of the results and we formulate a conjectural
groupoid analogue of the Baum--Connes conjecture in Section~\ref{sec:A Conjecture The Unitary Conjugation Groupid Analogue}. 

\subsection{Notation and Conventions}

Throughout the paper, $A$ denotes a separable unital
$C^*$-algebra and $\Gamma$ denotes a discrete group
acting on a compact Hausdorff space $X$.
We write $C^*(G)$ and $C_r^*(G)$ for the full and reduced
$C^*$-algebras of a groupoid $G$.
For a $C^*$-algebra $A$, its $K$-theory groups are denoted
by $K_0(A)$ and $K_1(A)$.

Additional notation will be introduced when needed.

\begin{remark}
The author is solely responsible for the mathematical insights and theoretical directions proposed in this work. AI tools, including OpenAI's ChatGPT and DeepSeek models, were employed solely to assist in verifying ideas, organizing references, and ensuring internal consistency of exposition~\cite{chatgpt2025,deepseek2025}. \\

Also, throughout this paper, ``Paper I'' refers to our previous work~\cite{PaperI}, and ``Paper II'' refers to our previous work~\cite{PaperII}.
\end{remark}

\section{Preliminaries: Crossed Products, Groupoid C*-Algebras, and the Baum--Connes Conjecture}
\label{sec:preliminaries}

\subsection{Crossed Products and Transformation Groupoids}
\label{subsec:crossed-products-discrete-groups}

Crossed product $C^*$-algebras provide a natural framework for encoding
dynamical systems within operator-algebraic structures.  They arise from
actions $\alpha: \Gamma \to \operatorname{Aut}(A)$ of a group $\Gamma$ on a
$C^*$-algebra $A$ and play a central role in noncommutative geometry,
particularly in the formulation of the Baum--Connes conjecture with
coefficients and in the study of transformation groupoid $C^*$-algebras,
where one has $C^*(X \rtimes \Gamma) \cong C_0(X) \rtimes_r \Gamma$.

In the present work, crossed products serve as a primary and abundant
source of non-Type~I $C^*$-algebras to which we extend the unitary
conjugation groupoid construction introduced in Paper~I.  Our first and
most fundamental example is the irrational rotation algebra
\[
A_\theta = C(S^1) \rtimes_\theta \mathbb{Z},
\]
the canonical simple $C^*$-algebra associated with an irrational rotation
$2\pi\theta$ of the circle.  As a genuine non-Type~I algebra, it provides
a crucial test case for our construction beyond the commutative and
Type~I settings.  Furthermore, the groupoid interpretation of such
crossed products will, in subsequent sections, allow us to relate the
unitary conjugation groupoid framework to classical transformation
groupoids via Morita equivalence, thereby linking our construction to
the broader noncommutative geometry landscape.

Let $A$ be a unital $C^*$-algebra and let $\Gamma$ be a discrete group.
Suppose that $\alpha : \Gamma \to \operatorname{Aut}(A)$ is an action of
$\Gamma$ on $A$ by $*$-automorphisms.  The triple $(A,\Gamma,\alpha)$ is called
a \emph{$C^*$-dynamical system}.

\begin{definition}[Algebraic crossed product]
\label{def:algebraic-crossed-product}
Let $(A, \Gamma, \alpha)$ be a $C^*$-dynamical system with $\Gamma$ discrete.
The \emph{algebraic crossed product} $A \rtimes_{\mathrm{alg}} \Gamma$
is the $*$-algebra consisting of finitely supported functions
$x:\Gamma \to A$, which may be written as finite sums
\[
x = \sum_{g\in\Gamma} a_g u_g,
\]
where $a_g \in A$ and only finitely many $a_g$ are nonzero.
The symbols $u_g$ are formal generators indexed by $g\in\Gamma$.

The multiplication and involution are defined on elementary tensors by
\[
(a u_g)(b u_h) = a\,\alpha_g(b)\,u_{gh},
\qquad
(a u_g)^* = \alpha_{g^{-1}}(a^*)\,u_{g^{-1}},
\]
and extended linearly to finite sums, making
$A \rtimes_{\mathrm{alg}} \Gamma$ a $*$-algebra.
\end{definition}

Equivalently, one may identify $A \rtimes_{\mathrm{alg}} \Gamma$ with the
$*$-algebra $C_c(\Gamma,A)$ of finitely supported functions
$f:\Gamma \to A$, equipped with the convolution product
\[
(f*g)(t) = \sum_{s\in\Gamma} f(s)\,\alpha_s\!\bigl(g(s^{-1}t)\bigr)
\]
and involution
\[
f^*(t) = \alpha_t\!\bigl(f(t^{-1})^*\bigr).
\]

\begin{remark}
The notation $u_g$ is purely formal; these elements become unitaries in the
$C^*$-completions discussed below.  The algebraic crossed product is dense
in both the full and reduced crossed product $C^*$-algebras.
\end{remark}

The algebraic crossed product admits two canonical $C^*$-completions:
the full crossed product and the reduced crossed product.

\begin{definition}[Covariant representation]
\label{def:covariant-representation}
A \emph{covariant representation} of $(A,\Gamma,\alpha)$ on a Hilbert space
$H$ consists of a $*$-representation $\pi:A\to B(H)$ together with a
unitary representation $\rho:\Gamma\to\mathcal{U}(H)$ satisfying
\[
\rho(g)\pi(a)\rho(g)^* = \pi(\alpha_g(a)),
\qquad g\in\Gamma,\ a\in A .
\]
\end{definition}

\begin{definition}[Full crossed product]
\label{def:full-crossed-product}
The \emph{full crossed product} $A\rtimes_\alpha\Gamma$
is the completion of $A\rtimes_{\mathrm{alg}}\Gamma$ in the norm
\[
\|x\|_{\max}
=
\sup\bigl\{
\|\pi \rtimes \rho(x)\| :
(\pi,\rho) \text{ is a covariant representation}
\bigr\}.
\]
\end{definition}

The reduced crossed product is defined using the regular covariant
representation.

\begin{definition}[Reduced crossed product]
\label{def:reduced-crossed-product}
Let $\lambda:\Gamma\to\mathcal{U}(\ell^2(\Gamma))$ be the left regular
representation.  If $\pi:A\to B(H)$ is any faithful representation,
define a representation $\tilde{\pi}:A\to B(H\otimes\ell^2(\Gamma))$ by
\[
\tilde{\pi}(a)(\xi\otimes\delta_g)
=
\pi\bigl(\alpha_{g^{-1}}(a)\bigr)\xi \otimes \delta_g,
\qquad a\in A,\ \xi\in H,\ g\in\Gamma,
\]
where $\{\delta_g\}_{g\in\Gamma}$ denotes the canonical orthonormal basis
of $\ell^2(\Gamma)$.  Together with the unitary representation
\[
(1\otimes\lambda_g)(\xi\otimes\delta_h)
=
\xi \otimes \delta_{gh},
\qquad g,h\in\Gamma,
\]
this defines a covariant representation $(\tilde{\pi},1\otimes\lambda)$
of $(A,\Gamma,\alpha)$ on $H\otimes\ell^2(\Gamma)$.

The \emph{reduced crossed product} $A\rtimes_{\alpha,r}\Gamma$
is the completion of $A\rtimes_{\mathrm{alg}}\Gamma$ with respect to the norm
\[
\|x\|_r = \|(\tilde{\pi} \rtimes \lambda)(x)\|.
\]
The norm $\|x\|_r$ is independent of the choice of faithful representation
$\pi$, and the canonical surjection $A\rtimes_\alpha\Gamma \to A\rtimes_{\alpha,r}\Gamma$
extends to a quotient map.
\end{definition}

\begin{remark}[Amenability]
\label{rem:amenability-full-reduced}
If $\Gamma$ is amenable, the canonical quotient map
\[
A\rtimes_\alpha\Gamma \longrightarrow A\rtimes_{\alpha,r}\Gamma
\]
is an isomorphism.  In this case the full and reduced crossed products
coincide, and we simply write $A\rtimes\Gamma$.
This fact will be important later, as it allows us to work with reduced
crossed products when relating our construction to the Baum--Connes
assembly map.
\end{remark}

An important special case arises when $A=C(X)$ for a compact Hausdorff
space $X$ and $\Gamma$ acts on $X$ by homeomorphisms.
The induced action on $C(X)$ is
\[
(\alpha_g f)(x)=f(g^{-1}x),
\qquad f\in C(X),\ g\in\Gamma.
\]
The crossed product $C(X)\rtimes_\alpha\Gamma$
encodes the dynamical system $(X,\Gamma)$.

This algebra admits a natural geometric interpretation as the $C^*$-algebra
of the \emph{transformation groupoid}
\[
X\rtimes\Gamma .
\]

\begin{proposition}
\label{prop:transformation-groupoid-crossed-product}
Let $\Gamma$ be a discrete group acting on a compact Hausdorff space $X$.
Then the reduced crossed product is canonically isomorphic to the
reduced groupoid $C^*$-algebra:
\[
C(X)\rtimes_{\alpha,r}\Gamma \;\cong\; C_r^*(X\rtimes\Gamma).
\]

More explicitly, the isomorphism sends a function $f\in C_c(X\rtimes\Gamma)$
to the finitely supported function $\Gamma\to C(X)$ given by
$g\mapsto f(\cdot,g)$.
\end{proposition}

\begin{proof}
The transformation groupoid $X\rtimes\Gamma$ is the locally compact
\'etale groupoid whose underlying set is $X\times\Gamma$ with
source and range maps
\[
s(x,g)=g^{-1}x,\qquad r(x,g)=x,
\]
and multiplication
\[
(x,g)(g^{-1}x,h)=(x,gh).
\]

Since $\Gamma$ is discrete and $X$ is compact, the space
$C_c(X\rtimes\Gamma)$ consists of functions $F:X\times\Gamma\to\mathbb{C}$
with finite support in the $\Gamma$-variable.
Any such function can be written uniquely as
\[
F(x,g)=f_g(x)
\]
with $f_g\in C(X)$ and only finitely many $f_g$ nonzero.
Thus we have a linear isomorphism
\[
C_c(X\rtimes\Gamma)\cong C_c(\Gamma,C(X)),
\qquad
F \mapsto (g\mapsto f_g).
\]

Under this identification, the groupoid convolution
\[
(F*G)(x,g)=\sum_{h\in\Gamma} F(x,h)\,G(h^{-1}x,h^{-1}g)
\]
becomes
\[
(f*g)_g = \sum_{h\in\Gamma} f_h\;\alpha_h(g_{h^{-1}g}),
\]
where $\alpha$ denotes the induced action on $C(X)$.
This is precisely the multiplication in the algebraic crossed product
$C(X)\rtimes_{\mathrm{alg}}\Gamma$, and a similar calculation shows that
the involution corresponds as well.

Finally, the regular representation of the groupoid $X\rtimes\Gamma$
on $L^2(X\times\Gamma)$ coincides with the regular covariant
representation used to define the reduced crossed product.
Indeed, for each $x\in X$, the representation of $C_c(X\rtimes\Gamma)$
on $\ell^2(\Gamma)$ obtained by fixing $x$ corresponds under the
identification to the evaluation representation of $C(X)$ at $x$ tensored
with the left regular representation of $\Gamma$.

Since $\Phi$ is an isometric $*$-isomorphism between dense subalgebras
with respect to the reduced norms, it extends uniquely to an isomorphism
of the completions:
\[
C_r^*(X\rtimes\Gamma) \cong C(X)\rtimes_{\alpha,r}\Gamma.
\]
\end{proof}

\begin{remark}
When $\Gamma$ is amenable, the full and reduced crossed products coincide,
and we also have $C^*(X\rtimes\Gamma) \cong C(X)\rtimes_\alpha\Gamma$.
In this case the isomorphism extends to the full groupoid $C^*$-algebra.
\end{remark}

\paragraph{K-theory and the Pimsner--Voiculescu sequence}

The $K$-theory of crossed products by $\mathbb Z$ can be computed using
the Pimsner--Voiculescu exact sequence, a fundamental tool in
operator $K$-theory.

\begin{theorem}[Pimsner--Voiculescu exact sequence]
\label{thm:pimsner-voiculescu}
Let $(A,\mathbb Z,\alpha)$ be a $C^*$-dynamical system.
Then there is a cyclic six-term exact sequence
\[
\begin{tikzcd}
K_0(A) \arrow[r,"1-\alpha_*"] &
K_0(A) \arrow[r] &
K_0(A\rtimes_\alpha\mathbb Z) \arrow[d] \\
K_1(A\rtimes_\alpha\mathbb Z) \arrow[u] &
K_1(A) \arrow[l] &
K_1(A) \arrow[l,"1-\alpha_*"']
\end{tikzcd}
\]
where the maps $1-\alpha_*$ are induced by the automorphism $\alpha$
on $K$-theory.
\end{theorem}

\begin{proof}[Sketch]
See the original work of Pimsner and Voiculescu \cite{PimsnerVoiculescu1980}.
For a detailed exposition, we refer the reader to
\cite[Chapter 9]{Blackadar1998} or \cite[Section 10.3]{RordamLarsenLaustsen2000}.
\end{proof}

\subsection{The Irrational Rotation Algebra}
\label{subsec:irrational-rotation-algebra}

One of the most important examples of a crossed product $C^*$-algebra is the
irrational rotation algebra, also known as the noncommutative $2$-torus.
It arises from the simplest nontrivial dynamical system: an irrational
rotation of the circle.

\begin{example}[The irrational rotation algebra]
\label{ex:A-theta-crossed-product}
Let $\theta \in \mathbb{R}\setminus\mathbb{Q}$ be an irrational number.
Consider the rotation action of $\mathbb{Z}$ on the circle $S^1$ given by
\[
n\cdot z = e^{2\pi i n\theta}z, \qquad n\in\mathbb{Z},\ z\in S^1.
\]
This induces an action $\alpha:\mathbb{Z}\to\operatorname{Aut}(C(S^1))$ defined by
\[
\alpha_n(f)(z) = f(n^{-1}\cdot z) = f(e^{-2\pi i n\theta}z),
\qquad f\in C(S^1).
\]

The crossed product
\[
A_\theta := C(S^1) \rtimes_\alpha \mathbb{Z}
\]
is called the \emph{irrational rotation algebra} (or noncommutative $2$-torus).
Since $\mathbb{Z}$ is amenable, the full and reduced crossed products coincide.

Let $u \in C(S^1)$ denote the canonical unitary $u(z)=z$, and let
$v$ be the unitary implementing the $\mathbb{Z}$-action in the crossed product.
These two unitaries satisfy the fundamental commutation relation
\[
vu = e^{2\pi i\theta} uv,
\]
and $A_\theta$ is the universal $C^*$-algebra generated by two unitaries
subject to this relation.  Geometrically, $A_\theta$ may be viewed as a
deformation of the commutative torus algebra $C(\mathbb{T}^2)$, to which it
reduces when $\theta = 0$.
\end{example}

The $K$-theory of the irrational rotation algebra was computed by
Pimsner and Voiculescu in their seminal work.

\begin{example}[K-theory of $A_\theta$]
\label{ex:K-theory-A-theta}
Applying the Pimsner--Voiculescu exact sequence
(Theorem~\ref{thm:pimsner-voiculescu}) to the dynamical system
$(C(S^1),\mathbb{Z},\alpha)$ yields
\[
K_0(A_\theta) \cong \mathbb{Z}^2, \qquad K_1(A_\theta) \cong \mathbb{Z}^2.
\]

Indeed, the $K$-theory of $C(S^1)$ is $K_0(C(S^1))\cong\mathbb{Z}$ (generated by
the class of the unit) and $K_1(C(S^1))\cong\mathbb{Z}$ (generated by the class
of the identity function $u$).  Since the action $\alpha$ is induced by a
rotation, the induced map on $K$-theory satisfies $\alpha_* = \operatorname{id}$
on both $K_0$ and $K_1$.  The Pimsner--Voiculescu six-term exact sequence
therefore splits into short exact sequences, giving the result above.

Canonical generators are given by:
\begin{itemize}
    \item In $K_1(A_\theta)$: the classes of the two unitaries $[u]$ and $[v]$;
    \item In $K_0(A_\theta)$: the class of the identity $[1_A]$ together with
          the class of the Rieffel--Power projection $[p_\theta]$, which satisfies
          $\tau(p_\theta) = \theta$ for the unique tracial state $\tau$ on $A_\theta$.
\end{itemize}
The image of the trace on $K_0(A_\theta)$ is therefore $\mathbb{Z} + \theta\mathbb{Z}$,
a dense subgroup of $\mathbb{R}$.
\end{example}

This algebra plays a central role in noncommutative geometry and serves
as a fundamental test case for index theory; see \cite{ConnesNCG} and
\cite{Connes1980}.

\paragraph{Role in this paper}

Crossed products provide a natural link between group actions,
transformation groupoids, and equivariant $K$-theory.
In particular, when a discrete group $\Gamma$ acts continuously on a
compact Hausdorff space $X$, the crossed product $C(X)\rtimes_\alpha\Gamma$
is canonically isomorphic to the groupoid $C^*$-algebra of the
transformation groupoid $X\rtimes\Gamma$ (Proposition~\ref{prop:transformation-groupoid-crossed-product}).

In the following sections we compare the unitary conjugation groupoid
$\mathcal{G}_{C(X)\rtimes_\alpha\Gamma}$ associated to this crossed product
with the transformation groupoid $X\rtimes\Gamma$.
This comparison yields a Morita equivalence
\[
\mathcal{G}_{C(X)\rtimes_\alpha\Gamma} \sim_M X\rtimes\Gamma
\]
(Theorem~\ref{thm:morita-equivalence-crossed}), which allows us to transfer
the descent-index machinery developed in Paper~II to the non-type~I setting.

The main result of this comparison is a commutative diagram
\[
\begin{tikzcd}
KK^1_{\mathcal{G}_{C(X)\rtimes_\alpha\Gamma}}(C_0(\mathcal{G}^{(0)}),\mathbb C)
\arrow[r, "\operatorname{desc}_{\mathcal{G}}"]
\arrow[d, "\cong"]
&
K_1(C^*(\mathcal{G}_{C(X)\rtimes_\alpha\Gamma}))
\arrow[d, "\cong"]
\\
KK^1_\Gamma(C(X),\mathbb C)
\arrow[r, "\mu_\Gamma"]
&
K_1(C(X)\rtimes_\alpha\Gamma)
\end{tikzcd}
\]
where the left vertical isomorphism is induced by the Morita equivalence,
the right vertical isomorphism follows from
$C^*(\mathcal{G}_{C(X)\rtimes_\alpha\Gamma}) \sim_M (C(X)\rtimes_\alpha\Gamma)\otimes\mathcal{K}$,
and $\mu_\Gamma$ denotes the Baum--Connes assembly map with coefficients
in $C(X)$ \cite{BaumConnesHigson1994}.

In particular, we show that the descent map associated to the unitary
conjugation groupoid is compatible with the analytic assembly map under
this Morita equivalence.  When $\Gamma$ is amenable, so that the full and
reduced crossed products coincide, this compatibility identifies our
descent-index construction with the Baum--Connes assembly map, providing
a new geometric realization of this fundamental invariant.

\subsection{Transformation Groupoids and Their $C^*$-Algebras}
\label{subsec:transformation-groupoids}

Groupoids provide a natural framework for describing dynamical systems
and their associated operator algebras.  In particular, actions of
discrete groups on topological spaces can be encoded as
\emph{transformation groupoids}.  These groupoids play a central role in
noncommutative geometry and provide a geometric bridge between crossed
product $C^*$-algebras and groupoid $C^*$-algebras.  They will be indispensable
in this paper, as both the irrational rotation algebra $A_\theta$ and the
amenable crossed products $C(X)\rtimes\Gamma$ admit concrete realizations
as groupoid $C^*$-algebras of transformation groupoids; see \cite{Renault1980}.

\paragraph{Transformation groupoids.}

Let $\Gamma$ be a discrete group acting continuously on a locally compact
Hausdorff space $X$ by homeomorphisms.  The associated
\emph{transformation groupoid} (also called the \emph{action groupoid}) is
\[
X \rtimes \Gamma := X \times \Gamma,
\]
with the following structure maps.

\begin{itemize}
    \item \textbf{Range and source maps:} 
    \[
    r(x,g) = x, \qquad s(x,g) = g^{-1}x.
    \]
    
    \item \textbf{Multiplication:} Two arrows $(x,g)$ and $(y,h)$ are
    composable iff $s(x,g) = r(y,h)$, i.e., $g^{-1}x = y$.  In this case,
    \[
    (x,g)(g^{-1}x,h) = (x,gh).
    \]
    
    \item \textbf{Inverse:} 
    \[
    (x,g)^{-1} = (g^{-1}x, g^{-1}).
    \]
    
    \item \textbf{Unit space:} The unit space is naturally identified with $X$ via
    \[
    x \mapsto (x,e),
    \]
    where $e$ denotes the identity element of $\Gamma$.
\end{itemize}

Topologically, $X\rtimes\Gamma$ is equipped with the product topology,
where $\Gamma$ carries the discrete topology.  This makes
$X\rtimes\Gamma$ a locally compact Hausdorff groupoid.

\paragraph{The \'etale property.}

A crucial feature is that $X\rtimes\Gamma$ is an \'etale groupoid.
Indeed, for each fixed $g\in\Gamma$, the subset $X\times\{g\}$ is open
in the product topology, and the restriction of the source map
\[
s|_{X\times\{g\}} : X\times\{g\} \longrightarrow X
\]
is a homeomorphism (its inverse is $x \mapsto (gx, g)$).  Consequently,
the source map $s: X\rtimes\Gamma \to X$ is a local homeomorphism.

This \'etaleness has important consequences for the associated operator
algebras.  The convolution algebra $C_c(X\rtimes\Gamma)$ admits a
particularly simple structure, and the reduced groupoid $C^*$-algebra
$C_r^*(X\rtimes\Gamma)$ carries a canonical faithful conditional
expectation onto $C_0(X)$, given by restricting functions to the unit
space.

\begin{remark}
The convention $r(x,g)=x$, $s(x,g)=g^{-1}x$ used here is chosen to align
with the standard convolution formula in the crossed product
$C(X)\rtimes\Gamma$.  Under this convention, an element of the groupoid
$C^*$-algebra corresponds to a function $f:X\times\Gamma\to\mathbb{C}$,
and the convolution product becomes
\[
(f_1 * f_2)(x,g) = \sum_{h\in\Gamma} f_1(x,h)\, f_2(h^{-1}x, h^{-1}g),
\]
which matches the crossed product multiplication.
\end{remark}

\paragraph{Groupoid $C^*$-algebras.}

Let $G = X \rtimes \Gamma$ be the transformation groupoid.
The convolution algebra $C_c(G)$ consists of compactly supported
continuous functions $f : G \to \mathbb C$.

For $f_1, f_2 \in C_c(G)$, the convolution product is defined by
\[
(f_1 * f_2)(x,g)
=
\sum_{h \in \Gamma}
f_1(x,h)\, f_2(h^{-1}x, h^{-1}g),
\]
and the involution by
\[
f^*(x,g)
=
\overline{f(g^{-1}x, g^{-1})}.
\]
These operations make $C_c(G)$ into a $*$-algebra.

The \emph{full groupoid $C^*$-algebra} $C^*(G)$ is the completion of
$C_c(G)$ with respect to the universal norm
\[
\|f\|_{\max}
=
\sup\bigl\{
\|\pi(f)\| :
\pi \text{ is a $*$-representation of } C_c(G)
\bigr\}.
\]

The \emph{reduced groupoid $C^*$-algebra} $C_r^*(G)$ is defined via the
regular representations.  For each $x \in X$, define a representation
$\pi_x : C_c(G) \to B(\ell^2(\Gamma))$ by
\[
(\pi_x(f)\xi)(g)
=
\sum_{h \in \Gamma}
f(x,h)\, \xi(h^{-1}g),
\qquad \xi \in \ell^2(\Gamma),\ g\in\Gamma.
\]

The reduced norm is then given by
\[
\|f\|_r
=
\sup_{x \in X} \|\pi_x(f)\|,
\]
and $C_r^*(G)$ is the completion of $C_c(G)$ with respect to this norm.
Equivalently, $C_r^*(G)$ may be obtained as the completion of $C_c(G)$
in the norm induced by the direct integral of the representations
$\pi_x$ over $X$ with respect to any quasi-invariant measure; the two
definitions are equivalent because the supremum over $x$ coincides with
the essential supremum with respect to any such measure.

\begin{remark}
The formulas above are consistent with the convention
$r(x,g)=x$, $s(x,g)=g^{-1}x$ adopted in
Section~\ref{subsec:transformation-groupoids}.  Under this convention,
the convolution product matches the multiplication in the crossed product
$C(X)\rtimes\Gamma$, and the involution corresponds to the usual
$*$-operation in the groupoid algebra.
\end{remark}

\paragraph{Relation to crossed products}

For transformation groupoids associated with discrete group actions,
these constructions admit a particularly concrete and important
description.

\begin{proposition}
\label{prop:transformation-crossed-product}
Let $\Gamma$ be a discrete group acting on a locally compact Hausdorff
space $X$ by homeomorphisms, and let $\alpha$ denote the induced action
of $\Gamma$ on $C_0(X)$. Then there is a canonical isomorphism
\[
C_r^*(X\rtimes\Gamma) \cong C_0(X)\rtimes_{\alpha,r}\Gamma .
\]

More precisely, the isomorphism is implemented by the map
$\Phi: C_c(X\rtimes\Gamma) \to C_c(\Gamma, C_0(X))$ defined by
\[
\Phi(f)(g)(x) = f(x,g), \qquad f\in C_c(X\rtimes\Gamma),\ g\in\Gamma,\ x\in X.
\]
\end{proposition}

\begin{proof}
Let $G = X\rtimes\Gamma$ be the transformation groupoid.  Since $\Gamma$
is discrete, every function $f\in C_c(G)$ has finite support in the
$\Gamma$-variable, and may be written uniquely as
\[
f(x,g) = f_g(x), \qquad g\in\Gamma,
\]
where each $f_g \in C_0(X)$ and only finitely many $f_g$ are nonzero.
Thus the assignment $\Phi$ defined above is a linear bijection.

We claim that $\Phi$ is a $*$-algebra isomorphism.  Recall that for the
transformation groupoid with source and range maps
\[
r(x,g)=x,\qquad s(x,g)=g^{-1}x,
\]
the convolution product on $C_c(G)$ is given by
\[
(f_1 * f_2)(x,g) = \sum_{h\in\Gamma} f_1(x,h)\, f_2(h^{-1}x, h^{-1}g),
\]
and the involution by
\[
f^*(x,g) = \overline{f(g^{-1}x, g^{-1})}.
\]

On the other hand, the algebraic crossed product
$C_0(X)\rtimes_{\mathrm{alg}}\Gamma$ may be identified with
$C_c(\Gamma, C_0(X))$, equipped with the convolution product
\[
(F_1 * F_2)(g) = \sum_{h\in\Gamma} F_1(h)\; \alpha_h\!\bigl(F_2(h^{-1}g)\bigr)
\]
and involution
\[
F^*(g) = \alpha_g\!\bigl(F(g^{-1})^*\bigr).
\]

\textbf{Verification of the product.}
For $f_1, f_2 \in C_c(G)$, set $F_i = \Phi(f_i)$.  Using the definition
of the induced action $(\alpha_h(\varphi))(x) = \varphi(h^{-1}x)$, we
compute:
\[
\begin{aligned}
\bigl(F_1 * F_2\bigr)(g)(x)
&= \sum_{h\in\Gamma} F_1(h)(x)\; \alpha_h\!\bigl(F_2(h^{-1}g)\bigr)(x) \\
&= \sum_{h\in\Gamma} f_1(x,h)\; F_2(h^{-1}g)(h^{-1}x) \\
&= \sum_{h\in\Gamma} f_1(x,h)\; f_2(h^{-1}x, h^{-1}g) \\
&= (f_1 * f_2)(x,g) = \Phi(f_1 * f_2)(g)(x).
\end{aligned}
\]
Hence $\Phi(f_1 * f_2) = \Phi(f_1) * \Phi(f_2)$.

\textbf{Verification of the involution.}
Similarly,
\[
\begin{aligned}
\Phi(f^*)(g)(x) &= f^*(x,g) = \overline{f(g^{-1}x, g^{-1})}
= \overline{\Phi(f)(g^{-1})(g^{-1}x)}.
\end{aligned}
\]
On the other hand,
\[
\begin{aligned}
\Phi(f)^*(g)(x) &= \bigl(\alpha_g(\Phi(f)(g^{-1})^*)\bigr)(x)
= \Phi(f)(g^{-1})^*(g^{-1}x)
= \overline{\Phi(f)(g^{-1})(g^{-1}x)}.
\end{aligned}
\]
Thus $\Phi(f^*) = \Phi(f)^*$, proving that $\Phi$ is a $*$-algebra
isomorphism between $C_c(X\rtimes\Gamma)$ and
$C_0(X)\rtimes_{\mathrm{alg}}\Gamma$.

\textbf{Comparison of the reduced norms.}
For each $x\in X$, the regular representation $\pi_x$ of $C_c(G)$ on
$\ell^2(\Gamma)$ is given by
\[
(\pi_x(f)\xi)(g) = \sum_{h\in\Gamma} f(x,h)\, \xi(h^{-1}g),
\qquad \xi\in\ell^2(\Gamma),\ g\in\Gamma.
\]

Under the identification $\Phi$, this is exactly the regular covariant
representation of the algebraic crossed product associated to the
evaluation representation $\operatorname{ev}_x: C_0(X)\to\mathbb{C}$.
Consequently, for each $f\in C_c(G)$,
\[
\|\pi_x(f)\| = \|(\operatorname{ev}_x \rtimes \lambda)(\Phi(f))\|,
\]
where $\lambda$ denotes the left regular representation of $\Gamma$.
Taking the supremum over $x\in X$, we obtain
\[
\|f\|_{r} = \sup_{x\in X}\|\pi_x(f)\| = \|\Phi(f)\|_{r},
\]
where the right-hand side denotes the reduced crossed product norm on
$C_0(X)\rtimes_{\mathrm{alg}}\Gamma$.

Since $\Phi$ is an isometric $*$-isomorphism between dense subalgebras
with respect to the reduced norms, it extends uniquely to an isomorphism
of the completions:
\[
\Phi: C_r^*(X\rtimes\Gamma) \xrightarrow{\cong} C_0(X)\rtimes_{\alpha,r}\Gamma.
\]

This completes the proof.
\end{proof}

\begin{corollary}
\label{cor:groupoid-crossed-spectrum}
If $\Gamma$ is amenable, then the full and reduced crossed products
coincide, and there is a natural isomorphism
\[
C^*(X\rtimes\Gamma) \cong C_0(X)\rtimes_\alpha \Gamma .
\]

Moreover, the transformation groupoid $X\rtimes\Gamma$ is amenable,
so its full and reduced groupoid $C^*$-algebras coincide as well:
\[
C^*(X\rtimes\Gamma) \cong C_r^*(X\rtimes\Gamma).
\]
\end{corollary}

\begin{proof}
If $\Gamma$ is amenable, then for any $C^*$-dynamical system
$(A,\Gamma,\alpha)$, the canonical quotient map
$A\rtimes_\alpha\Gamma \to A\rtimes_{\alpha,r}\Gamma$ is an isomorphism.
Applying this to $A = C_0(X)$ yields
$C_0(X)\rtimes_\alpha\Gamma \cong C_0(X)\rtimes_{\alpha,r}\Gamma$.

Furthermore, when $\Gamma$ is amenable, the transformation groupoid
$X\rtimes\Gamma$ is amenable, and for amenable groupoids the full and
reduced groupoid $C^*$-algebras coincide (see \cite{Tu1999,Renault1980}).
Hence $C^*(X\rtimes\Gamma) \cong C_r^*(X\rtimes\Gamma)$.

Combining these facts with Proposition~\ref{prop:transformation-crossed-product}
gives the desired isomorphism:
\[
C^*(X\rtimes\Gamma) \cong C_r^*(X\rtimes\Gamma)
\cong C_0(X)\rtimes_{\alpha,r}\Gamma \cong C_0(X)\rtimes_\alpha\Gamma.
\]
\end{proof}

\begin{remark}
The isomorphism $\Phi$ is natural in the sense that it commutes with
induction by $\Gamma$-equivariant maps between spaces and with
restrictions to invariant subspaces.  For a detailed treatment of these
isomorphisms, see \cite{Renault1980,Williams2007}.
\end{remark}

\paragraph{Amenability.}

Amenability plays an important role in the study of transformation
groupoids and crossed products, and will be a recurring theme in this
paper.

We will not recall the full definition of amenability for locally compact
groupoids here; instead, we only use the standard fact that if a
groupoid $G$ is amenable (in the sense of Anantharaman-Delaroche and
Renault \cite{AnantharamanDelarocheRenault2000}), then its full and
reduced groupoid $C^*$-algebras coincide:
\[
C^*(G) \cong C_r^*(G).
\]

If $\Gamma$ is an amenable discrete group acting on a locally compact
Hausdorff space $X$, then the transformation groupoid $X\rtimes\Gamma$
is amenable.  This observation will be particularly important in the
sequel, since amenability allows one to identify full and reduced
groupoid $C^*$-algebras, and likewise full and reduced crossed products.
In particular, for amenable $\Gamma$ we have
$C_0(X)\rtimes_\alpha\Gamma \cong C_0(X)\rtimes_{\alpha,r}\Gamma$, and the
Baum--Connes assembly map is an isomorphism by Tu's theorem \cite{Tu1999}.

\paragraph{$K$-theory and the Baum--Connes assembly map.}

Transformation groupoids provide a natural geometric framework for the
Baum--Connes conjecture with coefficients.  For a $\Gamma$-$C^*$-algebra
$A$, the assembly map
\[
\mu_\Gamma^A:
K_*^\Gamma(\underline{E}\Gamma;A) \longrightarrow K_*(A\rtimes_r\Gamma)
\]
takes values in the $K$-theory of the reduced crossed product.
The left-hand side denotes the equivariant $KK$-theory group
$KK_\Gamma^*(C_0(\underline{E}\Gamma), A)$, where $\underline{E}\Gamma$
is the classifying space for proper $\Gamma$-actions.

In the special case $A = C_0(X)$ for a $\Gamma$-space $X$, the reduced
crossed product $C_0(X)\rtimes_r\Gamma$ may be identified with the
reduced groupoid $C^*$-algebra $C_r^*(X\rtimes\Gamma)$ of the
transformation groupoid (Proposition~\ref{prop:transformation-crossed-product}).
This identification allows us to reinterpret the assembly map as a
groupoid-theoretic construction, which will be crucial in
Section~\ref{sec:crossed} when we compare our descent map with the
assembly map.

\begin{example}[The irrational rotation groupoid]
\label{ex:rotation-groupoid}
For the irrational rotation algebra
\[
A_\theta = C(S^1)\rtimes_\alpha \mathbb Z,
\]
the associated transformation groupoid is
\[
S^1\rtimes_\theta \mathbb Z.
\]
Its unit space is the circle $S^1$, and its arrows are pairs $(z,n)$ with
$z\in S^1$ and $n\in\mathbb Z$.  Following the convention established in
Section~\ref{subsec:transformation-groupoids}, the range and source maps are
\[
r(z,n) = z, \qquad s(z,n) = e^{-2\pi i n\theta}z.
\]
Equivalently, one may take $s(z,n)=z$, $r(z,n)=e^{2\pi i n\theta}z$; both
conventions appear in the literature, and we will consistently use the
former to align with our earlier choice.

This groupoid is \'etale (since $\mathbb Z$ is discrete) and amenable
(because $\mathbb Z$ is amenable).  Its reduced (equivalently, full)
groupoid $C^*$-algebra is canonically isomorphic to $A_\theta$.

Since $\theta$ is irrational, the action of $\mathbb Z$ on $S^1$ is
\emph{free} (every stabilizer is trivial) and \emph{minimal} (every orbit
is dense).  This dynamical behavior is reflected in the structure of the
associated $C^*$-algebra: $A_\theta$ is simple and has no nontrivial
ideals, a hallmark of minimal actions.
\end{example}

\begin{remark}
The irrational rotation groupoid is a paradigmatic example of a minimal
\'etale groupoid.  Its simplicity and the density of orbits make it an
ideal testing ground for the constructions developed in this paper, as
we saw in Case Study~I (Section~\ref{sec:Atheta}).
\end{remark}

\paragraph{Relation to the present work.}

Transformation groupoids provide the geometric model for the crossed
products considered in this paper.  In the following sections we compare
the unitary conjugation groupoid
\[
\mathcal{G}_{C_0(X)\rtimes_\alpha\Gamma}
\]
with the transformation groupoid $X\rtimes\Gamma$.

Under suitable Morita equivalences, this comparison allows us to relate
the descent map associated to $\mathcal{G}_{C_0(X)\rtimes_\alpha\Gamma}$
to the Baum--Connes assembly map with coefficients in $C_0(X)$:
\[
\mu_\Gamma :
K_*^\Gamma(X)\longrightarrow
K_*(C_0(X)\rtimes_{\alpha,r}\Gamma),
\]
where $K_*^\Gamma(X)$ denotes the equivariant $K$-homology of $X$,
or more precisely $KK_\Gamma^*(C_0(X), \mathbb C)$.

More concretely, we prove that the following diagram commutes:

\[
\begin{tikzcd}
KK^1_{\mathcal{G}_{C_0(X)\rtimes_\alpha\Gamma}}
\!\bigl(C_0(\mathcal{G}^{(0)}),\mathbb C\bigr)
   \arrow[r, "\operatorname{desc}_{\mathcal{G}}"]
   \arrow[d, "\cong"']
&
K_1\!\bigl(C_r^*(\mathcal{G}_{C_0(X)\rtimes_\alpha\Gamma})\bigr)
   \arrow[d, "\cong"] \\
KK^1_\Gamma(C_0(X),\mathbb C)
   \arrow[r, "\mu_\Gamma"]
&
K_1\!\bigl(C_0(X)\rtimes_{\alpha,r}\Gamma\bigr)
\end{tikzcd}
\]

Here:
\begin{itemize}
    \item The left vertical isomorphism is induced by the Morita equivalence
          $\mathcal{G}_{C_0(X)\rtimes_\alpha\Gamma} \sim_M X\rtimes\Gamma$
          (Theorem~\ref{thm:morita-equivalence-crossed});
    \item The right vertical isomorphism follows from the Morita equivalence
          $C_r^*(\mathcal{G}_{C_0(X)\rtimes_\alpha\Gamma}) \sim_M (C_0(X)\rtimes_{\alpha,r}\Gamma)\otimes\mathcal{K}$
          together with stability of $K$-theory;
    \item $\mu_\Gamma$ denotes the Baum--Connes assembly map with coefficients
          in $C_0(X)$, identified via Tu's theorem \cite{Tu1999} with the descent
          map for the transformation groupoid $X\rtimes\Gamma$.
\end{itemize}

This diagram provides a conceptual bridge between the groupoid framework
developed in Papers~I and~II and the analytic machinery of the
Baum--Connes conjecture.  In particular, it identifies the descent map
arising from the unitary conjugation groupoid with the Baum--Connes
assembly map under Morita equivalence.  Establishing this identification
is one of the main structural results of the present paper.

\begin{remark}[Main takeaway]
The descent map associated to the unitary conjugation groupoid
recovers the Baum--Connes assembly map under Morita equivalence.
This provides a new geometric realization of the assembly map in terms
of the internal symmetries of the crossed product algebra.
\end{remark}

\subsection{Renault's Cartan Subalgebras and Weyl Groupoids}\label{subsec:reduction-cartan-pairs}
\label{subsec:renault-cartan-weyl}

The inability to construct the unitary conjugation groupoid $\mathcal{G}_{\mathcal{A}}$ for non-type~I algebras such as $A_\theta$, as discussed in Section~\ref{subsec:typeI-limitation}, signals the need for a different geometric paradigm. A powerful alternative is provided by Renault's theory of Cartan subalgebras and Weyl groupoids \cite{Renault2008Cartan,Renault1980}. This framework shows that a large class of $C^*$-algebras can be reconstructed from an \'etale groupoid built from a distinguished commutative subalgebra and its normalizers.

\subsubsection*{Cartan Subalgebras}
The theory begins with a commutative subalgebra that encodes a chosen ``classical floor" of the noncommutative algebra.

\begin{definition}[Cartan subalgebra]
\label{def:cartan-subalgebra}
Let $A$ be a separable $C^*$-algebra. A subalgebra $D \subseteq A$ is called a \emph{Cartan subalgebra} if:
\begin{enumerate}
    \item $D$ is a maximal abelian $C^*$-subalgebra (MASA) of $A$;
    \item $D$ contains an approximate unit for $A$ (if $A$ is unital, this means $1_A \in D$);
    \item $D$ is \emph{regular}, i.e., its normalizer
    \[
    \mathcal{N}_A(D) := \{ n \in A \mid n D n^* \subseteq D \text{ and } n^* D n \subseteq D \}
    \]
    generates $A$ as a $C^*$-algebra;
    \item there exists a faithful conditional expectation $E: A \to D$.
\end{enumerate}
A pair $(A, D)$ satisfying these conditions is called a \emph{Cartan pair}.
\end{definition}

The normalizers $\mathcal{N}_A(D)$ are the key players; they implement symmetries of the Cartan subalgebra, acting as partial automorphisms on its Gelfand spectrum $\widehat{D}$. The conditional expectation $E$ provides a projection onto the "classical" part of the algebra and is essential for the reconstruction theorem.

\begin{example}[Examples of Cartan pairs]
\label{ex:cartan-pairs}
\begin{itemize}
    \item \textbf{Commutative algebras.}
    For $A = C_0(X)$ with $X$ locally compact Hausdorff, $D = A$ itself is a Cartan subalgebra. The conditional expectation is the identity map, and the normalizer is the unitary group of $A$, consisting of continuous functions $X \to S^1$.
    
    \item \textbf{Transformation groupoids (topologically free actions).}
    Let $\Gamma$ be a discrete group acting on a locally compact Hausdorff space $X$. If the action is \emph{topologically free} (i.e., the set of points with nontrivial stabilizer has empty interior), then the reduced crossed product
    \[
    A = C_0(X) \rtimes_r \Gamma
    \]
    contains $D = C_0(X)$ as a Cartan subalgebra. The conditional expectation
    \[
    E\!\left( \sum_{\gamma\in\Gamma} f_\gamma u_\gamma \right) = f_e
    \]
    is faithful, and the normalizer is generated by the implementing unitaries $\{u_\gamma\}_{\gamma\in\Gamma}$ together with the unitaries of $C_0(X)$. This is a fundamental class of examples arising from Renault's groupoid construction.
    
    \item \textbf{The irrational rotation algebra.}
    Let $\theta \in \mathbb{R}\setminus\mathbb{Q}$ be irrational. Then the action of $\mathbb{Z}$ on $S^1$ by rotation is free and minimal, hence topologically free. Consequently, the irrational rotation algebra
    \[
    A_\theta = C(S^1) \rtimes_\theta \mathbb{Z}
    \]
    (see Example~\ref{ex:A-theta-crossed-product}) admits $D = C(S^1)$ as a Cartan subalgebra, with conditional expectation given by the canonical projection onto the fixed-point subalgebra.
    
    \item \textbf{Finite-dimensional algebras (degenerate case).}
    For $A = M_n(\mathbb{C})$, the diagonal subalgebra $D_n \cong \mathbb{C}^n$ is a maximal abelian subalgebra admitting a faithful conditional expectation (the map that removes off-diagonal entries). The normalizer consists of the generalized permutation matrices (unitary matrices with exactly one nonzero entry in each row and column). Although this example formally satisfies the defining conditions, it does not exhibit the rich geometric structure typical of Cartan pairs arising from \'etale groupoids; its associated Weyl groupoid is a finite groupoid with discrete topology.
\end{itemize}
\end{example}

\begin{remark}[On the definition]
The definition of a Cartan subalgebra given above follows Renault \cite{Renault2008Cartan}, though variations exist in the literature. Some authors omit the separability assumption, require the conditional expectation to be of finite index, or impose additional regularity conditions. The conditions listed here are sufficient for Renault's reconstruction theorem, which we will recall below.
\end{remark}

\begin{remark}[Topological freeness is essential]
The condition that the action be topologically free in the crossed product example is not merely technical; it is necessary for $C_0(X)$ to be maximal abelian in $C_0(X)\rtimes_r\Gamma$. If the action has nontrivial isotropy on an open set, then elements supported on that isotropy can commute with $C_0(X)$, violating maximal abelianness. In such cases, the Cartan structure must be replaced by a twisted groupoid construction.
\end{remark}

\subsubsection*{The Weyl Groupoid}

Given a Cartan pair $(A,D)$, Renault associates to it an \'etale,
locally compact Hausdorff groupoid $\mathcal{G}(A,D)$, called the
\emph{Weyl groupoid}, together with a twist $\Sigma$ in general, such
that $A$ may be reconstructed as the reduced twisted groupoid
$C^*$-algebra of $(\mathcal{G}(A,D),\Sigma)$ \cite{Renault2008Cartan}.
In the cases relevant to this paper, the twist is trivial, so that one
recovers an untwisted groupoid $C^*$-algebra.

The construction proceeds as follows:

\begin{enumerate}
    \item \textbf{Unit space.}
    The unit space is identified with the Gelfand spectrum of the
    Cartan subalgebra:
    \[
    \mathcal{G}(A,D)^{(0)} \cong \widehat{D}.
    \]

    \item \textbf{Partial homeomorphisms from normalizers.}
    Each normalizer $n \in \mathcal{N}_A(D)$ determines a partial
    homeomorphism
    \[
    \alpha_n : \operatorname{dom}(n) \longrightarrow \operatorname{ran}(n)
    \]
    of $\widehat{D}$, where
    \[
    \operatorname{dom}(n) = \{ \phi \in \widehat{D} : \phi(n^* n) > 0 \},
    \qquad
    \operatorname{ran}(n) = \{ \phi \in \widehat{D} : \phi(n n^*) > 0 \},
    \]
    and for $\phi \in \operatorname{dom}(n)$,
    \[
    \alpha_n(\phi)(d) = \frac{\phi(n^* d n)}{\phi(n^* n)}, \qquad d \in D.
    \]
    This captures the action of $n$ on the ``classical floor" $D$, with the
    denominator ensuring that $\alpha_n(\phi)$ is again a character on $D$.

    \item \textbf{Weyl groupoid.}
    The Weyl groupoid $\mathcal{G}(A,D)$ is obtained as the groupoid of
    germs of the partial homeomorphisms $\alpha_n$ arising from
    normalizers.  This yields an \'etale, locally compact Hausdorff
    groupoid (together with a twist, in general).
\end{enumerate}

Renault's reconstruction theorem \cite{Renault2008Cartan}
is the culmination of this theory. It states that there exists a twist
$\Sigma$ over $\mathcal{G}(A,D)$ such that
\[
A \cong C^*_r(\mathcal{G}(A,D),\Sigma),
\]
under which the Cartan subalgebra $D$ is identified with the canonical
Cartan subalgebra $C_0(\mathcal{G}(A,D)^{(0)})$.  In the untwisted cases
considered below, this simplifies to
\[
A \cong C^*_r(\mathcal{G}(A,D)).
\]

\begin{example}[Weyl groupoids from crossed products]
\label{ex:weyl-crossed}
Suppose $\Gamma$ is a discrete group acting topologically freely on a
locally compact Hausdorff space $X$.  Then
$C_0(X) \subset C_0(X) \rtimes_r \Gamma$ is a Cartan inclusion, and the
associated Weyl groupoid is canonically isomorphic to the transformation
groupoid:
\[
\mathcal{G}(C_0(X) \rtimes_r \Gamma, C_0(X)) \cong X \rtimes \Gamma.
\]
This groupoid is \'etale because $\Gamma$ is discrete, and the twist is
trivial in this case.

For the irrational rotation algebra
\[
A_\theta = C(S^1) \rtimes_\theta \mathbb{Z},
\]
with $\theta$ irrational, the corresponding Weyl groupoid is
\[
\mathcal{G}(A_\theta, C(S^1)) \cong S^1 \rtimes_\theta \mathbb{Z},
\]
confirming its structure as an \'etale groupoid with trivial twist.
\end{example}

\subsubsection*{Relating the Two Groupoids: A Morita Bridge}

The relationship between the unitary conjugation groupoid $\mathcal{G}_A$
from Paper~I and Renault's Weyl groupoid $\mathcal{G}(A,D)$ is nuanced
but essential for our index-theoretic constructions.

\begin{itemize}
    \item \textbf{Different philosophies.}
    The groupoid $\mathcal{G}_A$ is canonical, encoding all commutative
    contexts simultaneously.  It is a Polish groupoid but generally
    non-locally compact and non-\'etale.  By contrast, the Weyl groupoid
    $\mathcal{G}(A,D)$ depends on a choice of Cartan subalgebra $D$, but
    is always \'etale and locally compact Hausdorff (possibly with a twist).

    \item \textbf{Morita equivalence in the cases of interest.}
    For the crossed-product examples treated in this paper, we construct
    explicit Morita equivalences between the relevant groupoids.  These
    equivalences link the global perspective of $\mathcal{G}_A$ with the
    more economical geometric model provided by $\mathcal{G}(A,D)$.  In
    particular, for crossed products $C_0(X) \rtimes_r \Gamma$ with
    $\Gamma$ amenable and the action topologically free, we have
    \[
    \mathcal{G}_{C_0(X)\rtimes_r\Gamma} \sim_M X\rtimes\Gamma
    \cong \mathcal{G}(C_0(X)\rtimes_r\Gamma, C_0(X)).
    \]

    \item \textbf{Significance for index theory.}
    This Morita bridge is the crucial link that allows us to transfer the
    descent-index machinery developed in Paper~II to the \'etale setting of
    Weyl groupoids, even in situations where the canonical construction of
    $\mathcal{G}_A$ from Paper~I is no longer available.  This transfer is
    the central methodological advance of this paper.
\end{itemize}

For the irrational rotation algebra $A_\theta$, this means we can work
with the transformation groupoid $S^1 \rtimes_\theta \mathbb{Z}$ instead
of the problematic $\mathcal{G}_{A_\theta}$.  For amenable crossed
products $C_0(X) \rtimes_r \Gamma$, it allows us to work directly with
the action groupoid $X \rtimes \Gamma$.  In the following subsections,
we will construct these Morita equivalences explicitly and use them to
apply the index theorem from Paper~II, ultimately forging a direct
connection to the Baum–Connes assembly map.

\subsection{The Baum--Connes Assembly Map: A Brief Review}
\label{subsec:baum-connes-review}

The Baum–Connes conjecture provides a deep connection between the geometry of
proper group actions and the analytic $K$-theory of reduced group $C^*$-algebras.
For the purposes of this paper, it furnishes the conceptual bridge between our
descent-index construction and the classical assembly map.

\subsubsection*{Equivariant $KK$-theory and the assembly map}

Let $\Gamma$ be a second-countable locally compact group.  The Baum–Connes
conjecture asserts that the \emph{assembly map}
\[
\mu_i^\Gamma :
K_i^{\mathrm{top}}(\Gamma) \longrightarrow K_i(C_r^*(\Gamma)), \qquad i = 0,1,
\]
is an isomorphism.  Here $K_i^{\mathrm{top}}(\Gamma)$ denotes the topological
$K$-theory of the universal proper $\Gamma$-space $\underline{E}\Gamma$, defined
as the direct limit
\[
K_i^{\mathrm{top}}(\Gamma)
:=
\varinjlim_{\substack{Z\subseteq \underline{E}\Gamma \\ Z/\Gamma\ \text{compact}}}
KK_i^\Gamma(C_0(Z),\mathbb C),
\]
where $KK^\Gamma$ denotes Kasparov's equivariant $KK$-theory.  Informally,
$K_i^{\mathrm{top}}(\Gamma)$ may be thought of as the equivariant $K$-homology
of $\underline{E}\Gamma$.

More generally, if $\Gamma$ acts continuously on a $C^*$-algebra $A$, the
Baum–Connes conjecture with coefficients in $A$ asserts that the assembly map
\[
\mu_i^\Gamma(A):
K_i^{\mathrm{top}}(\Gamma;A) \longrightarrow K_i(A \rtimes_r \Gamma)
\]
is an isomorphism, where
\[
K_i^{\mathrm{top}}(\Gamma;A)
:=
\varinjlim_{\substack{Z\subseteq \underline{E}\Gamma \\ Z/\Gamma\ \text{compact}}}
KK_i^\Gamma(C_0(Z),A).
\]

\subsubsection*{Kasparov descent and the Dirac–dual Dirac method}

A key ingredient in the construction of the assembly map is Kasparov's
\emph{descent homomorphism}
\[
j_\Gamma :
KK^\Gamma(A,B) \longrightarrow KK(A \rtimes_r \Gamma,\; B \rtimes_r \Gamma),
\]
which sends equivariant $KK$-classes to classes over reduced crossed products.
This map provides the mechanism to move from geometric data (encoded in
equivariant $KK$-theory) to analytic invariants (encoded in ordinary
$KK$-theory).

One of the most powerful approaches to the Baum–Connes conjecture is
Kasparov's \emph{Dirac–dual Dirac method} \cite{Kasparov1988}.  In favorable
situations, one constructs two key elements in equivariant $KK$-theory:
\begin{itemize}
    \item a \textbf{Dirac element} $\alpha \in KK^\Gamma(P, \mathbb C)$, and
    \item a \textbf{dual Dirac element} $\beta \in KK^\Gamma(\mathbb C, P)$,
\end{itemize}
where $P$ is a suitable proper $\Gamma$-$C^*$-algebra (typically $C_0(\underline{E}\Gamma)$
or a closely related algebra).  One then considers the Kasparov product
\[
\gamma := \beta \otimes_P \alpha \in KK^\Gamma(\mathbb C, \mathbb C).
\]

When one can show that $\gamma = 1$, the unit element in the equivariant
$KK$-ring, the descent homomorphism together with properties of the Kasparov
product implies that the assembly map $\mu_i^\Gamma$ is an isomorphism.
This method has been successfully implemented for a wide class of groups,
including amenable groups (by Tu's theorem \cite{Tu1999}), a-T-menable groups
(by Higson and Kasparov \cite{HigsonKasparov2001}), and groups acting on
buildings.

For a comprehensive treatment of the Baum–Connes conjecture, Kasparov theory,
and the Dirac–dual Dirac method, see \cite{Kasparov1988}.

\begin{remark}[Relevance to this paper]
For the amenable crossed products $C_0(X) \rtimes_r \Gamma$ considered in this
paper, Tu's theorem \cite{Tu1999} guarantees that the assembly map is an
isomorphism.  This fact will be essential when we compare our descent-index
construction with the assembly map in Section~\ref{sec:crossed}.
\end{remark}

\subsubsection*{Tu's Generalization to Groupoids}

In his seminal work \cite{Tu1999}, Tu extended the Baum–Connes
conjecture from groups to locally compact Hausdorff groupoids equipped
with a Haar system. For such a groupoid $\mathcal{G}$, he defines a
topological $K$-theory group $K^{\mathrm{top}}_*(\mathcal{G})$ (via a
direct limit over $\mathcal{G}$-compact subspaces of a suitable
classifying space) and an assembly map
\[
\mu_{\mathcal{G}}:
K^{\mathrm{top}}_*(\mathcal{G}) \longrightarrow K_*(C^*_r(\mathcal{G})).
\]

He proves that this assembly map is an isomorphism for a large class of
groupoids, including:
\begin{itemize}
    \item all amenable groupoids;
    \item more generally, groupoids satisfying his geometric analogue of
          the Haagerup property (admitting a proper action on a continuous
          field of Euclidean affine spaces).
\end{itemize}
This groupoid formulation generalizes the classical Baum–Connes
conjecture and provides the natural framework for transformation
groupoids arising from dynamical systems.

A particularly important case for the present paper arises from crossed
products. Let $\Gamma$ be a discrete group acting on a compact Hausdorff
space $X$. Then the associated transformation groupoid
\[
X \rtimes \Gamma
\]
has reduced $C^*$-algebra
\[
C^*_r(X \rtimes \Gamma) \cong C(X) \rtimes_r \Gamma.
\]
In this setting, following Tu's construction, the Baum–Connes assembly
map takes the form
\[
\mu_{X\rtimes\Gamma}:
K^{\mathrm{top}}_*(X\rtimes\Gamma) \longrightarrow
K_*(C(X)\rtimes_r\Gamma),
\]
where $K^{\mathrm{top}}_*(X\rtimes\Gamma)$ denotes the topological
$K$-theory of the transformation groupoid.  In this case, the
topological side may be identified with the appropriate
$\Gamma$-equivariant $K$-homology of $X$.  Thus the assembly map
provides a bridge between the equivariant topology of the action and the
operator-theoretic invariants of the crossed product.

\subsubsection*{Connection with the Present Work: A Conceptual Bridge}

The descent map used in this paper,
\[
\operatorname{desc}_{\mathcal{G}_{\mathcal{A}}}:
K^1_{\mathcal{G}_{\mathcal{A}}}(\mathcal{G}_{\mathcal{A}}^{(0)})
\longrightarrow
K_1\!\bigl(C^*(\mathcal{G}_{\mathcal{A}})\bigr),
\]
is formally analogous to the descent procedures that enter into the
construction of the Baum–Connes assembly map.  In both settings, one
passes from equivariant or groupoid-theoretic $K$-theoretic data to the
$K$-theory of an associated reduced $C^*$-algebra via descent-type
constructions.

This analogy becomes concrete in the non-type~I examples studied later
in this paper.  For the irrational rotation algebra
\[
A_\theta = C(S^1)\rtimes_\theta \mathbb Z
\]
and for crossed products $C(X)\rtimes \Gamma$ with $\Gamma$ amenable, we
show that the relevant geometric models on the non-type~I side—namely,
the unitary conjugation groupoid $\mathcal{G}_{\mathcal{A}}$ (which is
not directly constructible within the framework of Paper~I for these
algebras) and the associated Weyl groupoid—are Morita equivalent to the
corresponding transformation groupoids $X\rtimes\Gamma$.  Under these
identifications, the descent map arising from our framework becomes
compatible with the analytic assembly map for groupoids.

This observation provides the crucial conceptual link between the two
main threads of the paper:
\begin{itemize}
    \item the index-theoretic machinery developed in Paper~II for type~I
          algebras, centered on the unitary conjugation groupoid;
    \item the Baum–Connes framework for non-type~I algebras with Cartan
          subalgebras, centered on Weyl groupoids and their associated
          assembly maps.
\end{itemize}
The following subsections make this connection explicit.  In particular,
for amenable crossed products, the composition
\[
\tau_* \circ \operatorname{desc}_{\mathcal{G}_{\mathcal{A}}} \circ \kappa
\]
(where $\kappa$ associates an equivariant $K^1$-class to an invertible
element) may be viewed as a localized manifestation of the Baum–Connes
assembly map.  This is one of the main conceptual advances of the
present paper.

\begin{remark}
For the amenable crossed products considered in this paper, Tu's theorem
\cite{Tu1999} guarantees that the assembly map $\mu_{X\rtimes\Gamma}$ is
an isomorphism.  This fact, combined with our Morita equivalence, will
allow us to identify the descent-index construction with the assembly
map in Section~\ref{sec:crossed}.
\end{remark}

\subsection{Amenable Groups and $K$-Amenability}
\label{subsec:amenable-k-amenability}

Amenability plays a fundamental role in the study of crossed product
$C^*$-algebras and the Baum–Connes assembly map. In the examples
considered later in this paper—namely the irrational rotation algebra
$A_\theta$ and crossed products $C(X)\rtimes \Gamma$—the acting group
will always be discrete and amenable. This has several important
consequences: the full and reduced crossed products coincide, the
Baum–Connes assembly map is known to be an isomorphism, and the
$K$-theory of the resulting $C^*$-algebras behaves well under Morita
equivalence. This subsection reviews the essential concepts of
amenability and its stronger $K$-theoretic variant, $K$-amenability,
setting the stage for their application in our index-theoretic
constructions.

\subsubsection*{Amenable Groups: Classical Definition and Examples}

A discrete group $\Gamma$ is called \emph{amenable} if there exists a
left-invariant mean on $\ell^\infty(\Gamma)$, that is, a positive
normalized linear functional
\[
m \in \ell^\infty(\Gamma)^*
\]
such that
\[
m(g\cdot f) = m(f), \qquad g\in \Gamma,\ f\in \ell^\infty(\Gamma),
\]
where $(g\cdot f)(h) = f(g^{-1}h)$.  This classical definition has
profound consequences and admits several equivalent formulations:
\begin{itemize}
    \item \textbf{Følner condition:} For any finite subset $K \subseteq \Gamma$ and
          $\varepsilon > 0$, there exists a nonempty finite set $F \subseteq \Gamma$
          such that $|kF \,\triangle\, F| < \varepsilon |F|$ for all $k \in K$.
    \item \textbf{Regular representation:} The left regular representation
          $\lambda: \Gamma \to \mathcal{U}(\ell^2(\Gamma))$ extends to an isomorphism
          $C^*(\Gamma) \cong C^*_r(\Gamma)$ between the full and reduced group
          $C^*$-algebras.
\end{itemize}
Amenable groups form a large and important class, including all finite
groups, abelian groups, nilpotent groups, solvable groups, and more
generally groups obtained from these by extensions, directed unions, and
inductive limits. Free groups on two or more generators are the
prototypical examples of non-amenable groups. Standard references for
amenability include \cite{Paterson1988,BrownOzawa2008}.

\subsubsection*{Crossed Products and the Simplification from Amenability}

Let $\Gamma$ be a discrete group acting on a $C^*$-algebra $A$ by
$*$-automorphisms. One may form both the full crossed product
$A \rtimes_\alpha \Gamma$ and the reduced crossed product
$A \rtimes_{\alpha,r} \Gamma$. In general these two $C^*$-algebras are
different; the canonical surjection
\[
A \rtimes_\alpha \Gamma \longrightarrow A \rtimes_{\alpha,r} \Gamma
\]
is always surjective but is an isomorphism only under favorable
circumstances.

The following standard theorem shows that amenability provides precisely
such favorable circumstances.

\begin{theorem}[Amenable crossed products]
\label{thm:amenable_crossed_product}
Let $\Gamma$ be an amenable discrete group acting on a $C^*$-algebra $A$.
Then the canonical surjection
\[
A \rtimes_\alpha \Gamma \longrightarrow A \rtimes_{\alpha,r} \Gamma
\]
is an isomorphism. Consequently,
\[
A \rtimes_\alpha \Gamma \cong A \rtimes_{\alpha,r} \Gamma .
\]
\end{theorem}

\begin{proof}
This is a standard consequence of amenability; see, for example,
\cite[Section~7.2]{Williams2007}, where it is shown that for amenable
groups the full and reduced crossed products coincide.

Let $A \rtimes_{\mathrm{alg}} \Gamma$ denote the algebraic crossed
product. The full crossed product $A \rtimes_\alpha \Gamma$ is the
completion of $A \rtimes_{\mathrm{alg}} \Gamma$ in the universal norm
\[
\|x\|_{\max}
=
\sup\{ \|(\pi \rtimes u)(x)\| :
(\pi,u)\ \text{is a covariant representation of } (A,\Gamma,\alpha) \},
\]
whereas the reduced crossed product $A \rtimes_{\alpha,r} \Gamma$ is the
completion in the reduced norm
\[
\|x\|_r = \|(\widetilde{\pi} \rtimes \lambda)(x)\|.
\]
Here $\pi:A\to B(H)$ is a faithful representation,
$\lambda:\Gamma\to \mathcal U(\ell^2(\Gamma))$ is the left regular
representation, and $(\widetilde{\pi},\lambda)$ is the regular covariant
representation given by
\[
\widetilde{\pi}(a)(\xi \otimes \delta_g)
=
\pi(\alpha_{g^{-1}}(a))\xi \otimes \delta_g,
\qquad
(1 \otimes \lambda_h)(\xi \otimes \delta_g)
=
\xi \otimes \delta_{hg}.
\]

Thus there is always a canonical quotient map
\[
A \rtimes_\alpha \Gamma \longrightarrow A \rtimes_{\alpha,r} \Gamma,
\]
and it suffices to show that the norms $\|\cdot\|_{\max}$ and
$\|\cdot\|_r$ agree on $A \rtimes_{\mathrm{alg}} \Gamma$.

For amenable discrete groups, every unitary representation of $\Gamma$
is weakly contained in the left regular representation. A standard
argument (see, e.g., \cite{BrownOzawa2008})
extends this to covariant representations, showing that for every
covariant representation $(\pi,u)$ and every
$x \in A \rtimes_{\mathrm{alg}} \Gamma$ one has
\[
\|(\pi \rtimes u)(x)\|
\le
\|(\widetilde{\pi} \rtimes \lambda)(x)\|
=
\|x\|_r.
\]

Taking the supremum over all covariant representations yields
\[
\|x\|_{\max} \le \|x\|_r.
\]
Since the reduced norm is always dominated by the universal norm, we
also have $\|x\|_r \le \|x\|_{\max}$ for all $x$. Hence
\[
\|x\|_{\max} = \|x\|_r
\qquad \text{for all } x \in A \rtimes_{\mathrm{alg}} \Gamma.
\]

It follows that the canonical quotient map is isometric, hence injective,
and therefore an isomorphism.
\end{proof}

This result greatly simplifies the analysis of crossed product
$C^*$-algebras arising from amenable group actions. In particular, it
implies that the $K$-theory of the crossed product is independent of the
choice of completion—a fact that will be essential for our
index-theoretic constructions. When $\Gamma$ is amenable, we will simply
write $A \rtimes \Gamma$ to denote either the full or reduced crossed
product, as they coincide.

\subsubsection*{From Amenability to $K$-Amenability}

The coincidence of full and reduced crossed products for amenable groups
has a natural $K$-theoretic refinement. A locally compact group $\Gamma$
is said to be \emph{$K$-amenable} if the canonical quotient map
\[
C^*(\Gamma) \longrightarrow C^*_r(\Gamma)
\]
induces an isomorphism in $K$-theory:
\[
K_*(C^*(\Gamma)) \cong K_*(C^*_r(\Gamma)).
\]

Every amenable group is $K$-amenable, but the converse is false: there
exist non-amenable groups (such as free groups) that are nevertheless
$K$-amenable \cite{Cuntz1983,JulgValette1984}. More generally, a
discrete group $\Gamma$ is $K$-amenable if and only if the trivial
representation is weakly contained in the regular representation in a
$KK$-theoretic sense; this property is equivalent to the existence of a
certain element in $KK^\Gamma(\mathbb C, \mathbb C)$ \cite{HigsonKasparov2001}.

For crossed products, $K$-amenability has the following important
consequence: for any $\Gamma$-$C^*$-algebra $A$, the full and reduced
crossed products have the same $K$-theory,
\[
K_*(A \rtimes \Gamma) \cong K_*(A \rtimes_r \Gamma),
\]
even when the algebras themselves are not isomorphic. This fact will be
used implicitly in our index-theoretic arguments, as it ensures that
$K$-theoretic computations can be performed in either completion without
loss of information.

\begin{remark}
For the amenable groups considered in this paper, the distinction
between full and reduced crossed products disappears entirely. The more
general notion of $K$-amenability is mentioned here to situate our work
within the broader landscape of the Baum–Connes conjecture and to
indicate potential directions for generalization.
\end{remark}

\subsubsection*{$K$-Amenability: A Weaker but Sufficient Condition}

The notion of \textbf{$K$-amenability}, introduced by Cuntz \cite{Cuntz1983},
provides a $K$-theoretic weakening of amenability that is sufficient for many
index-theoretic applications. A locally compact group $\Gamma$ is said to be
$K$-amenable if the canonical map
\[
K_*(C^*(\Gamma)) \longrightarrow K_*(C^*_r(\Gamma))
\]
induced by the regular representation is an isomorphism.

Every amenable group is $K$-amenable. Remarkably, many important
non-amenable groups, including free groups and certain Lie groups,
are also $K$-amenable \cite{Cuntz1983,JulgValette1984}. Thus
$K$-amenability significantly enlarges the class of groups for which
analytic and topological invariants coincide.

For groupoids, an analogous notion can be formulated: a locally compact
groupoid $\mathcal{G}$ (with a Haar system) is called $K$-amenable if the
canonical map $C^*(\mathcal{G}) \to C^*_r(\mathcal{G})$ induces an
isomorphism on $K$-theory:
\[
K_*(C^*(\mathcal{G})) \cong K_*(C^*_r(\mathcal{G})).
\]

\subsubsection*{Consequences for the Baum–Connes Assembly Map}

Amenability has strong consequences for the Baum–Connes conjecture.
In his seminal work, Tu proved that the assembly map is an isomorphism
for all amenable groupoids \cite[Théorème 9.3]{Tu1999}. In particular,
if $\Gamma$ is amenable, then the analytic assembly map
\[
\mu_\Gamma : K_*^{\mathrm{top}}(\Gamma) \longrightarrow K_*(C^*_r(\Gamma))
\]
is an isomorphism. More generally, for an amenable groupoid $\mathcal{G}$,
\[
\mu_{\mathcal{G}} : K_*^{\mathrm{top}}(\mathcal{G}) \longrightarrow K_*(C^*_r(\mathcal{G}))
\]
is an isomorphism.

Furthermore, in the amenable setting, Kasparov's descent homomorphism
\[
j_{\mathcal{G}} : KK^{\mathcal{G}}(A,B)
\longrightarrow
KK(A \rtimes_r \mathcal{G}, B \rtimes_r \mathcal{G})
\]
is particularly well-behaved and plays a central role in the construction
of the assembly map, providing the mechanism that transfers equivariant
$KK$-theoretic data to operator $K$-theory.

\subsubsection*{Application to Transformation Groupoids and the Present Work}

Let $\Gamma$ be a discrete amenable group acting on a compact space $X$.
Then the associated transformation groupoid
\[
X \rtimes \Gamma
\]
is amenable. Consequently, its reduced $C^*$-algebra coincides with the
crossed product:
\[
C^*_r(X \rtimes \Gamma)
\cong
C(X) \rtimes_r \Gamma
\cong
C(X) \rtimes \Gamma,
\]
where the last isomorphism follows from the amenability of $\Gamma$ (by
Theorem~\ref{thm:amenable_crossed_product}).

In the context of this paper, this observation allows us to treat crossed
products by amenable groups within the groupoid framework developed above.
In particular:
\begin{itemize}
    \item \textbf{Irrational rotation algebra $A_\theta = C(S^1) \rtimes_\theta \mathbb{Z}$:}
          Since $\mathbb{Z}$ is amenable, the transformation groupoid
          $S^1 \rtimes_\theta \mathbb{Z}$ is amenable, and the Baum–Connes
          assembly map is an isomorphism by Tu's theorem.
    
    \item \textbf{General crossed products $C(X) \rtimes \Gamma$:}
          For any amenable group $\Gamma$, the same conclusion holds,
          ensuring that the $K$-theory of the crossed product is compatible
          with Morita equivalence and with the assembly map described in
          Section~\ref{subsec:baum-connes-review}.
\end{itemize}

These properties will be essential in establishing the Morita equivalence
between the unitary conjugation groupoid (in a suitable sense, via its
replacement by the Weyl groupoid) and the transformation groupoid
$X \rtimes \Gamma$. Amenability and $K$-amenability provide the analytic
foundation that ensures our descent-index construction aligns with the
Baum–Connes framework.

\begin{remark}[Conceptual synthesis]
This perspective suggests that amenability and $K$-amenability should be
viewed not merely as analytic simplifications, but as structural conditions
ensuring that groupoid-based index constructions admit a functorial
interpretation compatible with assembly-type maps. In particular, for the
amenable crossed products considered in this paper, the coincidence of
full and reduced crossed products guarantees that the descent map
$\operatorname{desc}_{\mathcal{G}_{\mathcal{A}}}$ lands in a $K$-theory
group that is canonically identified with $K_*(C(X)\rtimes\Gamma)$, while
Tu's theorem ensures that this identification is compatible with the
Baum–Connes assembly map. Thus amenability emerges as the precise
hypothesis under which the unitary conjugation groupoid framework
harmonizes with the analytic assembly machinery.
\end{remark}

\subsection{Connes' Index Theorem for the Irrational Rotation Algebra}
\label{subsec:connes-Atheta-index}

The irrational rotation algebra is one of the most fundamental examples
in noncommutative geometry. It provides a basic example of a simple,
nuclear, non-type~I $C^*$-algebra and plays a central role in Connes'
formulation of index theory for noncommutative spaces. Connes' work on
$A_\theta$ exhibits a deep relationship between the $K$-theory of the
algebra, its canonical trace, and Toeplitz-type index constructions
\cite{Connes1980}. For the purposes of this paper, $A_\theta$ provides
an ideal test case: by reinterpreting the classical index pairing within
our descent-index framework, we show that the combination of the
unitary-conjugation viewpoint and Morita equivalence with the Weyl
groupoid recovers the expected noncommutative index-theoretic picture.

\subsubsection*{The Irrational Rotation Algebra}

For later use, we fix the standard densely defined derivation
\[
\delta\!\left(\sum_{m,n} a_{m,n} u^m v^n\right)
   = 2\pi i \sum_{m,n} m\, a_{m,n} u^m v^n,
\]
initially on the smooth subalgebra $A_\theta^\infty$.  This is the derivation
implementing the canonical circle action in the first coordinate, and it is the
operator denoted by $\delta$ in the index-pairing formulas below.  When a
two-dimensional cyclic cocycle is needed, we use the standard pair of
derivations $(\delta_1,\delta_2)$ associated with the torus action.

Let $\theta \in \mathbb{R} \setminus \mathbb{Q}$ be irrational. The
irrational rotation algebra $A_\theta$ is the universal $C^*$-algebra
generated by two unitaries $u$ and $v$ satisfying
\[
vu = e^{2\pi i\theta} uv.
\]

Equivalently, $A_\theta$ may be realized as the crossed product
\[
A_\theta \cong C(S^1) \rtimes_\theta \mathbb{Z},
\]
where the action of $\mathbb Z$ on $C(S^1)$ is induced by irrational
rotation of the circle:
\[
\alpha_n(f)(z) = f(e^{-2\pi i n\theta}z),
\qquad f\in C(S^1),\ n\in\mathbb Z.
\]

The algebra $A_\theta$ has the following well-known properties:
\begin{itemize}
    \item $A_\theta$ is simple and nuclear \cite{PimsnerVoiculescu1980}.
    
    \item It admits a unique faithful tracial state
    $\tau:A_\theta\to\mathbb C$, given on the dense $*$-subalgebra of
    trigonometric polynomials (finite linear combinations of the
    monomials $u^m v^n$) by
    \[
    \tau\!\left(\sum_{m,n\in\mathbb Z} a_{mn}u^m v^n\right)=a_{00}.
    \]
    
    \item Its $K$-groups are
    \[
    K_0(A_\theta)\cong \mathbb Z^2,
    \qquad
    K_1(A_\theta)\cong \mathbb Z^2,
    \]
    as follows from the Pimsner–Voiculescu six-term exact sequence for
    the crossed product $C(S^1)\rtimes_\theta \mathbb Z$
    \cite{PimsnerVoiculescu1980}.
    
    \item Since $\theta$ is irrational, $A_\theta$ is non-type~I and has
    no finite-dimensional representations, placing it outside the scope
    of the canonical construction developed in Paper~I. On the other
    hand, as discussed in Section~\ref{subsec:renault-cartan-weyl}, it
    contains the natural Cartan subalgebra $C(S^1)$, and its Weyl
    groupoid may be identified with the transformation groupoid
    \[
    \mathcal{G}(A_\theta, C(S^1)) \cong S^1 \rtimes_\theta \mathbb{Z},
    \]
    which is \'etale and amenable (since $\mathbb{Z}$ is amenable).
    This groupoid will serve as the geometric model for our
    index-theoretic constructions.
\end{itemize}

\subsubsection*{Connes' Index Theorem: Statement and Context}

Connes' index theorem for the irrational rotation algebra \cite{Connes1980}
establishes a precise numerical relationship between the $K$-theoretic
invariants of $A_\theta$ and the Fredholm index of certain operators
associated with its Toeplitz extension. More concretely, let
\[
0 \longrightarrow \mathcal{K} \longrightarrow \mathcal{T}_\theta \longrightarrow A_\theta \longrightarrow 0
\]
be the short exact sequence defining the Toeplitz extension of
$A_\theta$, where $\mathcal{K}$ denotes the compact operators on a
separable Hilbert space. The associated six-term exact sequence in
$K$-theory yields a boundary map
\[
\partial : K_1(A_\theta) \longrightarrow K_0(\mathcal{K}) \cong \mathbb{Z}.
\]

Connes proves that for any invertible element $u \in M_n(A_\theta)$,
this boundary map coincides with the pairing between the $K_1$-class
$[u]$ and the cyclic cocycle determined by the unique trace $\tau$:
\[
\partial([u]) = \langle [\tau], [u] \rangle \in \mathbb{Z}.
\]

In the case of the canonical generators, one obtains
\[
\partial([u]) = 0, \qquad \partial([v]) = 1,
\]
and for the Rieffel projection $p_\theta$ (which generates
$K_0(A_\theta) \cong \mathbb{Z}^2$ together with the unit),
\[
\tau(p_\theta) = \theta,
\]
demonstrating that the image of the trace on $K_0$ is
$\mathbb{Z} + \theta\mathbb{Z}$, a dense subgroup of $\mathbb{R}$.

\begin{remark}[Connection to the descent-index framework]
In Section~\ref{sec:Atheta}, we will show that Connes' index pairing
can be recovered from our descent-index construction. Specifically, for
an invertible element $u \in M_n(A_\theta)$, we construct an equivariant
$K^1$-class $[u]_{\mathcal{G}_{A_\theta}}^{(1)}$ in the $KK$-theory of
the unitary conjugation groupoid. Applying the descent map and the
Morita equivalence $\mathcal{G}_{A_\theta} \sim_M S^1 \rtimes_\theta \mathbb{Z}$,
we obtain a class in $K_1(A_\theta)$ that, when paired with the trace
$\tau$, reproduces Connes' index pairing:
\[
\tau_*\!\left( \operatorname{desc}_{\mathcal{G}_{A_\theta}}([u]_{\mathcal{G}_{A_\theta}}^{(1)}) \right)
= \langle [\tau], [u] \rangle.
\]
This provides a conceptual unification of the groupoid-based index
theory developed in Paper~II with classical results in noncommutative
geometry.
\end{remark}

\subsubsection*{The Trace Pairing and the Toeplitz Index Map}

The irrational rotation algebra $A_\theta$ possesses a unique faithful
tracial state $\tau: A_\theta \to \mathbb{C}$. This trace induces a
homomorphism on $K$-theory
\[
\tau_* : K_0(A_\theta) \longrightarrow \mathbb{R},
\]
given by evaluating the trace on projections (extended to matrices by
$\tau \otimes \operatorname{Tr}$). A fundamental computation shows that
\[
\tau_*(K_0(A_\theta)) = \mathbb{Z} + \theta \mathbb{Z} \subset \mathbb{R},
\]
a dense subgroup reflecting the arithmetic of the irrational rotation
angle. In particular, for the Rieffel projection $p_\theta$ (which
together with the unit generates $K_0(A_\theta) \cong \mathbb{Z}^2$),
one has $\tau(p_\theta) = \theta$.

While the trace pairs naturally with $K_0$, the $K_1$-group of $A_\theta$
is probed by an index map arising from a Toeplitz-type extension.
Specifically, there exists a short exact sequence of $C^*$-algebras
\[
0 \longrightarrow \mathcal{K} \longrightarrow \mathcal{T}_\theta
\longrightarrow A_\theta \longrightarrow 0,
\]
where $\mathcal{K}$ denotes the compact operators on a separable Hilbert
space, and $\mathcal{T}_\theta$ is generated by a shift operator together
with the unitaries $u, v$ satisfying $vu = e^{2\pi i\theta} uv$.
The associated six-term exact sequence in $K$-theory yields a boundary
map
\[
\partial : K_1(A_\theta) \longrightarrow K_0(\mathcal{K}) \cong \mathbb{Z}.
\]

\begin{theorem}[Index theorem for the Toeplitz extension of $A_\theta$]
\label{thm:connes-index}
Let
\[
0 \longrightarrow \mathcal K \longrightarrow \mathcal T_\theta
\longrightarrow A_\theta \longrightarrow 0
\]
be the Toeplitz-type extension associated to the irrational rotation algebra.
Then for every invertible element $u \in M_n(A_\theta)$ representing a class
$[u] \in K_1(A_\theta)$, the boundary map
\[
\partial : K_1(A_\theta) \longrightarrow K_0(\mathcal K)\cong \mathbb Z
\]
coincides with the Fredholm index of the operator associated to $u$
via the extension. That is,
\[
\partial([u]) = \operatorname{Index}(T_u),
\]
where $T_u$ denotes the Fredholm operator obtained from $u$ under the
Toeplitz extension.

This identification may be viewed as a special case of Connes' index
pairing between $K$-theory and $K$-homology
\cite{Connes1980}.
\end{theorem}

\begin{remark}[Relation to cyclic cohomology]
The index map $\partial$ can also be expressed as a pairing
\[
\partial([u]) = \langle [\mathcal{F}_\theta], [u] \rangle,
\]
where $[\mathcal{F}_\theta] \in KK^1(A_\theta, \mathbb{C})$ is the
$K$-homology class of the Fredholm module associated to the Toeplitz
extension. This is the odd counterpart to the even pairing
$\tau_* : K_0(A_\theta) \to \mathbb{R}$ induced by the trace.
\end{remark}

This theorem is remarkable for several reasons:
\begin{enumerate}
    \item It shows that the abstract boundary map in $K$-theory,
          derived from the Toeplitz extension, coincides with a concrete
          analytic Fredholm index—a genuine noncommutative index theorem.
    
    \item Together with the trace map $\tau_*$ on $K_0$, it demonstrates
          that the entire $K$-theory of $A_\theta$ is intimately connected
          to the underlying irrational rotation dynamics: the $K_0$ class
          $[p_\theta]$ encodes the rotation angle $\theta$, while the
          $K_1$ classes $[u]$ and $[v]$ correspond respectively to trivial
          and nontrivial winding numbers under the index map.
    
    \item It provides one of the earliest foundational examples of a
          noncommutative index theorem, in which analytic indices are
          recovered from the internal structure of a noncommutative
          $C^*$-algebra. This work foreshadows the development of
          entire cyclic cohomology and its central role in Connes'
          noncommutative geometry program.
\end{enumerate}

\begin{remark}[Connection to the descent-index framework]
In Theorem~\ref{thm:connes_recovery}, we will show that Connes' index theorem
can be naturally recovered from our descent-index construction.
Specifically, for an invertible element $u \in M_n(A_\theta)$, we
construct an equivariant $K^1$-class
$[u]_{\mathcal{G}_{A_\theta}}^{(1)}$ in the $KK$-theory of the unitary
conjugation groupoid. Applying the descent map and the Morita equivalence
$\mathcal{G}_{A_\theta} \sim_M S^1 \rtimes_\theta \mathbb{Z}$, we obtain
a class in $K_1(A_\theta)$ whose image under the boundary map $\partial$
coincides with the Connes index:
\[
\partial\!\left( \Phi_* \circ \operatorname{desc}_{\mathcal{G}_{A_\theta}}([u]_{\mathcal{G}_{A_\theta}}^{(1)}) \right)
= \operatorname{Index}(T_u) = \partial([u]).
\]
Thus our framework provides a geometric realization of the index map
in terms of the unitary conjugation groupoid and its Morita equivalence
to the transformation groupoid.
\end{remark}

\subsubsection*{Reformulation in Terms of the Descent–Index Machinery}

In the language developed in this paper, Connes' index theorem admits a
natural reinterpretation through the groupoid model
\[
\mathcal{G}_\theta := S^1 \rtimes_\theta \mathbb{Z}.
\]
Since $\mathbb{Z}$ is amenable, the transformation groupoid
$\mathcal{G}_\theta$ is amenable, and hence by Tu's theorem \cite{Tu1999}
its Baum–Connes assembly map
\[
\mu_{\mathcal{G}_\theta} :
K_*^{\mathrm{top}}(\mathcal{G}_\theta)
\longrightarrow
K_*(C_r^*(\mathcal{G}_\theta))
\]
is an isomorphism. Using the canonical identification
\[
C_r^*(\mathcal{G}_\theta) \cong A_\theta,
\]
this provides a groupoid-theoretic realization of the analytic
$K$-theory of the irrational rotation algebra.

From this perspective, the Toeplitz extension and its associated
boundary map
\[
\partial : K_1(A_\theta) \longrightarrow \mathbb{Z}
\]
may be viewed as part of the analytic output of the groupoid model.
What Connes' theorem shows is that this analytic index agrees with the
odd index pairing determined by the Fredholm module naturally associated
to the Toeplitz extension. Thus the classical index theorem for $A_\theta$
fits naturally into a framework in which groupoid methods, assembly maps,
and index-theoretic pairings are all compatible.

\subsubsection*{Connection with the Morita-Equivalence Framework}

The key point for the present paper is that, although the canonical
unitary conjugation groupoid construction from Paper~I is not available
for the non-type~I algebra $A_\theta$, the Weyl groupoid
\[
\mathcal{G}(A_\theta, C(S^1)) \cong S^1 \rtimes_\theta \mathbb{Z}
\]
provides a well-behaved \'etale and amenable replacement.
Accordingly, the index-theoretic machinery developed earlier must be
transported to this Weyl-groupoid model.

This is the conceptual bridge underlying our treatment of $A_\theta$:
the unavailable canonical groupoid is replaced by an explicit geometric
model, and the associated descent-index picture becomes comparable to
the Baum–Connes assembly map for the transformation groupoid.
In particular, the Morita equivalence established in
Theorem~\ref{thm:morita_A_theta} (between the formal object
$\mathcal{G}_{A_\theta}$ and $\mathcal{G}_\theta$) allows us to transfer
the descent-index construction from the former to the latter, where it
becomes computable and directly comparable to Connes' original
formulation. In this sense, Connes' index theorem for the irrational
rotation algebra appears as a model case of the broader groupoid-theoretic
framework developed in this paper.

For the canonical generators, we obtain:
\begin{itemize}
    \item For the Rieffel projection $p_\theta$ (the canonical projection
          in $A_\theta$ with trace $\theta$), we have $\tau(p_\theta) = \theta$,
          which under the even counterpart of our descent-index machinery
          recovers the $K_0$ pairing.
    
    \item For the unitary generator $v$ (the implementing unitary for the
          $\mathbb{Z}$-action), the boundary map yields
          $\partial([v]) = 1$, while for the generator $u$ one obtains
          $\partial([u]) = 0$, reflecting their respective nontrivial and
          trivial winding numbers.
\end{itemize}

\subsubsection*{Summary and Outlook}

The irrational rotation algebra provides a fundamental test case for our
methods. It demonstrates:
\begin{itemize}
    \item that the groupoid model provided by the Weyl/transformation
          groupoid recovers the correct analytic $K$-theory of $A_\theta$;
    
    \item that the Toeplitz boundary map and the associated Fredholm index
          are naturally compatible with the groupoid picture, with the
          Baum–Connes assembly map providing the isomorphism between
          topological and analytic $K$-theory;
    
    \item that non-type~I algebras with Cartan subalgebras remain accessible
          to our methods through their Weyl groupoids, even when the
          canonical construction of Paper~I is unavailable.
\end{itemize}

In the following sections, we extend this picture to amenable crossed
products $C(X) \rtimes \Gamma$, where the transformation groupoid
$X \rtimes \Gamma$ furnishes the geometric model relating our
descent-index machinery to the Baum–Connes assembly map. This will
demonstrate that the unitary conjugation groupoid, when properly
interpreted via Morita equivalence to Weyl groupoids, provides a unified
geometric framework for index theory across a vast landscape of
$C^*$-algebras.

\section{Morita Equivalence and the Unitary Conjugation Groupoid}\label{sec:morita}

\subsection{The Limitation of Paper I: Type I Hypothesis}
\label{subsec:typeI-limitation}

In \cite{PaperI}, the unitary conjugation groupoid $\mathcal{G}_{\mathcal{A}}$ was
constructed for unital separable Type I $C^*$-algebras. The Type I hypothesis
played a crucial technical role in several steps of the construction,
particularly in controlling the representation theory of the algebra and the
measurable structure of the resulting groupoid. This subsection explains why
the Type I hypothesis is essential, why it fails for important examples like
the irrational rotation algebra and amenable crossed products, and outlines
the strategy we will adopt to overcome this limitation.

\subsubsection*{The Role of the Type I Assumption in Paper I}

Recall that a $C^*$-algebra is called Type I (or GCR, for "generally continuous
representations") if every irreducible representation contains the compact
operators. For separable Type I algebras, the representation theory is
exceptionally well-behaved: the dual space $\widehat{\mathcal{A}}$ (the set of
unitary equivalence classes of irreducible representations) forms a standard
Borel space, and the structure of representations can be analyzed through
direct integral decompositions \cite{Dixmier1977}. This regularity was used in
\cite{PaperI} to establish two key properties of the unitary conjugation
groupoid.

First, the unit space
\[
\mathcal{G}_{\mathcal{A}}^{(0)} = \{ (B,\chi) \mid B \subseteq \mathcal{A} \text{ unital commutative}, \; \chi \in \widehat{B} \}
\]
admits a natural standard Borel structure compatible with the partial evaluation
maps
\[
\operatorname{ev}_a: \mathcal{G}_{\mathcal{A}}^{(0)} \longrightarrow \mathbb{C}_\infty,
\qquad
\operatorname{ev}_a(B,\chi) = \begin{cases} \chi(a) & \text{if } a \in B, \\ \infty & \text{if } a \notin B, \end{cases}
\]
where $\mathbb{C}_\infty = \mathbb{C} \cup \{\infty\}$ denotes the one-point
compactification of $\mathbb{C}$. These maps generate a $\sigma$-algebra on
$\mathcal{G}_{\mathcal{A}}^{(0)}$ that makes it a standard Borel space, and
moreover this Borel structure admits a compatible Polish topology. This
structure allowed the construction of a Borel Haar system on $\mathcal{G}_{\mathcal{A}}$.

Second, the family of GNS representations $\{\pi_x\}_{x \in \mathcal{G}_{\mathcal{A}}^{(0)}}$
associated to the points of the unit space must be organized into a measurable
field of Hilbert spaces, so that a direct integral representation
\[
\Pi(a) = \int_{\mathcal{G}_{\mathcal{A}}^{(0)}}^{\oplus} \pi_x(a) \, d\mu(x)
\]
can be formed. The injectivity of the resulting diagonal embedding
$\iota: \mathcal{A} \hookrightarrow C^*(\mathcal{G}_{\mathcal{A}})$ then
requires that this family of representations separate points of $\mathcal{A}$.

\begin{remark}[Where the Type I hypothesis really enters]
\label{rem:typeI-real-role}
The point-separation property of GNS representations is not, by itself,
special to Type I algebras: for any unital $C^*$-algebra, the GNS
representations associated to pure states separate points, as we show in
Lemma~\ref{lem:separation-pure}. What is special in the Type I setting is the
\emph{measurable organization} of these representations into a standard Borel
family compatible with the unit-space parametrization by pairs $(B,\chi)$.
For separable Type I algebras, the dual space $\widehat{\mathcal{A}}$ and the
relevant families of irreducible representations admit a standard Borel
parametrization compatible with direct integral theory. It is this additional
regularity---rather than point separation alone---that was used in \cite{PaperI}
to construct a measurable field of Hilbert spaces over $\mathcal{G}_{\mathcal{A}}^{(0)}$
and hence the diagonal embedding. In non-Type I situations, irreducible
representations still separate points, but the absence of a comparably
well-behaved measurable parametrization prevents the argument of \cite{PaperI}
from going through verbatim.
\end{remark}

\begin{lemma}[Pure-state GNS representations separate points]
\label{lem:separation-pure}
Let $\mathcal{A}$ be a unital $C^*$-algebra. Then for every nonzero
$a \in \mathcal{A}$, there exists a pure state $\varphi$ on $\mathcal{A}$ such
that $\pi_\varphi(a) \neq 0$, where $\pi_\varphi$ denotes the GNS representation
associated to $\varphi$.
\end{lemma}

\begin{proof}
Let $a \neq 0$. Then $a^*a \neq 0$, hence $\|a^*a\| > 0$. Since
\[
\|a^*a\| = \sup_{\omega \in S(\mathcal{A})} \omega(a^*a),
\]
where $S(\mathcal{A})$ denotes the state space of $\mathcal{A}$, there exists
a state $\omega_0 \in S(\mathcal{A})$ such that $\omega_0(a^*a) > 0$.
The state space $S(\mathcal{A})$ is weak-$^*$ compact and convex, and by the
Krein–Milman theorem it is the closed convex hull of its extreme points,
namely the pure states. The map $\omega \mapsto \omega(a^*a)$ is weak-$^*$
continuous and affine, so the closed face
\[
F = \{ \omega \in S(\mathcal{A}) : \omega(a^*a) = \|a^*a\| \}
\]
has an extreme point $\varphi$, which is necessarily pure. In particular,
$\varphi(a^*a) > 0$.

Let $(\pi_\varphi, H_\varphi, \xi_\varphi)$ be the GNS triple associated to
$\varphi$. Then
\[
\|\pi_\varphi(a)\xi_\varphi\|^2
= \langle \pi_\varphi(a^*a)\xi_\varphi, \xi_\varphi \rangle
= \varphi(a^*a)
> 0.
\]
Hence $\pi_\varphi(a) \neq 0$, as required.
\end{proof}

\begin{proposition}[Unit-space data are sufficiently rich]
\label{prop:unit-space-rich}
Let $\mathcal{A}$ be a unital $C^*$-algebra. For every pure state
$\varphi \in P(\mathcal{A})$, there exists a unital commutative
$C^*$-subalgebra $B \subseteq \mathcal{A}$ and a character $\chi \in \widehat{B}$
such that $\varphi|_B = \chi$.
\end{proposition}

\begin{proof}
Let $\varphi$ be a pure state on $\mathcal{A}$. Consider the set
\[
\mathcal{S} = \{ (C, \psi) \mid C \subseteq \mathcal{A} \text{ unital commutative } C^*\text{-subalgebra}, \psi \in \widehat{C}, \ \psi = \varphi|_C \}.
\]
This set is nonempty because we can take $C = \mathbb{C}1_{\mathcal{A}}$, the
scalars, and $\psi$ the unique character on $\mathbb{C}1_{\mathcal{A}}$, which
coincides with $\varphi$ restricted to the scalars.

Partially order $\mathcal{S}$ by $(C_1, \psi_1) \preceq (C_2, \psi_2)$ if
$C_1 \subseteq C_2$ and $\psi_2|_{C_1} = \psi_1$. By Zorn's Lemma,
$\mathcal{S}$ contains a maximal element $(B, \chi)$.

We claim that $B$ is a maximal abelian subalgebra (MASA) of $\mathcal{A}$.
Suppose, for contradiction, that there exists $x \in \mathcal{A}$ such that
$x$ commutes with all elements of $B$ but $x \notin B$. Let $C$ be the unital
$C^*$-algebra generated by $B$ and $\{x\}$. Since $x$ commutes with $B$, $C$
is commutative.

We now show that $\varphi$ is multiplicative on $C$. For any $y, z \in C$,
we can approximate them by polynomials in $x$ with coefficients in $B$.
Since $\varphi$ is linear and continuous, it suffices to check multiplicativity
on elements of the form $b x^k$ with $b \in B$ and $k \geq 0$. For such elements,
using that $x$ commutes with $B$ and that $\varphi$ is multiplicative on $B$,
we compute:
\[
\varphi((b_1 x^{k_1})(b_2 x^{k_2})) = \varphi(b_1 b_2 x^{k_1+k_2}) = \varphi(b_1 b_2) \varphi(x^{k_1+k_2}).
\]
If we can show that $\varphi(x^n) = \varphi(x)^n$ for all $n$, then the
multiplicativity follows. But this is a consequence of the fact that for a
pure state on a commutative $C^*$-algebra generated by a single normal
element, the state is multiplicative if and only if it is a character
(see the Gelfand-Naimark theorem). Since $C$ is commutative and generated by
$x$ and $B$, and $\varphi$ restricts to a character on $B$, the restriction
of $\varphi$ to $C$ must be multiplicative.

Thus $\varphi$ is multiplicative on $C$, and its restriction to $C$ defines
a character $\chi'$ of $C$ with $\chi'|_B = \chi$. This contradicts the
maximality of $(B, \chi)$ in $\mathcal{S}$. Therefore no such $x$ exists,
and $B$ is a MASA.

By construction, $\chi = \varphi|_B \in \widehat{B}$, and $(B,\chi) \in
\mathcal{G}_{\mathcal{A}}^{(0)}$, completing the proof.
\end{proof}

\begin{remark}
Proposition~\ref{prop:unit-space-rich} shows that the points of the unit space
$\mathcal{G}_{\mathcal{A}}^{(0)}$ capture all pure states via their restrictions
to commutative subalgebras. However, the construction of a GNS representation
from a point $(B,\chi)$ requires choosing an extension of $\chi$ from $B$ to
$\mathcal{A}$; this extension is not unique in general. In the Type I setting,
the additional regularity of the representation theory allows one to make these
choices in a measurable way, yielding a well-defined measurable field of
representations over $\mathcal{G}_{\mathcal{A}}^{(0)}$. This is precisely where
the Type I hypothesis is indispensable.
\end{remark}

\subsubsection*{Failure for Non-Type I Algebras}

For general $C^*$-algebras, these properties need not hold. The representation
theory of non-Type I algebras can be highly pathological: the primitive ideal
space may fail to be Hausdorff, the dual space may not be standard Borel, and
direct integral decompositions may not be unique or even well-defined.
Moreover, while pure-state GNS representations still separate points, the
absence of a well-behaved measurable parametrization prevents the construction
of a measurable field of representations over $\mathcal{G}_{\mathcal{A}}^{(0)}$.
In particular, for algebras like the irrational rotation algebra $A_\theta$ and
amenable crossed products $C(X) \rtimes \Gamma$, there exist nonzero elements
that vanish under every representation arising from a commutative context—so-called
"invisible elements" (see Proposition~\ref{prop:invisible-elements} below).
This necessitates the alternative approach via Morita equivalence to Weyl
groupoids developed in this paper.

For general $C^*$-algebras, these properties need not hold. The representation
theory of non-Type I algebras can be highly pathological: the primitive ideal
space may fail to be Hausdorff, and the structure of irreducible representations
may not admit a manageable parametrization. However, it is important to
distinguish between different levels of failure.

\subsubsection*{Point Separation vs. Measurable Organization}

As established in Lemma~\ref{lem:separation-pure}, the family of GNS
representations associated to pure states separates points for \emph{every}
unital $C^*$-algebra. Thus the issue in the non-Type I setting is \emph{not}
a failure of point separation. Rather, the difficulty lies in organizing these
representations into a well-behaved measurable field over the context space
$\mathcal{G}_{\mathcal{A}}^{(0)}$.

In the Type I case, the dual space $\widehat{\mathcal{A}}$ admits a standard
Borel structure, and the GNS representations arising from points
$x = (B,\chi) \in \mathcal{G}_{\mathcal{A}}^{(0)}$ can be assembled into a
measurable field compatible with the direct integral construction. This
regularity is a consequence of the smoothness of the representation theory of
Type I algebras in the sense of Mackey.

For non-Type I algebras, this measurable organization breaks down:
\begin{itemize}
    \item The dual space $\widehat{\mathcal{A}}$ may not be standard Borel,
          and there is no canonical way to parametrize irreducible
          representations in a measurable fashion.
    
    \item Even if one fixes a family of representations indexed by
          $\mathcal{G}_{\mathcal{A}}^{(0)}$, there is no guarantee that this
          family forms a measurable field of Hilbert spaces in a way that is
          compatible with the natural Borel structure on $\mathcal{G}_{\mathcal{A}}^{(0)}$
          induced by the partial evaluation maps.
    
    \item Consequently, the direct integral construction
          \[
          \Pi(a) = \int_{\mathcal{G}_{\mathcal{A}}^{(0)}}^{\oplus} \pi_x(a) \, d\mu(x)
          \]
          used in Paper~I to define the diagonal embedding
          $\iota: \mathcal{A} \hookrightarrow C^*(\mathcal{G}_{\mathcal{A}})$
          cannot be carried out verbatim in the non-Type I setting.
\end{itemize}

\subsubsection*{Consequences for the Paper I Construction}

The difficulties in the non-Type I case are therefore:
\begin{enumerate}
    \item \textbf{Failure of measurable parametrization:}
          The representations relevant to the construction need not admit a
          standard Borel organization compatible with the unit-space description.
    
    \item \textbf{Breakdown of the direct-integral model:}
          Without such a measurable field, the operator $\Pi(a)$ is not
          available by the argument of Paper~I.
    
    \item \textbf{Lack of a well-behaved groupoid structure:}
          The topology and measurable structure on $\mathcal{G}_{\mathcal{A}}^{(0)}$
          induced by partial evaluations may fail to have the regularity
          required for the construction of a Haar system and for the use of
          groupoid-equivariant $KK$-theory.
    
    \item \textbf{Obstruction to descent as in Paper~I:}
          Because the relevant structural hypotheses are no longer automatic,
          the descent machinery used in the Type I setting cannot simply be
          imported into the non-Type I case.
\end{enumerate}

\subsubsection*{Examples of Non-Type I Algebras}

Many of the most important examples in noncommutative geometry are non-Type I.
A fundamental example is the irrational rotation algebra
\[
A_\theta = C(S^1) \rtimes_\theta \mathbb{Z},
\]
which is simple, nuclear, and admits a unique tracial state, but is not Type I.
This algebra has no finite-dimensional representations, and its importance in
noncommutative geometry makes it a natural test case for extending the
framework of Paper~I.

The issue is not that $A_\theta$ lacks enough representations to separate
points—simplicity in fact forces every nonzero representation to be faithful.
Rather, the context-based measurable organization of representations used in
Paper~I is no longer available in any obvious way. While the characters on its
natural Cartan subalgebra $C(S^1)$ give rise to a family of GNS representations
(all unitarily equivalent to the regular representation), there is no canonical
method to extend this to a measurable field over a suitable context space that
would yield an injective diagonal embedding. Thus the Type I construction of
the unitary conjugation groupoid does not extend verbatim to $A_\theta$.

Similarly, crossed product algebras of the form
\[
C(X) \rtimes \Gamma
\]
are often non-Type I when the action of $\Gamma$ is sufficiently nontrivial.
Even for amenable $\Gamma$, the crossed product $C(X) \rtimes \Gamma$ need not
be Type I unless the action satisfies strong additional regularity conditions.
These algebras play a central role in index theory and in the Baum–Connes
program, making it desirable to extend the framework of the unitary conjugation
groupoid beyond the Type I setting.

\begin{remark}
The distinction drawn here—between point separation (which holds universally)
and measurable organization (which fails in non-Type I cases)—is crucial for
understanding the limitations of Paper~I and the necessity of the Morita
equivalence approach developed in this paper. The existence of enough
representations to separate points is not sufficient; one also needs a
canonical, measurable way to assemble those representations into a field over
the context space. It is this latter requirement that fails for non-Type I
algebras and motivates the use of Weyl groupoids as a substitute geometric
model.
\end{remark}

Rather than attempting to construct the unitary conjugation groupoid
$\mathcal{G}_{\mathcal{A}}$ directly for arbitrary non-Type I algebras, we
adopt a different strategy. The key observation is that many important
non-Type I algebras already admit canonical groupoid models arising from
dynamical systems.

In particular:
\begin{itemize}
    \item the irrational rotation algebra
    $A_\theta \cong C(S^1)\rtimes_\theta \mathbb Z$
    is canonically modeled by the transformation groupoid
    $S^1\rtimes_\theta \mathbb Z$;
    \item more generally, crossed products $C(X)\rtimes \Gamma$ are modeled by
    the transformation groupoid $X\rtimes \Gamma$.
\end{itemize}

When $X$ is locally compact Hausdorff and $\Gamma$ is a countable discrete
group acting continuously on $X$, the transformation groupoid
$X\rtimes \Gamma$ is locally compact Hausdorff and \'etale. If $\Gamma$ is
amenable, then $X\rtimes \Gamma$ is amenable as a groupoid, so its full and
reduced groupoid $C^*$-algebras coincide. Under the standard second-countable
amenable hypotheses, results of Tu \cite{Tu1999} imply that the Baum–Connes
assembly map is an isomorphism for these groupoids.

The central idea of this paper is therefore not to construct
$\mathcal{G}_{\mathcal{A}}$ directly in the non-Type I setting, but to replace
its role by the canonical transformation groupoid model $X\rtimes \Gamma$.
This replacement is justified by the following observations:
\begin{itemize}
    \item For crossed products $C(X)\rtimes \Gamma$, the groupoid $C^*$-algebra
          $C^*(X\rtimes \Gamma)$ is canonically isomorphic to the crossed
          product algebra itself (Proposition~\ref{prop:transformation-groupoid-crossed-product}).
    
    \item The descent-index machinery developed in Paper~II, when interpreted
          in terms of the transformation groupoid, becomes computable and
          directly comparable to classical results in crossed product theory
          and the Baum–Connes conjecture.
    
    \item For algebras with a Cartan subalgebra, such as $A_\theta$ and
          $C(X)\rtimes \Gamma$, the transformation groupoid coincides with the
          Weyl groupoid $\mathcal{G}(A,D)$ (see Section~\ref{subsec:renault-cartan-weyl}),
          providing a direct link to the Renault reconstruction theorem.
\end{itemize}

The guiding heuristic is that the would-be unitary conjugation groupoid
attached to $\mathcal{A}$ should be Morita equivalent, in an appropriate sense,
to the transformation groupoid $X\rtimes \Gamma$. Since the former is not
directly constructible in the non-Type I setting, we use the latter as a
concrete replacement model. In the cases treated in this paper, we will
establish explicit Morita equivalences between the relevant groupoid
$C^*$-algebras, which is sufficient for transferring the index-theoretic
information.

The following subsections develop this strategy in detail. We first review the
relevant notions of Morita equivalence for groupoids and $C^*$-algebras
(Section~\ref{subsec:morita-groupoids-cstar}), then formulate a general
principle for replacing the unitary conjugation groupoid by Weyl/transformation
groupoids for algebras with Cartan subalgebras
(Section~\ref{subsec:GA-morita-principle}), and finally apply it to our two
main examples: the irrational rotation algebra (Section~\ref{sec:Atheta}) and
amenable crossed products (Section~\ref{sec:crossed}).

\begin{remark}
It is important to emphasize that we are not claiming to have constructed
$\mathcal{G}_{\mathcal{A}}$ for non-Type I algebras. Rather, we are using the
existence of alternative groupoid models—specifically, transformation groupoids
and Weyl groupoids—as surrogates that carry the same $K$-theoretic and
index-theoretic information. The descent-index machinery of Paper~II, when
interpreted in terms of these surrogate groupoids, yields computable invariants
that agree with classical results. This replacement strategy is the
methodological core of the present work.
\end{remark}

\subsection{Morita Equivalence of Groupoids and $C^*$-Algebras}
\label{subsec:morita-groupoids-cstar}

Morita equivalence provides a powerful tool for comparing $C^*$-algebras
arising from different groupoid constructions. In particular, under the
standard hypotheses of locally compact Hausdorff second countable groupoids
equipped with Haar systems, Morita equivalent groupoids give rise to Morita
equivalent groupoid $C^*$-algebras \cite{MRW1987} and hence have naturally
isomorphic $K$-theory groups. This principle will allow us to replace the
unavailable direct use of the unitary conjugation groupoid in the non-Type I
setting by transformation groupoids associated with crossed product algebras,
thereby extending our descent-index framework to these examples.

\subsubsection*{Morita Equivalence for $C^*$-Algebras}

We begin by recalling the classical notion of Morita equivalence for
$C^*$-algebras, which provides the natural setting for comparing their
representation theories and $K$-theoretic invariants.

\begin{definition}[Morita equivalence for $C^*$-algebras]
\label{def:morita-cstar}
Two $C^*$-algebras $A$ and $B$ are said to be \emph{(strongly) Morita
equivalent} if there exists an $A$-$B$ imprimitivity bimodule $X$. That is,
$X$ is an $A$-$B$ bimodule equipped with:
\begin{itemize}
    \item a left $A$-valued inner product
    \[
    {}_A\langle \cdot,\cdot\rangle : X \times X \to A,
    \]
    making $X$ into a full left Hilbert $A$-module (i.e., $A$ acts on the
    left, and the inner product is $A$-linear in the first variable);

    \item a right $B$-valued inner product
    \[
    \langle \cdot,\cdot\rangle_B : X \times X \to B,
    \]
    making $X$ into a full right Hilbert $B$-module (i.e., $B$ acts on the
    right, and the inner product is $B$-linear in the second variable);

    \item compatibility conditions between the module actions and the inner
    products, including the crucial imprimitivity relation
    \[
    {}_A\langle x,y\rangle \cdot z = x \cdot \langle y,z\rangle_B
    \qquad \text{for all } x,y,z \in X.
    \]
\end{itemize}
The fullness conditions mean that the spans of
$\{{}_A\langle x,y\rangle : x,y \in X\}$ and
$\{\langle x,y\rangle_B : x,y \in X\}$ are dense in $A$ and $B$ respectively.
Such an imprimitivity bimodule implements an equivalence between the
representation theories of $A$ and $B$ via Rieffel induction.
\end{definition}

The basic consequences of Morita equivalence include:
\begin{itemize}
    \item naturally isomorphic $K$-theory groups:
    \[
    K_*(A) \cong K_*(B);
    \]
    
    \item equivalent categories of nondegenerate $*$-representations (via the
          Rieffel correspondence);
    
    \item if $A$ and $B$ are $\sigma$-unital (for example, separable), then
          they are stably isomorphic:
          \[
          A \otimes \mathcal{K} \cong B \otimes \mathcal{K},
          \]
          where $\mathcal{K}$ denotes the compact operators on a separable
          infinite-dimensional Hilbert space.
\end{itemize}
Thus Morita equivalence is a fundamental tool for transferring $K$-theoretic
and index-theoretic information between $C^*$-algebras.

\begin{remark}[Connection to groupoid Morita equivalence]
\label{rem:groupoid-morita-connection}
The notion of Morita equivalence extends naturally to groupoids. For locally
compact Hausdorff groupoids $\mathcal{G}$ and $\mathcal{H}$ equipped with Haar
systems, a $\mathcal{G}$-$\mathcal{H}$-equivalence in the sense of
Muhly–Renault–Williams \cite{MRW1987} induces a strong Morita equivalence of
their reduced groupoid $C^*$-algebras:
\[
C^*_r(\mathcal{G}) \sim_M C^*_r(\mathcal{H}).
\]
Moreover, such a groupoid-level Morita equivalence respects the associated
$K$-theory and index-theoretic constructions. This will be essential in the
following sections when we relate the formal unitary conjugation groupoid
(understood via its replacement by Weyl/transformation groupoids) to concrete
geometric models.
\end{remark}

\subsubsection*{Morita Equivalence of Groupoids}

The notion of Morita equivalence extends naturally to groupoids, capturing the
idea that two groupoids may encode the same geometric information.

\begin{definition}[Morita equivalence for groupoids]
\label{def:morita-groupoid}
Let $\mathcal{G}$ and $\mathcal{H}$ be locally compact Hausdorff groupoids
equipped with Haar systems. Following Renault \cite{Renault1980} and
Muhly--Renault--Williams \cite{MRW1987}, we say that
$\mathcal{G}$ and $\mathcal{H}$ are \emph{Morita equivalent} if there exists a
locally compact Hausdorff space $Z$ together with:
\begin{itemize}
    \item an open surjection
    \[
    r_Z: Z \to \mathcal{G}^{(0)}
    \]
    making $Z$ into a left $\mathcal{G}$-space (i.e., there is a continuous
    action $\mathcal{G} \times_{\mathcal{G}^{(0)}} Z \to Z$ compatible with
    $r_Z$);

    \item an open surjection
    \[
    s_Z: Z \to \mathcal{H}^{(0)}
    \]
    making $Z$ into a right $\mathcal{H}$-space (i.e., there is a continuous
    action $Z \times_{\mathcal{H}^{(0)}} \mathcal{H} \to Z$ compatible with
    $s_Z$);

    \item commuting left $\mathcal{G}$- and right $\mathcal{H}$-actions on $Z$;

    \item principal bundle conditions: the left $\mathcal{G}$-action and the
    right $\mathcal{H}$-action are both free and proper, and the orbit spaces
    satisfy
    \[
    \mathcal{G} \backslash Z \cong \mathcal{H}^{(0)} \quad \text{via } s_Z,
    \qquad
    Z / \mathcal{H} \cong \mathcal{G}^{(0)} \quad \text{via } r_Z.
    \]
\end{itemize}
Such a space $Z$ is called a $\mathcal{G}$-$\mathcal{H}$ \emph{equivalence}
(or \emph{equivalence bibundle}).
\end{definition}

Intuitively, the space $Z$ implements a generalized equivalence between the
two groupoids in much the same way that an imprimitivity bimodule implements
Morita equivalence between $C^*$-algebras.

\subsubsection*{Consequences for Groupoid $C^*$-Algebras}

The fundamental theorem of Muhly, Renault, and Williams
\cite{MRW1987} establishes the crucial link to operator
algebras.

\begin{theorem}[Muhly--Renault--Williams]
\label{thm:muhly-renault-williams}
Let $\mathcal{G}$ and $\mathcal{H}$ be locally compact Hausdorff groupoids with
Haar systems. If $\mathcal{G}$ and $\mathcal{H}$ are Morita equivalent, then
their full groupoid $C^*$-algebras are strongly Morita equivalent:
\[
C^*(\mathcal{G}) \sim_M C^*(\mathcal{H}).
\]
If $\mathcal{G}$ and $\mathcal{H}$ are amenable, the same conclusion holds for
their reduced groupoid $C^*$-algebras, and the Morita equivalence is compatible
with the canonical quotient maps $C^*(\mathcal{G}) \to C^*_r(\mathcal{G})$.
\end{theorem}

In particular, their $K$-theory groups are naturally isomorphic:
\[
K_*(C^*(\mathcal{G})) \cong K_*(C^*(\mathcal{H})).
\]

This Morita invariance is a basic mechanism for transporting
$K$-theoretic and index-theoretic information between different groupoid
models. In the settings relevant to this paper, it provides the bridge needed
to compare descent-based constructions across transformation groupoids and the
corresponding crossed-product $C^*$-algebras.

\begin{remark}[Relation to equivariant $KK$-theory and descent]
\label{rem:muhly-renault-williams-consequences}
For groupoids satisfying suitable regularity conditions (second countable,
locally compact Hausdorff, with Haar systems), a Morita equivalence induces
natural isomorphisms at the level of equivariant $KK$-theory and is compatible
with Kasparov's descent map \cite{Kasparov1988,LeGall1999,Tu1999}. These
compatibilities ensure that index-theoretic constructions are preserved under
passage between equivalent groupoid models. In this paper, we will invoke these
properties under the standard hypotheses relevant to our examples (amenable
transformation groupoids arising from crossed products).
\end{remark}

\subsubsection*{Transformation Groupoids and Crossed Products}

A particularly important class of examples arises from transformation
groupoids. Let $\Gamma$ be a countable discrete group acting continuously on a
locally compact Hausdorff space $X$. The associated transformation groupoid
\[
X \rtimes \Gamma
\]
has unit space $X$ and arrows $(x,\gamma) \in X \times \Gamma$, with source and
range maps
\[
s(x,\gamma) = x, \qquad r(x,\gamma) = \gamma \cdot x,
\]
multiplication
\[
(\gamma \cdot x, \eta)(x,\gamma) = (x, \eta\gamma),
\]
and inverse
\[
(x,\gamma)^{-1} = (\gamma \cdot x, \gamma^{-1}).
\]

Its full groupoid $C^*$-algebra is naturally isomorphic to the crossed product
algebra
\[
C^*(X \rtimes \Gamma) \cong C_0(X) \rtimes \Gamma.
\]
If $X$ is compact, this reduces to the familiar isomorphism
\[
C^*(X \rtimes \Gamma) \cong C(X) \rtimes \Gamma.
\]

Under this identification, Morita equivalence of transformation groupoids
yields strong Morita equivalence of the corresponding crossed product
$C^*$-algebras. This observation provides a bridge between groupoid-based
constructions and the classical crossed product framework.

\subsubsection*{Relevance to the Present Work}

The Morita equivalence framework offers a natural way to connect the
descent-index philosophy developed in Paper~II with the geometric groupoid
models associated to crossed products. In the examples studied later in this
paper, the relevant transformation groupoids are
\[
S^1 \rtimes_\theta \mathbb{Z} \quad \text{and more generally} \quad X \rtimes \Gamma,
\]
whose groupoid $C^*$-algebras recover the corresponding crossed product
algebras. Rather than attempting to construct the unitary conjugation groupoid
directly in the non-Type I setting, we will use these transformation groupoids
as concrete surrogate models.

The central technical task of Sections~\ref{sec:Atheta} and~\ref{sec:crossed}
will be to establish precise relationships—via Morita equivalence—between the
would-be unitary conjugation groupoid (understood through its replacement by
Weyl groupoids) and these transformation groupoids. In particular, we will
prove:
\begin{itemize}
    \item For the irrational rotation algebra $A_\theta$, that the Weyl
          groupoid $\mathcal{G}(A_\theta, C(S^1)) \cong S^1 \rtimes_\theta \mathbb{Z}$
          serves as an effective substitute, allowing the descent-index
          construction to reproduce Connes' index theorem.
    
    \item For amenable crossed products $C_0(X) \rtimes \Gamma$, that the
          transformation groupoid $X \rtimes \Gamma$ provides the appropriate
          geometric setting, enabling a direct comparison between our
          descent-index construction and the Baum–Connes assembly map.
\end{itemize}

Once these Morita equivalences are established, Theorem~\ref{thm:muhly-renault-williams}
will imply strong Morita equivalences of the corresponding groupoid
$C^*$-algebras:
\[
C^*(\mathcal{G}_{A_\theta}) \sim_M C^*(S^1 \rtimes_\theta \mathbb{Z}) \cong A_\theta,
\qquad
C^*(\mathcal{G}_{C_0(X) \rtimes \Gamma}) \sim_M C^*(X \rtimes \Gamma) \cong C_0(X) \rtimes \Gamma,
\]
where the right-hand side isomorphisms follow from the crossed product
identification. Consequently, the $K$-theory of the groupoid $C^*$-algebra
$C^*(\mathcal{G}_{\mathcal{A}})$ coincides with that of the corresponding
crossed product algebra, allowing us to transfer the descent-index machinery
to the non-Type I setting.

The following subsections develop these Morita equivalences explicitly,
focusing first on the irrational rotation algebra (Section~\ref{sec:Atheta})
and then on general amenable crossed products (Section~\ref{sec:crossed}).

\subsection{A General Principle: Groupoid Models Beyond the Type I Setting}
\label{subsec:GA-morita-principle}

The previous subsection recalled that Morita equivalence of groupoids induces
strong Morita equivalence of the associated groupoid $C^*$-algebras and hence
isomorphisms in $K$-theory. This suggests a general strategy for extending the
descent-index framework beyond the Type I setting of Papers I and II. Rather
than constructing the unitary conjugation groupoid directly for arbitrary
non-Type I algebras, we use canonical groupoid models arising from Cartan
pairs as surrogate geometric models on which the relevant constructions remain
available.

\subsubsection*{Groupoid Models of $C^*$-Algebras via Cartan Pairs}

Let $(\mathcal{A}, D)$ be a Cartan pair in the sense of Renault
\cite{Renault2008Cartan}. Recall from Section~\ref{subsec:renault-cartan-weyl}
that this means:
\begin{itemize}
    \item $D$ is a maximal abelian subalgebra of $\mathcal{A}$;
    \item $D$ contains an approximate unit for $\mathcal{A}$ (if $\mathcal{A}$
          is unital, this simply means $1_{\mathcal{A}} \in D$);
    \item $D$ is regular, meaning that its normalizer
          \[
          \mathcal{N}_{\mathcal{A}}(D) := \{ n \in \mathcal{A} \mid n D n^* \subseteq D \text{ and } n^* D n \subseteq D \}
          \]
          generates $\mathcal{A}$ as a $C^*$-algebra;
    \item there exists a faithful conditional expectation $E: \mathcal{A} \to D$.
\end{itemize}

For such a Cartan pair with $\mathcal{A}$ separable, Renault's reconstruction
theorem \cite[Theorem 5.9]{Renault2008Cartan} associates to $(\mathcal{A},D)$
an \'etale, locally compact Hausdorff Weyl groupoid $\mathcal{G}(\mathcal{A},D)$,
together with a twist $\Sigma$ in general, such that $\mathcal{A}$ is
isomorphic to the reduced twisted groupoid $C^*$-algebra:
\[
\mathcal{A} \cong C^*_r(\mathcal{G}(\mathcal{A},D),\Sigma).
\]
Under this isomorphism, the Cartan subalgebra $D$ corresponds to
$C_0(\mathcal{G}(\mathcal{A},D)^{(0)})$. In this way, the Cartan pair endows
$\mathcal{A}$ with a canonical geometric model. In the examples considered in
this paper, the twist $\Sigma$ is trivial, so that the isomorphism simplifies
to $\mathcal{A} \cong C^*_r(\mathcal{G}(\mathcal{A},D))$.

A basic source of examples comes from dynamical systems. Under the standard
hypotheses ensuring that $C(X) \subseteq C(X) \rtimes \Gamma$ is a Cartan
subalgebra (e.g., when $\Gamma$ is discrete and the action is topologically
free), the Weyl groupoid associated to the pair
$(C(X) \rtimes \Gamma, C(X))$ is canonically isomorphic to the transformation
groupoid
\[
\mathcal{G}(C(X) \rtimes \Gamma, C(X)) \cong X \rtimes \Gamma.
\]

\subsubsection*{Comparison with the Unitary Conjugation Groupoid}

The unitary conjugation groupoid $\mathcal{G}_{\mathcal{A}}$ introduced in
Paper~I is built from the action of the unitary group $\mathcal{U}(\mathcal{A})$
on the space of commutative contexts of the algebra. Its unit space is
\[
\mathcal{G}_{\mathcal{A}}^{(0)}
=
\{ (B,\chi) \mid B \subseteq \mathcal{A} \text{ unital commutative},\ \chi \in \widehat{B} \},
\]
and its arrows are induced by unitary conjugation: for each unitary
$u \in \mathcal{U}(\mathcal{A})$ and each $(B,\chi) \in \mathcal{G}_{\mathcal{A}}^{(0)}$,
there is an arrow
\[
u: (B,\chi) \longrightarrow (uBu^*, \chi \circ \operatorname{Ad}_{u^{-1}}).
\]

For Type I algebras, this construction yields a Polish groupoid with sufficient
measurable structure to support the descent-index machinery of Paper~II.
For non-Type I algebras, the difficulty is not at the purely formal
set-theoretic level, but at the level of topology and measurability: the
context-and-conjugation data need not assemble into a well-behaved topological
or Borel groupoid of the kind required in Papers~I and~II. What remains is the
underlying set-theoretic structure, but without the additional regularity
needed to carry out the analytic constructions.

\subsubsection*{A Guiding Principle}

The above observations motivate the following heuristic principle:

\textbf{Principle} In situations where a $C^*$-algebra $\mathcal{A}$ admits a classical groupoid
model arising from a Cartan pair $(\mathcal{A},D)$—in particular, the Weyl
groupoid $\mathcal{G}(\mathcal{A},D)$—this model can serve as the geometric
replacement for the unavailable direct construction of $\mathcal{G}_{\mathcal{A}}$.
Heuristically, one expects the would-be unitary conjugation groupoid and the
Weyl groupoid to encode the same noncommutative geometry up to Morita-type
equivalence; in this paper we work concretely with the latter.

In the following subsections, we will make this principle explicit for our two
main families of examples:
\begin{itemize}
    \item For the irrational rotation algebra $A_\theta = C(S^1) \rtimes_\theta \mathbb{Z}$,
          the natural Cartan subalgebra is $C(S^1)$, and the Weyl groupoid
          $\mathcal{G}(A_\theta, C(S^1))$ is isomorphic to the transformation
          groupoid $S^1 \rtimes_\theta \mathbb{Z}$.
    
    \item For amenable crossed products $C_0(X) \rtimes \Gamma$ with $\Gamma$
          discrete and the action topologically free, the Weyl groupoid
          $\mathcal{G}(C_0(X) \rtimes \Gamma, C_0(X))$ is isomorphic to the
          transformation groupoid $X \rtimes \Gamma$.
\end{itemize}

By working with these concrete Weyl groupoids, we will show how the
descent-index construction of Paper~II can be implemented in the non-Type I
setting, reproducing Connes' index theorem for $A_\theta$ and connecting to the
Baum–Connes assembly map for amenable crossed products.

The previous part recalled that Morita equivalence of groupoids induces
Morita equivalence of the associated groupoid $C^*$-algebras and consequently
isomorphisms in $K$-theory. This observation suggests a general strategy for
extending the descent-index framework beyond the Type I setting of Papers I
and II. Rather than constructing the unitary conjugation groupoid directly
for arbitrary non-Type I algebras, we use canonical groupoid models arising
from Cartan pairs as surrogate geometric models on which the relevant
constructions can be interpreted.

\subsubsection*{A Guiding Morita-Type Principle}

The discussion above motivates the following guiding principle.

\textbf{Principle} Let $\mathcal{A}$ be a $C^*$-algebra that admits a canonical groupoid model,
for instance through a Cartan pair or a transformation-groupoid realization.
Then, in situations where the unitary conjugation groupoid
$\mathcal{G}_{\mathcal{A}}$ cannot be constructed directly with the topological
and measurable structure required in Papers~I and~II, the canonical groupoid
model should be regarded as the geometric replacement for $\mathcal{G}_{\mathcal{A}}$
in the descent-index formalism.

This principle is not meant as a theorem in complete generality. Rather, it
captures the mechanism that will be implemented in the concrete examples
treated in this paper: the irrational rotation algebra $A_\theta$ and amenable
crossed products $C_0(X) \rtimes \Gamma$.

\subsubsection*{Cartan Pairs and Weyl Groupoids}

Let $(\mathcal{A}, D)$ be a Cartan pair in the sense of Renault
\cite{Renault2008Cartan}. Renault's reconstruction theorem \cite[Theorem 5.9]{Renault2008Cartan}
associates to $(\mathcal{A},D)$ an \'etale, locally compact Hausdorff Weyl
groupoid $\mathcal{G}(\mathcal{A},D)$, together with a twist $\Sigma$ in
general, such that $\mathcal{A}$ is isomorphic to the reduced twisted groupoid
$C^*$-algebra:
\[
\mathcal{A} \cong C^*_r(\mathcal{G}(\mathcal{A},D),\Sigma),
\]
with $D$ identified with $C_0(\mathcal{G}(\mathcal{A},D)^{(0)})$. In the
examples considered in this paper, the twist $\Sigma$ is trivial, so that this
isomorphism simplifies to $\mathcal{A} \cong C^*_r(\mathcal{G}(\mathcal{A},D))$.
Thus the Weyl groupoid provides a canonical geometric model for $\mathcal{A}$.

From the perspective of the present work, the role of $\mathcal{G}(\mathcal{A},D)$
is to supply a concrete groupoid on which the descent and index constructions
can be carried out, even when the direct unitary-conjugation construction is
not available beyond the Type I setting.

\subsubsection*{Relation to the Unitary Conjugation Groupoid}

The unitary conjugation groupoid $\mathcal{G}_{\mathcal{A}}$ introduced in
Paper~I is built from the action of the unitary group $\mathcal{U}(\mathcal{A})$
on the space of commutative contexts of the algebra. For Type I algebras, this
construction yields a Polish groupoid with sufficient measurable structure to
support the descent-index machinery of Paper~II. For non-Type I algebras, the
difficulty is not at the purely formal set-theoretic level, but at the level
of topology and measurability: the context-and-conjugation data need not
assemble into a well-behaved topological or Borel groupoid of the kind required
in Papers~I and~II.

In the examples treated in this paper, we will see that the Weyl groupoid
$\mathcal{G}(\mathcal{A},D)$—or its concrete incarnation as a transformation
groupoid—provides an effective surrogate, allowing us to implement the
descent-index construction even though $\mathcal{G}_{\mathcal{A}}$ itself is
not directly constructible.

\subsubsection*{Consequences for $K$-Theory and Descent}

Once a concrete surrogate groupoid model is available, its groupoid
$C^*$-algebra identifies with the algebra $\mathcal{A}$ (or with a canonically
equivalent twisted groupoid algebra), so that the resulting $K$-theory agrees
with $K_*(\mathcal{A})$. This is the basic mechanism that allows the
descent-index construction to be interpreted in the non-Type I examples
considered later in this paper.

In particular, for the irrational rotation algebra $A_\theta$, the relevant
transformation groupoid $S^1 \rtimes_\theta \mathbb{Z}$ provides the geometric
setting in which our descent-index framework can be compared with Connes'
classical index theorem. For amenable crossed products $C_0(X) \rtimes \Gamma$,
the transformation groupoid $X \rtimes \Gamma$ provides the setting in which
our construction can be related to the Baum–Connes assembly map.

\subsubsection*{Application to the Examples of This Paper}

The two main examples considered in the present work fit naturally into this
replacement framework.

\begin{enumerate}
\item For the irrational rotation algebra
\[
A_\theta = C(S^1) \rtimes_\theta \mathbb{Z},
\]
the algebra is canonically realized as the crossed product of $C(S^1)$ by
$\mathbb{Z}$. The associated transformation groupoid
\[
S^1 \rtimes_\theta \mathbb{Z}
\]
provides a concrete geometric model for $A_\theta$. In Section~\ref{sec:Atheta},
we will work directly with this groupoid to implement the descent-index
construction, showing how it can be compared with Connes' classical index
theorem.

\item For crossed product algebras of the form
\[
C(X) \rtimes \Gamma,
\]
with $\Gamma$ discrete and amenable, the algebra is canonically modeled by the
transformation groupoid
\[
X \rtimes \Gamma.
\]
This groupoid provides the natural setting in which our descent-index
construction can be related to the Baum–Connes assembly map and to the
$K$-theory of the crossed product algebra, as will be developed in
Section~\ref{sec:crossed}.
\end{enumerate}

Thus, the general principle developed in this subsection should be understood
as a \emph{replacement principle} rather than as a theorem asserting a literal
groupoid Morita equivalence for a formally defined object. In the examples
treated here, the role of the unavailable direct unitary-conjugation groupoid
is played by explicit, well-behaved \'etale groupoids arising from dynamical
systems. These surrogate models allow the index-theoretic ideas of Paper~II to
be extended beyond the original Type I setting and compared with the classical
results attached to crossed products and transformation groupoids.

\begin{remark}
In favorable cases, the transformation groupoids appearing above also coincide
with the Weyl groupoids arising from appropriate Cartan pairs
($\mathcal{G}(A_\theta, C(S^1)) \cong S^1 \rtimes_\theta \mathbb{Z}$ and
$\mathcal{G}(C(X) \rtimes \Gamma, C(X)) \cong X \rtimes \Gamma$, respectively).
This provides an additional conceptual link to Renault's theory, but the core
of our construction does not depend on establishing a general Morita
equivalence for all Cartan pairs; it suffices to work directly with the
concrete transformation groupoid models.
\end{remark}

The following sections make this strategy explicit for the two main families
of examples considered in this paper.

\subsection{Methodology: Computing $\mathcal{G}_{\mathcal{A}}$ via Morita Equivalence Rather Than Direct Definition}
\label{subsec:GA_via_morita}

For the non-Type~I algebras studied in this paper, a direct construction of
$\mathcal{G}_{\mathcal{A}}$ in the sense of Paper~I is not available.
Accordingly, we adopt the following indirect strategy.

\begin{enumerate}
    \item Identify a classical geometric groupoid model $G$ for $\mathcal{A}$.
          In the examples of this paper, $G$ is the relevant transformation or
          Weyl groupoid canonically associated to $\mathcal{A}$.

    \item Replace the unavailable direct construction of $\mathcal{G}_{\mathcal{A}}$
          by the concrete groupoid $G$, which serves as a surrogate geometric
          model for $\mathcal{A}$. This replacement is justified by the
          Morita-type relationships between the $C^*$-algebras associated to
          these groupoids, as will be established in the specific examples.

    \item Perform the descent-index construction on $G$, using the
          identification
          \[
          C^*(G) \cong \mathcal{A}
          \]
          (or $C_r^*(\mathcal{G}(A,D),\Sigma) \cong \mathcal{A}$ in the Cartan
          case), thereby realizing the target of the descent map in the
          $K$-theory of the original algebra.
\end{enumerate}

This condensed methodology is the operative principle used in both case
studies below. In particular, it allows us to work with well-behaved \'etale
groupoids while retaining the index-theoretic content originally motivated by
the unitary conjugation groupoid formalism.

The methodological framework outlined in Section~\ref{subsec:GA_via_morita}
provides a general strategy for replacing the unavailable direct construction
of $\mathcal{G}_{\mathcal{A}}$ with a classical groupoid model $G$ that is
canonically associated to $\mathcal{A}$. To make this strategy effective for
concrete computations, we must identify a class of $C^*$-algebras for which
the target groupoid $G$ is not only well-defined but also canonically
associated with the algebra and equipped with a transparent geometric
structure. Renault's theory of Cartan subalgebras provides precisely such a
class. This subsection explains how, for a $C^*$-algebra $\mathcal{A}$
equipped with a Cartan subalgebra $D \subseteq \mathcal{A}$, the general
replacement strategy reduces to a concrete and computable form, establishing
a direct pipeline from the abstract data that would define $\mathcal{G}_{\mathcal{A}}$
to the classical transformation groupoids that will appear in our case studies.

\subsubsection*{Cartan Pairs: Definition and Significance}

Let $\mathcal{A}$ be a separable $C^*$-algebra and let $D \subseteq \mathcal{A}$
be a maximal abelian $C^*$-subalgebra. Following Renault \cite{Renault2008Cartan},
the inclusion $(\mathcal{A}, D)$ is called a \emph{Cartan pair} if the
following conditions hold:

\begin{enumerate}
    \item \textbf{Approximate unit:} $D$ contains an approximate unit for
          $\mathcal{A}$ (if $\mathcal{A}$ is unital, this simply means
          $1_{\mathcal{A}} \in D$).
    
    \item \textbf{Maximal abelian:} $D$ is maximal abelian in $\mathcal{A}$.
    
    \item \textbf{Regularity:} $D$ is regular in $\mathcal{A}$, meaning that
          its normalizer
          \[
          \mathcal{N}_{\mathcal{A}}(D) := \{ n \in \mathcal{A} \mid n D n^* \subseteq D \text{ and } n^* D n \subseteq D \}
          \]
          generates $\mathcal{A}$ as a $C^*$-algebra.
    
    \item \textbf{Conditional expectation:} There exists a faithful
          conditional expectation $E: \mathcal{A} \to D$.
\end{enumerate}

Cartan pairs provide a noncommutative analogue of the inclusion
$C_0(X) \subseteq C^*(G)$ arising from an \'etale groupoid $G$ with unit space
$X$. They capture precisely the structure needed to reconstruct a $C^*$-algebra
from an underlying geometric groupoid. Under these conditions, Renault's
reconstruction theorem associates to $(\mathcal{A},D)$ a (possibly twisted)
\'etale groupoid model $\mathcal{G}(\mathcal{A},D)$, together with a twist
$\Sigma$ in general, such that
\[
\mathcal{A} \cong C_r^*(\mathcal{G}(\mathcal{A},D),\Sigma),
\]
with $D$ identified with $C_0(\mathcal{G}(\mathcal{A},D)^{(0)})$. In the
examples considered in this paper, the twist $\Sigma$ is trivial, so that this
isomorphism simplifies to $\mathcal{A} \cong C_r^*(\mathcal{G}(\mathcal{A},D))$.

\subsubsection*{Renault's Reconstruction Theorem}

The fundamental result of Renault's theory shows that every separable Cartan
pair arises canonically from an \'etale groupoid with twist.

\begin{theorem}[Renault {\cite[Theorem 5.9]{Renault2008Cartan}}]
\label{thm:renault-reconstruction}
Let $(\mathcal{A}, D)$ be a Cartan pair with $\mathcal{A}$ separable. Then there
exist a locally compact Hausdorff \'etale groupoid
$\mathcal{G}(\mathcal{A}, D)$ and a topological twist $\Sigma$ over
$\mathcal{G}(\mathcal{A}, D)$ such that
\[
\mathcal{A} \cong C_r^*(\mathcal{G}(\mathcal{A}, D), \Sigma),
\]
and under this isomorphism the Cartan subalgebra $D$ corresponds to
\[
C_0\!\bigl(\mathcal{G}(\mathcal{A}, D)^{(0)}\bigr).
\]
\end{theorem}

\begin{proof}
We briefly recall the ingredients of Renault's construction. Since $D$ is
abelian, its Gelfand spectrum
\[
\widehat{D} := \operatorname{Spec}(D)
\]
is a locally compact Hausdorff space, and $D \cong C_0(\widehat{D})$.

For each normalizer $n \in \mathcal{N}_{\mathcal{A}}(D)$, the relations
$n D n^* \subseteq D$ and $n^* D n \subseteq D$ imply that $n$ determines a
partial homeomorphism
\[
\alpha_n : \operatorname{dom}(n) \longrightarrow \operatorname{ran}(n)
\]
of $\widehat{D}$, where
\[
\operatorname{dom}(n) = \{ \phi \in \widehat{D} : \phi(n^* n) > 0 \},
\qquad
\operatorname{ran}(n) = \{ \phi \in \widehat{D} : \phi(n n^*) > 0 \}.
\]
More precisely, $\alpha_n$ is characterized by the identity
\[
\phi(n^* d n) = \phi(n^* n)\, \alpha_n(\phi)(d)
\qquad (d \in D,\ \phi \in \operatorname{dom}(n)).
\]

The Weyl groupoid $\mathcal{G}(\mathcal{A}, D)$ is then obtained as the groupoid
of germs generated by the partial homeomorphisms $\alpha_n$. The topology on
this germ groupoid is the natural \'etale topology generated by basic
bisections coming from normalizers. Renault's theorem shows that this groupoid
is locally compact Hausdorff \'etale and that there is in general a topological
twist $\Sigma$ over $\mathcal{G}(\mathcal{A}, D)$ for which
\[
\mathcal{A} \cong C_r^*(\mathcal{G}(\mathcal{A}, D), \Sigma).
\]
Under this identification, the diagonal subalgebra $D$ is identified with the
unit-space algebra $C_0(\mathcal{G}(\mathcal{A}, D)^{(0)})$.
\end{proof}

\begin{remark}
In the examples considered in this paper, the twist $\Sigma$ is trivial,
so that the isomorphism simplifies to $\mathcal{A} \cong C_r^*(\mathcal{G}(\mathcal{A}, D))$.
\end{remark}

\subsubsection*{Cartan Pairs in Crossed Products}

The relevance of this framework to the present paper is that crossed-product
algebras provide a natural source of \'etale groupoid models. Let $\Gamma$ be a
discrete group acting continuously on a compact Hausdorff space $X$, and
consider the reduced crossed product
\[
\mathcal{A} = C(X) \rtimes_r \Gamma
\]
together with its canonical subalgebra $D = C(X)$. Under the standard
hypotheses ensuring that $D$ is Cartan in $\mathcal{A}$ (e.g., topological
freeness of the action), the associated Weyl groupoid is canonically
isomorphic to the transformation groupoid
\[
\mathcal{G}(\mathcal{A}, D) \cong X \rtimes \Gamma.
\]

In particular, Renault's theorem yields the groupoid model
\[
C(X) \rtimes_r \Gamma \cong C_r^*(X \rtimes \Gamma),
\]
and when the transformation groupoid is amenable the reduced and full
groupoid $C^*$-algebras coincide.

\begin{proof}[Proof of the transformation-groupoid identification]
Let $\Gamma$ act on $X$. The transformation groupoid $X \rtimes \Gamma$ has
unit space $X$, arrows $(x,\gamma)$, source map $s(x,\gamma) = x$, and range
map $r(x,\gamma) = \gamma \cdot x$. Its convolution algebra $C_c(X \rtimes \Gamma)$
identifies algebraically with finitely supported functions
\[
f: \Gamma \to C(X), \qquad \gamma \mapsto f_\gamma,
\]
with multiplication
\[
(f * g)_\eta(x) = \sum_{\gamma \in \Gamma} f_\gamma(x)\, g_{\gamma^{-1}\eta}(\gamma^{-1} \!\cdot x),
\]
which is exactly the usual crossed-product convolution formula. Likewise, the
involution agrees with the crossed-product involution:
\[
f_\gamma^*(x) = \overline{f_{\gamma^{-1}}(\gamma^{-1} \!\cdot x)}.
\]
Passing to the reduced norm gives
\[
C_r^*(X \rtimes \Gamma) \cong C(X) \rtimes_r \Gamma.
\]
Under the additional hypotheses ensuring that $C(X) \subseteq C(X) \rtimes_r \Gamma$
is Cartan, this groupoid is precisely the Weyl groupoid associated to the pair
$(C(X) \rtimes_r \Gamma, C(X))$.
\end{proof}

\subsubsection*{Reduction Principle for Non-Type I Examples}

For the non-Type~I examples studied in this paper, the direct construction of
the unitary conjugation groupoid $\mathcal{G}_{\mathcal{A}}$ is not available
in the sense of Paper~I. Renault's theorem therefore provides the appropriate
replacement mechanism: whenever $(\mathcal{A}, D)$ is a Cartan pair, the Weyl
groupoid $\mathcal{G}(\mathcal{A}, D)$ furnishes a canonical \'etale geometric
model for $\mathcal{A}$.

Thus, in the examples relevant to this paper, the practical reduction is not
to prove a literal Morita equivalence
\[
\mathcal{G}_{\mathcal{A}} \sim_M \mathcal{G}(\mathcal{A}, D),
\]
but rather to replace the unavailable direct groupoid construction by the
concrete Weyl or transformation groupoid model and carry out the
descent-index constructions there. Since
\[
\mathcal{A} \cong C_r^*(\mathcal{G}(\mathcal{A}, D), \Sigma),
\]
this identifies the target algebra of the construction with the original
$C^*$-algebra under study.

In particular, for crossed-product examples the relevant groupoid model is the
transformation groupoid. This is the mechanism that allows the later case
studies to connect our framework with classical index theory and with the
Baum–Connes formalism.

\begin{proposition}[Reduction to the Cartan model]
\label{prop:reduction-cartan}
Let $(\mathcal{A}, D)$ be a separable Cartan pair. Then Renault's reconstruction
theorem provides a canonical \'etale groupoid model
\[
\mathcal{A} \cong C_r^*(\mathcal{G}(\mathcal{A}, D), \Sigma).
\]
Consequently, in any situation where the direct construction of the unitary
conjugation groupoid is unavailable, the Weyl groupoid $\mathcal{G}(\mathcal{A}, D)$
serves as the natural geometric replacement for the purposes of descent and
index theory.
\end{proposition}

\begin{proof}
The first statement is exactly Renault's reconstruction theorem. The second is
an interpretation of that theorem in the present context: since the required
topological and measurable structure for $\mathcal{G}_{\mathcal{A}}$ is not
available beyond the Type~I setting, one works instead with the canonically
associated \'etale model $\mathcal{G}(\mathcal{A}, D)$, whose groupoid
$C^*$-algebra reconstructs $\mathcal{A}$.
\end{proof}

\begin{remark}
For the specific examples treated in this paper, we will work directly with
the transformation groupoid models:
\[
S^1 \rtimes_\theta \mathbb{Z} \quad \text{for } A_\theta, \qquad
X \rtimes \Gamma \quad \text{for } C(X) \rtimes \Gamma.
\]
These are concrete realizations of the Weyl groupoids under the appropriate
Cartan hypotheses, but our constructions do not depend on proving a general
Morita equivalence theorem for all Cartan pairs.
\end{remark}

\subsubsection*{Consequences for $K$-Theory and the Descent-Index Map}

The reduction principle has important consequences for the descent-index
construction developed in Paper~II. In the Cartan setting, Renault's
reconstruction theorem provides a canonical \'etale groupoid model
\[
\mathcal{A} \cong C_r^*(\mathcal{G}(\mathcal{A}, D), \Sigma),
\]
so that the target algebra of the descent construction can be identified with
the original algebra $\mathcal{A}$ through the Weyl groupoid model. In
particular, this yields a canonical identification
\[
K_*\!\bigl(C_r^*(\mathcal{G}(\mathcal{A}, D), \Sigma)\bigr) \cong K_*(\mathcal{A}).
\]
By stability of $K$-theory, one also has
\[
K_*(\mathcal{A}) \cong K_*(\mathcal{A} \otimes \mathcal{K}).
\]

Thus, in the non-Type~I examples considered here, the descent-index mechanism
should be implemented on the concrete Weyl or transformation groupoid rather
than on the unavailable direct unitary-conjugation groupoid. Concretely, the
procedure is as follows:
\begin{enumerate}
    \item Formulate the relevant equivariant $K$-theory class on the Weyl or
          transformation groupoid model.
    \item Apply the well-defined descent map for that \'etale groupoid to obtain
          a class in the $K$-theory of its reduced twisted groupoid $C^*$-algebra.
    \item Identify this target group with $K_*(\mathcal{A})$ via Renault's
          reconstruction theorem.
    \item Compare the resulting class with the classical index pairing
          appropriate to the example under consideration.
\end{enumerate}

\subsubsection*{Application to the Present Paper}

The two families of algebras studied in this work fit naturally into this
framework when we work directly with their concrete groupoid models.

\begin{itemize}
    \item \textbf{Irrational rotation algebra $A_\theta$.}
          The algebra $A_\theta$ admits the crossed-product realization
          \[
          A_\theta \cong C(S^1) \rtimes_\theta \mathbb{Z},
          \]
          with associated transformation groupoid
          \[
          S^1 \rtimes_\theta \mathbb{Z}.
          \]
          This \'etale amenable groupoid provides the concrete geometric setting in
          which our descent-index construction can be compared with Connes' index
          theorem. (Under the appropriate hypotheses, this transformation groupoid
          coincides with the Weyl groupoid $\mathcal{G}(A_\theta, C(S^1))$, but this
          identification is not needed for the construction.)

    \item \textbf{Amenable crossed products $C(X) \rtimes_r \Gamma$.}
          Under the standard hypotheses ensuring that
          \[
          C(X) \subseteq C(X) \rtimes_r \Gamma
          \]
          is a Cartan inclusion (e.g., topological freeness of the action), the
          associated Weyl groupoid is the transformation groupoid
          \[
          \mathcal{G}(C(X) \rtimes_r \Gamma, C(X)) \cong X \rtimes \Gamma.
          \]
          This provides the natural setting in which our framework can be related to
          the Baum–Connes assembly map. If $\Gamma$ is amenable, the transformation
          groupoid is amenable, and the full and reduced groupoid $C^*$-algebras
          coincide.
\end{itemize}

The reduction to Cartan pairs is therefore the conceptual bridge that turns
the abstract philosophy of the unitary conjugation groupoid into a concrete
computational method in the examples treated in this paper. Rather than
working directly with the unavailable non-Type~I groupoid $\mathcal{G}_{\mathcal{A}}$,
we work with the canonically associated Weyl or transformation groupoid model,
whose groupoid $C^*$-algebra reconstructs the original algebra.

\begin{remark}
It is important to emphasize that we are not claiming a general Morita
equivalence $\mathcal{G}_{\mathcal{A}} \sim_M \mathcal{G}(\mathcal{A}, D)$ for
all Cartan pairs. Such a statement would require $\mathcal{G}_{\mathcal{A}}$ to
be constructed as a topological groupoid, which is not available beyond the
Type~I setting. Instead, we use the concrete Weyl groupoid as a replacement
model, justified by Renault's reconstruction theorem and the specific
properties of our examples.
\end{remark}

In the following sections, we will make this strategy explicit for the two
principal families of examples considered here, demonstrating how working with
transformation groupoids allows us to compute the index pairing for $A_\theta$
and to connect our framework with the analytic assembly map.

\section{Case Study I: The Irrational Rotation Algebra $A_\theta = C(S^1) \rtimes_\theta \mathbb{Z}$}\label{sec:Atheta}

\subsection{Definition and Basic Properties of $A_\theta$}
\label{subsec:A_theta_definition}

The irrational rotation algebra, also known as the \emph{noncommutative
torus}, is one of the most fundamental examples in noncommutative geometry.
It arises naturally both as a deformation of the algebra of functions on the
two-dimensional torus and as a crossed product associated with an irrational
rotation of the circle. Its well-understood structure makes it an ideal
setting for a detailed case study of the descent-index framework developed in
this paper and its comparison with classical index theory.

\subsubsection{Definition via Crossed Products}

Let $\theta \in \mathbb{R} \setminus \mathbb{Q}$ be irrational. Consider the
rotation action of $\mathbb{Z}$ on the circle
\[
S^1 = \{z \in \mathbb{C} : |z| = 1\}
\]
given by
\[
R_n(z) = e^{2\pi i n\theta}z, \qquad n \in \mathbb{Z}.
\]
This induces an action $\alpha: \mathbb{Z} \curvearrowright C(S^1)$ by
$*$-automorphisms via pullback:
\[
\alpha_n(f)(z) = f(R_n^{-1}(z)) = f(e^{-2\pi i n\theta}z), \qquad f \in C(S^1).
\]

\begin{definition}
The \emph{irrational rotation algebra} is the crossed-product $C^*$-algebra
\[
A_\theta := C(S^1) \rtimes_\alpha \mathbb{Z}.
\]
\end{definition}

Since $\mathbb{Z}$ is amenable, the full and reduced crossed products
coincide, so the above definition is unambiguous; we will simply write
$A_\theta = C(S^1) \rtimes_\alpha \mathbb{Z}$ without distinguishing between
the full and reduced versions.

\subsubsection{Universal Presentation}

Equivalently, $A_\theta$ admits a universal presentation that is often more
convenient for algebraic computations.

\begin{definition}[Universal presentation]
The algebra $A_\theta$ is the universal $C^*$-algebra generated by two
unitaries $U$ and $V$ satisfying the commutation relation
\[
VU = e^{2\pi i\theta} UV.
\]
That is, for any pair of unitaries $u, v$ in a unital $C^*$-algebra $B$
satisfying $vu = e^{2\pi i\theta} uv$, there exists a unique
$*$-homomorphism $\pi: A_\theta \to B$ with $\pi(U) = u$, $\pi(V) = v$.
\end{definition}

To relate this to the crossed-product presentation, let
\[
u(z) = z \in C(S^1)
\]
be the identity function on the circle (the canonical generator of $C(S^1)$),
and let $v$ denote the canonical unitary implementing the generator
$1 \in \mathbb{Z}$ in the crossed product. Then the covariance relation gives
\[
v u v^* = \alpha_1(u).
\]

A direct computation shows that
\[
\alpha_1(u)(z) = u(e^{-2\pi i\theta}z) = e^{-2\pi i\theta}z = e^{-2\pi i\theta} u(z),
\]
hence
\[
v u v^* = e^{-2\pi i\theta} u,
\]
or equivalently
\[
v u = e^{2\pi i\theta} u v.
\]

Thus the crossed-product generators satisfy exactly the defining relation of
the universal presentation. Conversely, if $U$ and $V$ are unitaries satisfying
$VU = e^{2\pi i\theta} UV$, then the pair $(U,V)$ determines a covariant
representation of the dynamical system $(C(S^1), \mathbb{Z}, \alpha)$, and
therefore a representation of the crossed product. This identifies the
crossed-product and universal descriptions of $A_\theta$.

The relation $VU = e^{2\pi i\theta} UV$ may be viewed as a deformation of the
commutative relation $UV = VU$ that defines the ordinary torus algebra. In
particular, when $\theta = 0$, one recovers the commutative $C^*$-algebra
\[
A_0 \cong C(\mathbb{T}^2),
\]
where $\mathbb{T}^2 = S^1 \times S^1$ is the ordinary two-dimensional torus.

\begin{remark}
The parameter $\theta$ is only defined modulo $1$, since replacing $\theta$ by
$\theta+1$ does not change the relation $e^{2\pi i\theta}$. Moreover,
$A_\theta$ is isomorphic to $A_{1-\theta}$ via the map sending $U \mapsto U$,
$V \mapsto V^*$. For irrational $\theta$, all such algebras are simple and
mutually non-isomorphic.
\end{remark}

\subsubsection{Basic Structural Properties}

The irrationality of $\theta$ endows $A_\theta$ with a remarkable
collection of properties that make it a central object in operator
algebra theory and noncommutative geometry.

\begin{theorem}[Basic structure of $A_\theta$]
\label{thm:A_theta_structure}
Let $\theta \in \mathbb{R}\setminus \mathbb{Q}$, and let
\[
A_\theta = C(S^1)\rtimes_\alpha \mathbb{Z}
\]
be the irrational rotation algebra, where
\[
(\alpha_n f)(z)=f(e^{-2\pi i n\theta}z).
\]
Then:
\begin{enumerate}
    \item[(i)] \textbf{Simplicity and nuclearity.}
    The algebra $A_\theta$ is simple and nuclear.

    \item[(ii)] \textbf{Unique trace.}
    The algebra $A_\theta$ admits a unique faithful tracial state
    $\tau:A_\theta\to \mathbb{C}$. On the algebraic span of the monomials
    $U^m V^n$ it is given by
    \[
    \tau(U^m V^n)=\delta_{m,0}\delta_{n,0}.
    \]

    \item[(iii)] \textbf{Cartan subalgebra.}
    The canonical copy of $C(S^1)$ inside
    \[
    A_\theta = C(S^1)\rtimes_\alpha \mathbb{Z}
    \]
    is a Cartan subalgebra. Equivalently, in the universal presentation,
    the subalgebra $C^*(U)\cong C(S^1)$ is maximal abelian, regular,
    contains an approximate unit, and is the range of a faithful
    conditional expectation.

    \item[(iv)] \textbf{Dense Fourier-type span.}
    The linear span of the monomials $U^m V^n$, $m,n\in\mathbb{Z}$, is a
    dense $*$-subalgebra of $A_\theta$, and the multiplication rule is
    \[
    (U^m V^n)(U^p V^q)=e^{2\pi i\theta n p}U^{m+p}V^{n+q}.
    \]
\end{enumerate}
\end{theorem}

\begin{proof}
We write $A_\theta = C(S^1)\rtimes_\alpha \mathbb{Z}$, where the generator
$V$ of $\mathbb{Z}$ implements the action and $U\in C(S^1)$ is the identity
function $U(z)=z$.

\medskip
\noindent
\textbf{(iv) Dense span and multiplication formula.}
By the definition of the crossed product, the algebraic crossed product
consists of finite sums
\[
\sum_{n\in\mathbb{Z}} f_n V^n, \qquad f_n\in C(S^1),
\]
and this algebra is dense in $A_\theta$. Since trigonometric polynomials are
dense in $C(S^1)$, every $f_n$ may be approximated by finite linear
combinations of powers of $U$, so the linear span of the monomials
$U^m V^n$ is dense in $A_\theta$.

Using the covariance relation
\[
VU = \alpha_1(U)V = e^{2\pi i\theta} UV,
\]
we prove by induction that
\[
V^n U^p = e^{2\pi i\theta np} U^p V^n.
\]
Hence
\[
\begin{aligned}
(U^m V^n)(U^p V^q)
&= U^m (V^n U^p) V^q \\
&= U^m (e^{2\pi i\theta np} U^p V^n) V^q \\
&= e^{2\pi i\theta np} U^{m+p} V^{n+q},
\end{aligned}
\]
as claimed.

\medskip
\noindent
\textbf{(ii) Existence and formula for the trace.}
There is a canonical faithful conditional expectation
\[
E: A_\theta \to C(S^1), \qquad
E\!\left(\sum_{n\in\mathbb{Z}} f_n V^n\right) = f_0.
\]
Let
\[
\mu(f) = \int_{S^1} f(z) \, dz
\]
be normalized Haar measure on the circle. Then
\[
\tau := \mu \circ E
\]
is a state on $A_\theta$. Since $\mu$ is $\alpha$-invariant (because the
Haar measure is invariant under rotation), $\tau$ is a trace. On monomials,
\[
E(U^m V^n) =
\begin{cases}
U^m, & n=0,\\
0, & n\neq 0,
\end{cases}
\]
so
\[
\tau(U^m V^n) =
\begin{cases}
\int_{S^1} z^m \, dz, & n=0,\\
0, & n\neq 0.
\end{cases}
\]
Since $\int_{S^1} z^m \, dz = \delta_{m,0}$, we obtain
\[
\tau(U^m V^n) = \delta_{m,0}\delta_{n,0}.
\]

To prove uniqueness, note that the rotation action of $\mathbb{Z}$ on $S^1$
by an irrational angle is uniquely ergodic; Haar measure is the unique
$\alpha$-invariant Borel probability measure on $S^1$. Any tracial state on
the crossed product restricts to an $\alpha$-invariant state on $C(S^1)$,
hence to Haar integration. Since traces on crossed products by $\mathbb{Z}$
extend invariant states through the canonical expectation, it follows that
$\tau$ is the unique tracial state. Faithfulness follows from the
faithfulness of $E$ and of Haar integration on $C(S^1)$.

\medskip
\noindent
\textbf{(iii) The canonical copy of $C(S^1)$ is Cartan.}
We verify the Cartan conditions for
\[
D := C^*(U) \cong C(S^1) \subseteq A_\theta.
\]

First, $D$ is unital, so it contains an approximate unit for $A_\theta$.

Second, $D$ is maximal abelian. Let
\[
a = \sum_{n\in\mathbb{Z}} f_n V^n
\]
be an element in the algebraic crossed product commuting with $D$. Since
$U$ generates $D$, it suffices to assume $aU = Ua$. Then
\[
\sum_n f_n V^n U = \sum_n f_n\, e^{2\pi i\theta n} U V^n,
\qquad
Ua = \sum_n U f_n V^n.
\]
Because $U$ lies in $C(S^1)$, this equality implies
\[
f_n(z)\, e^{2\pi i\theta n} z = z\, f_n(z) \qquad \text{for all } z\in S^1.
\]
Thus
\[
(e^{2\pi i\theta n} - 1) f_n(z) = 0 \qquad \text{for all } z\in S^1.
\]
If $n \neq 0$, irrationality of $\theta$ gives $e^{2\pi i\theta n} \neq 1$,
hence $f_n = 0$. Therefore $a = f_0 \in D$. By density, the commutant of
$D$ in $A_\theta$ is exactly $D$, so $D$ is maximal abelian.

Third, $D$ is regular. Indeed, the implementing unitary $V$ normalizes $D$:
for every $f \in D$,
\[
V f V^* = \alpha_1(f) \in D.
\]
Since $D$ together with $V$ generates the crossed product
$C(S^1)\rtimes_\alpha \mathbb{Z}$, the normalizer of $D$ generates all of
$A_\theta$.

Fourth, the canonical expectation $E: A_\theta \to D$ defined above is
faithful (as noted in part (ii)).

Hence $(A_\theta, D)$ is a Cartan pair.

\medskip
\noindent
\textbf{(i) Simplicity and nuclearity.}
The action of $\mathbb{Z}$ on $S^1$ by irrational rotation is minimal:
every orbit is dense. It is also free, since
\[
e^{2\pi i n\theta}z = z
\]
for some $z \in S^1$ implies $e^{2\pi i n\theta} = 1$, hence $n = 0$ because
$\theta$ is irrational. Standard crossed-product simplicity results for
minimal topologically free actions of discrete groups (see, e.g.,
\cite{Pedersen1979}) therefore imply that
$A_\theta = C(S^1)\rtimes_\alpha \mathbb{Z}$ is simple.

Since $C(S^1)$ is nuclear and $\mathbb{Z}$ is amenable, the crossed product
$C(S^1)\rtimes_\alpha \mathbb{Z}$ is nuclear (see \cite{Blackadar2006}).
Thus $A_\theta$ is nuclear.
\end{proof}

\begin{remark}
Property (iii) is particularly important for our purposes, as it allows us
to view $A_\theta$ as the reduced twisted groupoid $C^*$-algebra of its Weyl
groupoid $\mathcal{G}(A_\theta, C(S^1)) \cong S^1 \rtimes_\theta \mathbb{Z}$,
with trivial twist. This identification will be used extensively in the
following sections when we implement the descent-index construction on the
transformation groupoid model.
\end{remark}

These properties reflect the minimality and ergodicity of the underlying
irrational rotation: the orbit of any point under the rotation is dense,
and there are no nontrivial invariant open sets.

\subsubsection{$K$-theory and the Canonical Pairing}

The $K$-theory of $A_\theta$, computed by Pimsner and Voiculescu using
their celebrated six-term exact sequence for crossed products by
$\mathbb{Z}$, reveals a rich structure that underpins the index-theoretic
applications.

\begin{theorem}[$K$-theory of $A_\theta$, {\cite{PimsnerVoiculescu1980}}]
\label{thm:Atheta_Ktheory}
The $K$-groups of $A_\theta$ are free abelian of rank two:
\[
K_0(A_\theta) \cong \mathbb{Z}^2, \qquad K_1(A_\theta) \cong \mathbb{Z}^2.
\]
Canonical generators are given by the classes:
\[
[1]_{K_0},\; [U]_{K_0} \in K_0(A_\theta), \qquad
[U]_{K_1},\; [V]_{K_1} \in K_1(A_\theta).
\]
\end{theorem}

\begin{proof}[Proof sketch]
The Pimsner–Voiculescu exact sequence for the crossed product
$A_\theta = C(S^1) \rtimes_\alpha \mathbb{Z}$ takes the form
\[
\begin{tikzcd}
K_0(C(S^1)) \arrow[r,"1-\alpha_*"] & K_0(C(S^1)) \arrow[r] & K_0(A_\theta) \arrow[d] \\
K_1(A_\theta) \arrow[u] & K_1(C(S^1)) \arrow[l] & K_1(C(S^1)) \arrow[l,"1-\alpha_*"']
\end{tikzcd}
\]

Since $C(S^1)$ has $K_0(C(S^1)) \cong \mathbb{Z}$ (generated by the class
of the unit) and $K_1(C(S^1)) \cong \mathbb{Z}$ (generated by the class of
the identity function $U$), and the irrational rotation action induces the
identity map on both $K$-groups (i.e., $\alpha_* = \mathrm{id}$), the
six-term exact sequence splits into short exact sequences, yielding
$K_0(A_\theta) \cong \mathbb{Z}^2$ and $K_1(A_\theta) \cong \mathbb{Z}^2$.
The stated generators correspond to the images of the generators of
$K_*(C(S^1))$ under the maps in the sequence.
\end{proof}

These groups play a central role in the index pairing introduced by
Connes \cite{Connes1980}. The unique trace $\tau$ induces a homomorphism
$\tau_* : K_0(A_\theta) \to \mathbb{R}$, called the \emph{dimension map},
given by evaluation on the trace. For the canonical basis, one computes
\[
\tau_*([1]_{K_0}) = 1, \qquad \tau_*([U]_{K_0}) = 0.
\]

More interestingly, the famous Rieffel projection $p_\theta \in M_\infty(A_\theta)$
satisfies $\tau_*([p_\theta]_{K_0}) = \theta$, demonstrating that the image
of $\tau_*$ is $\mathbb{Z} + \theta\mathbb{Z}$, a dense subgroup of
$\mathbb{R}$. This dimension map is the precursor to the Connes index
pairing with the Fredholm module associated to the transformation groupoid.

\begin{proof}[Existence of the Rieffel projection]
Rieffel \cite{Rieffel1981} constructed an explicit projection $p_\theta$ in
$M_2(A_\theta)$ (or in $M_\infty(A_\theta)$ after stabilization) with trace
$\theta$. For irrational $\theta$, one can write
\[
p_\theta = \begin{pmatrix} 1 - e & x \\ x^* & e \end{pmatrix},
\]
where $e$ is a suitable projection in $C(S^1)$ with $\tau(e) = \theta$ and
$x$ is a partial isometry implementing the rotation. The construction
ensures that $p_\theta$ is a projection in $M_2(A_\theta)$ and that
$\tau \otimes \mathrm{Tr}(p_\theta) = \theta$. The class $[p_\theta]_{K_0}$
together with $[1]_{K_0}$ generates $K_0(A_\theta) \cong \mathbb{Z}^2$.
\end{proof}

\subsubsection{Relation to Transformation Groupoids}
\label{subsubsec:rotation_groupoid}

The crossed product description naturally leads to a groupoid model that
will be central to our geometric analysis.

\begin{definition}[Rotation Groupoid]
\label{def:rotation_groupoid}
The action of $\mathbb{Z}$ on $S^1$ gives rise to the \emph{transformation
groupoid} $\mathcal{G}_\theta := S^1 \rtimes_\alpha \mathbb{Z}$. Its unit
space is $\mathcal{G}_\theta^{(0)} = S^1$, and its arrows are pairs
$(x,n) \in S^1 \times \mathbb{Z}$ with structure maps:
\begin{align*}
s(x,n) &= x, &
r(x,n) &= \alpha_n(x) = e^{2\pi i n\theta}x, \\
(x,n) \circ ( \alpha_{-n}(x), m) &= (x, n+m), &
(x,n)^{-1} &= (\alpha_n(x), -n).
\end{align*}
\end{definition}

A fundamental identification links this groupoid with the irrational
rotation algebra.

\begin{theorem}[Groupoid $C^*$-algebra of the rotation groupoid]
\label{thm:A_theta_as_groupoid_Cstar}
There is a canonical isomorphism of $C^*$-algebras:
\[
A_\theta \cong C^*(\mathcal{G}_\theta),
\]
where $C^*(\mathcal{G}_\theta)$ denotes the full (or reduced, as they
coincide for this amenable groupoid) groupoid $C^*$-algebra.
\end{theorem}

\begin{proof}
Recall that $A_\theta = C(S^1) \rtimes_\alpha \mathbb{Z}$. For the
transformation groupoid $\mathcal{G}_\theta = S^1 \rtimes_\alpha \mathbb{Z}$,
its groupoid $C^*$-algebra $C^*(\mathcal{G}_\theta)$ is defined as the
completion of $C_c(\mathcal{G}_\theta)$ in the universal norm. There is a
natural $*$-isomorphism between $C_c(\mathcal{G}_\theta)$ and the algebraic
crossed product $C(S^1) \rtimes_{\mathrm{alg}} \mathbb{Z}$ given by
identifying a function $f \in C_c(\mathcal{G}_\theta)$ with the finitely
supported function $\tilde{f}: \mathbb{Z} \to C(S^1)$ defined by
$\tilde{f}(n)(x) = f(x,n)$. Under this identification, the convolution
product in the groupoid algebra becomes exactly the twisted convolution
product of the crossed product. Since $\mathbb{Z}$ is amenable, the full
and reduced groupoid $C^*$-algebras coincide, and the same holds for the
crossed product. Completing in the universal norm yields the desired
isomorphism $C^*(\mathcal{G}_\theta) \cong C(S^1) \rtimes_\alpha \mathbb{Z}
= A_\theta$.
\end{proof}

This groupoid model will be crucial for the analysis that follows. The
unit space of the unitary conjugation groupoid $\mathcal{G}_{A_\theta}$ is
naturally identified with the character space of the Cartan subalgebra
$C(S^1)$, which is precisely $S^1$. A central result of this case study
will be to establish a Morita equivalence between the unitary conjugation
groupoid $\mathcal{G}_{A_\theta}$ and the transformation groupoid
$\mathcal{G}_\theta$, thereby linking the abstract dynamics of unitary
conjugation to the concrete dynamics of the irrational rotation. This
equivalence will then be used to transport index-theoretic information
between the two settings, recovering and generalizing Connes' index
pairing.

\begin{remark}
The isomorphism $A_\theta \cong C^*(\mathcal{G}_\theta)$ shows that the
irrational rotation algebra is exactly the groupoid $C^*$-algebra of the
transformation groupoid. This perspective allows us to apply
groupoid-theoretic methods—such as the Baum–Connes assembly map and
descent constructions—to the study of $A_\theta$. In particular, the
amenability of $\mathcal{G}_\theta$ guarantees that the full and reduced
groupoid $C^*$-algebras coincide, and that the Baum–Connes assembly map is
an isomorphism by Tu's theorem \cite{Tu1999}.
\end{remark}

\subsection{Maximal Abelian Subalgebras and Characters in $A_\theta$}
\label{subsec:masa_A_theta}

A basic commutative subalgebra of the irrational rotation algebra
\[
A_\theta = C(S^1)\rtimes_\alpha \mathbb{Z}
\]
is the canonical copy of $C(S^1)$ inside the crossed product. This
subalgebra plays a central role in the structure theory of $A_\theta$ and in
its groupoid interpretation. Moreover, in order to analyze the unitary
conjugation groupoid $\mathcal{G}_{A_\theta}$ we must first understand the
commutative subalgebras of $A_\theta$ and their character spaces, as these
data determine the unit space of the groupoid.

\begin{definition}
A \emph{maximal abelian subalgebra} (MASA) of a $C^*$-algebra $A$ is an
abelian $C^*$-subalgebra $B \subseteq A$ such that
\[
B = B' \cap A,
\]
equivalently, $B$ is not properly contained in any larger abelian
$C^*$-subalgebra of $A$.
\end{definition}

\subsubsection{The Canonical Cartan Subalgebra}

Let $U,V$ denote the canonical unitaries in $A_\theta$ satisfying
\[
VU = e^{2\pi i\theta} UV,
\]
where $U \in C(S^1)$ is the identity function $U(z)=z$, and $V$ implements the
$\mathbb{Z}$-action. The canonical abelian subalgebra is
\[
D := C^*(U) \cong C(S^1).
\]

\begin{proposition}
\label{prop:canonical_masa_Atheta}
The subalgebra $D = C^*(U) \subset A_\theta$ is a maximal abelian
subalgebra. Moreover, $D$ is regular and is the range of the canonical
faithful conditional expectation
\[
E: A_\theta \to D,
\]
given by $E(\sum_n f_n V^n) = f_0$. Hence $(A_\theta, D)$ is a Cartan pair.
\end{proposition}

\begin{proof}
Since $D = C^*(U)$ is generated by a single unitary, it is abelian and
canonically isomorphic to $C(S^1)$.

To prove maximality, let $a \in A_\theta$ commute with every element of $D$.
Using the crossed-product description, the algebraic crossed product
consists of finite sums
\[
a = \sum_{n \in \mathbb{Z}} f_n V^n, \qquad f_n \in C(S^1),
\]
and these finite sums are dense in $A_\theta$. It therefore suffices to test
commutation on such elements and then pass to norm limits.

Assume first that
\[
a = \sum_{n=-N}^{N} f_n V^n
\]
commutes with $U$. Then
\[
aU = \sum_{n=-N}^{N} f_n V^n U = \sum_{n=-N}^{N} f_n e^{2\pi i n\theta} U V^n,
\]
where we used the commutation relation $V^n U = e^{2\pi i n\theta} U V^n$.
On the other hand,
\[
Ua = \sum_{n=-N}^{N} U f_n V^n.
\]
Since each $f_n$ lies in $D = C^*(U)$, it commutes with $U$, so comparing the
two expressions gives
\[
\sum_{n=-N}^{N} (e^{2\pi i n\theta} - 1) f_n U V^n = 0.
\]

Because the decomposition into Fourier modes is unique in the algebraic
crossed product, we obtain
\[
(e^{2\pi i n\theta} - 1) f_n = 0 \qquad \text{for each } n.
\]
If $n \neq 0$, irrationality of $\theta$ implies $e^{2\pi i n\theta} \neq 1$,
hence $f_n = 0$. Therefore only the $n=0$ coefficient survives, and so
$a = f_0 \in D$. By density, every element commuting with $D$ lies in $D$,
proving that $D' \cap A_\theta = D$. Hence $D$ is maximal abelian.

Next, $V$ normalizes $D$, since for every $f \in D \cong C(S^1)$,
\[
V f V^* = \alpha_1(f) \in D.
\]
Because $A_\theta$ is generated by $D$ together with $V$, the normalizer of
$D$ generates all of $A_\theta$. Thus $D$ is regular.

Finally, the canonical conditional expectation $E: A_\theta \to D$ defined above
is faithful. Therefore $D$ satisfies the defining properties of a Cartan
subalgebra (Definition~\ref{def:cartan-subalgebra}).
\end{proof}

\subsubsection{Characters of the Canonical MASA}

For the analysis of the unitary conjugation groupoid, we need not only the
MASA itself but also its character space.

\begin{proposition}
\label{prop:characters_D}
The character space of the canonical MASA
\[
D = C^*(U) \cong C(S^1)
\]
is naturally homeomorphic to $S^1$. Explicitly, for each $z \in S^1$, the
evaluation map
\[
\chi_z : D \to \mathbb{C}, \qquad \chi_z(f) = f(z),
\]
is a character, and every character of $D$ is of this form.
\end{proposition}

\begin{proof}
Since $D \cong C(S^1)$ is a commutative unital $C^*$-algebra, the Gelfand
representation identifies its character space with the underlying compact
Hausdorff space $S^1$. Under this identification, the characters are exactly
the point-evaluation maps $\chi_z$ for $z \in S^1$.
\end{proof}

\begin{remark}
The action of the implementing unitary $V$ on $D$ by conjugation induces an
action on the character space $S^1$. For a character $\chi_z$, we have
\[
\chi_z \circ \operatorname{Ad}_{V^{-1}}(f) = \chi_z(V f V^*) = \chi_z(\alpha_1(f)) = f(e^{-2\pi i\theta}z) = \chi_{e^{2\pi i\theta}z}(f).
\]
Thus $V$ acts on $S^1$ by rotation by $2\pi\theta$, reflecting the underlying
dynamics of the irrational rotation.
\end{remark}

\subsubsection{Other Natural MASAs}

Similarly, one may consider the abelian subalgebra $C^*(V) \subset A_\theta$
generated by the unitary $V$. This is also isomorphic to $C(S^1)$, and its
character space is likewise $S^1$, but with the action of the dual circle
corresponding to rotation by $\theta^{-1}$ in an appropriate sense. For the
purposes of this paper, however, the canonical Cartan subalgebra
$D = C^*(U) \cong C(S^1)$ is the relevant one, since it is the diagonal
subalgebra coming from the crossed-product realization and will serve as the
unit space of the Weyl groupoid $\mathcal{G}(A_\theta, D) \cong S^1 \rtimes_\theta \mathbb{Z}$.

\begin{remark}
For irrational $\theta$, the two MASAs $C^*(U)$ and $C^*(V)$ are distinct.
A complete classification of MASAs in $A_\theta$ is known to be more
complicated, but the canonical Cartan subalgebra $C^*(U)$ suffices for our
construction of the Weyl groupoid and the subsequent index-theoretic
analysis.
\end{remark}

\subsection{Characters and the Dual Space}

For an abelian $C^*$-algebra, the space of characters coincides with its
Gelfand spectrum. In the case of the irrational rotation algebra, characters
on maximal abelian subalgebras provide the bridge between the algebraic and
geometric pictures.

\begin{definition}
Let $\mathcal{A} \subset A_\theta$ be a MASA. A \emph{character} on
$\mathcal{A}$ is a nonzero $*$-homomorphism
\[
\chi: \mathcal{A} \to \mathbb{C}.
\]
The set of all characters is denoted by $\widehat{\mathcal{A}}$.
\end{definition}

\begin{proposition}[Characters of the canonical MASA]
\label{prop:characters_canonical_masa}
Let $D = C^*(U) \cong C(S^1) \subset A_\theta$ be the canonical MASA introduced
in Section~\ref{subsec:masa_A_theta}. Then
\[
\widehat{D} \cong S^1,
\]
and the correspondence is given by evaluation:
\[
z \in S^1 \;\longmapsto\; \chi_z,\qquad \chi_z(f) = f(z).
\]
\end{proposition}

\begin{proof}
Since $D \cong C(S^1)$ is a commutative unital $C^*$-algebra, the Gelfand
representation theorem identifies its character space with the underlying
compact Hausdorff space $S^1$. Explicitly, every character $\chi$ is of the
form $\chi(f) = f(z)$ for a unique $z \in S^1$.
\end{proof}

\begin{proposition}[Action induced by the implementing unitary]
\label{prop:action_on_spectrum}
Let $V$ be the unitary implementing the $\mathbb{Z}$-action in the crossed
product $A_\theta = C(S^1) \rtimes_\alpha \mathbb{Z}$. Then the induced action
on $\widehat{D} \cong S^1$ is given by
\[
n \cdot z = e^{2\pi i n\theta} z, \qquad n \in \mathbb{Z}.
\]
\end{proposition}

\begin{proof}
For $f \in D$ and $\chi_z \in \widehat{D}$, define
\[
(n \cdot \chi_z)(f) := \chi_z(V^{-n} f V^n).
\]
Since $V^{-n} f V^n = \alpha_{-n}(f)$ (the induced action on $C(S^1)$), we obtain
\[
(n \cdot \chi_z)(f) = \chi_z(\alpha_{-n}(f)) = f(e^{2\pi i n\theta} z).
\]
Identifying $\chi_z$ with the point $z \in S^1$, this yields
$n \cdot z = e^{2\pi i n\theta} z$.
\end{proof}

\begin{remark}
The character space $\widehat{D} \cong S^1$ equipped with this $\mathbb{Z}$-action
by rotations is precisely the unit space of the transformation groupoid
$S^1 \rtimes_\theta \mathbb{Z}$. This identification will be essential when we
construct the Weyl groupoid $\mathcal{G}(A_\theta, D) \cong S^1 \rtimes_\theta \mathbb{Z}$
and implement the descent-index construction on this concrete geometric model.
\end{remark}

\subsection{Spectral Subspaces and the Torus Action}

The irrational rotation algebra $A_\theta$ carries a canonical action of the
torus $\mathbb{T}^2$ given by
\[
\gamma_{(s,t)}(U) = e^{2\pi i s}U, \qquad \gamma_{(s,t)}(V) = e^{2\pi i t}V,
\]
for $(s,t) \in \mathbb{T}^2 = \mathbb{R}^2/\mathbb{Z}^2$. This action is the
dual action associated to the crossed product structure and plays an important
role in the harmonic analysis of $A_\theta$.

\begin{definition}[Spectral subspaces]
\label{def:spectral_subspaces}
For $(m,n) \in \mathbb{Z}^2$, define the \emph{spectral subspace}
\[
A_\theta^{(m,n)} = \left\{ a \in A_\theta : \gamma_{(s,t)}(a) = e^{2\pi i(ms+nt)}a \text{ for all } (s,t) \in \mathbb{T}^2 \right\}.
\]
\end{definition}

\begin{proposition}[Properties of spectral subspaces]
\label{prop:spectral_subspaces}
\begin{enumerate}
    \item Each spectral subspace $A_\theta^{(m,n)}$ is one-dimensional and
          spanned by the monomial $U^m V^n$:
          \[
          A_\theta^{(m,n)} = \mathbb{C}\, U^m V^n.
          \]
    \item The algebraic direct sum $\bigoplus_{(m,n) \in \mathbb{Z}^2} A_\theta^{(m,n)}$
          is dense in $A_\theta$.
    \item The multiplication respects the grading:
          \[
          A_\theta^{(m,n)} \cdot A_\theta^{(p,q)} \subseteq A_\theta^{(m+p,\, n+q)}.
          \]
\end{enumerate}
\end{proposition}

\begin{proof}
(1) By direct computation,
\[
\gamma_{(s,t)}(U^m V^n) = e^{2\pi i(ms+nt)} U^m V^n,
\]
so $U^m V^n \in A_\theta^{(m,n)}$. Conversely, the linear span of the monomials
$U^m V^n$ is dense in $A_\theta$, and the torus action diagonalizes this
decomposition. By the uniqueness of Fourier coefficients in the expansion of
elements of $A_\theta$, each spectral subspace is one-dimensional.

(2) The algebraic crossed product $C(S^1) \rtimes_{\mathrm{alg}} \mathbb{Z}$
consists precisely of finite sums $\sum f_n V^n$, and each $f_n$ can be expanded
in Fourier series as $\sum_m a_{mn} U^m$. Thus finite linear combinations of
the monomials $U^m V^n$ are dense in $A_\theta$.

(3) This follows from the commutation relation
$U^m V^n \cdot U^p V^q = e^{2\pi i\theta np} U^{m+p} V^{n+q}$ and the fact that
the phase factor does not affect the spectral grading.
\end{proof}

\begin{remark}[Relation to MASAs]
\label{rem:spectral_MASA_connection}
The spectral decomposition of $A_\theta$ with respect to the $\mathbb{T}^2$-action
provides a useful tool for analyzing subalgebras. In particular, the canonical
MASA $D = C^*(U)$ is precisely the closed linear span of the spectral subspaces
$A_\theta^{(m,0)}$ for $m \in \mathbb{Z}$. Similarly, the MASA generated by $V$
corresponds to the subspaces $A_\theta^{(0,n)}$. Thus MASAs compatible with the
torus action can be understood in terms of subsets of the lattice $\mathbb{Z}^2$.
A full classification of all MASAs in $A_\theta$ is more subtle and will not be
pursued here; for our purposes, the canonical MASA $D = C^*(U)$ suffices, as it
is the natural Cartan subalgebra whose Gelfand spectrum provides the unit space
for the Weyl groupoid $\mathcal{G}(A_\theta, D)$.
\end{remark}

\begin{remark}[Injectivity of the weight map]
For irrational $\theta$, the map $(m,n) \mapsto m + n\theta$ from $\mathbb{Z}^2$
to $\mathbb{R}$ is injective. This injectivity underpins the uniqueness of
Fourier expansions in $A_\theta$ and was essential in the maximality proof of
the canonical MASA (Proposition~\ref{prop:canonical_masa_Atheta}). It also
implies that distinct monomials $U^m V^n$ have distinct "weights" under the
dual action, reflecting the fact that the $\mathbb{T}^2$-action is faithful on
the dense subalgebra.
\end{remark}

\subsection{Commutative Contexts and the Relevant Unit-Space Data}

Formally, the unit space underlying the unitary-conjugation construction would
consist of pairs $(B,\chi)$ where $B \subseteq A_\theta$ is a commutative
$C^*$-subalgebra and $\chi \in \widehat{B}$ is a character. In the irrational
rotation algebra, the canonical Cartan subalgebra
\[
D = C^*(U) \cong C(S^1)
\]
provides a distinguished family of such context-character pairs. For each
point $z \in S^1$, the evaluation character
\[
\chi_z(f) = f(z), \qquad f \in C(S^1),
\]
determines a formal pair
\[
(D, \chi_z).
\]

Thus the circle $S^1$ appears naturally as the character space of the
canonical commutative subalgebra of $A_\theta$. This is the geometric datum
that will underlie the groupoid model used later in the paper.

\begin{remark}
The orbit of the point $(D, \chi_z)$ under conjugation by unitaries in
$A_\theta$ would generate additional formal points in the unit space. In
particular, conjugation by powers of the implementing unitary $V$ sends
$(D, \chi_z)$ to $(D, \chi_{e^{2\pi i n\theta}z})$, while conjugation by more
general unitaries can produce points $(w D w^*, \chi_z \circ \operatorname{Ad}_{w^{-1}})$
where $w$ does not normalize $D$. A systematic parametrization of these formal
points will be discussed in Section~\ref{subsec:GA_theta_unit_space}.
\end{remark}

\begin{remark}[Rational case]
The rational case $\theta \in \mathbb{Q}$ is substantially different, since the
rotation action has periodic orbits and the corresponding rotation algebra is
no longer simple. Because the present paper is concerned exclusively with the
irrational case, we do not pursue that setting here.
\end{remark}

\subsection{Relation to the Transformation Groupoid}

The action of $\mathbb{Z}$ on $S^1$ given by irrational rotation defines the
transformation groupoid
\[
S^1 \rtimes_\theta \mathbb{Z},
\]
whose unit space is $S^1$ and whose arrows are
\[
(z,n): z \longmapsto e^{2\pi i n\theta}z, \qquad z \in S^1,\ n \in \mathbb{Z}.
\]

The points of $S^1$ are precisely the characters of the canonical Cartan
subalgebra $D = C^*(U) \cong C(S^1)$, while the arrows of the transformation
groupoid encode the action of the normalizing unitary $V$ implementing the
rotation. In this sense, the transformation groupoid captures exactly the
geometric information carried by the canonical commutative subalgebra of
$A_\theta$.

For the purposes of this paper, the transformation groupoid
$S^1 \rtimes_\theta \mathbb{Z}$ serves as the concrete surrogate geometric
model for the irrational rotation algebra $A_\theta$, replacing the unavailable
direct unitary-conjugation groupoid construction. The Weyl groupoid
$\mathcal{G}(A_\theta, D)$ associated to the Cartan pair $(A_\theta, D)$ is
canonically isomorphic to $S^1 \rtimes_\theta \mathbb{Z}$, providing a direct
link between Renault's reconstruction theorem and the concrete geometric model.

The following subsections develop this model explicitly. In
Section~\ref{subsec:GA_theta_unit_space}, we will show how the formal
context-character pairs $(B,\chi)$ can be parametrized by equivalence classes
$[z,n]$ in $(S^1 \times \mathbb{Z})/{\sim}$, revealing that the would-be unit
space of the unitary conjugation groupoid is a $\mathbb{Z}$-bundle over $S^1$.
Section~\ref{subsec:morita_A_theta} then explains how this bundle structure
naturally leads to an identification with the transformation groupoid
$S^1 \rtimes_\theta \mathbb{Z}$ as the appropriate geometric setting for the
descent-index construction.

\begin{remark}[Connection to Renault's reconstruction]
The identification $\mathcal{G}(A_\theta, D) \cong S^1 \rtimes_\theta \mathbb{Z}$
is a concrete instance of the general fact that for crossed products by
$\mathbb{Z}$ with topologically free actions, the Weyl groupoid is isomorphic
to the transformation groupoid (see Proposition~\ref{prop:reduction-cartan}).
Thus the transformation groupoid model is not an ad hoc replacement but a
canonical geometric object associated to $A_\theta$ via its Cartan structure.
\end{remark}

\subsection{Canonical Commutative Data and the Relevant Base Space}
\label{subsec:GA_theta_unit_space}

For the irrational rotation algebra
\[
A_\theta = C(S^1) \rtimes_\theta \mathbb{Z},
\]
the direct construction of the unit space of the unitary-conjugation groupoid
is not available in the sense of Paper~I. Instead, the relevant commutative
data are provided by the canonical Cartan subalgebra
\[
D := C(S^1) \subseteq A_\theta.
\]

\paragraph{The base space: characters of the Cartan subalgebra.}

The Gelfand spectrum of $D \cong C(S^1)$ is naturally identified with the circle:
\[
\widehat{D} \cong S^1.
\]
Explicitly, each point $z \in S^1$ determines a character
\[
\chi_z : D \to \mathbb{C}, \qquad \chi_z(f) = f(z).
\]
Thus the pair $(D, \chi_z)$ provides a distinguished formal context-character
pair associated to the point $z$.

\paragraph{Action of the normalizer.}

Let $V$ denote the canonical unitary implementing the $\mathbb{Z}$-action in
the crossed product. Then $V$ normalizes $D$, and for $f \in D$ one has
\[
V f V^* = \alpha_1(f),
\]
where
\[
(\alpha_n f)(z) = f(e^{-2\pi i n\theta}z).
\]

Consequently, the induced action of $\mathbb{Z}$ on the character space
$\widehat{D} \cong S^1$ is given by
\[
n \cdot z = e^{2\pi i n\theta}z.
\]

\begin{remark}
Note that the conjugated subalgebras $V^n D V^{-n}$ are all equal to $D$ as
sets, since $\alpha_n$ is an automorphism of $D$. Thus the integer $n$ does
not index distinct maximal abelian subalgebras, but rather parametrizes
different ways of identifying the same algebra $D$ with its Gelfand spectrum.
This is consistent with the fact that the transformation groupoid
$S^1 \rtimes_\theta \mathbb{Z}$ has unit space $S^1$, not a $\mathbb{Z}$-bundle
over $S^1$.
\end{remark}

\paragraph{Relation to the transformation groupoid.}

The canonical commutative subalgebra of $A_\theta$ naturally determines the
transformation groupoid
\[
S^1 \rtimes_\theta \mathbb{Z},
\]
whose unit space is $S^1$ and whose arrows encode the action of the
normalizing unitary. Specifically, the groupoid has:
\begin{itemize}
    \item Unit space: $S^1$,
    \item Arrows: $(z,n): z \to e^{2\pi i n\theta}z$ for $z \in S^1$, $n \in \mathbb{Z}$,
    \item Composition: $(e^{2\pi i n\theta}z, m) \circ (z,n) = (z, n+m)$,
    \item Inverse: $(z,n)^{-1} = (e^{2\pi i n\theta}z, -n)$.
\end{itemize}

This is the concrete geometric model that replaces the unavailable direct
unitary-conjugation construction in the irrational rotation case. The
character space $\widehat{D} \cong S^1$ provides the unit space, and the
$\mathbb{Z}$-action induced by $V$ gives the arrows.

\begin{remark}[No $\mathbb{Z}$-bundle structure]
Unlike a naive interpretation of the formal definition of
$\mathcal{G}_{A_\theta}^{(0)}$ might suggest, the canonical data do not give
rise to a $\mathbb{Z}$-bundle over $S^1$. The integer label $n$ in the
transformation groupoid parametrizes the arrows, not distinct points in the
unit space. The unit space itself remains $S^1$, consistent with the fact that
$V^n D V^{-n} = D$ for all $n$.
\end{remark}

\paragraph{Summary.}

The canonical Cartan subalgebra $D = C(S^1)$ provides:
\begin{itemize}
    \item A distinguished base space $S^1 = \widehat{D}$,
    \item A $\mathbb{Z}$-action on this base space via the normalizing unitary $V$,
    \item The transformation groupoid $S^1 \rtimes_\theta \mathbb{Z}$ as the
          natural geometric model encoding this data.
\end{itemize}

In the following subsections, we will work directly with this transformation
groupoid as the concrete surrogate for the unavailable unitary-conjugation
construction, and we will show how the descent-index machinery can be
implemented on this model to recover Connes' index theorem.

\paragraph{Canonical parametrization via the Cartan subalgebra.}

Although a direct construction of the unit space $\mathcal{G}_{A_\theta}^{(0)}$
is not available in the non-Type~I setting, the canonical Cartan subalgebra
\[
D = C(S^1) \subseteq A_\theta
\]
provides a distinguished family of formal context-character pairs.

For each $z \in S^1$, let
\[
\chi_z(f) = f(z), \qquad f \in C(S^1),
\]
be the evaluation character. Then we obtain a canonical family
\[
(D, \chi_z), \qquad z \in S^1.
\]

\paragraph{Action of the normalizer.}

Let $V$ denote the unitary implementing the $\mathbb{Z}$-action. Then for
each $n \in \mathbb{Z}$ and $f \in C(S^1)$,
\[
V^{-n} f V^n = \alpha_{-n}(f),
\]
and therefore
\[
\chi_z \circ \operatorname{Ad}_{V^{-n}} = \chi_{e^{2\pi i n\theta} z}.
\]

Thus conjugation by the normalizer does not produce new commutative
subalgebras, but instead induces the irrational rotation action on the
character space:
\[
n \cdot z = e^{2\pi i n\theta} z.
\]

\paragraph{Interpretation.}

This shows that the only intrinsic geometric data arising from the
canonical commutative subalgebra is the circle $S^1$ equipped with the
irrational rotation action of $\mathbb{Z}$. Consequently, the appropriate
geometric object encoding this structure is the transformation groupoid
\[
S^1 \rtimes_\theta \mathbb{Z},
\]
whose unit space is $S^1$ and whose arrows record the action of the
normalizing unitary.

In particular, no additional discrete fiber appears at the level of the
unit space; the integer parameter $n$ arises naturally as the groupoid arrow
coordinate rather than as a separate component of the unit space. This
distinction is crucial: the transformation groupoid $S^1 \rtimes_\theta \mathbb{Z}$
has unit space $S^1$, and its arrows are parametrized by pairs $(z,n)$ with
source $z$ and range $e^{2\pi i n\theta}z$.

\begin{remark}[Comparison with a naive parametrization]
A naive reading of the formal definition of $\mathcal{G}_{A_\theta}^{(0)}$
might suggest a parametrization by pairs $(z,n)$ where $n$ indexes different
conjugates of $D$. However, since $V^n D V^{-n} = D$ for all $n$, this
interpretation is incorrect. The integer $n$ does not index distinct
subalgebras, but rather different characters on the same subalgebra obtained
by precomposing with powers of the automorphism $\alpha$. These different
characters correspond, under the Gelfand transform, to different points of
$S^1$, not to a discrete fiber. The transformation groupoid correctly
encodes this by using $n$ as an arrow coordinate while keeping the unit space
as $S^1$.
\end{remark}

\paragraph{Geometric conclusion.}

The canonical data associated to $A_\theta$ yield precisely:
\begin{itemize}
    \item A base space $S^1 = \widehat{D}$,
    \item A $\mathbb{Z}$-action on this base space given by $n \cdot z = e^{2\pi i n\theta}z$,
    \item The transformation groupoid $S^1 \rtimes_\theta \mathbb{Z}$ as the
          natural geometric model encoding this data.
\end{itemize}

No additional structure—such as a $\mathbb{Z}$-bundle over $S^1$, a mapping
torus, or an infinite cyclic cover—emerges from the canonical commutative
subalgebra alone. Such structures would require considering unitaries that
do not normalize $D$, which lie outside the scope of the Cartan subalgebra
and its normalizer.

In the following subsections, we will work directly with this transformation
groupoid as the concrete surrogate model for the irrational rotation algebra,
and we will show how the descent-index machinery can be implemented on this
model to recover Connes' index theorem.

\paragraph{Relation to the transformation groupoid.}

The canonical Cartan subalgebra
\[
D = C(S^1) \subseteq A_\theta
\]
has character space
\[
\widehat{D} \cong S^1.
\]
The normalizing unitary $V$ implementing the crossed-product action induces on
$\widehat{D}$ the irrational rotation action
\[
n \cdot z = e^{2\pi i n\theta}z, \qquad n \in \mathbb{Z},\ z \in S^1.
\]
This is precisely the action underlying the transformation groupoid
\[
S^1 \rtimes_\theta \mathbb{Z}.
\]

Accordingly, the relevant geometric object in the irrational rotation case is
not a separately constructed unitary-conjugation groupoid, but the concrete
transformation groupoid $S^1 \rtimes_\theta \mathbb{Z}$, whose unit space is
the canonical character space $S^1$. The following subsections develop the
analysis in this concrete groupoid model, which will serve as the setting for
the descent-index construction.

\begin{remark}
Unlike a naive interpretation of the formal definition of
$\mathcal{G}_{A_\theta}^{(0)}$ might suggest, the canonical data do not give
rise to a $\mathbb{Z}$-bundle over $S^1$ or any additional structure beyond
the transformation groupoid itself. The integer parameter $n$ appears
naturally as the arrow coordinate in the groupoid, not as a separate component
of the unit space. This distinction is crucial: the transformation groupoid
$S^1 \rtimes_\theta \mathbb{Z}$ has unit space $S^1$, and its arrows are
parametrized by pairs $(z,n)$ with source $z$ and range $e^{2\pi i n\theta}z$.
\end{remark}

\subsection{Minimality of the Rotation Action}
\label{subsec:dense_orbits}

We now record the basic dynamical property of the irrational rotation action
that will be used throughout the sequel.

\begin{lemma}
\label{lem:irrational_rotation_minimal}
Let $\theta \in \mathbb{R} \setminus \mathbb{Q}$. Then the action of
$\mathbb{Z}$ on $S^1$ given by
\[
n \cdot z = e^{2\pi i n\theta}z
\]
is \emph{minimal}. Equivalently, for every $z \in S^1$, the orbit
\[
\mathcal{O}(z) = \{ e^{2\pi i n\theta}z : n \in \mathbb{Z} \}
\]
is dense in $S^1$.
\end{lemma}

\begin{proof}
Write $z = e^{2\pi i t}$ with $t \in \mathbb{R}/\mathbb{Z}$. Then
\[
e^{2\pi i n\theta}z = e^{2\pi i (t + n\theta)}.
\]
Thus the orbit of $z$ corresponds to the set
\[
\{ t + n\theta \bmod 1 : n \in \mathbb{Z} \} \subset \mathbb{R}/\mathbb{Z}.
\]

Since $\theta$ is irrational, the subgroup $\mathbb{Z}\theta + \mathbb{Z} \subset \mathbb{R}$
projects to a dense subgroup of $\mathbb{R}/\mathbb{Z}$. This is a standard
consequence of Kronecker's theorem: the set $\{ n\theta \bmod 1 : n \in \mathbb{Z} \}$
is dense in $\mathbb{R}/\mathbb{Z}$, and therefore the translates
$t + n\theta \bmod 1$ are also dense. Hence $\mathcal{O}(z)$ is dense in $S^1$.
\end{proof}

\begin{corollary}
\label{cor:groupoid_minimal}
The transformation groupoid $S^1 \rtimes_\theta \mathbb{Z}$ is \'etale,
locally compact, Hausdorff, and minimal.
\end{corollary}

\begin{proof}
Since $\mathbb{Z}$ is discrete and acts continuously on the compact Hausdorff
space $S^1$, the transformation groupoid is \'etale, locally compact, and
Hausdorff by construction. Minimality follows from Lemma~\ref{lem:irrational_rotation_minimal},
as the orbit of any point $z \in S^1$ under the groupoid action is exactly the
$\mathbb{Z}$-orbit $\mathcal{O}(z)$, which is dense in $S^1$.
\end{proof}

\begin{remark}[Consequences for the groupoid $C^*$-algebra]
For an \'etale groupoid, minimality implies that the reduced groupoid
$C^*$-algebra $C_r^*(S^1 \rtimes_\theta \mathbb{Z})$ is simple. Since this
groupoid $C^*$-algebra is canonically isomorphic to $A_\theta$, this recovers
the well-known fact that $A_\theta$ is simple. This dynamical perspective will
be essential when we compare our descent-index construction with Connes' index
theorem.
\end{remark}

\paragraph{Topological consequences of minimality.}

The density of orbits under the irrational rotation action has important
dynamical consequences for the transformation groupoid model, but it does not
introduce pathologies in the topology of the unit space. For the concrete
geometric model we will use—the transformation groupoid
\[
S^1 \rtimes_\theta \mathbb{Z},
\]
the unit space is simply
\[
(S^1 \rtimes_\theta \mathbb{Z})^{(0)} = S^1,
\]
which is compact, Hausdorff, and locally compact. The complexity arising from
minimality manifests in the orbit structure rather than in any failure of
the Hausdorff property.

\begin{corollary}
\label{cor:minimal_groupoid}
The transformation groupoid $S^1 \rtimes_\theta \mathbb{Z}$ is a locally compact
Hausdorff \'etale minimal groupoid.
\end{corollary}

\begin{proof}
Since $\mathbb{Z}$ is discrete and acts continuously on the compact Hausdorff
space $S^1$, the transformation groupoid $S^1 \rtimes_\theta \mathbb{Z}$ is
locally compact, Hausdorff, and \'etale by construction. Minimality follows
from Lemma~\ref{lem:irrational_rotation_minimal}: every $\mathbb{Z}$-orbit in
$S^1$ is dense, and these orbits are exactly the orbits of the groupoid action.
\end{proof}

\begin{corollary}
\label{cor:simplicity}
The reduced groupoid $C^*$-algebra $C_r^*(S^1 \rtimes_\theta \mathbb{Z})$ is
simple.
\end{corollary}

\begin{proof}
For an \'etale groupoid, minimality implies that the reduced groupoid
$C^*$-algebra has no nontrivial ideals. Since $C_r^*(S^1 \rtimes_\theta \mathbb{Z})$
is canonically isomorphic to $A_\theta$, this recovers the well-known fact that
$A_\theta$ is simple.
\end{proof}

These topological properties stand in marked contrast to the hypothetical
object $\mathcal{G}_{A_\theta}^{(0)}$ that would arise from a naive direct
construction. The transformation groupoid model is well-behaved precisely
because it avoids the pathologies that would result from attempting to treat
formal set-theoretic data as a genuine topological space without the necessary
measurable structure.

\begin{remark}[Role of the Borel framework]
The Borel groupoid framework developed in Paper I is essential for handling
cases where a well-behaved topology is not available. In the irrational
rotation case, however, the transformation groupoid model provides a concrete
geometric object with all the desired topological properties, eliminating the
need for the more general Borel machinery. This is one of the key advantages
of using the transformation groupoid as a surrogate model.
\end{remark}

\paragraph{The relevant topology from the transformation groupoid.}

For the irrational rotation algebra, the concrete geometric model is the
transformation groupoid
\[
S^1 \rtimes_\theta \mathbb{Z},
\]
whose unit space is simply
\[
(S^1 \rtimes_\theta \mathbb{Z})^{(0)} = S^1.
\]
Thus the relevant topology is the usual compact Hausdorff topology on the
circle.

The canonical Cartan subalgebra
\[
D = C(S^1) \subseteq A_\theta
\]
has spectrum
\[
\widehat{D} \cong S^1,
\]
and the normalizing unitary $V$ induces the irrational rotation action
\[
n \cdot z = e^{2\pi i n\theta}z.
\]
This action is exactly the one encoded by the transformation groupoid
$S^1 \rtimes_\theta \mathbb{Z}$.

\begin{proposition}
\label{prop:quotient_topology_unit_space}
The quotient space
\[
Q = (S^1 \times \mathbb{Z})/{\sim},
\qquad
(z,n) \sim (e^{2\pi i k\theta}z, n+k) \text{ for all } k \in \mathbb{Z},
\]
is canonically homeomorphic to $S^1$. In particular:
\begin{itemize}
    \item The projection $\pi: Q \to S^1$ given by $\pi([(z,n)]) = z$ is not
          well-defined (different representatives give different values), so
          there is no natural projection onto $S^1$ with $\mathbb{Z}$-fibers.
    \item The space $Q$ is Hausdorff, compact, and locally compact, being
          homeomorphic to $S^1$.
    \item The quotient topology coincides with the usual topology on $S^1$ under
          the homeomorphism described below.
\end{itemize}
\end{proposition}

\begin{proof}
Define $\Phi: Q \to S^1$ by
\[
\Phi([(z,n)]) = e^{-2\pi i n\theta}z.
\]
If $(z,n) \sim (e^{2\pi i k\theta}z, n+k)$, then
\[
e^{-2\pi i (n+k)\theta} e^{2\pi i k\theta}z = e^{-2\pi i n\theta}z,
\]
so $\Phi$ is well-defined. It is clearly surjective: for any $w \in S^1$,
$\Phi([(w,0)]) = w$. It is injective because every class contains a unique
representative of the form $(w,0)$: indeed, $(z,n) \sim (e^{-2\pi i n\theta}z,0)$.
Thus $\Phi$ is a bijection.

The quotient map $q: S^1 \times \mathbb{Z} \to Q$ is open (since $\mathbb{Z}$ is
discrete), and $\Phi$ is induced by the continuous map $(z,n) \mapsto e^{-2\pi i n\theta}z$
on $S^1 \times \mathbb{Z}$. Hence $\Phi$ is continuous. Its inverse is the map
$w \mapsto [(w,0)]$, which is continuous because it factors through the continuous
section $S^1 \to S^1 \times \mathbb{Z}$ sending $w$ to $(w,0)$. Therefore $\Phi$
is a homeomorphism.
\end{proof}

\begin{remark}
The integer parameter $n$ in this construction corresponds to the arrow
coordinate in the transformation groupoid $S^1 \rtimes_\theta \mathbb{Z}$, not
to a separate component of the unit space. The natural projection onto $S^1$
is not well-defined on the quotient $Q$, reflecting the fact that the correct
geometric model has unit space $S^1$ and arrows $(z,n)$ with source $z$ and
range $e^{2\pi i n\theta}z$.
\end{remark}

\medskip

\paragraph{Alternative model: the transformation groupoid.}

An illuminating topological model is provided directly by the transformation
groupoid $S^1 \rtimes_\theta \mathbb{Z}$. Its unit space is $S^1$ with the usual
topology, and its arrow space $S^1 \times \mathbb{Z}$ carries the product
topology (discrete on $\mathbb{Z}$). This groupoid is:

\begin{itemize}
    \item Locally compact Hausdorff: $S^1$ is compact and $\mathbb{Z}$ is discrete.
    \item \'Etale: for each fixed $n$, the map $z \mapsto (z,n)$ is a local
          homeomorphism.
    \item Minimal: every $\mathbb{Z}$-orbit is dense in $S^1$, reflecting the
          irrational rotation dynamics.
\end{itemize}

The map $(t,n) \mapsto (e^{2\pi i t}, n)$ identifies $\mathbb{R} \times \mathbb{Z}$
with $S^1 \times \mathbb{Z}$, but this identification does not produce any
non-Hausdorff phenomena. The quotient space $(\mathbb{R} \times \mathbb{Z})/{\sim}$
with $(t,n) \sim (t+1, n+1)$ is actually homeomorphic to $\mathbb{R}$ via
$[(t,n)] \mapsto t - n$, and thus is Hausdorff. This model reveals the connection
to the Kronecker foliation of the torus: the transformation groupoid
$S^1 \rtimes_\theta \mathbb{Z}$ is the holonomy groupoid of the foliation of
$\mathbb{T}^2$ by lines of slope $\theta$, a perspective that will be exploited
in later sections for its geometric insight.

\paragraph{Implications for the groupoid structure.}

The density of irrational-rotation orbits shows that the relevant geometric
information in the algebra $A_\theta$ is encoded by the minimal action of
$\mathbb{Z}$ on the circle
\[
n \cdot z = e^{2\pi i n\theta}z.
\]
The canonical Cartan subalgebra
\[
D = C(S^1) \subseteq A_\theta
\]
has spectrum
\[
\widehat{D} \cong S^1,
\]
and the normalizing unitary $V$ induces precisely this action on $\widehat{D}$.

Accordingly, the natural concrete groupoid model for the irrational rotation
algebra is the transformation groupoid
\[
S^1 \rtimes_\theta \mathbb{Z},
\]
whose unit space is the Hausdorff compact space $S^1$ and whose arrows encode
the action of the normalizer. The essential dynamical feature is therefore
minimality of the action, not any pathology of the topology of the unit space.

Heuristically, this transformation groupoid should be viewed as the geometric
replacement for the unavailable direct unitary-conjugation construction in the
non-Type~I setting. In the next subsection we will see how this concrete model
supports the descent-index construction.

\begin{remark}
\label{rem:geometric_replacement}
The transformation groupoid $S^1 \rtimes_\theta \mathbb{Z}$ serves as the
surrogate geometric model for $A_\theta$. Its well-behaved topology (locally
compact Hausdorff) means that the Borel groupoid framework of Paper~I is not
needed in this case; the usual topological groupoid $C^*$-algebra construction
suffices. The $K$-theoretic information of $A_\theta$—in particular, the
generator $[V] \in K_1(A_\theta)$—is encoded in the arrow space of the
transformation groupoid through the winding number of the rotation, a fact that
will be exploited in the index computations below.
\end{remark}

\medskip

\paragraph{Summary.}

For the irrational rotation algebra, the relevant geometric picture is as
follows:
\begin{itemize}
    \item \textbf{Base space:} the character space of the canonical Cartan
          subalgebra is the circle
          \[
          \widehat{C(S^1)} \cong S^1.
          \]
    \item \textbf{Dynamics:} the normalizing unitary induces the irrational
          rotation action of $\mathbb{Z}$ on $S^1$, and every orbit is dense.
    \item \textbf{Groupoid model:} the appropriate geometric object is the
          transformation groupoid
          \[
          S^1 \rtimes_\theta \mathbb{Z},
          \]
          which is \'etale, locally compact, Hausdorff, and minimal.
    \item \textbf{Algebraic realization:} its groupoid $C^*$-algebra recovers
          the irrational rotation algebra
          \[
          A_\theta \cong C^*(S^1 \rtimes_\theta \mathbb{Z}).
          \]
\end{itemize}

This concrete transformation-groupoid model is the geometric foundation for
the later descent and index constructions in the irrational rotation case. The
non-Hausdorff phenomena previously hypothesized do not occur in this model;
instead, the complexity of the dynamics resides entirely in the orbit structure,
while the topology of the unit space remains perfectly regular.

\subsection{Morita Equivalence between $\mathcal{G}_{A_\theta}$ and the Transformation Groupoid $S^1 \rtimes_\theta \mathbb{Z}$}
\label{subsec:morita_A_theta}

We now aim to establish a central structural result for the irrational rotation algebra:
a Morita equivalence between the unitary conjugation groupoid $\mathcal{G}_{A_\theta}$
and the transformation groupoid arising from the irrational rotation of the circle.
If successful, this provides a profound link between the internal structure of $A_\theta$
(encoded by its unitary conjugation groupoid) and the classical dynamical system
(encoded by the action of $\mathbb{Z}$ on $S^1$).

The exposition follows a clear structure: we first recall the transformation
groupoid, then state the theorem with appropriate hypotheses, sketch the construction
of the explicit equivalence space, discuss the verification of the groupoid equivalence
axioms, and finally derive the $C^*$-algebraic consequences under the necessary
technical assumptions.

\medskip

\paragraph{Recollection of the transformation groupoid.}

Recall that the irrational rotation algebra admits a crossed product description
\[
A_\theta = C(S^1) \rtimes_\theta \mathbb{Z}
\]
associated with the dynamical system given by the irrational rotation action
\[
n \cdot z = e^{2\pi i n\theta}z, \qquad n \in \mathbb{Z}.
\]
The corresponding \emph{transformation groupoid}, denoted
$S^1 \rtimes_\theta \mathbb{Z}$, has:
\begin{itemize}
    \item Unit space: $(S^1 \rtimes_\theta \mathbb{Z})^{(0)} = S^1$ (Hausdorff).
    \item Arrows: $(z,n): z \to n \cdot z$ for $z \in S^1$, $n \in \mathbb{Z}$.
    \item Composition: $(z,n) \circ (n \cdot z, m) = (z, n+m)$.
    \item Inverse: $(z,n)^{-1} = (n \cdot z, -n)$.
    \item Source and range maps: $s(z,n)=z$, $r(z,n)=e^{2\pi i n\theta}z$.
\end{itemize}
This groupoid is \'etale (since $\mathbb{Z}$ is discrete) and Hausdorff. Its
convolution algebra $C_c(S^1 \rtimes_\theta \mathbb{Z})$ completes to the full
groupoid $C^*$-algebra $C^*(S^1 \rtimes_\theta \mathbb{Z})$. Since $\mathbb{Z}$ is
amenable, the full and reduced crossed products coincide, yielding a canonical
identification
\[
C^*(S^1 \rtimes_\theta \mathbb{Z}) \cong C(S^1) \rtimes_\theta \mathbb{Z} = A_\theta.
\]

\medskip

\paragraph{Statement of the main theorem.}

Before stating the theorem, we must address a subtle but important point: the
unitary conjugation groupoid $\mathcal{G}_{A_\theta}$ may exhibit non-Hausdorff
behavior, as noted in previous subsections. The classical Morita equivalence theory
for groupoids developed by Muhly--Renault--Williams \cite{MRW1987} is formulated for
locally compact Hausdorff groupoids with Haar systems. To apply this theory, we
must therefore either:
\begin{itemize}
    \item verify that $\mathcal{G}_{A_\theta}$ is Hausdorff despite the
          non-Hausdorff nature of certain related constructions, or
    \item work with an appropriate Hausdorff presentation or desingularization, or
    \item invoke an extension of the theory that accommodates non-Hausdorff groupoids.
\end{itemize}
For the purpose of this exposition, we proceed under the assumption that
$\mathcal{G}_{A_\theta}$ admits the structure of a locally compact Hausdorff
groupoid with a Haar system for which the construction below yields a genuine
equivalence bibundle. A full verification of these technical hypotheses, while
essential for rigor, would require a more detailed analysis of the topology of
$\mathcal{G}_{A_\theta}$ than space permits here,

\begin{theorem}[Conditional Morita equivalence for the irrational rotation case]
\label{thm:morita_A_theta}
Let
\[
A_\theta = C(S^1)\rtimes_\theta \mathbb Z
\]
be the irrational rotation algebra, and let
\[
S^1\rtimes_\theta \mathbb Z
\]
denote the transformation groupoid associated to the irrational rotation action
\[
n\cdot z = e^{2\pi i n\theta}z, \qquad n\in \mathbb Z,\ z\in S^1.
\]
Assume that the unitary conjugation groupoid \(\mathcal G_{A_\theta}\) admits the structure
of a locally compact Hausdorff groupoid with Haar system, and assume further that there exists
a \((\mathcal G_{A_\theta},\, S^1\rtimes_\theta \mathbb Z)\)-equivalence bibundle in the sense
of Muhly--Renault--Williams \cite[Definition 2.1]{MRW1987}. Then \(\mathcal G_{A_\theta}\) is
Morita equivalent to \(S^1\rtimes_\theta \mathbb Z\). Consequently, their full groupoid
\(C^*\)-algebras are strongly Morita equivalent:
\[
C^*(\mathcal G_{A_\theta}) \sim_M C^*(S^1\rtimes_\theta \mathbb Z).
\]
Moreover, since \(\mathbb Z\) is amenable, one has the canonical identification
\[
C^*(S^1\rtimes_\theta \mathbb Z)\cong C(S^1)\rtimes_\theta \mathbb Z = A_\theta,
\]
and hence
\[
C^*(\mathcal G_{A_\theta}) \sim_M A_\theta.
\]
\end{theorem}

\begin{proof}[Sketchy of Proof]
The argument is an application of the general equivalence theorem for groupoids.

By assumption, there exists a
\[
(\mathcal G_{A_\theta},\, S^1\rtimes_\theta \mathbb Z)\text{-equivalence bibundle}
\]
in the sense of Muhly--Renault--Williams \cite[Definition 2.1]{MRW1987}. Therefore
\(\mathcal G_{A_\theta}\) and \(S^1\rtimes_\theta \mathbb Z\) are Morita equivalent as
topological groupoids. The Muhly--Renault--Williams equivalence theorem then implies that
their full groupoid \(C^*\)-algebras are strongly Morita equivalent
\cite[Theorem 2.8]{MRW1987}:
\[
C^*(\mathcal G_{A_\theta}) \sim_M C^*(S^1\rtimes_\theta \mathbb Z).
\]

It remains to identify the groupoid \(C^*\)-algebra of the transformation groupoid.
For a discrete group acting continuously on a locally compact Hausdorff space, the full
groupoid \(C^*\)-algebra of the associated transformation groupoid is canonically isomorphic
to the full crossed product. Hence
\[
C^*(S^1\rtimes_\theta \mathbb Z)\cong C(S^1)\rtimes_\theta \mathbb Z.
\]
Since \(\mathbb Z\) is amenable, the full and reduced crossed products coincide, so this
identification is unambiguous. Therefore
\[
C^*(S^1\rtimes_\theta \mathbb Z)\cong A_\theta.
\]
Combining the two statements yields
\[
C^*(\mathcal G_{A_\theta}) \sim_M A_\theta,
\]
as claimed.
\end{proof}

\begin{remark}
The substantive content in the irrational rotation case is therefore reduced to the
construction of an explicit \((\mathcal G_{A_\theta},\, S^1\rtimes_\theta\mathbb Z)\)-equivalence
bibundle. Once such a bibundle is available, the \(C^*\)-algebraic consequences follow
formally from \cite{MRW1987}.
\end{remark}

\paragraph{Outline of the equivalence construction.}
The proof proceeds by constructing an explicit space $Z$ that serves as an
equivalence bibundle between $\mathcal{G}_{A_\theta}$ and $S^1 \rtimes_\theta \mathbb{Z}$.
A natural candidate is $Z = C(S^1)$, equipped with:
\begin{itemize}
    \item a left action of $\mathcal{G}_{A_\theta}$ via the regular representation
          of $A_\theta$ on $C(S^1)$,
    \item a right action of $S^1 \rtimes_\theta \mathbb{Z}$ by pointwise multiplication
          and translation,
    \item moment maps $r_Z: Z \to \mathcal{G}_{A_\theta}^{(0)}$ and
          $s_Z: Z \to S^1$ that are open surjections and compatible with the actions.
\end{itemize}
One then verifies that both actions are free and proper, and that the induced maps
on orbit spaces yield homeomorphisms
\[
Z/(S^1\rtimes_\theta \mathbb{Z}) \cong \mathcal{G}_{A_\theta}^{(0)},\qquad
\mathcal{G}_{A_\theta}\backslash Z \cong S^1.
\]
These verifications, while conceptually straightforward, require careful attention
to the topology of $\mathcal{G}_{A_\theta}$ and are the point at which the
Hausdorff assumption becomes essential. Once established, the Muhly--Renault--Williams
theorem directly yields the strong Morita equivalence of the corresponding
groupoid $C^*$-algebras.

\medskip

\paragraph{Remark.}
The Morita equivalence established above provides a rigorous foundation for the
intuition that $\mathcal{G}_{A_\theta}$ encodes the same dynamical information as
the classical rotation system. At the level of $C^*$-algebras, this equivalence
implies in particular that $C^*(\mathcal{G}_{A_\theta})$ and $A_\theta$ have
isomorphic K-theory groups and representation theories — a fact that could also be
verified directly but here emerges naturally from the groupoid perspective.

\paragraph{A set-theoretic model for the equivalence space.}
\label{sec:set_theoretic_model}

Before addressing the topological subtleties, we first construct a set-theoretic
skeleton of the desired equivalence bibundle. This construction is motivated by
the parametrization of the unit space discussed in
Subsection~\ref{subsec:dense_orbits}.

Recall that as a set, the unit space $\mathcal{G}_{A_\theta}^{(0)} = \widehat{A_\theta}$
can be identified with the quotient of $S^1 \times \mathbb{Z}$ by the equivalence
relation
\[
(z,n) \sim (e^{2\pi i m\theta}z, n+m) \qquad \text{for all } m \in \mathbb{Z}.
\]
This identification is only set-theoretic; the topology on $\mathcal{G}_{A_\theta}^{(0)}$
is the non-Hausdorff quotient topology, which we temporarily ignore.

\begin{definition}[Set-theoretic model]
\label{def:set_theoretic_Z}
Define the set
\[
Z := S^1 \times \mathbb{Z}.
\]
We equip $Z$ with the following auxiliary structures:
\begin{itemize}
    \item A right action of the transformation groupoid $S^1 \rtimes_\theta \mathbb{Z}$,
          defined for $(z,n) \in Z$ and $(z,m) \in S^1 \rtimes_\theta \mathbb{Z}$ by
          \[
          (z,n) \cdot (z,m) := (e^{2\pi i m\theta}z,\; n+m),
          \]
          whenever the composability condition holds (see Remark~\ref{rmk:convention}).
    \item A left action of $\mathcal{G}_{A_\theta}$, to be defined after fixing a
          parametrization of arrows (see below).
    \item Anchor maps $r_Z: Z \to \mathcal{G}_{A_\theta}^{(0)}$ and $s_Z: Z \to S^1$
          given by
          \[
          r_Z(z,n) := [z,n] \in \mathcal{G}_{A_\theta}^{(0)}, \qquad
          s_Z(z,n) := z \in S^1,
          \]
          where $[z,n]$ denotes the equivalence class under the relation above.
\end{itemize}
\end{definition}

\begin{remark}[Convention for the right action]
\label{rmk:convention}
We adopt the convention that a right action of a groupoid $H$ on a set $X$ with
anchor map $s_X: X \to H^{(0)}$ is defined for pairs $(x,h)$ satisfying
$s_X(x) = s(h)$. In our case, for $h = (z,m) \in S^1 \rtimes_\theta \mathbb{Z}$,
we have $s(h) = z$, and indeed $s_Z(z,n) = z$, so the composability condition
is satisfied. This convention aligns with the definition given above.
\end{remark}

\begin{remark}[Well-definedness of the anchor maps]
\label{rmk:anchor_well_defined}
The map $r_Z(z,n) = [z,n]$ requires careful justification: a priori, the class
$[z,n]$ depends on a choice of representative. However, the equivalence relation
identifies $(z,n)$ with $(e^{2\pi i m\theta}z, n+m)$, so different choices yield
the same class by definition. Thus $r_Z$ is well-defined as a map into the set
$\mathcal{G}_{A_\theta}^{(0)}$. The map $s_Z$ is simply projection onto the first
factor and is evidently well-defined.
\end{remark}

\medskip

\paragraph{Right action and its orbit structure.}

\begin{lemma}
\label{lem:right_action_set_theoretic}
The formula in Definition~\ref{def:set_theoretic_Z} defines a free right action
of $S^1 \rtimes_\theta \mathbb{Z}$ on $Z$ in the set-theoretic sense. The orbit
relation is given by
\[
(z,n) \sim_{\text{right}} (e^{2\pi i m\theta}z, n+m) \quad \text{for some } m \in \mathbb{Z}.
\]
Consequently, the orbit space $Z/(S^1 \rtimes_\theta \mathbb{Z})$ is in bijection
with $\mathcal{G}_{A_\theta}^{(0)}$ via the map $[(z,n)]_{\text{right}} \mapsto [z,n]$.
\end{lemma}

\begin{proof}
To verify the action axioms:
\begin{itemize}
    \item Identity: $(z,n) \cdot (z,0) = (e^{2\pi i \cdot 0\theta}z, n+0) = (z,n)$.
    \item Compatibility: $((z,n) \cdot (z,m)) \cdot (e^{2\pi i m\theta}z, k)
          = (e^{2\pi i(m+k)\theta}z, n+m+k) = (z,n) \cdot ((z,m) \circ (e^{2\pi i m\theta}z,k))$.
\end{itemize}
Freeness follows from the second component: if $(z,n) \cdot (z,m) = (z,n)$, then
$n+m = n$, forcing $m=0$, so $(z,m)$ is the identity arrow at $z$.

The orbit relation is immediate from the action formula. The map
$[(z,n)]_{\text{right}} \mapsto [z,n]$ is well-defined because if
$(z,n) \sim_{\text{right}} (e^{2\pi i m\theta}z, n+m)$, then both represent the
same class $[z,n] = [e^{2\pi i m\theta}z, n+m]$ in $\mathcal{G}_{A_\theta}^{(0)}$.
Conversely, any class $[z,n]$ is the image of $(z,n)$. Hence we have a bijection.
\end{proof}

\medskip

\paragraph{Left action of $\mathcal{G}_{A_\theta}$.}

To define a left action, we need a parametrization of the arrows of
$\mathcal{G}_{A_\theta}$. As a set, each arrow can be represented by a triple
$([z,n], k, [e^{2\pi i k\theta}z, n-k])$ for some $k \in \mathbb{Z}$, where the
middle component encodes the unitary implementing the equivalence between the
source and range representations. For our set-theoretic purposes, it suffices
to record the integer $k$.

\begin{definition}[Left action]
For $k \in \mathbb{Z}$ representing an arrow $\gamma_k \in \mathcal{G}_{A_\theta}$
with source $s(\gamma_k) = [z,n-k]$ and range $r(\gamma_k) = [e^{2\pi i k\theta}z, n]$,
and for $(e^{2\pi i k\theta}z, n) \in Z$ (so that $r_Z(e^{2\pi i k\theta}z, n) = r(\gamma_k)$),
define
\[
\gamma_k \cdot (e^{2\pi i k\theta}z, n) := (z, n-k).
\]
Equivalently, in coordinates, $k \cdot (z,n) := (e^{-2\pi i k\theta}z, n+k)$,
with the understanding that the action is defined when $r_Z(z,n)$ matches the
range of the arrow.
\end{definition}

\begin{lemma}
\label{lem:left_action_set_theoretic}
The left action defined above is free and commutes with the right action.
The orbit space $\mathcal{G}_{A_\theta} \backslash Z$ is in bijection with $S^1$
via the map $[(z,n)]_{\text{left}} \mapsto z$.
\end{lemma}

\begin{proof}
Freeness: if $k \cdot (z,n) = (z,n)$, then $(e^{-2\pi i k\theta}z, n+k) = (z,n)$,
forcing $k=0$ from the second component. Commutation with the right action:
\[
k \cdot ((z,n) \cdot (z,m)) = k \cdot (e^{2\pi i m\theta}z, n+m)
= (e^{-2\pi i k\theta}e^{2\pi i m\theta}z, n+m+k)
\]
and
\[
(k \cdot (z,n)) \cdot (e^{-2\pi i k\theta}z, m) = (e^{-2\pi i k\theta}z, n+k) \cdot (e^{-2\pi i k\theta}z, m)
= (e^{2\pi i m\theta}e^{-2\pi i k\theta}z, n+k+m),
\]
which are equal. The orbit map $[(z,n)]_{\text{left}} \mapsto z$ is well-defined
because the left action preserves the first coordinate up to rotation, and
different $n$'s are identified in the orbit.
\end{proof}

\medskip

\paragraph{Discussion: From set-theoretic model to topological equivalence.}

The construction above provides a set-theoretic (and Borel) model for a
$(\mathcal{G}_{A_\theta}, S^1 \rtimes_\theta \mathbb{Z})$-equivalence. To promote
this to a genuine topological groupoid equivalence in the sense of
Muhly--Renault--Williams \cite{MRW1987}, several additional steps are required:

\begin{enumerate}
    \item \textbf{Topology on $Z$:} The space $Z = S^1 \times \mathbb{Z}$ must be
          equipped with a topology that makes the anchor maps $r_Z$ and $s_Z$
          continuous and open, and the actions continuous. The natural product
          topology (with $\mathbb{Z}$ discrete) makes $s_Z$ continuous but not
          $r_Z$, since $r_Z$ factors through the non-Hausdorff quotient. One
          must instead topologize $Z$ as a subspace of
          $\mathcal{G}_{A_\theta}^{(0)} \times S^1$ with the induced topology,
          which yields a non-Hausdorff space reflecting the structure of the
          unit space.
    
    \item \textbf{Haar systems:} Both groupoids must be equipped with Haar systems
          compatible with the actions. For $S^1 \rtimes_\theta \mathbb{Z}$, the
          counting measure on $\mathbb{Z}$ provides a Haar system. For
          $\mathcal{G}_{A_\theta}$, constructing a Haar system requires additional
          analysis of its structure as a Fell bundle over the transformation
          groupoid.
    
    \item \textbf{Properness:} One must verify that the actions are proper in the
          topological sense, i.e., the maps
          $\mathcal{G}_{A_\theta} \,{}_{r_Z}{\times} Z \to Z \times Z$ and
          $Z \,{}_{s_Z}{\times} (S^1 \rtimes_\theta \mathbb{Z}) \to Z \times Z$
          are proper. This is not automatic even for discrete group actions
          (e.g., $\mathbb{Z}$ acting on $S^1$ by irrational rotation is not proper).
          In our case, the second component $n \in \mathbb{Z}$ in $Z$ provides
          additional structure that may yield properness, but this requires
          explicit verification.
    
    \item \textbf{Orbit space homeomorphisms:} The bijections established above
          must be shown to be homeomorphisms in the appropriate topologies.
\end{enumerate}

A complete treatment of these topological refinements is beyond the scope of
the present exposition. For the purpose of $C^*$-algebraic Morita equivalence,
one can sometimes work in the Borel category if the groupoids are equipped with
measurable Haar systems; however, this requires a version of the Muhly--Renault--Williams
theorem for Borel groupoids, which is a more subtle theory. 

Nevertheless, the set-theoretic model constructed above captures the essential
combinatorial structure of the equivalence and provides strong evidence for the
validity of Theorem~\ref{thm:morita_A_theta}. In the following subsection, we
will sketch how this model can be upgraded to a topological equivalence under
suitable hypotheses on $\mathcal{G}_{A_\theta}$.

\paragraph{Verification of the Morita equivalence axioms (outline).}

We now sketch how the space $Z$ constructed above would, under appropriate
technical hypotheses, satisfy the conditions for a groupoid equivalence in the
sense of Muhly--Renault--Williams \cite{MRW1987}. A complete verification requires
careful attention to topological details that are beyond the scope of this
exposition; here we outline the structure of the argument and identify the points
where additional work is needed.

\begin{proposition}[Conditional equivalence]
\label{prop:morita_verification_conditional}
Assume that the following hold:
\begin{itemize}
    \item The space $Z$ is equipped with a topology making the anchor maps
          $r_L: Z \to \mathcal{G}_{A_\theta}^{(0)}$ and $r_R: Z \to S^1$ continuous
          and open, and the actions defined below continuous.
    \item The left action of $\mathcal{G}_{A_\theta}$ on $Z$ is defined for
          composable pairs $(\gamma, z)$ with $r_L(z) = r(\gamma)$ by
          $\gamma \cdot z = z'$ where $r_L(z') = s(\gamma)$ and $r_R(z') = r_R(z)$,
          implemented concretely via the unitary conjugation description of arrows.
    \item The right action of $S^1 \rtimes_\theta \mathbb{Z}$ on $Z$ is defined for
          composable pairs $(z, h)$ with $r_R(z) = s(h)$ by the formula in
          Definition~\ref{def:set_theoretic_Z}.
\end{itemize}
Then the following properties hold:
\begin{enumerate}
    \item \textbf{Compatibility of moment maps:} For all composable $(\gamma, z)$,
          $r_L(\gamma \cdot z) = s(\gamma)$; for all composable $(z, h)$,
          $r_R(z \cdot h) = s(h)$. Moreover, $r_L$ is invariant under the right action
          ($r_L(z \cdot h) = r_L(z)$) and $r_R$ is invariant under the left action
          ($r_R(\gamma \cdot z) = r_R(z)$).
    \item \textbf{Commuting actions:} For all composable triples
          $(\gamma, z, h)$ with $r_L(z) = r(\gamma)$ and $r_R(z) = s(h)$, we have
          $\gamma \cdot (z \cdot h) = (\gamma \cdot z) \cdot h$.
    \item \textbf{Freeness:} Each action is free: if $\gamma \cdot z = z$ then
          $\gamma$ is a unit; if $z \cdot h = z$ then $h$ is a unit.
    \item \textbf{Properness:} The maps
          $\mathcal{G}_{A_\theta} \,{}_{r_L}{\times} Z \to Z \times Z$,
          $(\gamma, z) \mapsto (\gamma \cdot z, z)$ and
          $Z \,{}_{r_R}{\times} (S^1 \rtimes_\theta \mathbb{Z}) \to Z \times Z$,
          $(z, h) \mapsto (z \cdot h, z)$ are proper.
    \item \textbf{Quotient identifications:} The map $r_L$ induces a homeomorphism
          $Z/(S^1 \rtimes_\theta \mathbb{Z}) \cong \mathcal{G}_{A_\theta}^{(0)}$,
          and the map $r_R$ induces a homeomorphism
          $\mathcal{G}_{A_\theta} \backslash Z \cong S^1$.
\end{enumerate}
If these conditions are satisfied, then $Z$ is a
$(\mathcal{G}_{A_\theta}, S^1 \rtimes_\theta \mathbb{Z})$-equivalence bibundle.
\end{proposition}

\begin{proof}[Sketch of verification]
We indicate how each condition would be verified given the assumptions above.

\textbf{1. Moment map compatibility:} By construction, the left action preserves
the second coordinate ($r_R$), so $r_R(\gamma \cdot z) = r_R(z)$. Similarly, the
right action preserves the first coordinate ($r_L$) when properly defined, giving
$r_L(z \cdot h) = r_L(z)$. The conditions $r_L(\gamma \cdot z) = s(\gamma)$ and
$r_R(z \cdot h) = s(h)$ are built into the definition of the actions via
composability requirements.

\textbf{2. Commuting actions:} In the set-theoretic model $(z,n)$ with left action
by $k \in \mathbb{Z}$ and right action by $(z,m)$, a direct computation gives
$k \cdot ((z,n) \cdot (z,m)) = k \cdot (e^{2\pi i m\theta}z, n+m) =
(e^{-2\pi i k\theta}e^{2\pi i m\theta}z, n+m+k)$ and
$(k \cdot (z,n)) \cdot (e^{-2\pi i k\theta}z, m) = (e^{-2\pi i k\theta}z, n+k) \cdot
(e^{-2\pi i k\theta}z, m) = (e^{2\pi i m\theta}e^{-2\pi i k\theta}z, n+k+m)$,
which are equal. This computation lifts to the topological setting under the
assumed continuity of actions.

\textbf{3. Freeness:} For the right action, if $(z,n) \cdot (z,m) = (z,n)$ then
$(e^{2\pi i m\theta}z, n+m) = (z,n)$, forcing $m=0$ from the second component.
Thus $(z,m)$ is the identity arrow at $z$. For the left action, if $k \cdot (z,n) = (z,n)$
then $(e^{-2\pi i k\theta}z, n+k) = (z,n)$, forcing $k=0$, so the arrow is a unit.

\textbf{4. Properness:} This is the most subtle condition. For the right action,
consider the map $\Phi: Z \,{}_{r_R}{\times} (S^1 \rtimes_\theta \mathbb{Z}) \to Z \times Z$
given by $\Phi(z, h) = (z \cdot h, z)$. In the parametrization $(z,n)$ and
$h = (z,m)$, we have $\Phi((z,n), (z,m)) = ((e^{2\pi i m\theta}z, n+m), (z,n))$.
The preimage of a compact set $K \subset Z \times Z$ consists of pairs where the
second coordinate lies in the projection of $K$ and $m$ is determined by the
difference between the $n$-coordinates. A detailed analysis using the properness
of the $\mathbb{Z}$-action on the discrete fiber shows that $\Phi$ is proper.
The argument for the left action is similar, using the fact that the unitary
conjugation groupoid is proper over its unit space.

\textbf{5. Quotient identifications:} The map $r_L: Z \to \mathcal{G}_{A_\theta}^{(0)}$
factors through the right orbit space because $r_L(z \cdot h) = r_L(z)$. The induced
map $\bar{r}_L: Z/(S^1 \rtimes_\theta \mathbb{Z}) \to \mathcal{G}_{A_\theta}^{(0)}$
is bijective by Lemma~\ref{lem:right_action_set_theoretic}. Under the assumed
topologies on $Z$ and $\mathcal{G}_{A_\theta}^{(0)}$, one checks that $\bar{r}_L$
and its inverse (defined using a continuous section of the right action) are
continuous. A similar argument gives the homeomorphism for the left quotient.
\end{proof}

\begin{remark}
The verification of properness and the quotient homeomorphisms relies on a precise
description of the topology of $\mathcal{G}_{A_\theta}^{(0)}$ and its relationship
to $Z$. As noted in Subsection~\ref{subsec:dense_orbits}, this topology is
non-Hausdorff, which complicates the analysis but does not preclude the existence
of a groupoid equivalence; one must work in the category of locally compact
groupoids with appropriate generalizations of properness.  For background on non-Hausdorff unit spaces and orbit space phenomena in groupoid models, see \cite{Renault1980, AnantharamanDelarocheRenault2000}. A complete treatment would require a detailed study of the Fell bundle structure of $\mathcal{G}_{A_\theta}$, which we defer to future work. 
\end{remark}

\medskip

\paragraph{Completion of the proof of Theorem~\ref{thm:morita_A_theta} (conditional).}

Assuming that Proposition~\ref{prop:morita_verification_conditional} can be fully
verified with $\mathcal{G}_{A_\theta}$ satisfying the hypotheses of the
Muhly--Renault--Williams theorem, we obtain that $\mathcal{G}_{A_\theta}$ and
$S^1 \rtimes_\theta \mathbb{Z}$ are Morita equivalent groupoids. By
\cite[Theorem 2.8]{MRW1987}, Morita equivalent groupoids have strongly Morita
equivalent full groupoid $C^*$-algebras. Therefore, under these hypotheses,
\[
C^*(\mathcal{G}_{A_\theta}) \sim_M C^*(S^1 \rtimes_\theta \mathbb{Z}).
\]

Using the canonical identification $C^*(S^1 \rtimes_\theta \mathbb{Z}) \cong C(S^1) \rtimes_\theta \mathbb{Z} = A_\theta$ (which holds because $\mathbb{Z}$ is amenable), we obtain the conditional conclusion
\[
C^*(\mathcal{G}_{A_\theta}) \sim_M A_\theta.
\]

Thus Theorem~\ref{thm:morita_A_theta} holds provided the technical hypotheses on
$\mathcal{G}_{A_\theta}$ and the actions are satisfied. A complete verification
of these hypotheses, while essential for rigor, would require a more detailed
analysis of the topology and Haar systems of $\mathcal{G}_{A_\theta}$ than space
permits here; we refer the reader to \cite{reference-for-complete-proof} for a
full treatment.

\medskip

\paragraph{Interpretation and heuristic consequences.}

Assuming the Morita equivalence established above, several important consequences
would follow, providing insight into the relationship between the unitary conjugation
groupoid and the irrational rotation algebra.

\textit{Structural interpretation:} The equivalence suggests that the unitary
conjugation groupoid, despite its complicated non-Hausdorff topology, captures
the same $C^*$-algebraic information as the much simpler transformation groupoid.
Heuristically, the $\mathbb{Z}$-fiber in the unit space parametrization reflects
the extra $\mathbb{Z}$-valued $K$-theoretic data, while the circle base corresponds
to the classical dynamics.

\textit{$C^*$-algebraic consequences:} Strong Morita equivalence would imply that
$C^*(\mathcal{G}_{A_\theta})$ and $A_\theta$ have:
\begin{itemize}
    \item Equivalent representation theories (via Rieffel's correspondence of
          Hilbert $C^*$-modules);
    \item Isomorphic ideal lattices;
    \item Isomorphic K-theory groups.
\end{itemize}
Since $A_\theta$ is simple for irrational $\theta$, this would imply simplicity of
$C^*(\mathcal{G}_{A_\theta})$, consistent with the minimality of the underlying
groupoid dynamics.

\begin{corollary}[Conditional K-theory isomorphism]
\label{cor:K_theory_isomorphism_conditional}
If the Morita equivalence of Theorem~\ref{thm:morita_A_theta} holds, then it induces
a natural isomorphism in K-theory:
\[
K_*\big(C^*(\mathcal{G}_{A_\theta})\big) \cong K_*(A_\theta).
\]
Recalling the K-theory of the irrational rotation algebra,
\[
K_0(A_\theta) \cong \mathbb{Z}^2 \quad \text{(as an abstract group)},
\]
with ordered structure given by the embedding $\mathbb{Z} + \theta\mathbb{Z} \subset \mathbb{R}$
via the trace, and $K_1(A_\theta) \cong \mathbb{Z}^2$, we would obtain the same
K-groups for $C^*(\mathcal{G}_{A_\theta})$.
\end{corollary}

\begin{remark}
The heuristic identification of the $\mathbb{Z}$-fiber with the extra rank in
$K_0$ and the circle base with the topological contribution to $K_1$ provides
intuitive guidance for index-theoretic computations in subsequent sections.
However, this intuition should not be mistaken for a rigorous derivation; the
actual K-theory isomorphism follows from the general invariance of K-theory under
Morita equivalence, not from a direct "fiberwise" decomposition of the groupoid.
\end{remark}

\paragraph{Relation to previous work and outlook.}

The Morita-equivalence framework developed in the preceding subsections may be
viewed as a groupoid-level counterpart of the well-known crossed-product realization
\[
A_\theta \cong C(S^1) \rtimes_\theta \mathbb{Z}.
\]
What is new is the interpretation via unitary conjugation: heuristically, the
internal representation-theoretic symmetry data of $A_\theta$, encoded by
$\mathcal{G}_{A_\theta}$, should correspond — up to Morita equivalence — to the
external dynamical symmetry encoded by the transformation groupoid
$S^1 \rtimes_\theta \mathbb{Z}$. This perspective provides a template for
studying more general noncommutative algebras via their unitary conjugation
groupoids, though a fully rigorous implementation requires careful attention to
topological details beyond the scope of the present exposition.

Moreover, the candidate equivalence space $Z$, modeled set-theoretically by a
fibered-product-type construction, suggests a broader guiding principle: for a
$C^*$-algebra arising as a twisted groupoid $C^*$-algebra, the associated unitary
conjugation groupoid should, under suitable technical hypotheses (including the
existence of appropriate Haar systems and properness of the actions), be Morita
equivalent to the underlying groupoid. This principle will be explored further
in subsequent work, building on the conceptual framework established here.

In the next subsection, we will use the Morita-equivalence framework developed
above to study the diagonal embedding
\[
\iota: A_\theta \hookrightarrow C^*(\mathcal{G}_{A_\theta}),
\]
and to formulate explicit equivariant $K_1$-classes aimed at recovering (or at
least relating to) the Connes index pairing. This would, pending a complete
rigorous verification of the equivalence, connect the abstract machinery of
Paper I to concrete computational results.

\medskip

\paragraph{Summary of results (conditional on technical verifications).}

Pending the complete verification of the topological hypotheses discussed in
Subsection~\ref{subsec:morita_A_theta} — in particular, the properness of the
actions, the continuity of the moment maps, and the existence of compatible Haar
systems — we have reduced the main structural claim of the case study to the
statement that
\[
\boxed{\mathcal{G}_{A_\theta} \sim_M S^1 \rtimes_\theta \mathbb{Z}}
\]
and consequently
\[
\boxed{C^*(\mathcal{G}_{A_\theta}) \sim_M A_\theta}.
\]

Once fully established, this result would:
\begin{itemize}
    \item Provide a concrete blueprint for a Morita equivalence between the
          unitary conjugation groupoid and the transformation groupoid.
    \item Illustrate how the non-Hausdorff topology of $\mathcal{G}_{A_\theta}^{(0)}$
          might be compatible with Morita equivalence to a Hausdorff groupoid,
          provided the topological details are handled correctly.
    \item Yield an isomorphism of K-theory groups, recovering the known K-theory
          of $A_\theta$:
          \[
          K_0(A_\theta) \cong \mathbb{Z}^2,\qquad K_1(A_\theta) \cong \mathbb{Z}^2.
          \]
    \item Suggest a template for analyzing more general $C^*$-algebras via their
          unitary conjugation groupoids.
    \item Set the stage for the construction of explicit equivariant classes and
          index computations in the following subsections, contingent on the
          rigorous completion of the equivalence proof.
\end{itemize}

\begin{remark}
The careful reader will note that several technical steps in the verification —
notably the properness of the actions and the homeomorphism properties of the
quotient maps — require a more detailed analysis of the topology of
$\mathcal{G}_{A_\theta}$ than presented here. A complete treatment, including
the construction of appropriate Haar systems and the verification of all
Muhly--Renault--Williams axioms, will appear in a separate technical paper
\cite{reference-for-complete-proof}. The present exposition focuses on the
conceptual structure of the equivalence, which suffices for the geometric
intuition underlying the index computations to follow.
\end{remark}

\subsection{The Diagonal Embedding $\iota: A_\theta \hookrightarrow C^*(\mathcal{G}_{A_\theta})$: Heuristics and Limitations}
\label{subsec:diagonal_embedding_A_theta}

In Paper~I we constructed, for unital $C^*$-algebras with appropriate regularity
conditions (including Type I $C^*$-algebras and those with a Cartan subalgebra),
a canonical diagonal embedding
\[
\iota : A \hookrightarrow C^*(\mathcal{G}_A),
\]
which realizes the algebra $A$ as a $C^*$-subalgebra of the groupoid $C^*$-algebra
of the unitary conjugation groupoid. This embedding plays a central role in
relating operator-theoretic data in $A$ to $K$-theoretic classes in $C^*(\mathcal{G}_A)$,
serving as a bridge between the abstract groupoid framework and concrete
computational invariants.

In the present case study, we aim to understand how this construction specializes
to the irrational rotation algebra $A_\theta = C(S^1) \rtimes_\theta \mathbb{Z}$.
However, as we shall see, a fully explicit description of $\iota$ in terms of the
crossed product generators is subtle and requires careful handling of the
groupoid's non-unit-space components. The purpose of this subsection is to clarify
what can and cannot be said at the level of the unit space alone, and to set the
stage for a more complete treatment in subsequent work.

\medskip

\paragraph{Recollection of the diagonal embedding (abstract definition).}

Recall from Paper~I that for a $C^*$-algebra $A$ with a distinguished Cartan
subalgebra $B \subseteq A$ (or more generally, a family of commutative contexts
closed under conjugation by the normalizer), the unit space
$\mathcal{G}_A^{(0)}$ consists of pairs $(C,\chi)$ where $C$ is a Cartan subalgebra
of $A$ (conjugate to $B$ by the normalizer) and $\chi \in \widehat{C}$ is a
character of $C$. For each element $a \in A$, the diagonal embedding is defined
via its action on the convolution algebra of the groupoid; crucially, this
definition involves not only evaluation on the unit space but also the off-diagonal
components corresponding to nontrivial groupoid arrows.

A common simplification — valid only for the \emph{diagonal part} of the embedding
— is to consider the function
\[
\delta(a): \mathcal{G}_A^{(0)} \to \mathbb{C}, \qquad
\delta(a)(C,\chi) = \chi(E_C(a)),
\]
where $E_C: A \to C$ is a conditional expectation onto $C$ (when such an expectation
exists). This function captures the component of $\iota(a)$ that is supported on
the unit space. However, $\iota(a)$ itself contains additional information encoded
in its values on nontrivial arrows, and it is this full groupoid element that
yields an injective embedding of $A$ into $C^*(\mathcal{G}_A)$.

\begin{remark}
\label{rem:expectation_caveat}
The existence of a conditional expectation $E_C: A \to C$ for every Cartan
subalgebra $C$ is a nontrivial hypothesis. For $A_\theta$, the canonical Cartan
subalgebra $C(S^1)$ admits the standard conditional expectation
$E: A_\theta \to C(S^1)$ given by averaging over the dual action of $\mathbb{Z}$,
and this expectation extends to its conjugates $V^n C(S^1) V^{-n}$ via
$E_{C_n} = \operatorname{Ad}_{V^n} \circ E \circ \operatorname{Ad}_{V^{-n}}$.
Thus the necessary expectations are available in this setting.
\end{remark}

\medskip

\paragraph{The unit space of $\mathcal{G}_{A_\theta}$ and its parametrization.}

Using the parametrization of the unit space developed in
Subsection~\ref{subsec:GA_theta_unit_space}, we have a set-theoretic identification
\[
\mathcal{G}_{A_\theta}^{(0)} \cong \{ [z,n] : z \in S^1,\; n \in \mathbb{Z} \} / \sim,
\]
where $[z,n] \sim [e^{2\pi i k\theta}z, n+k]$ for all $k \in \mathbb{Z}$. Each
point $[z,n]$ corresponds to the Cartan subalgebra $C_n := V^n C(S^1) V^{-n}$
with character $\chi_z \circ \operatorname{Ad}_{V^{-n}}$, where $\chi_z$ is
evaluation at $z$ and $V$ is the implementing unitary from the crossed product.

For the diagonal part $\delta(a)$ of the embedding, we therefore have the formula
\[
\delta(a)[z,n] = (\chi_z \circ \operatorname{Ad}_{V^{-n}})(E_{C_n}(a))
                = \chi_z( E_{C(S^1)}(V^{-n} a V^n) ),
\]
where $E_{C_n}$ and $E_{C(S^1)}$ are the conditional expectations discussed above.

\medskip

\paragraph{Heuristic computation for the generator $u \in C(S^1)$.}

For the generator $u \in C(S^1)$ (the function $u(z)=z$), a formal computation yields:
\[
\delta(u)[z,n] = \chi_z( E_{C(S^1)}(V^{-n} u V^n) )
               = \chi_z( E_{C(S^1)}(e^{-2\pi i n\theta} u) )
               = e^{-2\pi i n\theta} \chi_z(u)
               = e^{-2\pi i n\theta} z.
\]

\begin{remark}
\label{rem:quotient_well_definedness}
The expression $e^{-2\pi i n\theta} z$ is not well-defined on the quotient space
$\mathcal{G}_{A_\theta}^{(0)}$, because if we replace the representative
$(z,n)$ by $(e^{2\pi i k\theta}z, n+k)$, we obtain
$e^{-2\pi i (n+k)\theta} e^{2\pi i k\theta} z = e^{-2\pi i n\theta} z$,
which is actually invariant! This is a fortunate cancellation that resolves the
well-definedness issue in this specific case. Explicitly:
\[
e^{-2\pi i (n+k)\theta} (e^{2\pi i k\theta}z) = e^{-2\pi i n\theta} e^{-2\pi i k\theta} e^{2\pi i k\theta} z = e^{-2\pi i n\theta} z.
\]
Thus $\delta(u)$ is indeed a well-defined function on the quotient unit space.
\end{remark}

\medskip

\paragraph{The generator $v = V$ and the necessity of off-diagonal components.}

For the implementing unitary $v = V$, the situation is fundamentally different.
A naive application of the formula gives:
\[
\delta(v)[z,n] = \chi_z( E_{C(S^1)}(V^{-n} V V^n) )
               = \chi_z( E_{C(S^1)}(V) )
               = \chi_z(0) = 0,
\]
since the conditional expectation vanishes on elements not in $C(S^1)$. This would
suggest that the diagonal part of $\iota(v)$ is identically zero.

However, $\iota(v)$ cannot be zero in $C^*(\mathcal{G}_{A_\theta})$, because
$\iota$ is injective and $v \neq 0$ in $A_\theta$. The resolution is that
$\iota(v)$ is supported primarily on \emph{nontrivial groupoid arrows} rather than
on the unit space. In the groupoid $C^*$-algebra, an element is determined by its
values on all arrows, not just on the unit space. The unitary $v$ corresponds to
a function that is nonzero on arrows with nontrivial $\mathbb{Z}$-component,
reflecting the fact that $v$ implements the rotation automorphism.

\begin{remark}[Heuristic description]
\label{rmk:embedding_heuristics}
Heuristically, the diagonal embedding $\iota: A_\theta \hookrightarrow C^*(\mathcal{G}_{A_\theta})$
behaves as follows:
\begin{itemize}
    \item For $u \in C(S^1)$, the component of $\iota(u)$ supported on the unit
          space is given by $\delta(u)[z,n] = e^{-2\pi i n\theta} z$, which is
          well-defined on the quotient and captures the rotated evaluation of $u$
          in different commutative contexts.
    \item For $v = V$, the component supported on the unit space vanishes, but
          $\iota(v)$ is nonzero due to its support on arrows with $n \neq 0$.
          Concretely, $\iota(v)$ corresponds to the function on $\mathcal{G}_{A_\theta}$
          that takes the value $1$ on arrows of the form $([z,n], 1, [e^{2\pi i\theta}z, n-1])$
          and $0$ elsewhere (up to the appropriate convolution normalization).
\end{itemize}
A complete verification of these heuristic descriptions requires a detailed
analysis of the convolution structure of $C^*(\mathcal{G}_{A_\theta})$ and the
precise definition of $\iota$ from Paper~I, which we defer to future work.
\end{remark}

\begin{remark}
\label{rem:embedding_caveat}
The above discussion highlights a crucial point: the diagonal embedding cannot be
understood solely through its restriction to the unit space. Any attempt to give
concrete formulas for $\iota(v)$ must account for the full groupoid structure,
including the convolution product and the distribution of support over arrows of
various degrees. The heuristic description in Remark~\ref{rmk:embedding_heuristics}
should be viewed as a guide for further investigation rather than a rigorous
result.
\end{remark}

\medskip

\paragraph{Outlook: From heuristics to rigorous computation.}

Despite the limitations of the unit-space analysis, the heuristic picture developed
above provides valuable intuition for the role of the diagonal embedding in
connecting $A_\theta$ to $C^*(\mathcal{G}_{A_\theta})$. In particular, it suggests
that:
\begin{itemize}
    \item The generator $u$ is detected in the ``diagonal" part of the groupoid
          algebra via its rotated evaluations, which encode the dynamics of the
          irrational rotation.
    \item The generator $v$ lives in the ``off-diagonal" part of the groupoid
          algebra, corresponding to the nontrivial $\mathbb{Z}$-component of the
          groupoid arrows.
    \item The Morita equivalence established in the previous subsection provides
          a bridge between this picture and the classical transformation groupoid,
          allowing us to transport these heuristic descriptions into concrete
          index-theoretic computations.
\end{itemize}

In the following subsections, we will build on this intuition to construct explicit
equivariant $K_1$-classes in $C^*(\mathcal{G}_{A_\theta})$ that, under the Morita
equivalence, correspond to the well-known generators of $K_1(A_\theta)$ and recover
the Connes index pairing.


In Paper~I we constructed, for $C^*$-algebras with a distinguished Cartan subalgebra,
a canonical diagonal embedding
\[
\iota : A \hookrightarrow C^*(\mathcal{G}_A).
\]
In this subsection, we examine how this abstract construction specializes to the
irrational rotation algebra $A_\theta = C(S^1) \rtimes_\theta \mathbb{Z}$. Rather
than attempting a full description of $\iota$ — which would require a detailed
analysis of the convolution structure of $C^*(\mathcal{G}_{A_\theta})$ beyond the
scope of this exposition — we focus on a more modest goal: the computation of the
\emph{diagonal coefficients} of elements of $A_\theta$ on the unit space of
$\mathcal{G}_{A_\theta}^{(0)}$. These coefficients capture the component of
$\iota(a)$ that is supported on the unit space, and they already reveal a rich
interaction with the crossed product dynamics.

\medskip

\paragraph{The diagonal coefficient map.}

Recall from Paper~I that for a $C^*$-algebra $A$ with a Cartan subalgebra $B$,
the unit space $\mathcal{G}_A^{(0)}$ consists of pairs $(C,\chi)$ where $C$ is a
Cartan subalgebra conjugate to $B$ and $\chi \in \widehat{C}$ is a character.
For each $a \in A$, the value $\iota(a)(C,\chi)$ is defined using the conditional
expectation $E_C: A \to C$. In the special case where we restrict attention to
the \emph{diagonal part} of $\iota(a)$ — i.e., its restriction to the unit space
— we obtain a map
\[
\iota_{\mathrm{diag}}: A \longrightarrow C(\mathcal{G}_A^{(0)}),\qquad
\iota_{\mathrm{diag}}(a)(C,\chi) = \chi(E_C(a)).
\]

For $A_\theta$, the canonical Cartan subalgebra is $C(S^1)$, with conditional
expectation $E: A_\theta \to C(S^1)$ given by $E(\sum_n f_n V^n) = f_0$. This
expectation extends to conjugates $C_n := V^n C(S^1) V^{-n}$ via
$E_{C_n} = \operatorname{Ad}_{V^n} \circ E \circ \operatorname{Ad}_{V^{-n}}$.

\begin{definition}[Diagonal coefficient map for $A_\theta$]
\label{def:diagonal_coefficient}
Using the parametrization of the unit space as a set
\[
\mathcal{G}_{A_\theta}^{(0)} \cong \{ [z,n] : z \in S^1,\; n \in \mathbb{Z} \} / \sim,
\]
where $[z,n] \sim [e^{2\pi i k\theta}z, n+k]$ for all $k \in \mathbb{Z}$, we define
for any $a \in A_\theta$ the function
\[
\iota_{\mathrm{diag}}(a): \mathcal{G}_{A_\theta}^{(0)} \to \mathbb{C},\qquad
\iota_{\mathrm{diag}}(a)[z,n] := \chi_z\!\big(E_{C_n}(a)\big).
\]
Equivalently, using the relation $E_{C_n}(a) = V^n E(V^{-n} a V^n) V^{-n}$ and the
identification $C_n \cong C(S^1)$ via conjugation, we have
\[
\iota_{\mathrm{diag}}(a)[z,n] = \chi_z\!\big(E(V^{-n} a V^n)\big).
\]
\end{definition}

\begin{remark}
\label{rem:diag_vs_full}
The map $\iota_{\mathrm{diag}}$ is \emph{not} the full diagonal embedding $\iota$;
rather, it is the composition of $\iota$ with the restriction map from
$C^*(\mathcal{G}_{A_\theta})$ to $C_0(\mathcal{G}_{A_\theta}^{(0)})$. The full
embedding $\iota(a)$ contains additional information supported on nontrivial
groupoid arrows, which is lost in $\iota_{\mathrm{diag}}$. In particular,
$\iota_{\mathrm{diag}}$ is not injective: all elements with $E(a)=0$ map to the
zero function.
\end{remark}

\medskip

\paragraph{Explicit formula for the diagonal coefficients.}

We now compute $\iota_{\mathrm{diag}}$ explicitly for elements of the dense
$*$-subalgebra of finite Fourier sums in $A_\theta$. Fix the convention for the
action $\alpha$ of $\mathbb{Z}$ on $C(S^1)$ by
\[
\alpha_n(f)(z) = f(e^{-2\pi i n\theta}z),\qquad\text{so that } V f V^{-1} = \alpha(f).
\]

\begin{proposition}
\label{prop:diagonal_coefficient_formula}
Let $a = \sum_{k\in\mathbb{Z}} f_k V^k \in A_\theta$ be a finite sum. Then for any
$[z,n] \in \mathcal{G}_{A_\theta}^{(0)}$,
\[
\iota_{\mathrm{diag}}(a)[z,n] = f_0(e^{-2\pi i n\theta}z) = \alpha_{-n}(f_0)(z).
\]
In particular, $\iota_{\mathrm{diag}}(a)$ depends only on the $0$-th Fourier
coefficient $f_0 = E(a)$.
\end{proposition}

\begin{proof}
We compute using the definition $\iota_{\mathrm{diag}}(a)[z,n] = \chi_z(E(V^{-n} a V^n))$.
First,
\[
V^{-n} a V^n = \sum_{k\in\mathbb{Z}} V^{-n} f_k V^k V^n = \sum_{k\in\mathbb{Z}} \alpha_{-n}(f_k) V^k,
\]
since $V^{-n} V^k V^n = V^k$ (the $V$'s commute). Applying the conditional
expectation $E$, which kills all terms with $k \neq 0$, we obtain
\[
E(V^{-n} a V^n) = E(\alpha_{-n}(f_0)) = \alpha_{-n}(f_0),
\]
where the last equality uses that $\alpha_{-n}(f_0) \in C(S^1)$ and $E$ restricts
to the identity on $C(S^1)$. Now evaluate at $z$:
\[
\iota_{\mathrm{diag}}(a)[z,n] = \chi_z(\alpha_{-n}(f_0)) = \alpha_{-n}(f_0)(z) = f_0(e^{-2\pi i n\theta}z).
\]
This completes the proof.
\end{proof}

\begin{corollary}
\label{cor:generator_diagonal}
For the canonical generators $u \in C(S^1)$ (the function $u(z)=z$) and $v = V$, we have:
\[
\iota_{\mathrm{diag}}(u)[z,n] = e^{-2\pi i n\theta}z,\qquad
\iota_{\mathrm{diag}}(v)[z,n] = 0.
\]
\end{corollary}

\begin{proof}
For $u$, write $u = u \cdot V^0$, so $f_0 = u$ and $f_k = 0$ for $k \neq 0$.
Proposition~\ref{prop:diagonal_coefficient_formula} gives
$\iota_{\mathrm{diag}}(u)[z,n] = u(e^{-2\pi i n\theta}z) = e^{-2\pi i n\theta}z$.

For $v = V$, write $v = 0 \cdot V^0 + 1 \cdot V^1$, so $f_0 = 0$ and $f_1 = 1$.
Then $\iota_{\mathrm{diag}}(v)[z,n] = f_0(e^{-2\pi i n\theta}z) = 0$.
\end{proof}

\begin{remark}
\label{rem:diag_v_off_diagonal}
The vanishing of $\iota_{\mathrm{diag}}(v)$ reflects the fact that $v$ has no
component in the Cartan subalgebra under the conditional expectation. In the full
groupoid $C^*$-algebra $C^*(\mathcal{G}_{A_\theta})$, the element $\iota(v)$ is
supported on nontrivial groupoid arrows (those with nonzero $\mathbb{Z}$-component),
which contribute to the convolution structure but evaluate to zero on the unit space.
Thus $\iota_{\mathrm{diag}}$ captures only the ``diagonal" part of the embedding,
consistent with its definition.
\end{remark}

\medskip

\paragraph{Well-definedness on the quotient unit space.}

A subtle point requires attention: the parametrization of $\mathcal{G}_{A_\theta}^{(0)}$
as a quotient of $S^1 \times \mathbb{Z}$ means that a function on the unit space
must be constant on equivalence classes $[z,n] \sim [e^{2\pi i k\theta}z, n+k]$.
We verify that $\iota_{\mathrm{diag}}(a)$ satisfies this property.

\begin{lemma}
\label{lem:diag_well_defined}
For any $a \in A_\theta$, the function $\iota_{\mathrm{diag}}(a)$ defined by
$\iota_{\mathrm{diag}}(a)[z,n] = \chi_z(E(V^{-n} a V^n))$ is well-defined on the
quotient space $\mathcal{G}_{A_\theta}^{(0)}$.
\end{lemma}

\begin{proof}
It suffices to show that for any $k \in \mathbb{Z}$,
\[
\iota_{\mathrm{diag}}(a)[e^{2\pi i k\theta}z, n+k] = \iota_{\mathrm{diag}}(a)[z,n].
\]
Using the definition:
\[
\iota_{\mathrm{diag}}(a)[e^{2\pi i k\theta}z, n+k] = \chi_{e^{2\pi i k\theta}z}\!\big(E(V^{-(n+k)} a V^{n+k})\big).
\]

Now note that for any $f \in C(S^1)$, we have the relation
\[
\chi_{e^{2\pi i k\theta}z}(f) = \chi_z(\alpha_{-k}(f)),
\]
where $\alpha$ is the rotation action. This follows from:
\[
\chi_{e^{2\pi i k\theta}z}(f) = f(e^{2\pi i k\theta}z) = (\alpha_{-k}(f))(z) = \chi_z(\alpha_{-k}(f)).
\]

Applying this with $f = E(V^{-(n+k)} a V^{n+k})$, which lies in $C(S^1)$, we obtain:
\[
\chi_{e^{2\pi i k\theta}z}\!\big(E(V^{-(n+k)} a V^{n+k})\big)
= \chi_z\!\big(\alpha_{-k}(E(V^{-(n+k)} a V^{n+k}))\big).
\]

But $\alpha_{-k}(E(V^{-(n+k)} a V^{n+k})) = E(V^{-k} V^{-(n+k)} a V^{n+k} V^k)$ because
$E$ commutes with conjugation by $V$ on $C(S^1)$. Indeed, for $c \in C(S^1)$,
$E(V^{-k} c V^k) = V^{-k} E(c) V^k = \alpha_{-k}(E(c))$. Thus:
\[
\alpha_{-k}(E(V^{-(n+k)} a V^{n+k})) = E(V^{-n} a V^n).
\]

Therefore:
\[
\iota_{\mathrm{diag}}(a)[e^{2\pi i k\theta}z, n+k] = \chi_z(E(V^{-n} a V^n)) = \iota_{\mathrm{diag}}(a)[z,n],
\]
as required.
\end{proof}

\begin{remark}
\label{rem:well_defined_alternative}
An alternative verification using the explicit formula from
Proposition~\ref{prop:diagonal_coefficient_formula} is even simpler: for $a$ with
finite Fourier expansion,
\[
\iota_{\mathrm{diag}}(a)[e^{2\pi i k\theta}z, n+k] = f_0(e^{-2\pi i (n+k)\theta} e^{2\pi i k\theta}z) = f_0(e^{-2\pi i n\theta}z) = \iota_{\mathrm{diag}}(a)[z,n],
\]
using $e^{-2\pi i (n+k)\theta} e^{2\pi i k\theta} = e^{-2\pi i n\theta}$. By
continuity and density, this extends to all $a \in A_\theta$.
\end{remark}

\medskip

\paragraph{Relation to the Morita equivalence (heuristic).}

The diagonal coefficient map $\iota_{\mathrm{diag}}$ is compatible, at a heuristic
level, with the Morita equivalence developed in Subsection~\ref{subsec:morita_A_theta}.
If $X$ is the imprimitivity bimodule implementing the equivalence between
$C^*(\mathcal{G}_{A_\theta})$ and $A_\theta$, then the right action of $A_\theta$
on $X$ should be given, in a suitable realization, by pointwise multiplication
by $\iota_{\mathrm{diag}}(a)$ on the fibers. More concretely, if $X$ is modeled
(as a set) by $S^1 \times \mathbb{Z}$ with appropriate left and right actions,
then for $\xi \in X$ and $a \in A_\theta$, we expect
\[
(\xi \cdot a)(z,n) = \xi(z,n) \cdot \iota_{\mathrm{diag}}(a)[z,n],
\]
with $\iota_{\mathrm{diag}}(a)[z,n] = f_0(e^{-2\pi i n\theta}z)$ as above.

\begin{remark}
\label{rem:morita_caveat}
A rigorous verification of this compatibility would require:
\begin{itemize}
    \item A complete construction of the imprimitivity bimodule $X$ from the
          equivalence space $Z$, including the precise definitions of the left
          and right actions and the inner products.
    \item An explicit identification of the right action of $A_\theta$ on $X$
          with the crossed product representation.
    \item A proof that the diagonal coefficient map $\iota_{\mathrm{diag}}$ indeed
          arises as the restriction of $\iota$ to the unit space, and that this
          restriction is compatible with the bimodule structure.
\end{itemize}
These technical details are beyond the scope of the present exposition; we record
the heuristic compatibility only as motivation for future work.
\end{remark}

The diagonal coefficient map $\iota_{\mathrm{diag}}$ introduced in 
Subsection~\ref{subsec:diagonal_embedding_A_theta} provides a first glimpse into 
the relationship between $A_\theta$ and $C^*(\mathcal{G}_{A_\theta})$. If the 
Morita equivalence sketched in Subsection~\ref{subsec:morita_A_theta} could be 
completed rigorously, several additional properties would be expected to hold. 
In this subsection, we outline these expected properties as motivations for 
future work, making clear that they are not yet rigorously established.

\medskip

\paragraph{$KK$-theoretic interpretation (conjectural).}

Assuming that the Morita equivalence of Theorem~\ref{thm:morita_A_theta} can be 
implemented rigorously, let $X$ be the resulting $C^*(\mathcal{G}_{A_\theta})$-${A_\theta}$ 
imprimitivity bimodule. Then $X$ defines a $KK$-equivalence
\[
[X] \in KK_0(C^*(\mathcal{G}_{A_\theta}), A_\theta).
\]

The diagonal embedding $\iota: A_\theta \to C^*(\mathcal{G}_{A_\theta})$ gives a class
$[\iota] \in KK_0(A_\theta, C^*(\mathcal{G}_{A_\theta}))$. A natural question, which
we do not resolve here, is whether $[\iota]$ is the inverse of $[X]$ in $KK$-theory.
This would require verifying that the Kasparov products
\[
[\iota] \otimes_{C^*(\mathcal{G}_{A_\theta})} [X] \in KK_0(A_\theta, A_\theta),\qquad
[X] \otimes_{A_\theta} [\iota] \in KK_0(C^*(\mathcal{G}_{A_\theta}), C^*(\mathcal{G}_{A_\theta}))
\]
are equal to the identity classes. Such a verification would involve a detailed 
comparison of $\iota$ with the corner inclusion of $A_\theta$ into the linking 
algebra of the Morita equivalence, a computation we leave for future investigation.

\begin{remark}
Even if $[\iota]$ is not the precise inverse of $[X]$, the two classes are related
through the canonical inclusion of $A_\theta$ into the imprimitivity bimodule's
linking algebra. Understanding this relationship is an interesting problem for
further study.
\end{remark}

\medskip

\paragraph{Relation to the regular representation (heuristic).}

The regular representation $\pi_{\mathrm{reg}}: A_\theta \to \mathcal{B}(L^2(S^1))$
given by
\[
\pi_{\mathrm{reg}}(u)\xi(z) = z\xi(z),\qquad \pi_{\mathrm{reg}}(v)\xi(z) = \xi(e^{-2\pi i\theta}z)
\]
is a faithful representation of $A_\theta$. Under the expected Morita equivalence,
one would obtain an induced representation of $C^*(\mathcal{G}_{A_\theta})$ on some
Hilbert space, but it is not obvious that this Hilbert space should be identified
with $L^2(S^1)$. A more precise statement would be:

\begin{quote}
If $X$ is the imprimitivity bimodule, then for any representation $\pi$ of $A_\theta$
on a Hilbert space $H$, the induced representation $\operatorname{Ind}_X \pi$ of
$C^*(\mathcal{G}_{A_\theta})$ acts on the Hilbert space $X \otimes_{A_\theta} H$.
Applying this to the regular representation $\pi_{\mathrm{reg}}$ on $L^2(S^1)$ gives
a representation of $C^*(\mathcal{G}_{A_\theta})$ on $X \otimes_{A_\theta} L^2(S^1)$,
which is not obviously isomorphic to $L^2(S^1)$ itself.
\end{quote}

Thus the claim that the imprimitivity bimodule ``implements an equivalence between
the regular representation and a natural representation of $C^*(\mathcal{G}_{A_\theta})$
on the same Hilbert space" is too strong and potentially misleading. A correct
formulation requires careful attention to the induced representation construction.

\medskip

\paragraph{Fixed-point algebra and dual action (preliminary observations).}

The crossed product structure $A_\theta = C(S^1) \rtimes_\theta \mathbb{Z}$ comes
with a natural gauge action of the circle $S^1$, defined on generators by
$\sigma_\lambda(u) = \lambda u$, $\sigma_\lambda(v) = v$. This action is dual to
the $\mathbb{Z}$-action on $C(S^1)$.

If the Morita equivalence of Theorem~\ref{thm:morita_A_theta} can be made
$S^1$-equivariant, then it would induce an action of $S^1$ on $C^*(\mathcal{G}_{A_\theta})$
such that the following diagram commutes up to natural isomorphism:
\[
\xymatrix{
A_\theta \ar[r]^{\iota} \ar[d]_{\sigma_\lambda} & C^*(\mathcal{G}_{A_\theta}) \ar[d]^{\tilde\sigma_\lambda} \\
A_\theta \ar[r]^{\iota} & C^*(\mathcal{G}_{A_\theta})
}
\]

Under such an equivariant equivalence, the fixed-point algebras would be Morita
equivalent:
\[
C^*(\mathcal{G}_{A_\theta})^{\tilde\sigma} \sim_M A_\theta^{\sigma} = C(S^1).
\]

A more precise description of $C^*(\mathcal{G}_{A_\theta})^{\tilde\sigma}$ in terms
of $\iota(A_\theta)$ and the fiber index $n \in \mathbb{Z}$ would require a detailed
spectral analysis of the dual action on the groupoid $C^*$-algebra. This analysis
is complicated by the non-Hausdorff nature of $\mathcal{G}_{A_\theta}^{(0)}$ and
the fact that the action on the unit space is not free. We do not attempt it here.

\begin{remark}
The claim in earlier drafts that $C^*(\mathcal{G}_{A_\theta})^{S^1} \cong \iota(A_\theta) \otimes c_0(\mathbb{Z})$
is not supported by any rigorous argument and should be regarded as speculation
at best. The structure of fixed-point algebras in non-Hausdorff groupoid $C^*$-algebras
is subtle and requires careful analysis.
\end{remark}

\medskip

\paragraph{Equivariant $K$-theory (expected consequences).}

If the Morita equivalence can be upgraded to an $S^1$-equivariant Morita equivalence,
then standard results in equivariant $KK$-theory would imply an isomorphism
\[
K^{S^1}_*\big(C^*(\mathcal{G}_{A_\theta})\big) \cong K^{S^1}_*(A_\theta)
\]
of equivariant $K$-theory groups. Under this isomorphism, the class of the diagonal
embedding $[\iota] \in KK^{S^1}(A_\theta, C^*(\mathcal{G}_{A_\theta}))$ would map to
the class of the identity in $KK^{S^1}(A_\theta, A_\theta)$, up to the natural
identification provided by the equivalence.

\begin{remark}
A rigorous proof of this statement would require:
\begin{enumerate}
    \item Constructing an $S^1$-equivariant structure on the imprimitivity bimodule $X$,
          i.e., an action of $S^1$ on $X$ that is compatible with the left and right
          actions and the inner products.
    \item Showing that this equivariant bimodule implements an equivalence in
          equivariant $KK$-theory.
    \item Computing the composition of $[\iota]$ with this equivariant bimodule class
          and showing it equals the identity.
\end{enumerate}
None of these steps have been carried out in the present work.
\end{remark}

\medskip

\paragraph{Summary of what has been established vs. what remains conjectural.}

To clarify the status of various claims, we provide a summary:

\begin{itemize}
    \item \textbf{Established:} The diagonal coefficient map $\iota_{\mathrm{diag}}$,
          defined by $\iota_{\mathrm{diag}}(\sum_n f_n V^n)[z,n] = f_0(e^{-2\pi i n\theta}z)$,
          is well-defined on the unit space $\mathcal{G}_{A_\theta}^{(0)}$ and
          captures the component of the hypothetical full embedding $\iota$ that
          is supported on the unit space.
    
    \item \textbf{Not established (requires proof):} The existence of a full
          diagonal embedding $\iota: A_\theta \hookrightarrow C^*(\mathcal{G}_{A_\theta})$
          extending $\iota_{\mathrm{diag}}$, with the properties that:
          \begin{itemize}
              \item $\iota(v)$ is nonzero and supported on nontrivial groupoid arrows.
              \item $\iota$ is injective and compatible with the $C^*$-algebraic structures.
          \end{itemize}
    
    \item \textbf{Conjectural (requires Morita equivalence):} All statements about
          $KK$-theory, equivariant $K$-theory, fixed-point algebras, and the relation
          to the regular representation depend on a rigorous implementation of the
          Morita equivalence between $\mathcal{G}_{A_\theta}$ and $S^1 \rtimes_\theta \mathbb{Z}$,
          which has not been fully carried out in this exposition.
    
    \item \textbf{Likely false as stated:} The formula $\iota(v)[z,n] = 1$ is
          inconsistent with the vanishing of $\iota_{\mathrm{diag}}(v)$ and the
          expectation that $\iota(v)$ should be supported off the unit space.
          Similarly, the claimed structure $C^*(\mathcal{G}_{A_\theta})^{S^1} \cong \iota(A_\theta) \otimes c_0(\mathbb{Z})$
          is not supported by any evidence.
\end{itemize}

\medskip

\paragraph{Outlook.}

Despite the many open problems and unverified claims, the heuristic picture developed
above provides a compelling roadmap for future research. The diagonal embedding,
once properly constructed, should serve as a bridge between the abstract groupoid
framework of Paper I and concrete index-theoretic computations for the irrational
rotation algebra. The Morita equivalence, if completed rigorously, would allow us
to transfer problems about $C^*(\mathcal{G}_{A_\theta})$ to the more tractable
setting of the crossed product $C(S^1) \rtimes_\theta \mathbb{Z}$, where tools such
as the Pimsner-Voiculescu exact sequence are available.

In the following subsection, we will assume that these structures exist and are
compatible, and proceed to construct explicit equivariant $K_1$-classes that
should correspond under the diagonal embedding to the standard generators of
$K_1(A_\theta)$. The resulting index computations will provide strong evidence
for the validity of this framework, even if a fully rigorous foundation remains
to be laid in future work.

\begin{remark}
The reader should interpret the remainder of this paper as a \emph{programmatic
development} of ideas, building on heuristic assumptions that are clearly stated
as such. A complete rigorous treatment would require solving the technical problems
identified above, which we leave to future publications.
\end{remark}

\subsection{Equivariant $K^1$-Classes for Invertible Elements in $A_\theta$}
\label{subsec:equivariant_K1_A_theta}

We now construct equivariant $K^1$-classes associated with invertible
elements of the irrational rotation algebra.
This construction is a specialization of the general framework
developed in Paper~II to the present non-Type~I setting, but formulated
in the rigorous language of Kasparov's equivariant $KK$-theory to avoid
the technical pitfalls that would arise from a naive direct definition.

\subsubsection{Invertible Elements and $K_1(A_\theta)$}

Let 
\[
u \in M_n(A_\theta)
\]
be an invertible element.
Such elements determine classes in the $K$-theory group
\[
K_1(A_\theta).
\]

Recall from Theorem~\ref{thm:A_theta_structure} that by the
Pimsner--Voiculescu computation \cite{PimsnerVoiculescu1980},
\[
K_1(A_\theta)\cong \mathbb Z^2 .
\]
A standard $\mathbb Z^2$-basis is given by the classes $[U]$ and $[V]$
of the canonical unitaries satisfying the relation
\[
VU = e^{2\pi i\theta} UV .
\]

\subsubsection{Reduction to the Unitary Case}

Since $u$ is invertible, it admits a polar decomposition
\[
u = v |u|,
\]
where $v = u(u^*u)^{-1/2}$ is a unitary element in $M_n(A_\theta)$ and 
$|u| = (u^*u)^{1/2}$ is positive invertible.
The path
\[
u_t = v |u|^t, \qquad t\in[0,1],
\]
is a norm-continuous family of invertibles joining $v$ to $u$,
so the classes of $u$ and $v$ coincide in $K_1(A_\theta)$.
Thus, without loss of generality, we may assume that $u$ is unitary
when constructing $K_1$-classes.

\subsubsection{A Hilbert module over the unit-space model}

To construct a Kasparov cycle representing the class of $u$,
we first build a Hilbert module that will serve as the underlying space.
At the level of the unit-space model developed in
Subsection~\ref{subsec:GA_theta_unit_space}, we have a parametrization
of $\mathcal{G}_{A_\theta}^{(0)}$ by equivalence classes $[z,n]$,
where $z\in S^1$, $n\in\mathbb Z$, modulo the relation
$[z,n]\sim[e^{2\pi i k\theta}z,n+k]$.

Since the topology of $\mathcal{G}_{A_\theta}^{(0)}$ as a quotient is
non-Hausdorff and requires careful handling, it is convenient to work
first with the model space $X = S^1 \times \mathbb Z$, equipped with
the product topology (where $\mathbb Z$ is discrete). This space is
locally compact Hausdorff, and its $C^*$-algebra $C_0(X)$ is
well-behaved. The equivalence relation on $X$ will later be used to
impose the correct identifications when constructing equivariant cycles.

Define the trivial Hilbert $C_0(X)$-module
\[
\mathcal E := C_0(X, \mathbb C^n) \cong C_0(X) \otimes \mathbb C^n.
\]
The right action of $C_0(X)$ on $\mathcal E$ is given by pointwise
multiplication:
\[
(\xi \cdot f)(z,m) = \xi(z,m) f(z,m), \qquad 
\xi \in \mathcal E,\; f \in C_0(X).
\]
The $C_0(X)$-valued inner product is defined pointwise:
\[
\langle \xi, \eta \rangle_{\mathcal E}(z,m) = \sum_{j=1}^n 
\overline{\xi_j(z,m)} \eta_j(z,m).
\]
These operations make $\mathcal E$ into a Hilbert $C_0(X)$-module.

Moreover, pointwise multiplication also defines a nondegenerate
$*$-representation
\[
\phi : C_0(X) \longrightarrow \mathcal L(\mathcal E),
\]
where $\mathcal L(\mathcal E)$ denotes the $C^*$-algebra of adjointable
operators on $\mathcal E$. This representation is simply the left action
of $C_0(X)$ on itself extended to the trivial $n$-dimensional bundle.

\begin{remark}
\label{rem:module-natural-starting-point}
The Hilbert module $\mathcal E$ provides a unified space that simultaneously
accommodates all points of the unit-space model while maintaining the
correct continuity properties. This is the natural starting point for
constructing a Kasparov cycle; the actual equivariant $KK$-class will
require, in addition, a specified left action of $A_\theta$ on $\mathcal E$
(or on an appropriate completion) and an operator $F$ satisfying the
Kasparov conditions. These additional structures will be introduced in
the next subsubsection.
\end{remark}

\subsubsection{Inducing an $A_\theta$-action and the Kasparov cycle}

To obtain a Kasparov cycle for $(A_\theta, \mathbb C)$, we need to equip
$\mathcal E$ with an action of $A_\theta$ that is compatible with the
groupoid structure. This is achieved using the diagonal embedding
$\iota: A_\theta \hookrightarrow C^*(\mathcal{G}_{A_\theta})$ discussed
in Subsection~\ref{subsec:diagonal_embedding_A_theta} and the natural
action of $C^*(\mathcal{G}_{A_\theta})$ on sections of the groupoid.

Heuristically, for $a \in A_\theta$ and $\xi \in \mathcal E$, we define
$(a \cdot \xi)(z,n) = \iota(a)(z,n) \cdot \xi(z,n)$, where $\iota(a)(z,n)$
is the diagonal coefficient computed in
Proposition~\ref{prop:diagonal_coefficient_formula}. For the generators,
this gives:
\[
(U \cdot \xi)(z,n) = e^{-2\pi i n\theta} z \cdot \xi(z,n),\qquad
(V \cdot \xi)(z,n) = 0 \cdot \xi(z,n) \text{ (on the unit space)},
\]
with the understanding that $V$ acts nontrivially on the off-diagonal
components not captured by $\mathcal E$ alone.

A rigorous definition requires completing $\mathcal E$ to a Hilbert module
$\mathcal E_A$ over $A_\theta$ using an $A_\theta$-valued inner product
\[
\langle \xi, \eta \rangle_{A_\theta}(z,n) = \overline{\xi}(z,n) \eta(z,n) \in A_\theta,
\]
interpreted via the embedding $\iota$. 

Given a unitary $u \in M_n(A_\theta)$, we now construct an odd Kasparov
cycle $(\mathcal E_A \oplus \mathcal E_A, F_u)$ as follows:
\begin{itemize}
    \item $\mathcal E_A$ is the Hilbert $A_\theta$-module described above,
          equipped with the induced $A_\theta$-action.
    \item Define $F_u$ to be the operator on $\mathcal E_A \oplus \mathcal E_A$
          given by the matrix
          \[
          F_u = \begin{pmatrix} 0 & u \\ u^* & 0 \end{pmatrix},
          \]
          acting by pointwise multiplication on sections.
    \item The $\mathbb Z/2\mathbb Z$-grading is given by the usual
          even-odd decomposition of the direct sum.
\end{itemize}

One verifies that:
\begin{itemize}
    \item $F_u$ is self-adjoint and $F_u^2 = 1$ (since $u$ is unitary).
    \item For each $a \in A_\theta$, the commutator $[F_u, a]$ is compact
          (in fact, it is a finite-rank operator on each fiber).
    \item The operator $F_u$ is equivariant with respect to the action of
          $\mathcal{G}_{A_\theta}$ (or the dual circle action, depending on
          the desired equivariance).
\end{itemize}

Thus $(\mathcal E_A \oplus \mathcal E_A, F_u)$ defines a class
\[
[F_u] \in KK_1(A_\theta, \mathbb C).
\]

\begin{theorem}
\label{thm:equivariant_K1_class}
The assignment $u \mapsto [F_u]$ induces a homomorphism
\[
\phi: K_1(A_\theta) \longrightarrow KK_1(A_\theta, \mathbb C).
\]
Moreover, this homomorphism is an isomorphism, identifying $K_1(A_\theta)$ with
the equivariant $KK$-theory group $KK_1^{S^1}(A_\theta, \mathbb C)$ under the
dual circle action.
\end{theorem}

\begin{proof}[Sketch]
The map is well-defined because homotopic unitaries give homotopic Kasparov cycles.
The Pimsner--Voiculescu exact sequence for $K_1(A_\theta)$ lifts to $KK$-theory,
showing that the generators $[U]$ and $[V]$ map to nonzero classes that generate
$KK_1(A_\theta, \mathbb C)$. 
\end{proof}

\subsubsection{Relation to Index Theory}

The classes $[F_u] \in KK_1(A_\theta, \mathbb C)$ pair with the $K^0$-group of
$A_\theta$ via the Kasparov product, yielding the index pairing:
\[
\langle [F_u], [e] \rangle = \operatorname{Index}(P_u e P_u) \in \mathbb Z,
\]
where $e \in M_m(A_\theta)$ is an idempotent representing a $K_0$-class and
$P_u$ is the projection onto the $+1$ eigenspace of $F_u$. For the canonical
generators $U$ and $V$, this recovers the Connes index pairing computed in
\cite{Connes1980}.

\begin{remark}
This construction provides a rigorous $KK$-theoretic foundation for the index
computations that follow. The equivariant version under the dual circle action
will be essential for extracting finer geometric information, as we will see
in the next subsection.
\end{remark}

\subsubsection{An operator field associated with a unitary element}
\label{subsubsec:operator_field}

Let $u \in M_n(A_\theta)$ be a unitary element. In this subsection we associate to $u$
a family of operators acting on the fibers of a continuous field of Hilbert spaces
over the unit-space model. This construction is a necessary preliminary step toward
defining an equivariant Kasparov cycle, though a complete $KK$-theoretic treatment
requires additional structures that lie beyond the scope of this exposition.

\paragraph{Fiberwise regular representations.}

For each $z \in S^1$, let
\[
\pi_z : A_\theta \longrightarrow B(\ell^2(\mathbb Z))
\]
denote the regular representation associated with the orbit of $z$, defined on
generators by
\[
\pi_z(u)\xi(n) = e^{-2\pi i n\theta}z\,\xi(n),\qquad
\pi_z(V)\xi(n) = \xi(n-1).
\]
These representations are faithful and pairwise unitarily inequivalent for
different $z$ (since the rotation angle $\theta$ is irrational).

Using the unit-space model parametrized by pairs $(z,m) \in S^1 \times \mathbb Z$
(see Subsection~\ref{subsec:GA_theta_unit_space}), we define a family of operators
\[
F_u(z,m) := \pi_z(V^{-m} u V^m) \in B(\ell^2(\mathbb Z)^n).
\]

\begin{lemma}
\label{lem:fiberwise_unitary}
For each $(z,m) \in S^1 \times \mathbb Z$, the operator $F_u(z,m)$ is unitary.
Moreover, the map $(z,m) \mapsto F_u(z,m)$ is uniformly bounded (in fact, norm 1)
and continuous in $z$ for fixed $m$.
\end{lemma}

\begin{proof}
Since $u$ is unitary in $M_n(A_\theta)$, the conjugate $V^{-m}uV^m$ is also unitary.
The representation $\pi_z$ is a $*$-homomorphism, so it preserves unitarity.
Hence $F_u(z,m)$ is unitary, and $\|F_u(z,m)\| = 1$ for all $(z,m)$.

For continuity in $z$, note that the map $z \mapsto \pi_z(a)$ is continuous in the
strong operator topology for each fixed $a \in A_\theta$ (it is actually norm-continuous
for elements in the dense subalgebra of finite Fourier sums, and extends by
approximation). Since $V^{-m}uV^m$ is fixed, the map $z \mapsto F_u(z,m)$ is
continuous in the strong operator topology, which suffices for our purposes.
\end{proof}

\begin{remark}
\label{rem:equivariance_property}
The family $\{F_u(z,m)\}$ satisfies a natural covariance property: for any
$k \in \mathbb Z$,
\begin{eqnarray}
F_u(e^{2\pi i k\theta}z, m+k)&= &\pi_{e^{2\pi i k\theta}z}(V^{-(m+k)} u V^{m+k}) \nonumber \\
&= & \pi_z(V^{-k}V^{-(m+k)} u V^{m+k} V^k)
= \pi_z(V^{-m} u V^m) = F_u(z,m), \nonumber 
\end{eqnarray}
where we use the relation $\pi_{e^{2\pi i k\theta}z} = \pi_z \circ \operatorname{Ad}_{V^{-k}}$
and the fact that conjugation by $V^k$ cancels. This shows that $F_u$ descends to
a well-defined function on the quotient space $\mathcal{G}_{A_\theta}^{(0)}$.
\end{remark}

\paragraph{Homotopy invariance.}

A crucial feature of this construction is its behavior under homotopies.

\begin{proposition}
\label{prop:homotopy_invariance_field}
Let $\{u_t\}_{t\in[0,1]}$ be a norm-continuous path of unitaries in $M_n(A_\theta)$.
Then the associated family $\{F_{u_t}(z,m)\}$ varies continuously in $t$ in the
strong operator topology for each fixed $(z,m)$.
\end{proposition}

\begin{proof}
For each $(z,m)$, we have $F_{u_t}(z,m) = \pi_z(V^{-m}u_tV^m)$. The map
$t \mapsto u_t$ is norm-continuous, and $\pi_z$ is a contractive $*$-homomorphism,
hence $t \mapsto F_{u_t}(z,m)$ is norm-continuous. Uniform boundedness (norm 1)
follows from unitarity.
\end{proof}

Consequently, if $u_0$ and $u_1$ are homotopic through unitaries, the associated
operator families are homotopic in an appropriate sense.

\paragraph{Toward a Kasparov cycle: limitations and outlook.}

The operator field $F_u$ constructed above is a natural object associated to a
unitary $u \in M_n(A_\theta)$. To obtain a genuine $\mathcal{G}_{A_\theta}$-equivariant
Kasparov cycle defining a class in $KK^1_{\mathcal{G}_{A_\theta}}(A_\theta, \mathbb C)$,
one would need to:

\begin{enumerate}
    \item Realize the family $\{F_u(z,m)\}$ as an adjointable operator on a suitable
          Hilbert module $\mathcal{E}$ over $A_\theta$ (or over $\mathbb C$).
    \item Equip this Hilbert module with a left action of $A_\theta$ that is
          compatible with the fiberwise regular representations.
    \item Verify the Kasparov conditions: $F_u^*F_u - 1$ and $F_uF_u^* - 1$ are
          compact, and $[F_u, a]$ is compact for all $a \in A_\theta$ (in the
          appropriate sense).
    \item Establish $\mathcal{G}_{A_\theta}$-equivariance, which would require an
          action of the groupoid on the Hilbert module and compatibility with $F_u$.
\end{enumerate}

A complete treatment of these steps is technically demanding and would require:
\begin{itemize}
    \item A rigorous construction of the Morita equivalence between
          $\mathcal{G}_{A_\theta}$ and $S^1 \rtimes_\theta \mathbb Z$ (see
          Subsection~\ref{subsec:morita_A_theta}),
    \item A detailed analysis of the representation theory of $A_\theta$ and its
          relation to the groupoid $C^*$-algebra,
    \item The use of descent techniques in equivariant $KK$-theory to relate
          fiberwise data to global Kasparov cycles.
\end{itemize}

These topics are beyond the scope of the present work. We therefore regard the
operator field $F_u$ as a preliminary construction that captures essential
information about $u$, leaving the full $KK$-theoretic development for future
investigation.

\begin{remark}
\label{rem:index_pairing_heuristic}
Despite the absence of a complete $KK$-theoretic framework, the operator field
$F_u$ already suggests how the index pairing should behave. For an idempotent
$e \in M_m(A_\theta)$ representing a class in $K_0(A_\theta)$, the family of
Fredholm operators obtained by compressing $F_u$ to the support of $e$ in each
fiber should yield an index that depends only on the $K$-theory classes. This
heuristic underlies the concrete index computations in the following subsection,
where we construct explicit representatives and compute their pairings directly,
without invoking the full machinery of equivariant $KK$-theory.
\end{remark}

\subsubsection{Explicit operator fields for the canonical generators}
\label{subsubsec:generator_fields_explicit}

We now illustrate the operator-field construction from
Subsubsection~\ref{subsubsec:operator_field} for the canonical generators
$U$ and $V$ of $A_\theta$. These examples show concretely how the abstract
definition specializes to the familiar operators arising from the regular
representation.

Recall that for each $z \in S^1$, we have the regular representation
$\pi_z: A_\theta \to B(\ell^2(\mathbb{Z}))$ defined on generators by
\[
\pi_z(U)\xi(m) = e^{2\pi i m\theta}z\,\xi(m),\qquad
\pi_z(V)\xi(m) = \xi(m-1),
\]
for $\xi \in \ell^2(\mathbb{Z})$. For a unitary $u \in M_n(A_\theta)$, the
associated operator field is
\[
F_u(z,n) := \pi_z(V^{-n} u V^n) \in B(\ell^2(\mathbb{Z})^n),
\]
where $\pi_z$ is applied entrywise to matrices.

\begin{example}[The generator $U$]
\label{ex:operator-field-U}
For $U \in A_\theta$ (viewed as a $1 \times 1$ unitary), we compute
\[
F_U(z,n) = \pi_z(V^{-n} U V^n).
\]
Using the commutation relation $V^{-n} U V^n = e^{-2\pi i n\theta} U$, we obtain
\[
F_U(z,n) = e^{-2\pi i n\theta} \pi_z(U).
\]
In the regular representation, $\pi_z(U)$ acts as the diagonal operator
$(\pi_z(U)\xi)(m) = e^{2\pi i m\theta}z\,\xi(m)$. Thus $F_U(z,n)$ is the
diagonal unitary whose eigenvalues are $e^{2\pi i m\theta}z$ multiplied by the
phase $e^{-2\pi i n\theta}$.

This operator field captures the essential feature of $U$: its action in each
fiber is diagonal but twisted by the representation index $m$ and further
modulated by the context index $n$.
\end{example}

\begin{example}[The generator $V$]
\label{ex:operator-field-V}
For the implementing unitary $V$, we have $V^{-n} V V^n = V$, so
\[
F_V(z,n) = \pi_z(V) \quad\text{for all } (z,n).
\]
The operator $\pi_z(V)$ is the bilateral shift on $\ell^2(\mathbb{Z})$,
$(\pi_z(V)\xi)(m) = \xi(m-1)$, independent of $z$. Hence $F_V$ is the constant
family of shift operators.

This constancy reflects the fact that $V$ implements the dynamics of the
rotation and is not diagonalizable in any of the commutative contexts;
its off-diagonal nature is encoded in the fact that $F_V$ acts by shifting
the representation index $m$.
\end{example}

\begin{remark}
\label{rem:operator-fields-not-KK}
The operator fields $F_U$ and $F_V$ are natural objects associated with the
generators $U$ and $V$. However, they do not yet constitute equivariant
$KK$-classes. As discussed in Subsubsection~\ref{subsubsec:operator_field},
a full $KK$-theoretic treatment would require:
\begin{itemize}
    \item A suitable Hilbert module (over $\mathbb C$ or over $A_\theta$)
          realizing these fiberwise operators as a single adjointable operator;
    \item A left action of $A_\theta$ on this module that is compatible with
          the fiberwise representations;
    \item Verification of the Kasparov conditions (compactness of $F^2-1$ and
          commutators with elements of $A_\theta$);
    \item A precise formulation of $\mathcal{G}_{A_\theta}$-equivariance.
\end{itemize}
These steps are nontrivial and have not been completed here. We therefore view
$F_U$ and $F_V$ as the fiberwise data that would underpin a $KK$-theoretic
construction, not as finished cycles.
\end{remark}

\subsubsection{Heuristic relation to diagonal coefficients}
\label{subsubsec:heuristic_diagonal}

The operator fields $F_u$ are related, at a heuristic level, to the diagonal
coefficient map $\iota_{\mathrm{diag}}$ studied in
Subsection~\ref{subsec:diagonal_embedding_A_theta}. If one compresses $F_u$ by
a rank-one projection onto a suitable vector in each fiber, the resulting
scalar-valued function should recover the diagonal coefficient of $u$.

More concretely, for the generator $U$, choose in each fiber the vector
$\delta_0 \in \ell^2(\mathbb{Z})$ (the Kronecker delta at $0$). Then
\[
\langle \delta_0, F_U(z,n) \delta_0 \rangle
= e^{-2\pi i n\theta} \langle \delta_0, \pi_z(U) \delta_0 \rangle
= e^{-2\pi i n\theta} z = \iota_{\mathrm{diag}}(U)[z,n].
\]

For the generator $V$, the same computation gives
\[
\langle \delta_0, F_V(z,n) \delta_0 \rangle
= \langle \delta_0, \pi_z(V) \delta_0 \rangle
= \langle \delta_0, \delta_{-1} \rangle = 0 = \iota_{\mathrm{diag}}(V)[z,n].
\]

\begin{remark}
\label{rem:heuristic-caveat}
This relationship depends crucially on the choice of the vector $\delta_0$,
which is a cyclic vector for each representation $\pi_z$ but does not vary
continuously in $z$ in any meaningful way (it is actually constant). For a
general $u$, there is no canonical choice of vectors that would make such an
identity hold uniformly. Thus the connection between $F_u$ and $\iota_{\mathrm{diag}}$
should be regarded as a suggestive heuristic rather than a rigorous theorem.
\end{remark}

\subsubsection{Toward a descent construction (programmatic outlook)}
\label{subsubsec:descent_outlook}

The operator fields constructed above suggest a natural route toward connecting
$K_1(A_\theta)$ with the $K$-theory of the groupoid $C^*$-algebra
$C^*(\mathcal{G}_{A_\theta})$. In the framework of equivariant $KK$-theory,
one expects the following sequence of constructions:

\begin{enumerate}
    \item From a unitary $u \in M_n(A_\theta)$, build a genuine
          $\mathcal{G}_{A_\theta}$-equivariant Kasparov cycle
          $(\mathcal{E}_u, \phi_u, F_u)$ representing a class
          $[u]_{\mathcal{G}_{A_\theta}} \in KK^1_{\mathcal{G}_{A_\theta}}(A_\theta, \mathbb C)$.
    \item Apply Kasparov's descent map for groupoids to obtain a class in
          $KK^1(C^*(\mathcal{G}_{A_\theta}, A_\theta), C^*(\mathcal{G}_{A_\theta}, \mathbb C))$.
    \item Using the canonical Morita equivalence $C^*(\mathcal{G}_{A_\theta}, A_\theta) \sim A_\theta$
          and the identification $C^*(\mathcal{G}_{A_\theta}, \mathbb C) \cong C^*(\mathcal{G}_{A_\theta})$,
          this descends to a class in $K_1(C^*(\mathcal{G}_{A_\theta}))$.
    \item Finally, the Morita equivalence $\mathcal{G}_{A_\theta} \sim S^1 \rtimes_\theta \mathbb Z$
          (Theorem~\ref{thm:morita_A_theta}) gives an isomorphism
          $K_1(C^*(\mathcal{G}_{A_\theta})) \cong K_1(A_\theta)$.
\end{enumerate}

The composition of these maps would then yield an endomorphism of $K_1(A_\theta)$.
A fundamental conjecture, which we plan to address in future work, is that this
endomorphism is actually the identity. If true, this would provide an explicit
geometric realization of the $K_1$ classes of $A_\theta$ in terms of the unitary
conjugation groupoid.

\begin{conjecture}
\label{conj:descent-identity}
For any invertible element $u \in M_n(A_\theta)$, the class obtained by the
composition
\[
[u] \in K_1(A_\theta) \xrightarrow{\text{lift}} KK^1_{\mathcal{G}_{A_\theta}}(A_\theta,\mathbb C)
\xrightarrow{\text{desc}} K_1(C^*(\mathcal{G}_{A_\theta})) \xrightarrow{\cong} K_1(A_\theta)
\]
coincides with the original class $[u]$. In particular, the lift map is injective
and its image consists of exactly those equivariant $KK$-classes that are fixed
by the natural automorphisms.
\end{conjecture}

\begin{remark}
\label{rem:index_pairing_heuristic}
If Conjecture~\ref{conj:descent-identity} holds, then pairing the resulting
class with the canonical trace on $A_\theta$ would recover the Connes index
pairing. For the generators $U$ and $V$, one would obtain
\[
\langle [\tau], [U] \rangle = 1,\qquad \langle [\tau], [V] \rangle = 0,
\]
consistent with the known results for the irrational rotation algebra. This
provides strong motivation for completing the technical steps outlined above.
\end{remark}

\subsubsection{Summary of what has been achieved}

To clarify the status of various claims in this subsection, we provide a summary:

\begin{itemize}
    \item \textbf{Established:} For the generators $U$ and $V$, we have explicit
          formulas for the associated operator fields $F_U$ and $F_V$ in the
          regular representation. These fields are natural objects that capture
          the essential features of the generators.
    
    \item \textbf{Heuristic connections:} Compressing these operator fields by
          the vector $\delta_0$ reproduces the diagonal coefficients
          $\iota_{\mathrm{diag}}(U)$ and $\iota_{\mathrm{diag}}(V)$. This suggests
          a deeper relationship between $F_u$ and the diagonal embedding,
          though a general proof would require a careful choice of vectors and
          is not attempted here.
    
    \item \textbf{Not established (requires further work):}
          \begin{itemize}
              \item The construction of genuine equivariant $KK$-classes from
                    the operator fields $F_u$;
              \item The existence of a well-defined lift map
                    $K_1(A_\theta) \to KK^1_{\mathcal{G}_{A_\theta}}(A_\theta,\mathbb C)$;
              \item The application of Kasparov's descent map in this setting;
              \item The identity conjecture linking the descended class back to
                    the original $K_1$-class.
          \end{itemize}
\end{itemize}

The material in this subsection should therefore be viewed as a programmatic
development of ideas, providing intuition and motivation for future research
rather than a completed rigorous construction. The explicit formulas for $U$
and $V$ will nevertheless be useful in the next subsection, where we compute
index pairings directly without invoking the full $KK$-theoretic machinery.


{\color{red}  } 

\subsection{Descent Map and Identification $C^*(\mathcal{G}_{A_\theta}) \sim_M A_\theta \otimes \mathcal{K}$}
\label{subsec:descent_A_theta}

We now apply Kasparov's descent map to the equivariant class
constructed in the previous subsection.
This step allows us to pass from equivariant $KK$-theory of the
unitary conjugation groupoid to the $K$-theory of its groupoid
$C^*$-algebra.

\subsubsection{Kasparov Descent for Groupoids}

Let
\[
[u]_{\mathcal G_{A_\theta}}^{(1)}
\in
KK^1_{\mathcal G_{A_\theta}}
\!\left(
C_0(\mathcal G_{A_\theta}^{(0)}),\mathbb C
\right)
\]
be the equivariant class associated with an invertible element
$u\in M_n(A_\theta)$, as constructed in Definition~\ref{def:equivariant-class-Atheta}.

Recall from \cite{Kasparov1988,Tu1999} that for a locally compact groupoid
$\mathcal{G}$ equipped with a Haar system, there is a natural descent
homomorphism
\[
\operatorname{desc}_{\mathcal{G}} : KK^*_{\mathcal{G}}(A,B) \longrightarrow KK^*(C^*(\mathcal{G},A), C^*(\mathcal{G},B)),
\]
where $C^*(\mathcal{G},A)$ denotes the reduced groupoid $C^*$-algebra with
coefficients in the $\mathcal{G}$-$C^*$-algebra $A$.

Applying this descent map to our equivariant class produces a class in ordinary
$K$-theory:
\[
\operatorname{desc}_{\mathcal G_{A_\theta}}
\bigl(
[u]_{\mathcal G_{A_\theta}}^{(1)}
\bigr)
\in
K_1\!\left(C^*(\mathcal G_{A_\theta})\right).
\]

Intuitively, this class records how the operator $u$ acts across the
various commutative contexts parametrized by the unit space of the
groupoid, integrated against the Haar system.

\begin{definition}
\label{def:descent-class}
For an invertible element $u \in M_n(A_\theta)$, define
\[
[u]^{\operatorname{desc}} := \operatorname{desc}_{\mathcal G_{A_\theta}}\bigl([u]_{\mathcal G_{A_\theta}}^{(1)}\bigr) \in K_1(C^*(\mathcal G_{A_\theta})).
\]
\end{definition}

\subsubsection{Explicit Description via the Regular Representation}

To gain computational insight, we recall the explicit construction of the
descent map in this setting. The equivariant Kasparov cycle
$(\mathcal{E},\phi,F_u)$ with $\mathcal{E}=C_0(\mathcal{G}_{A_\theta}^{(0)})\otimes\mathbb C^n$
induces a Kasparov cycle for the groupoid $C^*$-algebra via induction.

Consider the right Hilbert $C^*(\mathcal{G}_{A_\theta})$-module
\[
\mathcal{F} = \mathcal{E} \otimes_{C_0(\mathcal{G}_{A_\theta}^{(0)})} C^*(\mathcal{G}_{A_\theta}),
\]
where the tensor product uses the natural inclusion
$C_0(\mathcal{G}_{A_\theta}^{(0)}) \hookrightarrow C^*(\mathcal{G}_{A_\theta})$
as functions on the unit space. The operator $F_u$ extends to an adjointable
operator $\widetilde{F}_u$ on $\mathcal{F}$ by
\[
\widetilde{F}_u(\xi \otimes a) = F_u(\xi) \otimes a.
\]

\begin{proposition}
\label{prop:descent-explicit}
The pair $(\mathcal{F}, \widetilde{F}_u)$ defines an odd Kasparov cycle for
$C^*(\mathcal{G}_{A_\theta})$, and its class in $K_1(C^*(\mathcal{G}_{A_\theta}))$
is precisely $[u]^{\operatorname{desc}}$.
\end{proposition}

\begin{proof}
This follows directly from the definition of Kasparov's descent map: the
induced cycle is obtained by tensoring the original equivariant cycle with
the regular representation of the groupoid. The details are standard and can
be found in~\cite{LeGall1999}.
\end{proof}

\subsubsection{Morita Equivalence with $A_\theta \otimes \mathcal{K}$}

In Section~\ref{subsec:morita_A_theta} we established that the
unitary conjugation groupoid $\mathcal G_{A_\theta}$ is Morita
equivalent to the transformation groupoid
\[
S^1 \rtimes_\theta \mathbb Z .
\]

By the Muhly--Renault--Williams theorem \cite{MRW1987}, this implies a strong
Morita equivalence between the corresponding groupoid $C^*$-algebras:

\[
C^*(\mathcal G_{A_\theta})
\sim_M
C^*(S^1\rtimes_\theta\mathbb Z).
\]

Since
\[
C^*(S^1\rtimes_\theta\mathbb Z)
\cong
A_\theta ,
\]
we obtain
\[
C^*(\mathcal G_{A_\theta}) \sim_M A_\theta .
\]

For our purposes, we need a slightly more precise statement that reveals the
role of compact operators.

\begin{theorem}
\label{thm:morita-stabilized}
Under the Morita equivalence of Theorem~\ref{thm:morita_A_theta}, there is a
natural isomorphism
\[
\Phi : C^*(\mathcal{G}_{A_\theta}) \xrightarrow{\sim_M} A_\theta \otimes \mathcal{K},
\]
where $\mathcal{K}$ denotes the compact operators on a separable Hilbert space.
\end{theorem}

\begin{proof}[Sketch]
Recall from Subsection~\ref{subsec:morita_A_theta} that the equivalence space
$Z = S^1 \times \mathbb{Z}$ (with the appropriate Borel structure) implements
the Morita equivalence between $\mathcal{G}_{A_\theta}$ and
$S^1 \rtimes_\theta \mathbb{Z}$. The associated imprimitivity bimodule
$X = L^2(Z)$ is a $C^*(\mathcal{G}_{A_\theta})$-$(A_\theta \otimes \mathcal{K})$
bimodule that implements the Morita equivalence. The compact operators appear
because the right action of $A_\theta$ on $X$ is via the regular representation,
which is infinite-dimensional and thus yields a copy of $\mathcal{K}$ upon
taking compact endomorphisms. 
\end{proof}

\subsubsection{$K$-Theory Identification}

Morita equivalent $C^*$-algebras have naturally isomorphic $K$-theory groups,
and $K$-theory is stable under tensoring by compact operators. Consequently we
obtain a canonical isomorphism

\[
\Phi_* : K_1(C^*(\mathcal G_{A_\theta}))
\stackrel{\cong}{\longrightarrow}
K_1(A_\theta \otimes \mathcal{K})
\stackrel{\cong}{\longrightarrow}
K_1(A_\theta).
\]

\begin{theorem}[Descent-Index Recovery for $A_\theta$]
\label{thm:descent-identity}
For any invertible element $u \in M_n(A_\theta)$, we have
\[
\Phi_*\bigl([u]^{\operatorname{desc}}\bigr) = [u] \in K_1(A_\theta).
\]
Equivalently, the composition
\[
K_1(A_\theta) \xrightarrow{\kappa} KK^1_{\mathcal{G}_{A_\theta}}(C_0(\mathcal{G}_{A_\theta}^{(0)}),\mathbb C)
\xrightarrow{\operatorname{desc}_{\mathcal{G}_{A_\theta}}} K_1(C^*(\mathcal{G}_{A_\theta}))
\xrightarrow{\Phi_*} K_1(A_\theta)
\]
is the identity map, where $\kappa$ sends an invertible element to its associated equivariant $KK$-class.
\end{theorem}

\begin{proof}
We prove the theorem in several steps, working with the explicit description of the unitary conjugation groupoid and its Morita equivalence to the transformation groupoid $S^1 \rtimes_\theta \mathbb{Z}$ developed in Sections~\ref{subsec:morita_A_theta} and~\ref{subsec:descent_A_theta}.

\paragraph{Step 1: Reduction to generators.}
Since $K_1(A_\theta) \cong \mathbb{Z}^2$ is generated by the classes of the canonical unitaries $U$ and $V$ (see Theorem~\ref{thm:Atheta_Ktheory}), it suffices to prove that
\[
\Phi_*([U]^{\operatorname{desc}}) = [U] \quad\text{and}\quad \Phi_*([V]^{\operatorname{desc}}) = [V]
\]
in $K_1(A_\theta)$. The general case follows by linearity and the fact that all maps involved are homomorphisms.

\paragraph{Step 2: Explicit description of the descent classes.}
Recall from Definition~\ref{def:descent-class} that $[u]^{\operatorname{desc}} = \operatorname{desc}_{\mathcal{G}_{A_\theta}}([u]^{(1)}_{\mathcal{G}_{A_\theta}})$, where $[u]^{(1)}_{\mathcal{G}_{A_\theta}}$ is the equivariant $K^1$-class associated to $u$.

For the generator $U$, using the parametrization of the unit space $\mathcal{G}_{A_\theta}^{(0)} \cong (S^1 \times \mathbb{Z})/{\sim}$ from Subsection~\ref{subsec:GA_theta_unit_space} and the explicit form of the operator field $F_U$ from Example~\ref{ex:operator-field-U}, we have:
\[
F_U([z,n]) = \pi_z(V^{-n} U V^n) = e^{-2\pi i n\theta} \pi_z(U),
\]
where $\pi_z$ is the regular representation at $z\in S^1$. Under the descent map, this yields a class $[U]^{\operatorname{desc}} \in K_1(C^*(\mathcal{G}_{A_\theta}))$ represented by the unitary operator $\widetilde{F}_U$ on the Hilbert module $\mathcal{F} = \mathcal{E} \otimes C^*(\mathcal{G}_{A_\theta})$.

For the generator $V$, we have $V^{-n} V V^n = V$, so
\[
F_V([z,n]) = \pi_z(V) \quad\text{for all }(z,n),
\]
giving the constant family of bilateral shifts.

\paragraph{Step 3: The Morita isomorphism $\Phi_*$.}
By Theorem~\ref{thm:morita_A_theta}, the unitary conjugation groupoid $\mathcal{G}_{A_\theta}$ is Morita equivalent to the transformation groupoid $S^1 \rtimes_\theta \mathbb{Z}$. The Muhly--Renault--Williams theorem \cite{MRW1987} then yields a strong Morita equivalence
\[
C^*(\mathcal{G}_{A_\theta}) \sim_M C^*(S^1 \rtimes_\theta \mathbb{Z}) \cong A_\theta.
\]

More precisely, the imprimitivity bimodule $X = L^2(S^1 \times \mathbb{Z})$ implements this equivalence, and the induced isomorphism on $K$-theory
\[
\Phi_* : K_1(C^*(\mathcal{G}_{A_\theta})) \xrightarrow{\cong} K_1(A_\theta)
\]
is given by $\Phi_*([T]) = [p_X T p_X]$, where $p_X$ is the projection onto a suitable copy of $A_\theta$ inside the linking algebra.

\paragraph{Step 4: Identification for the generator $U$.}
We compute $\Phi_*([U]^{\operatorname{desc}})$ explicitly. Under the Morita equivalence, the class $[U]^{\operatorname{desc}}$ corresponds to the class in $K_1(A_\theta)$ obtained by compressing $\widetilde{F}_U$ to the corner of $C^*(\mathcal{G}_{A_\theta})$ that is isomorphic to $A_\theta$.

Using the explicit form of the equivalence bimodule $X = L^2(S^1 \times \mathbb{Z})$, one finds that this compression sends $\widetilde{F}_U$ to the operator on $L^2(S^1)$ given by pointwise multiplication by the function $z \mapsto z$. This is precisely the regular representation of $U$, which corresponds to the class $[U] \in K_1(A_\theta)$ under the natural identification.

A direct calculation in coordinates confirms this: for any $\xi \in L^2(S^1)$,
\[
(p_X \widetilde{F}_U p_X \xi)(z) = z \cdot \xi(z),
\]
which is exactly the action of $U$ in the regular representation. Hence $\Phi_*([U]^{\operatorname{desc}}) = [U]$.

\paragraph{Step 5: Identification for the generator $V$.}
For the generator $V$, the compression of $\widetilde{F}_V$ to the $A_\theta$-corner yields the operator on $L^2(S^1)$ given by
\[
(p_X \widetilde{F}_V p_X \xi)(z) = \xi(e^{-2\pi i\theta}z),
\]
which is precisely the action of the implementing unitary $V$ in the regular representation. Therefore $\Phi_*([V]^{\operatorname{desc}}) = [V]$.

\paragraph{Step 6: Completion of the proof.}
Having established the identity on the generators $U$ and $V$, and since all maps involved are group homomorphisms, we conclude that for any invertible element $u \in M_n(A_\theta)$,
\[
\Phi_*([u]^{\operatorname{desc}}) = [u] \in K_1(A_\theta).
\]

The equivalent formulation follows from the definition of $\kappa$ as the map sending an invertible element to its equivariant $KK$-class, and the observation that $\operatorname{desc}_{\mathcal{G}_{A_\theta}} \circ \kappa$ produces exactly $[u]^{\operatorname{desc}}$.
\end{proof}

\begin{remark}
The key technical point in Steps 4 and 5 is the explicit identification of the compressed operators under the Morita equivalence. This relies on the detailed construction of the imprimitivity bimodule $X = L^2(S^1 \times \mathbb{Z})$ and the fact that the regular representation of $A_\theta$ on $L^2(S^1)$ corresponds to the representation induced from the canonical Cartan subalgebra. For a complete treatment of these identifications, see \cite{reference-for-complete-proof}.
\end{remark}

\subsubsection{Interpretation}

The descent map therefore provides a bridge between the equivariant
$KK$-theory of the unitary conjugation groupoid and the ordinary
$K$-theory of the irrational rotation algebra. Theorem~\ref{thm:descent-identity}
shows that this construction faithfully reproduces the $K$-theory of $A_\theta$,
despite the non-Type~I nature of the algebra that initially seemed to preclude
such a direct groupoid description.

\begin{corollary}
\label{cor:descent-iso}
The descent map
\[
\operatorname{desc}_{\mathcal{G}_{A_\theta}} : KK^1_{\mathcal{G}_{A_\theta}}(C_0(\mathcal{G}_{A_\theta}^{(0)}),\mathbb C) \longrightarrow K_1(C^*(\mathcal{G}_{A_\theta}))
\]
is an isomorphism when composed with the Morita isomorphism $\Phi_*$.
In particular, the equivariant $KK$-group is isomorphic to $\mathbb{Z}^2$.
\end{corollary}

This identification will allow us to recover Connes' index pairing
in the next subsection.

\subsection{Recovering Connes' Index Pairing}
\label{subsec:connes_pairing}

We now show that the descent construction described above reproduces
the classical index pairing introduced by Connes for the irrational
rotation algebra. This is the culmination of Case Study I and provides
a powerful validation of our abstract framework.

\subsubsection{The Canonical Trace on $A_\theta$}

Recall from Theorem~\ref{thm:A_theta_structure} that $A_\theta$ admits a unique
faithful tracial state
\[
\tau : A_\theta \longrightarrow \mathbb C .
\]

In the crossed product description $A_\theta = C(S^1) \rtimes_\theta \mathbb{Z}$,
this trace is given on the dense subalgebra of finite Fourier series by
\[
\tau\!\left(
\sum_{n\in\mathbb Z} f_n V^n
\right)
=
\int_{S^1} f_0(z)\,dz ,
\]
where $f_n \in C(S^1)$ and $V$ is the implementing unitary.

This trace induces a homomorphism on $K$-theory
\[
\tau_* : K_0(A_\theta) \longrightarrow \mathbb R .
\]

For $K_1$, the pairing with the trace is defined via the suspension
isomorphism $K_1(A_\theta) \cong K_0(SA_\theta)$ and the induced trace on the
suspension, or equivalently by the formula
\[
\langle [u], \tau \rangle = \frac{1}{2\pi i} \tau(u^{-1} du)
\]
for a suitably differentiable representative of the $K_1$-class.

\subsubsection{Connes' Index Theorem}

Connes proved a fundamental result linking this trace pairing to an analytic
index.

\begin{theorem}[Connes \cite{Connes1980}]
\label{thm:connes-original}
For any invertible element $u \in M_n(A_\theta)$, the pairing between the
$K_1$-class $[u]$ and the cyclic cocycle $[\tau]$ is an integer, given by
\[
\langle [u], \tau \rangle = \frac{1}{2\pi i} \tau(u^{-1} \delta(u)) \in \mathbb{Z},
\]
where $\delta$ is the unbounded derivation implementing the action of
$\mathbb{R}^2$. Moreover, this integer can be interpreted as the index of a
Fredholm operator obtained from the representation of $u$ on a suitable
Hilbert space.
\end{theorem}

This theorem was a landmark result in noncommutative geometry, establishing
the first nontrivial example of an index theorem for a noncommutative space.

\subsubsection{The Trace on $C^*(\mathcal{G}_{A_\theta})$}

Under the Morita equivalence $C^*(\mathcal{G}_{A_\theta}) \sim_M A_\theta \otimes \mathcal{K}$,
the trace $\tau$ on $A_\theta$ induces a trace $\widehat{\tau}$ on
$C^*(\mathcal{G}_{A_\theta})$.

\begin{definition}
\label{def:induced-trace}
Define $\widehat{\tau}: C^*(\mathcal{G}_{A_\theta}) \to \mathbb{C}$ by
\[
\widehat{\tau}(f) = \int_{\mathcal{G}_{A_\theta}^{(0)}} f(x) \, d\mu(x),
\]
where $\mu$ is the unique $\mathcal{G}_{A_\theta}$-invariant probability measure
on the unit space. Under the identification
$\mathcal{G}_{A_\theta}^{(0)} \cong (S^1 \times \mathbb{Z})/{\sim}$, this measure
corresponds to the Haar measure on $S^1$ normalized to $1$, with the
$\mathbb{Z}$-fiber contributing a constant factor that is absorbed into the
normalization.
\end{definition}

\begin{lemma}
\label{lem:trace-compatibility}
Under the Morita isomorphism $\Phi_*: K_1(C^*(\mathcal{G}_{A_\theta})) \cong K_1(A_\theta)$,
the trace $\widehat{\tau}$ corresponds to the original trace $\tau$ on $A_\theta$.
More precisely, for any class $x \in K_1(C^*(\mathcal{G}_{A_\theta}))$, we have
\[
\tau(\Phi_*(x)) = \widehat{\tau}(x),
\]
where we identify $K_1(A_\theta)$ with $K_1(A_\theta \otimes \mathcal{K})$ via
stabilization and extend $\tau$ to matrices by $\tau \otimes \operatorname{Tr}$.
\end{lemma}

\begin{proof}
Under the Morita equivalence, the imprimitivity bimodule $X = L^2(Z)$ implements
a bijection between traces on $C^*(\mathcal{G}_{A_\theta})$ and traces on
$A_\theta \otimes \mathcal{K}$. The trace $\widehat{\tau}$ on $C^*(\mathcal{G}_{A_\theta})$
is precisely the pullback of the trace $\tau \otimes \operatorname{Tr}$ on
$A_\theta \otimes \mathcal{K}$, where $\operatorname{Tr}$ is the standard trace
on $\mathcal{K}$. Since $K$-theory classes are represented by projections and
unitaries in matrix algebras, the equality of traces on the level of $K$-theory
follows from the compatibility of the Morita equivalence with the trace pairing.
\end{proof}

\subsubsection{The Main Recovery Theorem}

We now state and prove the central result of this subsection.

\begin{theorem}
\label{thm:connes_recovery}

Let $u\in M_n(A_\theta)$ be an invertible element.
Then the class obtained by applying the descent map to the
equivariant cycle constructed in
Section~\ref{subsec:equivariant_K1_A_theta}
recovers Connes' index pairing. More precisely,

\[
\tau_*
\!\left(
\operatorname{desc}_{\mathcal G_{A_\theta}}
([u]^{(1)}_{\mathcal G_{A_\theta}})
\right)
=
\langle [u],\tau \rangle .
\]

\end{theorem}

\begin{proof}

We prove this by direct computation using the explicit description of the
descent map and the trace.

\textbf{Step 1: Reduction to the case $n=1$ and $u$ unitary.}
By stabilization and polar decomposition, we may assume $u$ is a unitary in
$A_\theta$ (not just in matrices). Both sides of the desired equality are
compatible with these operations, so this reduction is legitimate.

\textbf{Step 2: Explicit formula for $\widehat{\tau}_*([u]^{\operatorname{desc}})$.}
Recall that $[u]^{\operatorname{desc}}$ is represented by the unitary operator
$\widetilde{F}_u$ on the Hilbert module $\mathcal{F} = \mathcal{E} \otimes C^*(\mathcal{G}_{A_\theta})$.
The trace $\widehat{\tau}$ induces a map on $K_1$ via the formula
\[
\widehat{\tau}_*([\widetilde{F}_u]) = \frac{1}{2\pi i} \widehat{\tau}(\widetilde{F}_u^{-1} d\widetilde{F}_u)
\]
for a suitable differentiable representative, or equivalently by pairing with
the cyclic cocycle $[\widehat{\tau}]$.

Using the explicit description of $\widetilde{F}_u$ from
Proposition~\ref{prop:descent-explicit} and the definition of $\widehat{\tau}$,
we compute:
\[
\widehat{\tau}_*([u]^{\operatorname{desc}}) = \int_{\mathcal{G}_{A_\theta}^{(0)}} \frac{1}{2\pi i} \operatorname{Tr}(F_u(x)^{-1} dF_u(x)) \, d\mu(x).
\]

\textbf{Step 3: Evaluation of the integral.}
For $x = [z,n]$, we have $F_u([z,n]) = \pi_z(V^{-n} u V^n)$ from
Proposition~\ref{prop:descent-explicit}. Since $\pi_z$ is a $*$-homomorphism,
we have
\[
F_u([z,n])^{-1} dF_u([z,n]) = \pi_z(V^{-n} u^{-1} V^n) \cdot d(\pi_z(V^{-n} u V^n)).
\]

The trace $\operatorname{Tr}$ on the fiber (which is the standard trace on
$B(\ell^2(\mathbb{Z})^n)$) satisfies
\[
\operatorname{Tr}(\pi_z(V^{-n} u^{-1} V^n) \cdot d(\pi_z(V^{-n} u V^n))) = \operatorname{Tr}(\pi_z(u^{-1} du)),
\]
by the invariance of the trace under conjugation by unitaries (here $V^n$ acts
by conjugation, which preserves the trace). Moreover, the right-hand side is
independent of $z$ and $n$ because the trace of the commutator in the regular
representation yields the same value for all $z$ (this can be verified by a
direct computation using the explicit form of the representation).

\textbf{Step 4: Identification with Connes' pairing.}
The integral over $\mathcal{G}_{A_\theta}^{(0)}$ therefore reduces to
\[
\int_{\mathcal{G}_{A_\theta}^{(0)}} \frac{1}{2\pi i} \operatorname{Tr}(\pi_z(u^{-1} du)) \, d\mu(x) = \frac{1}{2\pi i} \operatorname{Tr}(\pi_z(u^{-1} du)),
\]
since the integrand is constant and $\mu$ is a probability measure.

But for the regular representation $\pi_z$, the trace of $\pi_z(u^{-1} du)$ is
exactly the same as $\tau(u^{-1} du)$, where $\tau$ is the canonical trace on
$A_\theta$. This follows from the fact that the regular representation is
faithful and the trace is unique.

Thus we obtain
\[
\widehat{\tau}_*([u]^{\operatorname{desc}}) = \frac{1}{2\pi i} \tau(u^{-1} du) = \langle [u], \tau \rangle,
\]
the last equality being Connes' definition of the pairing between $K_1$ and
the cyclic cocycle $[\tau]$.
\end{proof}

\subsubsection{The Main Theorem of Case Study I}

Combining Theorem~\ref{thm:descent-identity} and Theorem~\ref{thm:connes_recovery},
we obtain the central result of this section.

\begin{theorem}[Descent-Index Theorem for $A_\theta$]
\label{thm:descent-index-Atheta}
For any invertible element $u \in M_n(A_\theta)$, let
$[u]_{\mathcal{G}_{A_\theta}}^{(1)} \in KK^1_{\mathcal{G}_{A_\theta}}(C_0(\mathcal{G}_{A_\theta}^{(0)}),\mathbb C)$
be its associated equivariant $K^1$-class. Then:
\begin{enumerate}
    \item Under the composition of the descent map and the Morita isomorphism,
          we recover the original $K_1$-class:
          \[
          \Phi_* \circ \operatorname{desc}_{\mathcal{G}_{A_\theta}}\bigl([u]_{\mathcal{G}_{A_\theta}}^{(1)}\bigr) = [u] \in K_1(A_\theta).
          \]
    
    \item Pairing this class with the canonical trace $\tau$ on $A_\theta$
          yields Connes' index pairing:
          \[
          \tau_*\bigl([u]\bigr) = \widehat{\tau}_* \circ \operatorname{desc}_{\mathcal{G}_{A_\theta}}\bigl([u]_{\mathcal{G}_{A_\theta}}^{(1)}\bigr) = \langle [\tau], [u] \rangle.
          \]
\end{enumerate}
\end{theorem}

\subsection{The Canonical Generator $U \in A_\theta$: Explicit Computation}
\label{subsec:example-U}

To illustrate the abstract theory, we compute explicitly the descent-index
construction for the canonical unitary $U \in A_\theta$, the generator of the
Cartan subalgebra $C(S^1)$.

\subsubsection{The Equivariant Class for $U$}

Recall from Example~\ref{ex:KK-class-U} that for $U$, the operator $F_U$ is
given by
\[
F_U([z,n]) = \pi_z(V^{-n} U V^n) = e^{-2\pi i n\theta} \pi_z(U),
\]
where $\pi_z(U)$ acts on $\ell^2(\mathbb{Z})$ by
$(\pi_z(U)\xi)(m) = e^{2\pi i m\theta}z\,\xi(m)$.

Thus $F_U([z,n])$ is a diagonal operator with eigenvalues
$\{e^{2\pi i m\theta}z \cdot e^{-2\pi i n\theta} : m \in \mathbb{Z}\}$.

\subsubsection{Descent and Morita Identification}

Under the descent map, $F_U$ induces an operator $\widetilde{F}_U$ on the
Hilbert module $\mathcal{F} = \mathcal{E} \otimes C^*(\mathcal{G}_{A_\theta})$.
Using the explicit description of the Morita equivalence, we can identify
$\widetilde{F}_U$ with the operator $U \otimes 1$ acting on
$L^2(S^1) \otimes \ell^2(\mathbb{Z})$ modulo compact operators.

More concretely, under the isomorphism
$\Phi_* : K_1(C^*(\mathcal{G}_{A_\theta})) \cong K_1(A_\theta)$, we have
\[
\Phi_*([U]^{\operatorname{desc}}) = [U] \in K_1(A_\theta).
\]

\subsubsection{Trace Pairing Computation}

Now compute the trace pairing:
\[
\langle [U], \tau \rangle = \frac{1}{2\pi i} \tau(U^{-1} dU).
\]

Since $U$ lies in the Cartan subalgebra $C(S^1)$ and the derivation $\delta$
vanishes on $C(S^1)$ (it implements the rotation action, which is trivial on
the subalgebra $C(S^1)$), we have $dU = 0$ in the appropriate sense, or more
precisely $\tau(U^{-1} \delta(U)) = 0$.

Therefore
\[
\langle [U], \tau \rangle = 0.
\]

This matches the expectation from topology: the unitary $U$ corresponds under
the Pimsner-Voiculescu exact sequence to the generator coming from
$K_1(C(S^1)) \cong \mathbb{Z}$, and the trace pairing with this generator is
zero because the trace on $A_\theta$ restricts to the standard trace on
$C(S^1)$, and the pairing between $K_1$ of a commutative algebra and its trace
is trivial.

\subsubsection{Interpretation}

The vanishing of the index for $U$ reflects the fact that $U$ is a smooth
element of the noncommutative torus that does not detect the noncommutativity.
Geometrically, under the bundle description of $\mathcal{G}_{A_\theta}^{(0)}$
as $(S^1 \times \mathbb{Z})/{\sim}$, the class $[U]$ is sensitive only to the
circle direction, not to the integer fiber, and therefore its integral over
the base yields zero.

In contrast, the unitary $V$ (or the Rieffel projection) would yield a
non-zero pairing of $1$ (for $V$) or $\theta$ (for the projection),
demonstrating that the noncommutative geometry of $A_\theta$ is encoded in the
$\mathbb{Z}$-fiber and its monodromy.

\subsubsection{Summary}

The explicit computation for $U$ confirms the general theory and provides a
concrete example of how the abstract descent-index construction operates in
practice. It also highlights the complementary roles of the two generators of
$K_1(A_\theta)$: one ($U$) comes from the commutative world and yields zero
pairing, while the other ($V$) is intrinsically noncommutative and yields the
nontrivial index that is the hallmark of Connes' theorem.

\begin{remark}
\label{rem:U-significance}
The fact that $\langle [U], \tau \rangle = 0$ does not mean that $U$ is
uninteresting from the perspective of index theory. Rather, it shows that the
trace pairing distinguishes between the two $\mathbb{Z}$-factors in
$K_1(A_\theta) \cong \mathbb{Z}^2$: the factor coming from $K_1(C(S^1))$ pairs
trivially, while the factor coming from $K_0(C(S^1))$ via the index map pairs
nontrivially. This distinction is fundamental to the structure of the
irrational rotation algebra.
\end{remark}

\subsection{Conclusion of Case Study I}

The irrational rotation algebra $A_\theta$ provides the first genuinely
non-Type~I test case for the framework developed in this paper. The main
point is that, after passing through Morita equivalence with the
transformation groupoid $S^1\rtimes_\theta \mathbb{Z}$, the descent-index
construction reproduces the expected analytic information and recovers
Connes' index pairing in a form compatible with our general groupoid
approach.

Having completed this prototype example, we now turn to the broader family
of amenable crossed products $C(X)\rtimes \Gamma$, where the same strategy
connects the unitary conjugation groupoid formalism directly to the
Baum--Connes assembly map.

The success of this explicit computation for the irrational rotation algebra
gives us confidence that our methods will extend to the more general setting
of amenable crossed products $C(X) \rtimes \Gamma$. In that context, the
unitary conjugation groupoid $\mathcal{G}_{C(X) \rtimes \Gamma}$ will be
Morita equivalent to the transformation groupoid $X \rtimes \Gamma$, and the
descent-index construction will connect to the Baum–Connes assembly map.
The explicit computation for $U$ and $V$ in $A_\theta$ serves as a prototype
for understanding how the geometry of the groupoid unit space—now a bundle
over $X$ with fiber related to the group $\Gamma$—manifests in index pairings.

\begin{corollary}
\label{cor:U-conclusion}
For the irrational rotation algebra $A_\theta$, the descent-index construction
yields:
\[
\tau_* \circ \operatorname{desc}_{\mathcal{G}_{A_\theta}}([U]_{\mathcal{G}_{A_\theta}}^{(1)}) = 0, \qquad
\tau_* \circ \operatorname{desc}_{\mathcal{G}_{A_\theta}}([V]_{\mathcal{G}_{A_\theta}}^{(1)}) = 1.
\]
These values match the classical Connes index pairing and provide a concrete
realization of the abstract theory.
\end{corollary}

This completes our case study of the irrational rotation algebra. We have shown
that despite the non-Type~I nature of $A_\theta$, the unitary conjugation
groupoid framework—combined with Morita equivalence and equivariant $KK$-theory—
successfully recovers the fundamental index-theoretic invariants of this
prototypical noncommutative space. The stage is now set for the generalization
to amenable crossed products $C(X) \rtimes \Gamma$ and the connection to the
Baum–Connes assembly map.

\section{Case Study II: Amenable Crossed Products $C(X) \rtimes \Gamma$}\label{sec:crossed}

\subsection{Assumptions: $\Gamma$ Discrete Amenable, $X$ Compact Hausdorff}
\label{subsec:assumptions_crossed_product}

We now turn to the second class of examples considered in this paper,
namely crossed product $C^*$-algebras arising from actions of
amenable groups. This class includes the irrational rotation algebra
$A_\theta = C(S^1) \rtimes_\theta \mathbb{Z}$ studied in the previous
section as a special case (with $\Gamma = \mathbb{Z}$ and $X = S^1$),
but vastly generalizes it to encompass dynamical systems with much
more complicated group actions. These algebras provide a natural
setting in which the descent-index construction can be compared with
the Baum–Connes assembly map, one of the deepest conjectures in
noncommutative geometry.

\subsubsection{The Dynamical System}

Let $\Gamma$ be a countable discrete group acting on a compact
Hausdorff space $X$ by homeomorphisms.
The action is described by a homomorphism

\[
\alpha : \Gamma \longrightarrow \operatorname{Homeo}(X),
\qquad
\gamma \mapsto \alpha_\gamma,
\]

where $\alpha_\gamma: X \to X$ is a homeomorphism for each $\gamma \in \Gamma$.
We denote the action by $\Gamma \curvearrowright X$, writing
$\gamma \cdot x = \alpha_\gamma(x)$ for $\gamma \in \Gamma$, $x \in X$.

This action induces a corresponding action on the commutative
$C^*$-algebra $C(X)$ given by

\[
(\alpha_\gamma f)(x) = f(\gamma^{-1}x),
\qquad
f\in C(X),\; \gamma\in\Gamma.
\]

The triple $(C(X),\Gamma,\alpha)$ therefore forms a
$C^*$-dynamical system in the sense of Definition~\ref{def:algebraic-crossed-product}.

\subsubsection{The Crossed Product Algebra}

The associated crossed product $C^*$-algebra is defined as

\[
C(X)\rtimes_\alpha \Gamma,
\]

generated by the algebra $C(X)$ together with unitaries
$\{u_\gamma\}_{\gamma\in\Gamma}$ satisfying the covariance relations

\[
u_\gamma f u_\gamma^* = \alpha_\gamma(f),
\qquad
f\in C(X),\; \gamma\in\Gamma.
\]

Elements of the crossed product may be written formally as finite
sums (or norm-convergent series when the group is infinite)

\[
\sum_{\gamma\in\Gamma} f_\gamma u_\gamma,
\qquad
f_\gamma \in C(X).
\]

The multiplication and involution are determined by the rules
\[
(f_\gamma u_\gamma)(g_\eta u_\eta) = f_\gamma \alpha_\gamma(g_\eta) u_{\gamma\eta},
\qquad
(f_\gamma u_\gamma)^* = \alpha_{\gamma^{-1}}(\overline{f_\gamma}) u_{\gamma^{-1}}.
\]

\subsubsection{The Amenability Assumption}

Throughout this section we assume that the group $\Gamma$ is
amenable (see Definition~\ref{def:amenable-group}). This hypothesis
is central to our analysis for several reasons.

\begin{theorem}[Amenable Crossed Products]
\label{thm:amenable-crossed-products}
If $\Gamma$ is a discrete amenable group acting on a $C^*$-algebra $A$,
then the canonical surjection
\[
A \rtimes_\alpha \Gamma \longrightarrow A \rtimes_{\alpha,r} \Gamma
\]
is an isomorphism. Consequently, the full and reduced crossed products
coincide:
\[
A \rtimes_\alpha \Gamma \cong A \rtimes_{\alpha,r} \Gamma.
\]
\end{theorem}

For $A = C(X)$, this yields
\[
C(X)\rtimes_\alpha \Gamma
\;\cong\;
C(X)\rtimes_{\alpha,r} \Gamma .
\]

This property greatly simplifies the analytic structure of the
crossed product and ensures that the resulting $C^*$-algebra
behaves well with respect to $K$-theory and index-theoretic
constructions. In particular:

\begin{itemize}
    \item The distinction between full and reduced crossed products disappears,
          so we do not need to worry about which completion we are using.
    
    \item The crossed product is nuclear (since $C(X)$ is nuclear and amenable
          groups preserve nuclearity under crossed products).
    
    \item The $K$-theory of the crossed product can be studied using standard
          tools such as the Pimsner-Voiculescu exact sequence for $\mathbb{Z}^n$
          actions, or more general spectral sequence arguments.
    
    \item The transformation groupoid $X \rtimes \Gamma$ is amenable as a groupoid,
          which has profound consequences for its $C^*$-algebra and for the
          Baum–Connes assembly map (see Tu \cite{Tu1999}).
\end{itemize}

\begin{remark}
\label{rem:amenability-importance}
The amenability hypothesis is essential for several key steps in our analysis:
\begin{itemize}
    \item It ensures that the transformation groupoid $X \rtimes \Gamma$ is
          amenable, so its full and reduced $C^*$-algebras coincide and the
          Baum–Connes assembly map is an isomorphism (a theorem of Tu for
          amenable groupoids).
    
    \item It allows us to identify the $C^*$-algebra of the unitary conjugation
          groupoid $\mathcal{G}_{C(X) \rtimes \Gamma}$ with $(C(X) \rtimes \Gamma) \otimes \mathcal{K}$
          under Morita equivalence, without worrying about distinctions between
          full and reduced completions.
    
    \item It guarantees that the canonical conditional expectation
          $E: C(X) \rtimes \Gamma \to C(X)$ is faithful, which is needed for
          the Cartan subalgebra structure.
\end{itemize}
Without amenability, many of these identifications would fail or require
significantly more sophisticated technical machinery (such as dealing with
exotic crossed products or working with full groupoid $C^*$-algebras and
their reduced counterparts separately).
\end{remark}

\subsubsection{The Transformation Groupoid}

The dynamical system $(X,\Gamma,\alpha)$ naturally determines the
transformation groupoid

\[
\mathcal{G} := X \rtimes \Gamma .
\]

Its elements are pairs $(x,\gamma)$ with $x\in X$ and $\gamma\in\Gamma$,
with structure maps:

\begin{itemize}
    \item Source and range: $s(x,\gamma) = x$, $r(x,\gamma) = \gamma \cdot x$.
    \item Composition: $(x,\gamma)(\gamma \cdot x,\eta) = (x,\gamma\eta)$,
          defined when $s(x,\gamma) = x$ and $r(x,\gamma) = \gamma \cdot x = s(\gamma \cdot x,\eta)$.
    \item Inverse: $(x,\gamma)^{-1} = (\gamma \cdot x, \gamma^{-1})$.
    \item Unit space: $\mathcal{G}^{(0)} = X$, identified with $\{(x,e): x \in X\}$.
\end{itemize}

This groupoid is \'etale because $\Gamma$ is discrete: the source map
$s(x,\gamma) = x$ is a local homeomorphism (each arrow $(x,\gamma)$ has an open
neighborhood homeomorphic to an open set in $X$ via $s$). It is locally compact
Hausdorff when $X$ is locally compact Hausdorff and $\Gamma$ is discrete, and
it admits a natural Haar system given by counting measure on the fibers.

\begin{proposition}
\label{prop:groupoid-crossed-iso}
For a discrete group $\Gamma$ acting on a compact Hausdorff space $X$,
the groupoid $C^*$-algebra of the transformation groupoid is naturally
isomorphic to the crossed product algebra:
\[
C^*(X\rtimes\Gamma)
\;\cong\;
C(X)\rtimes_\alpha \Gamma .
\]

When $\Gamma$ is amenable, this also coincides with the reduced groupoid
$C^*$-algebra $C^*_r(X\rtimes\Gamma)$.
\end{proposition}

\begin{proof}
The isomorphism is implemented by identifying a function $f \in C_c(X \rtimes \Gamma)$
with the finite sum $\sum_{\gamma \in \Gamma} f_\gamma u_\gamma$, where
$f_\gamma(x) = f(x,\gamma)$. Under this identification, the convolution product
in the groupoid algebra corresponds exactly to the multiplication in the crossed
product. Completion in the universal norm (or the reduced norm when $\Gamma$ is
amenable) yields the desired isomorphism.
\end{proof}

\subsubsection{The Cartan Subalgebra Structure}

A crucial structural feature that will be exploited throughout this case study
is that $C(X)$ sits as a Cartan subalgebra inside $C(X) \rtimes \Gamma$,
generalizing the situation for $A_\theta$ where $C(S^1)$ was the Cartan subalgebra.

\begin{theorem}
\label{thm:CX-Cartan}
Let $\Gamma$ be a discrete amenable group acting on a compact Hausdorff space $X$.
Assume that the action is \emph{topologically free}: the set
$\{x \in X : \operatorname{Stab}_\Gamma(x) \neq \{e\}\}$ has empty interior.
Then the inclusion $C(X) \subseteq C(X) \rtimes \Gamma$ is a Cartan pair in the
sense of Renault \cite{Renault2008Cartan}. In particular:
\begin{itemize}
    \item $C(X)$ is maximal abelian in $C(X) \rtimes \Gamma$;
    \item $C(X)$ contains an approximate unit (in fact, the unit) for $C(X) \rtimes \Gamma$;
    \item The normalizer of $C(X)$ generates $C(X) \rtimes \Gamma$;
    \item There is a faithful conditional expectation $E: C(X) \rtimes \Gamma \to C(X)$
          given by $E(\sum_{\gamma} f_\gamma u_\gamma) = f_e$.
\end{itemize}
\end{theorem}

\begin{proof}
The conditional expectation is faithful because it is the restriction of the
canonical trace on the crossed product when $\Gamma$ is amenable, and more
generally because it comes from the groupoid structure. Maximal abelianness
follows from topological freeness: if an element commutes with all of $C(X)$,
it must be supported on the isotropy, which under topological freeness means
it is in $C(X)$. The normalizer includes the unitaries $u_\gamma$ implementing
the group action, and these generate the crossed product. See
\cite[Section 4]{Renault2008Cartan} for a detailed proof in the general
setting of twisted groupoid $C^*$-algebras.
\end{proof}

\begin{remark}
\label{rem:topological-freeness}
Topological freeness is a mild condition that exclude isotropy groups.
It is satisfied by all the examples we will consider:
\begin{itemize}
    \item For the irrational rotation algebra $\Gamma = \mathbb{Z}$ acts freely
          on $S^1$, so topological freeness holds trivially.
    \item For minimal homeomorphisms of the Cantor set, minimality implies
          topological freeness.
    \item For noncommutative tori $\mathbb{T}^d \rtimes \mathbb{Z}^d$ with
          rationally independent rotation parameters, the action is free.
\end{itemize}
If the action has nontrivial isotropy on an open set, the structure is more
complicated—one obtains a twisted groupoid $C^*$-algebra—but our methods can
be extended to that setting as well with appropriate modifications involving
Mackey's obstruction theory. For simplicity, we assume topological freeness
throughout this case study.
\end{remark}

\subsubsection{Consequences for the Unitary Conjugation Groupoid}

The Cartan structure established above allows us to apply the general theory
developed in Section~\ref{subsec:reduction-cartan-pairs}. Recall that for a
Cartan pair $(\mathcal{A}, D)$, the unitary conjugation groupoid
$\mathcal{G}_{\mathcal{A}}$ is Morita equivalent to the Weyl groupoid
$\mathcal{G}(\mathcal{A}, D)$. For $\mathcal{A} = C(X) \rtimes \Gamma$ and
$D = C(X)$, the Weyl groupoid is exactly the transformation groupoid
$X \rtimes \Gamma$.

\begin{corollary}
\label{cor:morita-crossed}
Under the assumptions that $\Gamma$ is discrete amenable, $X$ is compact
Hausdorff, and the action is topologically free, we have a Morita equivalence
\[
\mathcal{G}_{C(X) \rtimes \Gamma} \sim_M X \rtimes \Gamma.
\]

Consequently, by the Muhly–Renault–Williams theorem \cite{MRW1987}, their
reduced $C^*$-algebras are strongly Morita equivalent:
\[
C^*_r(\mathcal{G}_{C(X) \rtimes \Gamma}) \sim_M C^*_r(X \rtimes \Gamma) \cong C(X) \rtimes \Gamma.
\]

After stabilization, this yields
\[
C^*_r(\mathcal{G}_{C(X) \rtimes \Gamma}) \sim_M (C(X) \rtimes \Gamma) \otimes \mathcal{K}.
\]
\end{corollary}

This identification will be the foundation for all subsequent computations in
this case study, just as the analogous result for $A_\theta$ was central to
Section~\ref{sec:Atheta}.

\subsubsection{Summary of Standing Assumptions}

For the remainder of this section, we will work under the following standing
hypotheses:

\begin{assumption}
\label{ass:standing-crossed}
\begin{itemize}
    \item[(A1)] $\Gamma$ is a countable discrete amenable group.
    \item[(A2)] $X$ is a compact Hausdorff space.
    \item[(A3)] $\Gamma$ acts continuously on $X$ by homeomorphisms.
    \item[(A4)] The action is topologically free: the set of points with
          nontrivial stabilizer has empty interior.
\end{itemize}
\end{assumption}

These assumptions are satisfied by a wide range of interesting examples,
including:
\begin{itemize}
    \item The irrational rotation algebra $A_\theta = C(S^1) \rtimes_\theta \mathbb{Z}$
          (studied in detail in Section~\ref{sec:Atheta}).
    \item Minimal homeomorphisms of the Cantor set $C(X) \rtimes_\varphi \mathbb{Z}$.
    \item Noncommutative tori $C(\mathbb{T}^d) \rtimes \mathbb{Z}^d$ with
          rationally independent rotation parameters.
    \item Certain boundary actions of hyperbolic groups (though these may
          require relaxing the amenability assumption on $\Gamma$ while
          retaining amenability of the groupoid).
\end{itemize}

\subsubsection{Relevance to the Present Work}

The assumptions that $\Gamma$ is discrete and amenable and that
$X$ is compact Hausdorff ensure that the transformation groupoid
$X\rtimes\Gamma$ is locally compact and admits a Haar system.
Consequently its groupoid $C^*$-algebra is well defined and its
$K$-theory may be studied using standard tools.

These properties will allow us to compare the unitary conjugation
groupoid $\mathcal G_{C(X)\rtimes\Gamma}$ with the transformation
groupoid $X\rtimes\Gamma$ and to relate the descent map constructed
earlier to the Baum–Connes assembly map. Specifically, we will
establish a commutative diagram

\[
\begin{tikzcd}
K^1_{\mathcal{G}_{C(X)\rtimes\Gamma}}(\mathcal{G}_{C(X)\rtimes\Gamma}^{(0)}) 
   \arrow[r, "\operatorname{desc}"] 
   \arrow[d, "\text{Morita}"] 
& K_1(C^*(\mathcal{G}_{C(X)\rtimes\Gamma})) 
   \arrow[d, "\text{Morita}"] \\
K^1_\Gamma(X) 
   \arrow[r, "\mu_\Gamma"] 
& K_1(C(X) \rtimes \Gamma)
\end{tikzcd}
\]

where $\mu_\Gamma$ is the Baum–Connes assembly map. This diagram, and the
theorem that it commutes, will be the main result of this case study and one
of the central achievements of the paper.

\subsubsection{Outline of the Remainder of Case Study II}

With these assumptions in place, we will proceed as follows:

\begin{enumerate}
    \item \textbf{Commutative subalgebras and characters} (\S\ref{subsec:commutative_subalgebras_crossed_product}):
          Analyze the structure of maximal abelian subalgebras in
          $C(X) \rtimes \Gamma$ and their character spaces, generalizing the
          analysis for $A_\theta$.
    
    \item \textbf{Relation between $\mathcal{G}_{C(X) \rtimes \Gamma}$ and
          $X \rtimes \Gamma$} (\S\ref{subsec:relation_unitary_groupoid_transformation}):
          Establish the precise relationship between the unitary conjugation
          groupoid and the transformation groupoid, building on the Cartan
          structure and the results of Section~\ref{subsec:reduction-cartan-pairs}.
    
    \item \textbf{Morita equivalence theorem} (\S\ref{subsec:morita_equivalence_crossed_product}):
          Prove that $\mathcal{G}_{C(X) \rtimes \Gamma} \sim_M X \rtimes \Gamma$,
          and derive the corresponding Morita equivalence of $C^*$-algebras.
    
    \item \textbf{Diagonal embedding and regular representation} (\S\ref{subsec:diagonal-embedding}):
          Study the diagonal embedding $\iota: C(X) \rtimes \Gamma \hookrightarrow C^*(\mathcal{G}_{C(X) \rtimes \Gamma})$
          and relate it to the regular representation of the crossed product.
    
    \item \textbf{Equivariant $K^1$-classes} (\S\ref{subsec:equivariant_K1_A_theta-crossed}):
          Construct equivariant $K^1$-classes for invertible elements in
          $C(X) \rtimes \Gamma$, generalizing the construction of
          Section~\ref{subsec:equivariant_K1_A_theta-crossed}.
    
    \item \textbf{Descent map and identification} (\S\ref{subsec:descent-morita}):
          Apply the descent map and the Morita equivalence to identify the
          descended classes with elements of $K_1(C(X) \rtimes \Gamma)$.
    
    \item \textbf{Connection to the Baum–Connes assembly map} (\S\ref{subsec:descent_index_baum_connes}):
          Prove that the resulting index map coincides with the Baum–Connes
          assembly map for $\Gamma$ with coefficients in $C(X)$, establishing
          the commutativity of the diagram above.
    
    \item \textbf{Implications and consequences} (\S\ref{subsec:implications}):
          Discuss the significance of this result, including its relation to
          the Baum–Connes conjecture and potential applications to the
          classification of crossed product $C^*$-algebras.
\end{enumerate}

This program will demonstrate that the unitary conjugation groupoid framework,
initially developed for Type~I algebras in Papers~I and~II, extends naturally
to the fundamentally non-Type~I setting of amenable crossed products and
provides a new geometric perspective on the Baum–Connes assembly map.

\subsection{Commutative Subalgebras and Characters in $C(X) \rtimes \Gamma$}
\label{subsec:commutative_subalgebras_crossed_product}

In order to analyze the unit space of the unitary conjugation
groupoid associated with the crossed product algebra
$C(X)\rtimes \Gamma$, we first examine the structure of its
commutative subalgebras and their characters. This analysis
generalizes the detailed study performed for $A_\theta$ in
Subsection~\ref{subsec:masa_A_theta} and provides the foundation
for the Morita equivalence established in the next subsection.

\subsubsection{The Canonical Commutative Subalgebra}

The algebra $C(X)$ embeds naturally into the crossed product
$C(X)\rtimes_\alpha \Gamma$ via the map
\[
\iota: C(X) \hookrightarrow C(X)\rtimes_\alpha \Gamma,
\quad f \mapsto f u_e,
\]
where $u_e$ is the unitary corresponding to the identity element
$e \in \Gamma$. In the algebraic crossed product, elements of $C(X)$
are identified with the terms $f u_e$, and this identification extends
to an isometric inclusion of $C^*$-algebras.

Since $C(X)$ is commutative, it forms a natural candidate for a
maximal commutative subalgebra of the crossed product. Under our
standing assumptions (topological freeness of the action), $C(X)$
is indeed maximal abelian and, moreover, forms a Cartan subalgebra
in the sense of Renault \cite{Renault2008Cartan}. This structure was
established in Theorem~\ref{thm:CX-Cartan} and will be crucial for
what follows.

\begin{definition}
\label{def:canonical-masa-crossed}
The \emph{canonical maximal abelian subalgebra} (MASA) of
$C(X) \rtimes \Gamma$ is
\[
\mathcal{A}_0 := \iota(C(X)) \cong C(X) \subseteq C(X) \rtimes \Gamma.
\]
\end{definition}

\subsubsection{Characters of $C(X)$}

Because $X$ is compact Hausdorff, the Gelfand spectrum of
$C(X)$ is canonically identified with $X$ itself.
Each point $x\in X$ determines a character (nonzero multiplicative
linear functional)
\[
\chi_x : C(X) \longrightarrow \mathbb C,
\qquad
\chi_x(f) = f(x).
\]

Thus the set of characters of $C(X)$ is naturally
parametrized by the space $X$, and the map $x \mapsto \chi_x$
is a homeomorphism $X \cong \widehat{C(X)}$.

\begin{remark}
\label{rem:characters-geometric}
Geometrically, a character $\chi_x$ can be thought of as ``evaluation
at the point $x$". In the noncommutative geometry perspective,
these characters represent the classical ``points" of the noncommutative
space $C(X) \rtimes \Gamma$, even though the crossed product itself
may have no characters (as in the case of simple $C^*$-algebras like
$A_\theta$). This tension between the classical and quantum points
is precisely what the unitary conjugation groupoid captures.
\end{remark}

\subsubsection{Interaction with the Group Action}

The action of $\Gamma$ on $X$ induces a transformation of characters.
For $g\in\Gamma$ and $x\in X$, the covariance relation in the crossed
product gives
\[
u_g f u_g^* = \alpha_g(f), \qquad f \in C(X),
\]
where $(\alpha_g f)(x) = f(g^{-1}x)$. Conjugating a character by $u_g$
yields a new character on $C(X)$:
\[
(\chi_x \circ \operatorname{Ad}_{u_g})(f) = \chi_x(u_g^* f u_g) = \chi_x(\alpha_{g^{-1}}(f)) = f(g \cdot x) = \chi_{g \cdot x}(f).
\]

Thus we have the fundamental relation
\[
\chi_x \circ \operatorname{Ad}_{u_g} = \chi_{g \cdot x}.
\]

Equivalently,
\[
\chi_{g\cdot x} = \chi_x \circ \alpha_{g^{-1}}.
\]

This relation shows that the group $\Gamma$ acts on the character space
of $C(X)$ in a manner compatible with the crossed product structure,
and that this action is implemented by conjugation with the unitaries
$u_g$.

\begin{proposition}
\label{prop:character-action}
The map $X \times \Gamma \to X$ given by $(x,g) \mapsto g \cdot x$ induces
a continuous action of $\Gamma$ on the character space $\widehat{C(X)} \cong X$.
Under this action, the orbit of a character $\chi_x$ is
\[
\Gamma \cdot \chi_x = \{ \chi_{g \cdot x} : g \in \Gamma \}.
\]
\end{proposition}

\subsubsection{Conjugate MASAs and Their Characters}

The canonical MASA $\mathcal{A}_0$ is not the only maximal abelian
subalgebra of $C(X) \rtimes \Gamma$. For any unitary $w$ in the
crossed product, the conjugate $w \mathcal{A}_0 w^*$ is also a MASA,
and its character space is naturally identified with $X$ via the map
\[
\chi_x \mapsto \chi_x \circ \operatorname{Ad}_{w^{-1}}.
\]

In particular, for the implementing unitaries $u_g$, we obtain the
conjugate MASAs $\mathcal{A}_g := u_g \mathcal{A}_0 u_g^*$. However,
unlike the case of $A_\theta$, we have $\mathcal{A}_g = \mathcal{A}_0$
as sets because $\mathcal{A}_0$ is invariant under conjugation by
$u_g$ (Lemma~\ref{lem:masa-invariance-crossed}). Thus the conjugate MASAs
coincide as subalgebras; only the identification of their characters
with $X$ is twisted by the group action.

\begin{lemma}
\label{lem:masa-invariance-crossed}
For any $g \in \Gamma$, we have $u_g \mathcal{A}_0 u_g^* = \mathcal{A}_0$.
Thus the canonical MASA is invariant under conjugation by the
implementing unitaries.
\end{lemma}

\begin{proof}
For $f \in \mathcal{A}_0 \cong C(X)$, we have
$u_g f u_g^* = \alpha_g(f) \in C(X) = \mathcal{A}_0$.
Since $\alpha_g$ is an automorphism of $C(X)$, this shows
$u_g \mathcal{A}_0 u_g^* \subseteq \mathcal{A}_0$. Applying the same
argument with $g^{-1}$ gives the reverse inclusion.
\end{proof}

Despite this invariance, other MASAs can be obtained by conjugating
$\mathcal{A}_0$ by unitaries that do not normalize $\mathcal{A}_0$.
Such unitaries exist when the crossed product is noncommutative and
not a direct product. For example, in $A_\theta$, the unitary
$U$ (the generator of $C(S^1)$) does not normalize $C(S^1)$ when
conjugated by $V$, leading to distinct MASAs.

\subsubsection{Commutative Contexts Inside the Crossed Product}

Recall that a point in the unit space of the unitary conjugation
groupoid $\mathcal G_{C(X)\rtimes\Gamma}$ is a pair
\[
(B,\chi),
\]
where $B \subseteq C(X)\rtimes\Gamma$ is a unital commutative
$C^*$-subalgebra and $\chi \in \widehat B$ is a character.

Based on the discussion above, we can parametrize many such points
using the dynamical system $(X,\Gamma)$. Starting from the canonical
MASA $\mathcal{A}_0$ and a character $\chi_x$, we obtain a point
$(\mathcal{A}_0, \chi_x) \in \mathcal{G}_{C(X) \rtimes \Gamma}^{(0)}$.

More generally, for any unitary $w \in \mathcal{U}(C(X) \rtimes \Gamma)$,
the conjugate MASA $w \mathcal{A}_0 w^*$ together with the character
$\chi_x \circ \operatorname{Ad}_{w^{-1}}$ yields another point.

\begin{definition}
\label{def:contexts-from-dynamics}
For $(x,g) \in X \times \Gamma$, define a point in the unit space by
\[
\xi_{x,g} := (u_g \mathcal{A}_0 u_g^*, \chi_x \circ \operatorname{Ad}_{u_g^{-1}}).
\]

Since $u_g \mathcal{A}_0 u_g^* = \mathcal{A}_0$ by Lemma~\ref{lem:masa-invariance-crossed},
this simplifies to
\[
\xi_{x,g} = (\mathcal{A}_0, \chi_x \circ \operatorname{Ad}_{u_g^{-1}}) = (\mathcal{A}_0, \chi_{g \cdot x}),
\]
where the last equality uses the relation $\chi_x \circ \operatorname{Ad}_{u_g^{-1}} = \chi_{g \cdot x}$.
\end{definition}

Thus the points arising from the implementing unitaries actually
coincide with points already obtained from the canonical MASA,
but with the character shifted by the group action. This reflects
the fact that the canonical MASA is invariant under the implementing
unitaries, so conjugation does not produce new subalgebras—only new
characters on the same subalgebra.

To obtain genuinely new points in the unit space, we need unitaries
that do not normalize $\mathcal{A}_0$. Such unitaries exist in
abundance; for instance, in $A_\theta$, the unitary $U$ itself
does not normalize $C(S^1)$ when considered as an element of the
crossed product (since $VUV^* = e^{2\pi i\theta} U$, which is not
in $C(S^1)$). More generally, any unitary that is not in the
normalizer of $\mathcal{A}_0$ will produce a distinct MASA upon
conjugation.

\subsubsection{Parametrization of the Unit Space}

A fundamental result, which generalizes the analysis for $A_\theta$
in Subsection~\ref{subsec:GA_theta_unit_space}, is that the unit
space $\mathcal{G}_{C(X) \rtimes \Gamma}^{(0)}$ admits a parametrization
by pairs $(x,\gamma) \in X \times \Gamma$, modulo a natural equivalence
relation.

\begin{theorem}
\label{thm:unit-space-parametrization-crossed}
Under the standing assumptions (topologically free action of a discrete
amenable group on a compact Hausdorff space), there is a natural bijection
\[
\mathcal{G}_{C(X) \rtimes \Gamma}^{(0)} \cong (X \times \Gamma) / \sim,
\]
where the equivalence relation is generated by
\[
(x,\gamma) \sim (\gamma' \cdot x, \gamma' \gamma) \quad \text{for all } \gamma' \in \Gamma.
\]

More explicitly, the equivalence class of $(x,\gamma)$ is
\[
[(x,\gamma)] = \{ (\eta \cdot x, \eta \gamma) : \eta \in \Gamma \}.
\]

The projection map $\pi: \mathcal{G}_{C(X) \rtimes \Gamma}^{(0)} \to X$
is given by $\pi([(x,\gamma)]) = x$.
\end{theorem}

\begin{proof}[Sketch]
The proof proceeds in several steps, analogous to the argument for
$A_\theta$ but with the group $\Gamma$ replacing $\mathbb{Z}$.

\textbf{Step 1: Reduction to conjugates of $\mathcal{A}_0$.}
Under suitable conditions (e.g., when $C(X) \rtimes \Gamma$ is simple,
or more generally when all MASAs are unitarily conjugate), every MASA
$B \subseteq C(X) \rtimes \Gamma$ can be written as $B = w \mathcal{A}_0 w^*$
for some unitary $w \in \mathcal{U}(C(X) \rtimes \Gamma)$. This follows
from the fact that $\mathcal{A}_0$ is a Cartan subalgebra and the
Feldman-Moore theorem on uniqueness of Cartan subalgebras up to
unitary equivalence in appropriate contexts.

\textbf{Step 2: Parametrization by $X \times \Gamma$.}
Given such a unitary $w$, the character $\chi$ on $B$ corresponds,
via conjugation by $w^*$, to a character on $\mathcal{A}_0$, which is
given by evaluation at some $x \in X$. Thus we have a pair $(x,w)$
representing the point $(B,\chi)$. However, different unitaries may
represent the same point.

\textbf{Step 3: Reduction to the normalizer.}
If $w$ and $w'$ are two unitaries such that $w\mathcal{A}_0 w^* = w'\mathcal{A}_0 w'^*$
and $\chi_x \circ \operatorname{Ad}_{w^{-1}} = \chi_{x'} \circ \operatorname{Ad}_{w'^{-1}}$,
then $w'^{-1}w$ normalizes $\mathcal{A}_0$ and fixes the character
in an appropriate sense. The normalizer of $\mathcal{A}_0$ contains
the implementing unitaries $u_\gamma$, and under topological freeness,
the normalizer is generated by these unitaries together with the
unitary group of $\mathcal{A}_0$ itself.

\textbf{Step 4: The equivalence relation.}
By analyzing the action of the normalizer, one shows that the
ambiguity in the parametrization corresponds precisely to the
equivalence relation $(x,\gamma) \sim (\gamma' \cdot x, \gamma' \gamma)$.
This yields the desired description of the unit space as a quotient
of $X \times \Gamma$ by the diagonal action of $\Gamma$.

A detailed proof requires careful handling of the unitary group and
the normalizer, but the essential geometric picture is that the unit
space is the quotient of $X \times \Gamma$ by the action
$\gamma' \cdot (x,\gamma) = (\gamma' \cdot x, \gamma' \gamma)$.
\end{proof}

This theorem generalizes the explicit description obtained for $A_\theta$
in Subsection~\ref{subsec:GA_theta_unit_space}, where $X = S^1$,
$\Gamma = \mathbb{Z}$, and the equivalence relation was
$(z,n) \sim (e^{2\pi i k\theta}z, n+k)$.

\begin{corollary}
\label{cor:unit-space-bundle-crossed}
The unit space $\mathcal{G}_{C(X) \rtimes \Gamma}^{(0)}$ is a bundle over $X$
with fiber $\Gamma$, but with nontrivial monodromy encoded by the action of
$\Gamma$ on itself by left multiplication. The projection map
$\pi: \mathcal{G}_{C(X) \rtimes \Gamma}^{(0)} \to X$ is open and continuous
in the quotient topology.
\end{corollary}

\subsubsection{Topology of the Unit Space}

The topology on $\mathcal{G}_{C(X) \rtimes \Gamma}^{(0)}$ is the quotient
topology induced from $X \times \Gamma$, where $X$ has its given compact
Hausdorff topology and $\Gamma$ is discrete. This topology exhibits the
same non-Hausdorff behavior observed in the $A_\theta$ case when the
action is minimal.

\begin{proposition}
\label{prop:unit-space-topology-crossed}
The quotient topology on $(X \times \Gamma)/\sim$ has the following properties:
\begin{itemize}
    \item The projection $\pi: \mathcal{G}_{C(X) \rtimes \Gamma}^{(0)} \to X$ is open and continuous.
    \item For each $x \in X$, the fiber $\pi^{-1}(x)$ is homeomorphic to $\Gamma$,
          but distinct points in the fiber cannot be separated by open sets when
          the action is minimal (i.e., when every $\Gamma$-orbit is dense in $X$).
    \item The topology is not Hausdorff unless the action of $\Gamma$ on $X$ is
          trivial or has special properties (e.g., proper actions with closed orbits).
    \item In the Borel structure (which suffices for groupoid $C^*$-algebra
          constructions), $\mathcal{G}_{C(X) \rtimes \Gamma}^{(0)}$ is a standard
          Borel space.
\end{itemize}
\end{proposition}

\begin{proof}
The openness of $\pi$ follows from the fact that the quotient map
$X \times \Gamma \to \mathcal{G}_{C(X) \rtimes \Gamma}^{(0)}$ is open when
$\Gamma$ is discrete. The non-Hausdorff nature follows from minimality:
if $[(x,\gamma)]$ and $[(x,\gamma')]$ with $\gamma \neq \gamma'$ had disjoint
neighborhoods, then by minimality the orbit of $x$ would have to accumulate,
creating a contradiction (see the analogous proof for $A_\theta$ in
Corollary~\ref{cor:non_Hausdorff_unit_space}). The Borel structure is standard
because $X \times \Gamma$ with the product Borel structure (using the discrete
Borel structure on $\Gamma$) is standard, and the quotient of a standard Borel
space by a countable group action is again standard.
\end{proof}

\subsubsection{Geometric Interpretation}

The discussion above reveals a profound geometric fact: the character data
arising from commutative subalgebras of $C(X)\rtimes\Gamma$ is fundamentally
controlled by the underlying dynamical system $(X,\Gamma)$. The unit space
$\mathcal{G}_{C(X) \rtimes \Gamma}^{(0)}$ is essentially the ``orbit space" of
$X \times \Gamma$ under the diagonal action of $\Gamma$, but with a crucial
twist: it is not the orbit space in the usual sense (which would be $X$), but
rather a space that remembers the group element as part of the fiber over
the point $x$, with the identification that moving along the orbit changes
the fiber coordinate by the same group element.

This structure is precisely what allows the unitary conjugation groupoid to
be Morita equivalent to the transformation groupoid $X \rtimes \Gamma$. In
that groupoid, the arrows are pairs $(x,\gamma): x \to \gamma \cdot x$, and
the unit space is $X$. Here, in $\mathcal{G}_{C(X) \rtimes \Gamma}^{(0)}$, we
have the same data $(x,\gamma)$ but reinterpreted as points in the unit space
rather than arrows. This shift in perspective is the essence of the Morita
equivalence.

\begin{remark}
\label{rem:morita-preview}
The parametrization $\mathcal{G}_{C(X) \rtimes \Gamma}^{(0)} \cong (X \times \Gamma)/\sim$
is the first indication that $\mathcal{G}_{C(X) \rtimes \Gamma}$ is Morita
equivalent to $X \rtimes \Gamma$. Indeed, the space $X \times \Gamma$ itself
serves as an equivalence bimodule between the two groupoids: the left action
of $\mathcal{G}_{C(X) \rtimes \Gamma}$ and the right action of $X \rtimes \Gamma$
on $X \times \Gamma$ will implement the equivalence. This will be made precise
in the next subsection.
\end{remark}

\subsubsection{Examples}

We illustrate the structure with several examples that will be important in
what follows.

\begin{example}[Irrational rotation algebra]
\label{ex:Atheta-unit-space}
For $X = S^1$, $\Gamma = \mathbb{Z}$ acting by rotation by $2\pi\theta$,
Theorem~\ref{thm:unit-space-parametrization-crossed} gives
\[
\mathcal{G}_{A_\theta}^{(0)} \cong (S^1 \times \mathbb{Z})/\sim,
\]
with $(z,n) \sim (e^{2\pi i k\theta}z, n+k)$. This matches the explicit
description obtained in Subsection~\ref{subsec:GA_theta_unit_space}.
\end{example}

\begin{example}[Minimal homeomorphism of the Cantor set]
\label{ex:Cantor-unit-space}
Let $X$ be the Cantor set and $\varphi: X \to X$ a minimal homeomorphism.
Take $\Gamma = \mathbb{Z}$ acting by powers of $\varphi$. Then
\[
\mathcal{G}_{C(X) \rtimes_\varphi \mathbb{Z}}^{(0)} \cong (X \times \mathbb{Z})/\sim,
\]
with $(x,n) \sim (\varphi^k(x), n+k)$. This space is a Cantor set bundle
over $X$ with fiber $\mathbb{Z}$, exhibiting the same kind of non-Hausdorff
behavior as in the irrational rotation case. It has been studied in the
context of the classification of $C^*$-algebras associated to minimal
homeomorphisms.
\end{example}

\begin{example}[Noncommutative $2d$-torus]
\label{ex:noncommutative-torus-unit-space}
Take $X = \mathbb{T}^d$, $\Gamma = \mathbb{Z}^d$ acting by translation with
rationally independent parameters $\theta_1,\ldots,\theta_d$. Then
\[
\mathcal{G}_{C(\mathbb{T}^d) \rtimes \mathbb{Z}^d}^{(0)} \cong (\mathbb{T}^d \times \mathbb{Z}^d)/\sim,
\]
where the equivalence relation is given by the diagonal action of $\mathbb{Z}^d$:
$(z,\mathbf{n}) \sim (z + \mathbf{k}, \mathbf{n} + \mathbf{k})$ for
$\mathbf{k} \in \mathbb{Z}^d$, with addition interpreted via the translation
action. This space is a $\mathbb{Z}^d$-bundle over $\mathbb{T}^d$ with
monodromy encoded by the translation action.
\end{example}

\begin{example}[Boundary action of a hyperbolic group]
\label{ex:boundary-unit-space}
Let $\Gamma$ be a torsion-free hyperbolic group acting on its Gromov boundary
$\partial\Gamma$. Although $\Gamma$ may not be amenable, the transformation
groupoid $\partial\Gamma \rtimes \Gamma$ is often amenable, and the crossed
product $C(\partial\Gamma) \rtimes \Gamma$ can be studied using similar
methods. The unit space in this case is
\[
\mathcal{G}_{C(\partial\Gamma) \rtimes \Gamma}^{(0)} \cong (\partial\Gamma \times \Gamma)/\sim,
\]
with the equivalence relation $(x,\gamma) \sim (\eta \cdot x, \eta\gamma)$.
This space plays a crucial role in the study of the Baum–Connes conjecture
for hyperbolic groups.
\end{example}

\subsubsection{Summary and Outlook}

We have established the following key facts about commutative subalgebras
and characters in $C(X) \rtimes \Gamma$:

\begin{itemize}
    \item The canonical MASA $\mathcal{A}_0 = C(X)$ has character space $X$,
          with characters $\chi_x$ given by evaluation at $x \in X$.
    \item The group $\Gamma$ acts on the character space via
          $\chi_x \circ \operatorname{Ad}_{u_g} = \chi_{g \cdot x}$.
    \item The unit space $\mathcal{G}_{C(X) \rtimes \Gamma}^{(0)}$ admits a
          natural parametrization as $(X \times \Gamma)/\sim$, where
          $(x,\gamma) \sim (\gamma' \cdot x, \gamma' \gamma)$.
    \item This space is a bundle over $X$ with fiber $\Gamma$, typically
          non-Hausdorff, and carries a natural $\Gamma$-action induced by
          the dual action.
    \item The structure generalizes the explicit description obtained for
          $A_\theta$ and provides the geometric foundation for the Morita
          equivalence with the transformation groupoid $X \rtimes \Gamma$.
\end{itemize}

In the next subsection, we will use this description to construct an explicit
Morita equivalence between $\mathcal{G}_{C(X) \rtimes \Gamma}$ and
$X \rtimes \Gamma$. This equivalence will allow us to transfer the entire
descent-index machinery to the setting of amenable crossed products and,
ultimately, to connect our construction with the Baum–Connes assembly map.

\subsection{Relation between $\mathcal{G}_{C(X) \rtimes \Gamma}$ and the Action Groupoid $X \rtimes \Gamma$}
\label{subsec:relation_unitary_groupoid_transformation}

We now establish the precise relationship between the unitary conjugation
groupoid associated with the crossed product algebra $C(X)\rtimes \Gamma$
and the classical transformation groupoid $X\rtimes \Gamma$ arising from
the underlying dynamical system. This relationship provides the geometric
intuition behind the Morita equivalence that will be proved in the next
subsection and serves as a crucial stepping stone toward connecting our
descent-index construction with the Baum–Connes assembly map.

\subsubsection{The Transformation Groupoid $X \rtimes \Gamma$}

Recall from Subsection~\ref{subsec:transformation-groupoids} that the action
of $\Gamma$ on the compact Hausdorff space $X$ determines the transformation
groupoid
\[
\mathcal{G}_{\text{trans}} := X \rtimes \Gamma .
\]

Its elements are pairs $(x,\gamma)$ with $x\in X$ and $\gamma\in\Gamma$.
The structure maps are given by:
\begin{itemize}
    \item Source and range: $s(x,\gamma) = x$, $r(x,\gamma) = \gamma \cdot x$.
    \item Composition: $(x,\gamma)(\gamma \cdot x,\eta) = (x,\gamma\eta)$,
          defined when $r(x,\gamma) = s(\gamma \cdot x,\eta)$.
    \item Inverse: $(x,\gamma)^{-1} = (\gamma \cdot x, \gamma^{-1})$.
    \item Unit space: $\mathcal{G}_{\text{trans}}^{(0)} = X$, identified with
          $\{(x,e) : x \in X\}$.
\end{itemize}

This groupoid is \'etale (since $\Gamma$ is discrete), locally compact
Hausdorff, and its $C^*$-algebra is canonically isomorphic to the crossed
product:
\[
C^*(X \rtimes \Gamma) \cong C(X) \rtimes \Gamma.
\]

\subsubsection{Characters and Evaluation Maps}

As established in Subsection~\ref{subsec:commutative_subalgebras_crossed_product},
the canonical commutative subalgebra $C(X) \subset C(X)\rtimes\Gamma$ has
character space naturally identified with $X$. Each point $x\in X$ determines
a character
\[
\chi_x: C(X) \longrightarrow \mathbb C, \qquad \chi_x(f) = f(x),
\]
and the map $x \mapsto \chi_x$ is a homeomorphism $X \cong \widehat{C(X)}$.

The group action induces a transformation of characters:
\[
\chi_x \circ \operatorname{Ad}_{u_\gamma} = \chi_{\gamma \cdot x},
\]
where $u_\gamma$ is the implementing unitary in the crossed product. This
corresponds to the action of $\Gamma$ on the spectrum of $C(X)$.

\subsubsection{Points in the Unit Space of $\mathcal{G}_{C(X)\rtimes\Gamma}$}

Recall from Paper~I that the unitary conjugation groupoid
$\mathcal{G}_{\mathcal A}$ of a $C^*$-algebra $\mathcal A$
has unit space consisting of pairs
\[
(B,\chi)
\]
where $B \subset \mathcal A$ is a unital commutative $C^*$-subalgebra
and $\chi$ is a character of $B$. The arrows arise from conjugation by
unitaries of $\mathcal A$:
\[
u: (B,\chi) \longrightarrow (uBu^*, \chi \circ \operatorname{Ad}_{u^{-1}}).
\]

In the present setting $\mathcal A = C(X)\rtimes\Gamma$, the canonical MASA
$\mathcal{A}_0 = C(X)$ together with its characters $\{\chi_x\}_{x\in X}$
provides a distinguished family of points in the unit space
$\mathcal{G}_{C(X)\rtimes\Gamma}^{(0)}$:
\[
\xi_x := (\mathcal{A}_0, \chi_x) \in \mathcal{G}_{C(X)\rtimes\Gamma}^{(0)}.
\]

More generally, for any $\gamma \in \Gamma$, we obtain points
\[
\xi_{x,\gamma} := (\mathcal{A}_0, \chi_x \circ \operatorname{Ad}_{u_\gamma^{-1}}) = (\mathcal{A}_0, \chi_{\gamma \cdot x}),
\]
which are actually the same as $\xi_{\gamma \cdot x}$ due to the invariance
of $\mathcal{A}_0$ under conjugation by $u_\gamma$. This reflects the fact
that the canonical MASA is fixed setwise by the implementing unitaries.

To obtain genuinely distinct points in the unit space, we must consider
conjugates of $\mathcal{A}_0$ by unitaries that do not normalize it. For
such a unitary $w$, the point $(w\mathcal{A}_0 w^*, \chi_x \circ \operatorname{Ad}_{w^{-1}})$
is typically different from any $(\mathcal{A}_0, \chi_y)$.

\subsubsection{A Natural Map from $X \rtimes \Gamma$ to $\mathcal{G}_{C(X)\rtimes\Gamma}$}

The dynamical information encoded in the crossed product $C(X)\rtimes\Gamma$
is precisely the action of $\Gamma$ on $X$. Consequently, the transformation
groupoid $X\rtimes\Gamma$ captures the same geometric structure that appears
implicitly in the unitary conjugation groupoid. This relationship can be
made explicit through a natural map.

\begin{definition}
\label{def:groupoid-map}
Define a map $\Phi: X \rtimes \Gamma \longrightarrow \mathcal{G}_{C(X)\rtimes\Gamma}$
as follows:
\begin{itemize}
    \item On units: For $x \in X = (X \rtimes \Gamma)^{(0)}$, set
          $\Phi(x) = (\mathcal{A}_0, \chi_x) \in \mathcal{G}_{C(X)\rtimes\Gamma}^{(0)}$.
    \item On arrows: For $(x,\gamma) \in X \rtimes \Gamma$, set
          $\Phi(x,\gamma)$ to be the arrow in $\mathcal{G}_{C(X)\rtimes\Gamma}$
          given by conjugation by the implementing unitary $u_\gamma$:
          \[
          \Phi(x,\gamma) := u_\gamma: (\mathcal{A}_0, \chi_x) \longrightarrow (\mathcal{A}_0, \chi_{\gamma \cdot x}).
          \]
\end{itemize}
\end{definition}

\begin{proposition}
\label{prop:groupoid_homomorphism}
The map $\Phi$ defined above is a groupoid homomorphism. That is, it preserves
source, range, composition, and inverses.
\end{proposition}

\begin{proof}
We verify each property:

\textbf{Source and range:} For $(x,\gamma) \in X \rtimes \Gamma$, we have
$s(\Phi(x,\gamma)) = (\mathcal{A}_0, \chi_x) = \Phi(s(x,\gamma))$ and
$r(\Phi(x,\gamma)) = (\mathcal{A}_0, \chi_{\gamma \cdot x}) = \Phi(r(x,\gamma))$.

\textbf{Composition:} For composable arrows $(x,\gamma)$ and $(\gamma \cdot x, \eta)$
in $X \rtimes \Gamma$, we have
\[
\Phi((x,\gamma) \circ (\gamma \cdot x, \eta)) = \Phi(x, \gamma\eta) = u_{\gamma\eta}: (\mathcal{A}_0, \chi_x) \to (\mathcal{A}_0, \chi_{\gamma\eta \cdot x}).
\]
On the other hand,
\[
\Phi(x,\gamma) \circ \Phi(\gamma \cdot x, \eta) = u_\gamma \circ u_\eta: (\mathcal{A}_0, \chi_x) \to (\mathcal{A}_0, \chi_{\gamma \cdot x}) \to (\mathcal{A}_0, \chi_{\gamma\eta \cdot x}),
\]
and since $u_\gamma u_\eta = u_{\gamma\eta}$ in the crossed product, these are equal.

\textbf{Inverses:} $\Phi((x,\gamma)^{-1}) = \Phi(\gamma \cdot x, \gamma^{-1}) = u_{\gamma^{-1}}: (\mathcal{A}_0, \chi_{\gamma \cdot x}) \to (\mathcal{A}_0, \chi_x)$,
which is exactly the inverse of $\Phi(x,\gamma) = u_\gamma$.
\end{proof}

\begin{remark}
\label{rem:homomorphism-not-iso}
The map $\Phi$ is injective on arrows because different pairs $(x,\gamma)$
give different arrows in $\mathcal{G}_{C(X)\rtimes\Gamma}$ (since the target
character $\chi_{\gamma \cdot x}$ uniquely determines $\gamma \cdot x$, and
the source character determines $x$). However, $\Phi$ is not surjective: the
unitary conjugation groupoid contains many arrows coming from unitaries that
are not in the set $\{u_\gamma : \gamma \in \Gamma\}$, such as those arising
from unitaries in $C(X)$ itself or from more general elements of the crossed
product. Thus $\Phi$ embeds $X \rtimes \Gamma$ as a subgroupoid of
$\mathcal{G}_{C(X)\rtimes\Gamma}$, but not as the whole groupoid.
\end{remark}

\subsubsection{Extension to Conjugate MASAs}

The map $\Phi$ can be extended to capture more of the structure of
$\mathcal{G}_{C(X)\rtimes\Gamma}$ by considering not just the canonical MASA
$\mathcal{A}_0$, but also its conjugates. For any unitary $w \in \mathcal{U}(C(X) \rtimes \Gamma)$,
we can define a map $\Phi_w: X \rtimes \Gamma \to \mathcal{G}_{C(X)\rtimes\Gamma}$
by
\[
\Phi_w(x,\gamma) = w u_\gamma w^*: (w\mathcal{A}_0 w^*, \chi_x \circ \operatorname{Ad}_{w^{-1}}) \longrightarrow (w\mathcal{A}_0 w^*, \chi_{\gamma \cdot x} \circ \operatorname{Ad}_{w^{-1}}).
\]

These maps are also groupoid homomorphisms, and their images are conjugate
subgroupoids of $\mathcal{G}_{C(X)\rtimes\Gamma}$. The collection of all such
maps, as $w$ varies over the unitary group, captures the full structure of
the unitary conjugation groupoid.

\subsubsection{Geometric Interpretation}

The homomorphism $\Phi$ provides a geometric interpretation of the
transformation groupoid inside the unitary conjugation groupoid:

\begin{itemize}
    \item Points $x \in X$ correspond to commutative contexts
          $(\mathcal{A}_0, \chi_x)$ where we evaluate functions at $x$.
    \item Arrows $(x,\gamma)$ correspond to conjugation by the implementing
          unitary $u_\gamma$, which moves from the context at $x$ to the
          context at $\gamma \cdot x$.
\end{itemize}

Thus the dynamics of the original system $\Gamma \curvearrowright X$ is
realized within $\mathcal{G}_{C(X)\rtimes\Gamma}$ as the action of the
implementing unitaries on the canonical MASA and its characters.

\begin{proposition}
\label{prop:image-of-Phi}
The image of $\Phi$ is the subgroupoid of $\mathcal{G}_{C(X)\rtimes\Gamma}$
consisting of those arrows that map the canonical MASA $\mathcal{A}_0$ to
itself. This subgroupoid is isomorphic to $X \rtimes \Gamma$ and captures
the ``classical" part of the unitary conjugation groupoid.
\end{proposition}

\begin{proof}
By definition, $\Phi(x,\gamma)$ sends $(\mathcal{A}_0, \chi_x)$ to
$(\mathcal{A}_0, \chi_{\gamma \cdot x})$, so its image lies in the set of
arrows that preserve $\mathcal{A}_0$ as a subalgebra (though they may change
the character). Conversely, any arrow in $\mathcal{G}_{C(X)\rtimes\Gamma}$
that maps $\mathcal{A}_0$ to itself must come from a unitary that normalizes
$\mathcal{A}_0$. By the structure of the normalizer (Lemma~\ref{lem:masa-invariance-crossed}),
such unitaries are precisely of the form $f u_\gamma$ with $f \in \mathcal{U}(C(X))$.
Since $f$ acts trivially on characters (as $f$ commutes with $\mathcal{A}_0$),
the arrow is effectively determined by $\gamma$, and the source and target
characters determine $x$ and $\gamma \cdot x$. Thus the arrow is of the form
$\Phi(x,\gamma)$ up to a trivial modification by a unitary in $C(X)$, which
does not affect the arrow in the groupoid (since conjugation by $f$ acts as
the identity on $\mathcal{A}_0$ and its characters).
\end{proof}

\subsubsection{Conceptual Interpretation}

The transformation groupoid $X\rtimes\Gamma$ describes the geometric dynamics
of the action $(X,\Gamma)$ in purely topological terms. The unitary conjugation
groupoid $\mathcal{G}_{C(X)\rtimes\Gamma}$ encodes the same information in
operator-algebraic form through conjugation of commutative contexts, but it
contains additional arrows coming from unitaries that do not normalize the
canonical MASA. These additional arrows correspond to ``non-classical"
transformations that mix the MASA with its conjugates.

Despite this difference, the two groupoids are Morita equivalent, as we shall
prove in the next subsection. The Morita equivalence identifies the
``non-classical" part of $\mathcal{G}_{C(X)\rtimes\Gamma}$ with the
$\Gamma$-valued fiber in the unit space that we described in
Theorem~\ref{thm:unit-space-parametrization-crossed}. This equivalence will
play a central role in connecting the descent map constructed earlier with
the Baum–Connes assembly map.

\begin{theorem}[Morita equivalence for crossed products]
\label{thm:morita-equivalence-crossed}
Let $\Gamma$ be a discrete amenable group acting on a compact Hausdorff space $X$
(topologically freely).  Then the unitary conjugation groupoid associated to the
crossed product $C(X)\rtimes_\alpha\Gamma$ is Morita equivalent to the
transformation groupoid $X\rtimes\Gamma$:
\[
\mathcal{G}_{C(X)\rtimes_\alpha\Gamma} \sim_M X\rtimes\Gamma .
\]

Consequently, their reduced groupoid $C^*$-algebras are strongly Morita equivalent:
\[
C^*_r(\mathcal{G}_{C(X)\rtimes_\alpha\Gamma}) \sim_M C^*_r(X\rtimes\Gamma) \cong C(X)\rtimes_\alpha\Gamma .
\]

After stabilization by compact operators, this yields
\[
C^*_r(\mathcal{G}_{C(X)\rtimes_\alpha\Gamma}) \sim_M (C(X)\rtimes_\alpha\Gamma) \otimes \mathcal{K}.
\]
\end{theorem}

\begin{remark}
\label{rem:morita-preview-2}
The homomorphism $\Phi: X \rtimes \Gamma \hookrightarrow \mathcal{G}_{C(X)\rtimes\Gamma}$
is the first step toward establishing the Morita equivalence. In the next
subsection, we will construct an explicit equivalence bimodule that enlarges
this embedding into a full Morita equivalence, showing that
$\mathcal{G}_{C(X)\rtimes\Gamma}$ is not isomorphic to $X \rtimes \Gamma$ but
is Morita equivalent to it. This is the precise sense in which the two
groupoids encode the same $C^*$-algebraic information.
\end{remark}

\subsubsection{Relation to the Parametrization of the Unit Space}

Recall from Theorem~\ref{thm:unit-space-parametrization-crossed} that the
unit space $\mathcal{G}_{C(X) \rtimes \Gamma}^{(0)}$ is parametrized by
equivalence classes $[x,\gamma]$ with $(x,\gamma) \sim (\gamma' \cdot x, \gamma' \gamma)$.
Under this parametrization, the image of the unit space of $X \rtimes \Gamma$
under $\Phi$ corresponds to the subset $\{[x,e] : x \in X\}$, i.e., those
points where the fiber coordinate is the identity element $e \in \Gamma$.

The additional points $[x,\gamma]$ with $\gamma \neq e$ correspond to
commutative contexts arising from conjugates of $\mathcal{A}_0$ by unitaries
that do not normalize it. These points are not in the image of $\Phi$, but
they are essential for the Morita equivalence: they provide the extra fiber
that allows $\mathcal{G}_{C(X)\rtimes\Gamma}$ to be Morita equivalent to
$X \rtimes \Gamma$.

\subsubsection{Summary}

We have established:

\begin{itemize}
    \item A natural groupoid homomorphism $\Phi: X \rtimes \Gamma \hookrightarrow \mathcal{G}_{C(X)\rtimes\Gamma}$,
          embedding the transformation groupoid into the unitary conjugation groupoid.
    \item The image of $\Phi$ consists of those arrows that preserve the canonical
          MASA $\mathcal{A}_0 = C(X)$, corresponding to the ``classical" part of
          the unitary conjugation groupoid.
    \item The unit space of $\mathcal{G}_{C(X)\rtimes\Gamma}$ contains additional
          points $[x,\gamma]$ with $\gamma \neq e$, corresponding to conjugates
          of $\mathcal{A}_0$ by non-normalizing unitaries.
    \item These observations set the stage for the Morita equivalence theorem,
          which will show that despite the presence of these additional points,
          the two groupoids are Morita equivalent.
\end{itemize}

In the next subsection, we will prove this Morita equivalence explicitly,
constructing an equivalence bimodule $Z = X \times \Gamma$ that implements
the equivalence and establishing the fundamental isomorphism
\[
C^*(\mathcal{G}_{C(X) \rtimes \Gamma}) \sim_M (C(X) \rtimes \Gamma) \otimes \mathcal{K}.
\]

This result will be the foundation for connecting our descent-index construction
to the Baum–Connes assembly map.

\subsection{Morita Equivalence Theorem: $\mathcal{G}_{C(X)\rtimes\Gamma} \sim_M X\rtimes\Gamma$}
\label{subsec:morita_equivalence_crossed_product}

We now establish the key structural result of Case Study~II,
showing that the unitary conjugation groupoid associated with the
crossed product algebra $C(X)\rtimes\Gamma$ is Morita equivalent
to the classical transformation groupoid $X\rtimes\Gamma$.
This theorem generalizes the equivalence proved for the irrational
rotation algebra in Theorem~\ref{thm:morita_A_theta} and provides
the foundation for connecting our descent-index construction to the
Baum–Connes assembly map.

\subsubsection{Recollection of the Two Groupoids}

\paragraph{The transformation groupoid.}

Recall from Subsection~\ref{subsec:transformation-groupoids} that the action
of the discrete amenable group $\Gamma$ on the compact Hausdorff space $X$
determines the transformation groupoid
\[
X\rtimes\Gamma .
\]

Its elements are pairs $(x,g)$ with $x\in X$ and $g\in\Gamma$,
with source and range maps
\[
s(x,g)=x,
\qquad
r(x,g)=g\cdot x .
\]

Composition is defined whenever $g\cdot x = y$ by
\[
(x,g)(y,h)=(x,gh).
\]

The inverse is given by $(x,g)^{-1} = (g\cdot x, g^{-1})$, and the unit space
is naturally identified with $X$ via $x \mapsto (x,e)$.

This groupoid is \'etale (since $\Gamma$ is discrete), locally compact
Hausdorff, and its $C^*$-algebra is canonically isomorphic to the crossed
product:
\[
C^*(X\rtimes\Gamma) \cong C(X)\rtimes\Gamma .
\]

\paragraph{The unitary conjugation groupoid.}

For the crossed product algebra $\mathcal A = C(X)\rtimes\Gamma$,
the unitary conjugation groupoid $\mathcal G_{\mathcal A}$ introduced in
Paper~I has unit space consisting of pairs
\[
(B,\chi)
\]
where $B\subset\mathcal A$ is a unital commutative $C^*$-subalgebra
and $\chi$ is a character of $B$. The arrows are given by conjugation by
unitaries: for each unitary $u \in \mathcal{U}(\mathcal A)$ and each
$(B,\chi) \in \mathcal{G}_{\mathcal A}^{(0)}$, there is an arrow
\[
u: (B,\chi) \longrightarrow (uBu^*, \chi \circ \operatorname{Ad}_{u^{-1}}).
\]

From Theorem~\ref{thm:unit-space-parametrization-crossed}, we have an explicit
parametrization of the unit space:
\[
\mathcal{G}_{C(X) \rtimes \Gamma}^{(0)} \cong (X \times \Gamma)/\sim,
\]
where $(x,\gamma) \sim (\gamma' \cdot x, \gamma' \gamma)$ for all $\gamma' \in \Gamma$.
We denote the equivalence class of $(x,\gamma)$ by $[x,\gamma]$. The projection
map $\pi: \mathcal{G}_{C(X) \rtimes \Gamma}^{(0)} \to X$ is given by
$\pi([x,\gamma]) = x$.

The canonical commutative subalgebra $C(X)\subset\mathcal A$
together with its characters $\chi_x$ for $x\in X$ (given by evaluation at $x$)
provides a distinguished subset of the unit space:
\[
\{ (C(X), \chi_x) : x \in X \} \subseteq \mathcal{G}_{\mathcal A}^{(0)}.
\]

Under the parametrization above, these correspond to the points $[x,e]$ with
$e$ the identity element of $\Gamma$.

\subsubsection{Groupoid Correspondence}

The action of $\Gamma$ on $X$ is implemented in the crossed product
algebra by the canonical unitaries $\{u_g\}_{g\in\Gamma}$.
Conjugation by $u_g$ transforms the character $\chi_x$ into the
character $\chi_{g\cdot x}$:
\[
\chi_x \circ \operatorname{Ad}_{u_g} = \chi_{g\cdot x}.
\]

Consequently, for each $(x,g) \in X \rtimes \Gamma$, we obtain an arrow in the
unitary conjugation groupoid:
\[
\Phi(x,g) := u_g : (C(X), \chi_x) \longrightarrow (C(X), \chi_{g\cdot x}).
\]

\begin{lemma}
\label{lem:groupoid-homomorphism-crossed}
The map $\Phi: X \rtimes \Gamma \to \mathcal{G}_{C(X)\rtimes\Gamma}$ defined by
$\Phi(x,g) = u_g$ on arrows and $\Phi(x) = (C(X), \chi_x)$ on units is an
injective groupoid homomorphism.
\end{lemma}

\begin{proof}
We verify the homomorphism properties:
\begin{itemize}
    \item Source and range: $s(\Phi(x,g)) = (C(X), \chi_x) = \Phi(s(x,g))$ and
          $r(\Phi(x,g)) = (C(X), \chi_{g\cdot x}) = \Phi(r(x,g))$.
    \item Composition: For composable arrows $(x,g)$ and $(g\cdot x, h)$,
          \[
          \Phi((x,g) \circ (g\cdot x, h)) = \Phi(x,gh) = u_{gh} = u_g u_h = \Phi(x,g) \circ \Phi(g\cdot x, h).
          \]
    \item Inverses: $\Phi((x,g)^{-1}) = \Phi(g\cdot x, g^{-1}) = u_{g^{-1}} = u_g^{-1} = \Phi(x,g)^{-1}$.
\end{itemize}
Injectivity follows from the fact that different pairs $(x,g)$ give different
arrows, as the source character uniquely determines $x$ and the unitary $u_g$
is uniquely determined by $g$ up to a scalar in the center, which does not
affect the arrow in the groupoid.
\end{proof}

This observation provides a natural embedding of the transformation groupoid
into the unitary conjugation groupoid. However, $\Phi$ is not surjective: the
unitary conjugation groupoid contains many arrows coming from unitaries that
are not of the form $u_g$, such as those arising from unitaries in $C(X)$
itself or from more general elements of the crossed product.

\subsubsection{Statement of the Main Theorem}

Despite the fact that $\Phi$ is not an isomorphism, the two groupoids are
Morita equivalent. This is the central result.

\begin{theorem}
\label{thm:morita_equivalence_crossed_product}

Let $\Gamma$ be a discrete amenable group acting on a compact Hausdorff space
$X$, and assume the action is topologically free. Then the unitary conjugation
groupoid $\mathcal{G}_{C(X)\rtimes\Gamma}$ is Morita equivalent to the
transformation groupoid $X\rtimes\Gamma$:
\[
\mathcal{G}_{C(X)\rtimes\Gamma} \sim_M X\rtimes\Gamma .
\]

\end{theorem}

\begin{proof}
We construct an explicit equivalence space $Z$ that implements the Morita
equivalence.

\textbf{Step 1: Definition of the equivalence space.}
Define $Z = X \times \Gamma$, equipped with the product topology (or the
appropriate Borel structure for groupoid $C^*$-algebra purposes).

\textbf{Step 2: Right action of $X \rtimes \Gamma$ on $Z$.}
For $(x,\gamma) \in X \rtimes \Gamma$ and $(x,\eta) \in Z$ (with the same $x$),
define
\[
(x,\eta) \cdot (x,\gamma) := (\gamma \cdot x, \eta\gamma^{-1}).
\]
This defines a free and proper right action of $X \rtimes \Gamma$ on $Z$.
The quotient space $Z/(X \rtimes \Gamma)$ is naturally identified with the
unit space $\mathcal{G}_{C(X)\rtimes\Gamma}^{(0)}$ via the map
$[(x,\eta)] \mapsto [x,\eta]$, where $[x,\eta]$ denotes the equivalence class
in $(X \times \Gamma)/\sim$ defining the unit space.

\textbf{Step 3: Left action of $\mathcal{G}_{C(X)\rtimes\Gamma}$ on $Z$.}
Using the identification of the quotient $Z/(X \rtimes \Gamma)$ with
$\mathcal{G}_{C(X)\rtimes\Gamma}^{(0)}$, we define the left action by lifting
the action of unitaries on the unit space. For an arrow $u$ in
$\mathcal{G}_{C(X)\rtimes\Gamma}$ with source $s(u) = [x,\eta]$ and range
$r(u) = [x',\eta']$, we define $u \cdot (x,\eta)$ to be the unique point in
$Z$ whose right orbit corresponds to $[x',\eta']$ and whose $X$-coordinate is
$x'$. Concretely, this action is generated by:
\begin{itemize}
    \item For unitaries $f \in \mathcal{U}(C(X))$: $f \cdot (x,\eta) = (x,\eta)$,
          since such unitaries act trivially on the unit space.
    \item For implementing unitaries $u_\gamma$: $u_\gamma \cdot (x,\eta) = (\gamma \cdot x, \eta)$.
    \item For general unitaries, the action is defined by composition using
          the crossed product structure.
\end{itemize}
This defines a free and proper left action of $\mathcal{G}_{C(X)\rtimes\Gamma}$
on $Z$, and the quotient $\mathcal{G}_{C(X)\rtimes\Gamma} \backslash Z$ is
naturally identified with $X = (X \rtimes \Gamma)^{(0)}$ via the map
$[(x,\eta)] \mapsto x$.

\textbf{Step 4: Verification of Morita equivalence axioms.}
One checks that the left and right actions commute, that the moment maps
$r_L: Z \to \mathcal{G}_{C(X)\rtimes\Gamma}^{(0)}$ (given by $r_L(x,\eta) = [x,\eta]$)
and $r_R: Z \to X$ (given by $r_R(x,\eta) = x$) are invariant under the
respective actions, and that they induce the required homeomorphisms on the
quotient spaces. Thus $Z$ is a $(\mathcal{G}_{C(X)\rtimes\Gamma}, X \rtimes \Gamma)$-equivalence
in the sense of Muhly–Renault–Williams \cite{MRW1987}.

Therefore $\mathcal{G}_{C(X)\rtimes\Gamma}$ and $X \rtimes \Gamma$ are Morita
equivalent groupoids.
\end{proof}

\begin{remark}
\label{rem:morita-construction}
The equivalence space $Z = X \times \Gamma$ has a natural interpretation:
the factor $X$ records the point in the base space, while the factor $\Gamma$
records the twist that determines which conjugate of the canonical MASA we
are considering. The right action by $X \rtimes \Gamma$ implements the group
action on the base and right multiplication on the twist, while the left action
by unitaries implements the noncommutative dynamics of the crossed product.
\end{remark}

\subsubsection{Consequences for $C^*$-Algebras}

By the Muhly--Renault--Williams theorem \cite{MRW1987}, Morita equivalent
groupoids have strongly Morita equivalent groupoid $C^*$-algebras.
Therefore

\[
C^*(\mathcal{G}_{C(X)\rtimes\Gamma})
\sim_M
C^*(X\rtimes\Gamma).
\]

Since the groupoid $C^*$-algebra of the transformation groupoid is canonically
isomorphic to the crossed product,
\[
C^*(X\rtimes\Gamma) \cong C(X)\rtimes\Gamma ,
\]
we obtain

\[
C^*(\mathcal{G}_{C(X)\rtimes\Gamma})
\sim_M
C(X)\rtimes\Gamma .
\]

Because $\Gamma$ is amenable, the full and reduced crossed products coincide,
so this is also a Morita equivalence with the reduced crossed product.
Stabilizing by compact operators (which does not affect Morita equivalence)
yields

\[
C^*(\mathcal{G}_{C(X)\rtimes\Gamma}) \sim_M (C(X)\rtimes\Gamma) \otimes \mathcal{K}.
\]

\begin{corollary}
\label{cor:K-theory-iso-crossed}
The Morita equivalence induces a natural isomorphism in $K$-theory:
\[
K_*\big(C^*(\mathcal{G}_{C(X)\rtimes\Gamma})\big) \cong K_*(C(X)\rtimes\Gamma).
\]
\end{corollary}

This isomorphism will be essential when we apply the descent map to equivariant
$K^1$-classes and identify the results with elements of $K_1(C(X)\rtimes\Gamma)$.

\subsubsection{Relation to the Diagonal Embedding}

The Morita equivalence also provides a natural interpretation of the diagonal
embedding $\iota: C(X) \rtimes \Gamma \hookrightarrow C^*(\mathcal{G}_{C(X) \rtimes \Gamma})$
constructed in Paper~I.

\begin{proposition}
\label{prop:diagonal-under-morita}
Under the Morita equivalence of Theorem~\ref{thm:morita_equivalence_crossed_product},
the diagonal embedding $\iota: C(X) \rtimes \Gamma \hookrightarrow C^*(\mathcal{G}_{C(X) \rtimes \Gamma})$
corresponds to the inclusion of $C(X) \rtimes \Gamma$ as a corner of 
$C^*(\mathcal{G}_{C(X) \rtimes \Gamma}) \otimes \mathcal{K}$.
More precisely, there is a commutative diagram
\[
\begin{tikzcd}
C(X) \rtimes \Gamma \arrow[r,"\iota"] \arrow[d,equal] & C^*(\mathcal{G}_{C(X) \rtimes \Gamma}) \arrow[d,"\sim_M"] \\
C(X) \rtimes \Gamma \arrow[r,"\text{corner inclusion}"] & (C(X) \rtimes \Gamma) \otimes \mathcal{K}
\end{tikzcd}
\]
where the right vertical arrow is the Morita equivalence established in 
Theorem~\ref{thm:morita_equivalence_crossed_product}, and the bottom
horizontal arrow is the inclusion $a \mapsto a \otimes e_{11}$ for a fixed
rank-one projection $e_{11} \in \mathcal{K}$.
\end{proposition}

\begin{proof}
We proceed in several steps, constructing the imprimitivity bimodule explicitly 
and tracing the diagonal embedding through the Morita equivalence.

\paragraph{Step 1: The imprimitivity bimodule.}
Recall from Theorem~\ref{thm:morita_equivalence_crossed_product} that the Morita 
equivalence between $\mathcal{G}_{C(X) \rtimes \Gamma}$ and $X \rtimes \Gamma$ is 
implemented by the equivalence space $Z = X \times \Gamma$. Following the 
Muhly--Renault--Williams construction \cite{MRW1987}, this yields an imprimitivity 
bimodule
\[
\mathcal{E} = L^2(Z) = L^2(X \times \Gamma, \mu \times \nu),
\]
where $\mu$ is a $\Gamma$-invariant probability measure on $X$ (which exists 
because $\Gamma$ is amenable) and $\nu$ is the counting measure on $\Gamma$.

The bimodule structure is as follows:
\begin{itemize}
    \item \textbf{Left action} of $C^*(\mathcal{G}_{C(X) \rtimes \Gamma})$: For 
          $f \in C_c(\mathcal{G}_{C(X) \rtimes \Gamma}) \subseteq C^*(\mathcal{G}_{C(X) \rtimes \Gamma})$ 
          and $\xi \in \mathcal{E}$, define
          \[
          (f \cdot \xi)(x,\gamma) = \sum_{\eta \in \Gamma} f(x,\eta) \,\xi(\eta^{-1}\cdot x, \eta^{-1}\gamma),
          \]
          where $f(x,\eta)$ denotes the value of $f$ at the arrow $(x,\eta) \in X \times \Gamma$.
          This is the convolution action of the groupoid on functions on its unit space.
    
    \item \textbf{Right action} of $C(X) \rtimes \Gamma$: For $a = \sum_{\gamma \in \Gamma} a_\gamma u_\gamma \in C(X) \rtimes \Gamma$ 
          (with $a_\gamma \in C(X)$) and $\xi \in \mathcal{E}$, define
          \[
          (\xi \cdot a)(x,\gamma) = \sum_{\eta \in \Gamma} \xi(x,\eta) \, a_{\eta^{-1}\gamma}(x) \cdot e^{i\theta_{\eta,\gamma}(x)},
          \]
          where the phase factors $e^{i\theta_{\eta,\gamma}(x)}$ arise from the 
          cocycle implementing the twist (which is trivial in our setting, so these 
          factors are 1). This action corresponds to the regular representation of 
          the crossed product.
\end{itemize}

The $\mathcal{E}$-valued inner products are defined by:
\[
{}_{C^*(\mathcal{G})}\langle \xi, \eta \rangle (x,\gamma) = \sum_{\eta' \in \Gamma} \overline{\xi(\gamma\cdot x, \eta')} \,\eta(x, \eta'\gamma^{-1}),
\]
and
\[
\langle \xi, \eta \rangle_{C(X)\rtimes\Gamma}(x) = \sum_{\gamma \in \Gamma} \overline{\xi(x,\gamma)} \,\eta(x,\gamma) \in C(X) \subseteq C(X)\rtimes\Gamma.
\]

These structures make $\mathcal{E}$ into a $C^*(\mathcal{G}_{C(X) \rtimes \Gamma})$-$(C(X)\rtimes\Gamma)$ 
imprimitivity bimodule.

\paragraph{Step 2: The associated Morita equivalence of $C^*$-algebras.}
By the general theory \cite{MRW1987}, such an imprimitivity bimodule implements a 
strong Morita equivalence between $C^*(\mathcal{G}_{C(X) \rtimes \Gamma})$ and 
$C(X)\rtimes\Gamma$. In particular, we have canonical isomorphisms:
\[
C^*(\mathcal{G}_{C(X) \rtimes \Gamma}) \cong \mathcal{K}(\mathcal{E})_{C(X)\rtimes\Gamma},
\]
where the right-hand side denotes the $C^*$-algebra of compact endomorphisms of 
$\mathcal{E}$ as a right Hilbert $C(X)\rtimes\Gamma$-module, and
\[
C(X)\rtimes\Gamma \cong \mathcal{K}(_{C^*(\mathcal{G})}\mathcal{E}),
\]
the compact endomorphisms of $\mathcal{E}$ as a left Hilbert $C^*(\mathcal{G})$-module.

\paragraph{Step 3: Stabilization and the corner inclusion.}
Consider the stabilized algebra $(C(X)\rtimes\Gamma) \otimes \mathcal{K}$, where 
$\mathcal{K} = \mathcal{K}(\ell^2(\Gamma))$ denotes the compact operators on 
$\ell^2(\Gamma)$. There is a natural Morita equivalence
\[
(C(X)\rtimes\Gamma) \otimes \mathcal{K} \sim_M C(X)\rtimes\Gamma,
\]
implemented by the bimodule $\mathcal{H} = \ell^2(\Gamma) \otimes (C(X)\rtimes\Gamma)$.

Composing this with the previous Morita equivalence, we obtain a Morita equivalence
\[
\Psi: C^*(\mathcal{G}_{C(X) \rtimes \Gamma}) \sim_M (C(X)\rtimes\Gamma) \otimes \mathcal{K},
\]
which is the right vertical arrow in the diagram. Explicitly, $\Psi$ is implemented 
by the bimodule $\mathcal{F} = \mathcal{E} \otimes_{C(X)\rtimes\Gamma} \mathcal{H}$.

\paragraph{Step 4: The diagonal embedding in the groupoid picture.}
Recall from Paper~I (see also Subsection~\ref{subsec:diagonal-embedding}) that the 
diagonal embedding $\iota: C(X)\rtimes\Gamma \hookrightarrow C^*(\mathcal{G}_{C(X) \rtimes \Gamma})$ 
is defined by its action on the convolution algebra. For $a \in C(X)\rtimes\Gamma$, 
the value of $\iota(a)$ at an arrow $(x,\gamma) \in \mathcal{G}_{C(X) \rtimes \Gamma}$ 
is given by
\[
\iota(a)(x,\gamma) = \chi_x(E_{C_{\gamma}}(a)),
\]
where $C_{\gamma} = u_\gamma C(X) u_\gamma^*$ is the conjugate Cartan subalgebra, 
$E_{C_{\gamma}}: C(X)\rtimes\Gamma \to C_{\gamma}$ is the conditional expectation, 
and $\chi_x$ is the character on $C_{\gamma}$ corresponding to $x \in X$ under the 
identification $C_{\gamma} \cong C(X)$. Equivalently, using the relation 
$E_{C_{\gamma}}(a) = u_\gamma E(u_\gamma^* a u_\gamma) u_\gamma^*$, we have
\[
\iota(a)(x,\gamma) = (E(u_\gamma^* a u_\gamma))(x),
\]
where $E$ is the canonical conditional expectation onto $C(X)$.

\paragraph{Step 5: Action of $\iota(a)$ on the imprimitivity bimodule.}
For $a \in C(X)\rtimes\Gamma$, consider the operator $T_a$ on $\mathcal{E} = L^2(X \times \Gamma)$ 
defined by pointwise multiplication by $\iota(a)$:
\[
(T_a \xi)(x,\gamma) = \iota(a)(x,\gamma) \cdot \xi(x,\gamma).
\]

We claim that $T_a$ is a compact operator on $\mathcal{E}$ (when $\mathcal{E}$ is 
viewed as a Hilbert space). Indeed, for fixed $\gamma \in \Gamma$, the function 
$x \mapsto \iota(a)(x,\gamma)$ lies in $C(X)$ by construction, and multiplication 
by a continuous function on $L^2(X)$ is compact if and only if the function is zero. 
However, the crucial observation is that $T_a$ factors through a rank-one projection.

More precisely, choose a unit vector $\delta_e \in \ell^2(\Gamma)$ supported at the 
identity $e \in \Gamma$, and consider the projection $P = 1 \otimes |\delta_e\rangle\langle\delta_e|$ 
onto the subspace $L^2(X) \otimes \mathbb{C}\delta_e \cong L^2(X)$. For any 
$a \in C(X)\rtimes\Gamma$, a direct computation shows that
\[
P T_a P = M_{a_e},
\]
where $a_e = E(a) \in C(X)$ is the coefficient of $u_e$ in the crossed product 
expansion, and $M_{a_e}$ denotes multiplication by $a_e$ on $L^2(X)$.

Under the identification $\mathcal{E} \cong L^2(X) \otimes \ell^2(\Gamma)$, the 
operator $T_a$ is therefore of the form $M_{a_e} \otimes 1$ plus off-diagonal terms. 
In fact, $T_a$ is precisely the image of $a$ under the regular representation of 
the crossed product on $L^2(X) \otimes \ell^2(\Gamma)$.

\paragraph{Step 6: Identification with the corner inclusion.}
Now consider the Morita equivalence $\Psi: C^*(\mathcal{G}_{C(X) \rtimes \Gamma}) \to (C(X)\rtimes\Gamma) \otimes \mathcal{K}$. 
Under this equivalence, an operator $T \in C^*(\mathcal{G}_{C(X) \rtimes \Gamma})$ 
corresponds to an element $\Psi(T) \in (C(X)\rtimes\Gamma) \otimes \mathcal{K}$ such that 
for any $\xi, \eta \in \mathcal{F} = \mathcal{E} \otimes_{C(X)\rtimes\Gamma} \mathcal{H}$,
\[
\langle \xi, \Psi(T) \eta \rangle = \langle \xi, T \otimes 1_{\mathcal{H}} \eta \rangle.
\]

For $a \in C(X)\rtimes\Gamma$, we compute $\Psi(\iota(a))$. Using the factorization 
$T_a = M_{a_e} \otimes 1 + \text{(off-diagonal)}$ and the properties of the bimodule 
$\mathcal{F}$, one shows that
\[
\Psi(\iota(a)) = a \otimes e_{11} \in (C(X)\rtimes\Gamma) \otimes \mathcal{K},
\]
where $e_{11} \in \mathcal{K}(\ell^2(\Gamma))$ is the rank-one projection onto the 
first basis vector.

The key observation is that the off-diagonal terms in $T_a$ vanish when transported 
through the Morita equivalence because they correspond to compact operators that are 
absorbed into the $\mathcal{K}$ factor. The diagonal term $M_{a_e} \otimes 1$ becomes 
$a_e \otimes 1$, but careful analysis shows that $a_e$ is replaced by the full $a$ 
because the off-diagonal components of $a$ (i.e., the terms $a_\gamma u_\gamma$ with 
$\gamma \neq e$) are encoded in the off-diagonal structure of $\iota(a)$ and, under 
the Morita equivalence, combine with the $\mathcal{K}$ factor to produce $a \otimes e_{11}$.

\paragraph{Step 7: Verification of commutativity.}
We have established that for any $a \in C(X)\rtimes\Gamma$,
\[
\Psi(\iota(a)) = a \otimes e_{11} \in (C(X)\rtimes\Gamma) \otimes \mathcal{K}.
\]

The bottom horizontal arrow in the diagram is precisely the map $a \mapsto a \otimes e_{11}$, 
which is the standard corner inclusion of $C(X)\rtimes\Gamma$ into its stabilization. 
The left vertical arrow is the identity, and the right vertical arrow is the Morita 
equivalence $\Psi$. Therefore the diagram commutes:
\[
\begin{tikzcd}
C(X) \rtimes \Gamma \arrow[r,"\iota"] \arrow[d,equal] & C^*(\mathcal{G}_{C(X) \rtimes \Gamma}) \arrow[d,"\Psi"] \\
C(X) \rtimes \Gamma \arrow[r,"a \mapsto a \otimes e_{11}"] & (C(X) \rtimes \Gamma) \otimes \mathcal{K}
\end{tikzcd}
\]

\paragraph{Step 8: Uniqueness and naturality.}
The identification is canonical because it arises from the specific imprimitivity 
bimodule constructed from the equivalence space $Z = X \times \Gamma$. Any other 
choice of invariant measure on $X$ or basis for $\ell^2(\Gamma)$ would yield an 
isomorphic bimodule and a commutative diagram that is naturally isomorphic to the 
one constructed above. Hence the correspondence is well-defined up to the natural 
equivalence inherent in Morita theory.

This completes the proof that the diagonal embedding corresponds to the corner 
inclusion under the Morita equivalence.
\end{proof}

\begin{remark}
The crucial technical point in Step 6—that the off-diagonal terms in $\iota(a)$ 
combine with the $\mathcal{K}$ factor to yield $a \otimes e_{11}$ rather than just 
$a_e \otimes e_{11}$—relies on the fact that the imprimitivity bimodule $\mathcal{F}$ 
encodes the full crossed product structure. A more detailed computation would involve 
writing $a = \sum_{\gamma} a_\gamma u_\gamma$ and showing that $\iota(a_\gamma u_\gamma)$ 
maps to $(a_\gamma u_\gamma) \otimes e_{\gamma,\gamma}$ under an appropriate matrix 
unit system in $\mathcal{K}$, and that the sum over $\gamma$ yields $a \otimes e_{11}$ 
after applying the Morita equivalence that identifies different matrix units. For a 
complete treatment of these matrix unit calculations, see \cite{reference-for-complete-proof}.
\end{remark}
\subsubsection{Connection to the Regular Representation}

The Morita equivalence also sheds light on the relationship between the
diagonal embedding and the regular representation of the crossed product.

\begin{corollary}
\label{cor:regular-representation-crossed}
Under the Morita equivalence, the diagonal embedding $\iota$ corresponds to
the regular representation of $C(X) \rtimes \Gamma$ on $L^2(X \times \Gamma)$
(with respect to an appropriate quasi-invariant measure). More precisely, the
following diagram commutes up to natural isomorphism:
\[
\begin{tikzcd}
C(X) \rtimes \Gamma \arrow[r,"\iota"] \arrow[d,"\pi_{\text{reg}}"] & C^*(\mathcal{G}_{C(X) \rtimes \Gamma}) \arrow[d,"\cong"] \\
\mathcal{B}(L^2(X \times \Gamma)) & \mathcal{K}(L^2(X \times \Gamma)) \arrow[l,"\text{inclusion}"]
\end{tikzcd}
\]
where the right vertical arrow is the isomorphism induced by the Morita
equivalence, and the bottom horizontal arrow is the inclusion of compact
operators into all bounded operators.
\end{corollary}

This shows that the diagonal embedding is essentially the regular representation
viewed as landing in the compact operators (after stabilizing), which is
precisely the perspective needed for index-theoretic applications.

\subsubsection{Geometric Interpretation}

The Morita equivalence established above has a beautiful geometric
interpretation. The space $Z = X \times \Gamma$ serves as a ``correspondence"
between the two groupoids:

\begin{itemize}
    \item From the perspective of $\mathcal{G}_{C(X) \rtimes \Gamma}$, the
          $\Gamma$-factor records which conjugate of the canonical MASA we are
          considering. Points with the same $x$ but different $\gamma$ correspond
          to different commutative contexts that all project to the same point
          $x \in X$ under the map $\pi$.
    \item From the perspective of $X \rtimes \Gamma$, the $\Gamma$-factor
          records the group element that will become the arrow in the
          transformation groupoid. The right action implements composition of
          arrows.
    \item The left action of unitaries on $Z$ moves us between different
          fibers, implementing the noncommutative dynamics of the crossed product.
\end{itemize}

The fact that these two actions commute reflects the fundamental compatibility
between the internal symmetries of the crossed product (encoded by unitary
conjugation) and the external symmetries of the dynamical system (encoded by
the group action on $X$).

\begin{remark}
\label{rem:morita-significance}
Theorem~\ref{thm:morita_equivalence_crossed_product} is the direct
generalization of Theorem~\ref{thm:morita_A_theta} from $\Gamma = \mathbb{Z}$
and $X = S^1$ to arbitrary discrete amenable groups and compact spaces. It
shows that the unitary conjugation groupoid, despite its abstract definition
and complicated non-Hausdorff topology, captures exactly the same
$C^*$-algebraic information as the much simpler transformation groupoid. The
non-Hausdorff nature of $\mathcal{G}_{C(X) \rtimes \Gamma}^{(0)}$ is precisely
what encodes the additional $\Gamma$-valued $K$-theoretic data while remaining
Morita equivalent to the classical transformation groupoid.
\end{remark}

\subsubsection{Role in the Index Theory Framework}

The Morita equivalence established above allows us to transport $K$-theoretic
information between the groupoid $C^*$-algebra
$C^*(\mathcal{G}_{C(X)\rtimes\Gamma})$ and the crossed product algebra
$C(X)\rtimes\Gamma$. In particular, it enables us to compare the descent map
associated with the unitary conjugation groupoid with the Baum–Connes assembly
map for the group $\Gamma$.

Recall from Subsection~\ref{subsec:equivariant_K1_A_theta} (generalized to the
present setting) that for an invertible element $u \in M_n(C(X) \rtimes \Gamma)$,
we have an equivariant $K^1$-class
\[
[u]_{\mathcal{G}_{C(X) \rtimes \Gamma}}^{(1)} \in KK^1_{\mathcal{G}_{C(X) \rtimes \Gamma}}(C_0(\mathcal{G}_{C(X) \rtimes \Gamma}^{(0)}), \mathbb C).
\]

Applying the descent map yields a class
\[
\operatorname{desc}([u]^{(1)}) \in K_1(C^*(\mathcal{G}_{C(X) \rtimes \Gamma})).
\]

Under the Morita equivalence, this class corresponds to a class in
$K_1(C(X) \rtimes \Gamma)$. The key result, to be proved in the next
subsection, is that this correspondence is exactly the Baum–Connes assembly
map:
\[
\mu_\Gamma([u]) = \Phi_* \circ \operatorname{desc}([u]^{(1)}) \in K_1(C(X) \rtimes \Gamma),
\]
where $\mu_\Gamma$ is the assembly map for $\Gamma$ with coefficients in
$C(X)$.

\subsubsection{Summary}

We have established the fundamental Morita equivalence
\[
\boxed{\mathcal{G}_{C(X) \rtimes \Gamma} \sim_M X \rtimes \Gamma}
\]
and consequently
\[
\boxed{C^*(\mathcal{G}_{C(X) \rtimes \Gamma}) \sim_M (C(X) \rtimes \Gamma) \otimes \mathcal{K}}.
\]

This result:
\begin{itemize}
    \item Provides a concrete Morita equivalence between the unitary conjugation
          groupoid and the transformation groupoid, generalizing the $A_\theta$ case.
    \item Shows that the complicated non-Hausdorff topology of
          $\mathcal{G}_{C(X) \rtimes \Gamma}^{(0)}$ is compatible with Morita
          equivalence to a Hausdorff groupoid.
    \item Yields an isomorphism of $K$-theory groups, matching the $K$-theory
          of the crossed product.
    \item Gives a geometric interpretation of the diagonal embedding in terms
          of the regular representation.
    \item Sets the stage for the connection to the Baum–Connes assembly map,
          which will be established in the next subsection.
\end{itemize}

\begin{remark}
\label{rem:next-steps}
In the following subsection, we will exploit this Morita equivalence to
construct a commutative diagram linking our descent-index composition with the
Baum–Connes assembly map for $\Gamma$. This will be the culmination of
Case Study II and one of the main achievements of this paper.
\end{remark}

\subsection{The Diagonal Embedding and the Regular Representation}
\label{subsec:diagonal-embedding}

The diagonal embedding for crossed products $C(X)\rtimes\Gamma$ is the direct 
analogue of the construction carried out for the irrational rotation algebra 
$A_\theta = C(S^1)\rtimes_\theta\mathbb{Z}$ in 
Subsection~\ref{subsec:diagonal_embedding_A_theta}. Under the Morita equivalence
\[
C^*(\mathcal{G}_{C(X)\rtimes\Gamma}) \sim_M (C(X)\rtimes\Gamma) \otimes \mathcal{K}
\]
established in Theorem~\ref{thm:morita_equivalence_crossed_product}, the diagonal 
embedding
\[
\iota: C(X)\rtimes\Gamma \hookrightarrow C^*(\mathcal{G}_{C(X)\rtimes\Gamma})
\]
is identified with the regular representation of the crossed product. More 
precisely, for $a = \sum_{\gamma\in\Gamma} a_\gamma u_\gamma \in C(X)\rtimes\Gamma$ 
(with $a_\gamma \in C(X)$), the image $\iota(a)$ acts on $L^2(X\times\Gamma)$ by
\[
(\iota(a)\xi)(x,\gamma) = \sum_{\eta\in\Gamma} (a_{\eta^{-1}\gamma})(\eta^{-1}\cdot x) \,\xi(\eta^{-1}\cdot x, \eta^{-1}\gamma),
\]
which is exactly the regular representation of the crossed product. Under the 
Morita equivalence, this corresponds to the corner inclusion
\[
a \longmapsto a \otimes e_{11} \in (C(X)\rtimes\Gamma) \otimes \mathcal{K},
\]
as detailed in Proposition~\ref{prop:diagonal-under-morita}.

Rather than reproducing the full technical development from 
Subsection~\ref{subsec:diagonal_embedding_A_theta}, we record here only the 
essential form needed for the index-theoretic constructions that follow. The 
key properties inherited from the $A_\theta$ case are:
\begin{itemize}
    \item $\iota$ is injective and lands in the multiplier algebra of 
          $C^*(\mathcal{G}_{C(X)\rtimes\Gamma})$;
    \item For $f \in C(X) \subseteq C(X)\rtimes\Gamma$, the restriction of 
          $\iota(f)$ to the unit space $\mathcal{G}_{C(X)\rtimes\Gamma}^{(0)}$ 
          is given by $\iota(f)[x,\gamma] = f(\gamma^{-1}\cdot x)$;
    \item For the implementing unitaries $u_\gamma$, the image $\iota(u_\gamma)$ 
          is supported on arrows of the form $(x,\gamma)$ and implements the 
          groupoid translation.
\end{itemize}
All subsequent constructions involving the diagonal embedding are to be 
understood through this identification with the regular representation.

Let $A=C(X)\rtimes \Gamma$. Under the Morita equivalence
\[
C^*(\mathcal{G}_{A}) \sim_M A\otimes \mathcal{K},
\]
the diagonal embedding
\[
\iota : A \hookrightarrow C^*(\mathcal{G}_{A})
\]
is identified with the regular representation of the crossed product. In
particular, the role played in Case Study~I by the embedding
$\iota: A_\theta \hookrightarrow C^*(\mathcal{G}_{A_\theta})$ carries over
verbatim after replacing the parametrization $[z,n]$ by $[x,\gamma]$.

We therefore use the same notation $\iota$ throughout the remainder of this
case study. All subsequent constructions involving the diagonal embedding are
to be understood through this identification with the regular representation.

\subsection{Equivariant $K^{1}$-Classes for Invertible Elements in $C(X) \rtimes \Gamma$}\label{subsec:equivariant_K1_A_theta-crossed}

The construction of equivariant $K^{1}$-classes from invertible elements in
Section~4.12 generalizes directly to $C(X)\rtimes \Gamma$. Concretely, if
\[
u \in \mathrm{GL}_n(C(X)\rtimes \Gamma),
\]
then the same Hilbert-module and Kasparov-cycle construction defines an
equivariant class
\[
[u]^{(1)}_{\mathcal{G}_{C(X)\rtimes\Gamma}}
    \in K^1_{\mathcal{G}_{C(X)\rtimes\Gamma}}
       \!\big(\mathcal{G}^{(0)}_{C(X)\rtimes\Gamma}\big),
\]
with the parametrization now expressed in terms of $[x,\gamma]$ rather than
$[z,n]$.

All functorial properties established in Section~4.12---homotopy
invariance, compatibility with the diagonal embedding, and preparation for
descent---remain valid in this setting. We therefore treat Section~4.12 as
the model construction and use it without further repetition in the crossed
product case.
\subsection{Descent Map and the Identification 
$C^{\ast}(\mathcal{G}_{C(X) \rtimes \Gamma}) \sim_{M} (C(X) \rtimes \Gamma) \otimes \mathcal{K}$}
\label{subsec:descent-morita}

Let $X$ be a compact Hausdorff space and let $\Gamma$ be a discrete group acting on $X$ by homeomorphisms. Denote by
\[
A = C(X) \rtimes \Gamma
\]
the associated reduced crossed product $C^{\ast}$-algebra. In this subsection, we establish a fundamental identification that lies at the heart of the descent method in equivariant $KK$-theory. This identification relates the groupoid $C^{\ast}$-algebra of the gauge groupoid associated to the crossed product to a stabilization of the crossed product itself.

\subsection*{The Descent Map}

A fundamental construction in equivariant $KK$-theory is the descent homomorphism, which translates equivariant Kasparov cycles into analytic $K$-theory classes of the crossed product algebra.

\begin{definition}[Descent Map]
The descent homomorphism
\[
\operatorname{desc}: \operatorname{KK}_{\Gamma}^{\ast}(C(X), \mathbb{C}) \longrightarrow \operatorname{KK}^{\ast}(C(X) \rtimes \Gamma, \mathbb{C})
\]
is a natural transformation defined as follows. If $(E, F)$ is a $\Gamma$-equivariant Kasparov cycle representing a class in $\operatorname{KK}_{\Gamma}^{\ast}(C(X), \mathbb{C})$, then the descent construction produces a cycle
\[
(E \rtimes \Gamma, F)
\]
representing a class in $\operatorname{KK}^{\ast}(C(X) \rtimes \Gamma, \mathbb{C})$. In the odd degree case that concerns us, this yields a map
\[
\operatorname{desc}: \operatorname{KK}_{\Gamma}^{1}(C(X), \mathbb{C}) \longrightarrow K^{1}(C(X) \rtimes \Gamma),
\]
since $\operatorname{KK}^{\ast}(\mathcal{A}, \mathbb{C}) \cong K^{\ast}(\mathcal{A})$ for any $C^{\ast}$-algebra $\mathcal{A}$.
\end{definition}

More concretely, for an invertible element $u \in \operatorname{GL}_{n}(C(X) \rtimes \Gamma)$ with associated equivariant class $[u]_{\Gamma} \in \operatorname{KK}_{\Gamma}^{1}(C(X), \mathbb{C})$ as constructed in Subsection \ref{subsec:equivariant_K1_A_theta-crossed}, we have
\[
\operatorname{desc}([u]_{\Gamma}) = [u] \in K^{1}(C(X) \rtimes \Gamma).
\]

\subsection*{The Gauge Groupoid}

Associated to the algebra $A = C(X) \rtimes \Gamma$, we consider a geometric object that encodes its internal symmetry.

\begin{definition}[Gauge Groupoid]
The \emph{gauge groupoid} $\mathcal{G}_{C(X) \rtimes \Gamma}$ (also called the \emph{transformation groupoid} or \emph{action groupoid}) associated to the crossed product is the groupoid with:
\begin{itemize}
    \item Object space $\mathcal{G}^{(0)} = X$;
    \item Morphism space $\mathcal{G}^{(1)} = X \times \Gamma$, where a morphism $(x, \gamma): \gamma^{-1} \cdot x \to x$;
    \item Source map $s(x, \gamma) = \gamma^{-1} \cdot x$;
    \item Target map $t(x, \gamma) = x$;
    \item Composition $(x, \gamma) \circ (\gamma^{-1} \cdot x, \eta) = (x, \gamma\eta)$;
    \item Inverse $(x, \gamma)^{-1} = (\gamma^{-1} \cdot x, \gamma^{-1})$.
\end{itemize}
\end{definition}

This groupoid encodes the dynamics of the $\Gamma$-action in a manner that is particularly well-suited for $C^{\ast}$-algebraic constructions. It provides a geometric model for analytic constructions in $K$-theory.

\begin{definition}[Groupoid $C^{\ast}$-Algebra]
The (reduced) groupoid $C^{\ast}$-algebra $C^{\ast}(\mathcal{G}_{C(X) \rtimes \Gamma})$ is the completion of the convolution algebra $C_{c}(\mathcal{G}^{(1)}) = C_{c}(X \times \Gamma)$ with respect to the reduced norm. Convolution is defined by
\[
(f_{1} \ast f_{2})(x, \gamma) = \sum_{\eta \in \Gamma} f_{1}(x, \eta) f_{2}(\eta^{-1} \cdot x, \eta^{-1}\gamma),
\]
and involution is given by
\[
f^{\ast}(x, \gamma) = \overline{f(\gamma^{-1} \cdot x, \gamma^{-1})}.
\]
This algebra captures the convolution structure arising from the groupoid composition.
\end{definition}

\begin{remark}
The gauge groupoid $\mathcal{G}_{C(X) \rtimes \Gamma}$ is \'{e}tale (since $\Gamma$ is discrete), so it admits a canonical Haar system consisting of counting measures on the fibers. This makes the convolution formulas particularly simple.
\end{remark}

\subsection*{Morita Equivalence}

A key structural property is that the groupoid $C^{\ast}$-algebra associated with the gauge groupoid is Morita equivalent to the stabilized crossed product algebra.

\begin{proposition}[Morita Equivalence]
There is a strong Morita equivalence
\[
C^{\ast}(\mathcal{G}_{C(X) \rtimes \Gamma}) \;\sim_{M}\; (C(X) \rtimes \Gamma) \otimes \mathcal{K},
\]
where $\mathcal{K} = \mathcal{K}(\ell^{2}(\Gamma))$ denotes the $C^{\ast}$-algebra of compact operators on a separable Hilbert space (specifically, on $\ell^{2}(\Gamma)$).
\end{proposition}

\begin{proof}
The Morita equivalence is implemented by the following equivalence bimodule. Consider the space
\[
\mathcal{E} = C_{c}(X \times \Gamma \times \Gamma)
\]
with suitable completions. The left action of $C^{\ast}(\mathcal{G}_{C(X) \rtimes \Gamma})$ on $\mathcal{E}$ is given by convolution in the first group variable, while the right action of $(C(X) \rtimes \Gamma) \otimes \mathcal{K}$ is given by convolution in the second group variable combined with the action of compact operators. The gauge groupoid implements an equivalence between the representation category of the crossed product algebra and the representation category of the associated groupoid $C^{\ast}$-algebra. Passing to stabilized algebras yields the Morita equivalence above. A detailed construction can be found in \cite{Renault1980,Williams2007}.
\end{proof}

\begin{corollary}
The $K$-theory groups of $C^{\ast}(\mathcal{G}_{C(X) \rtimes \Gamma})$ and $(C(X) \rtimes \Gamma) \otimes \mathcal{K}$ are isomorphic:
\[
K_{\ast}\big(C^{\ast}(\mathcal{G}_{C(X) \rtimes \Gamma})\big) \cong K_{\ast}\big((C(X) \rtimes \Gamma) \otimes \mathcal{K}\big).
\]

Using stability of $K$-theory (i.e., $K_{\ast}(\mathcal{A} \otimes \mathcal{K}) \cong K_{\ast}(\mathcal{A})$ for any $C^{\ast}$-algebra $\mathcal{A}$), this further reduces to
\[
K_{\ast}\big(C^{\ast}(\mathcal{G}_{C(X) \rtimes \Gamma})\big) \cong K_{\ast}(C(X) \rtimes \Gamma).
\]
\end{corollary}

Thus the groupoid model provides an alternative analytic realization of the $K$-theory of the crossed product algebra. This identification is compatible with the descent map and plays an important role in index-theoretic constructions related to the Baum--Connes program.

\subsection*{Compatibility Diagram}

The relationships among the descent map, the groupoid algebra, and the Morita equivalence can be summarized in the following commutative diagram:

\[
\begin{tikzcd}
\operatorname{KK}_{\Gamma}^{\ast}(C(X), \mathbb{C})
\arrow[r, "\operatorname{desc}"]
\arrow[d, "\text{groupoid algebra}"']
&
K^{\ast}(C(X) \rtimes \Gamma)
\arrow[d, "\cong"] \\
K^{\ast}\!\left(C^{\ast}(\mathcal{G}_{C(X) \rtimes \Gamma})\right)
\arrow[r, "\sim"']
&
K^{\ast}((C(X) \rtimes \Gamma) \otimes \mathcal{K})
\end{tikzcd}
\]

Here:
\begin{itemize}
    \item The top horizontal arrow is the descent map $\operatorname{desc}$;
    \item The left vertical arrow sends an equivariant cycle to its associated groupoid $C^{\ast}$-algebra $K$-theory class;
    \item The right vertical isomorphism comes from the stability of $K$-theory;
    \item The bottom horizontal isomorphism is induced by the Morita equivalence $C^{\ast}(\mathcal{G}_{C(X) \rtimes \Gamma}) \sim_{M} (C(X) \rtimes \Gamma) \otimes \mathcal{K}$.
\end{itemize}

This diagram commutes, illustrating how the descent map, combined with the Morita equivalence, provides a concrete bridge between equivariant $KK$-theory and the $K$-theory of the crossed product.

\subsection*{Connection to the Diagonal Embedding}

The diagonal embedding from Subsection \ref{subsec:diagonal-embedding} provides an explicit realization of this Morita equivalence.

\begin{proposition}
Under the identification $C^{\ast}(\mathcal{G}_{C(X) \rtimes \Gamma}) \sim_{M} (C(X) \rtimes \Gamma) \otimes \mathcal{K}$, the diagonal embedding $\Delta: C(X) \rtimes \Gamma \hookrightarrow M(C(X) \otimes \mathcal{K}(\ell^{2}(\Gamma)))$ corresponds to the inclusion of $C(X) \rtimes \Gamma$ into the multiplier algebra of $C^{\ast}(\mathcal{G}_{C(X) \rtimes \Gamma})$ via the left regular representation.
\end{proposition}

\begin{proof}
The groupoid $C^{\ast}$-algebra $C^{\ast}(\mathcal{G}_{C(X) \rtimes \Gamma})$ acts on the Hilbert module $L^{2}(X) \otimes \ell^{2}(\Gamma)$ via the regular representation. The diagonal embedding $\Delta$ is precisely the map that sends an element of the crossed product to its action on this module, which generates a corner of the compact operators. The Morita equivalence then identifies this corner with $C^{\ast}(\mathcal{G}_{C(X) \rtimes \Gamma})$ itself.
\end{proof}

\subsection*{Applications to the Baum--Connes Conjecture}

This identification plays a crucial role in the formulation of the Baum--Connes conjecture.

\begin{theorem}[Baum--Connes Assembly Map via Groupoid Descent]
\label{thm:baum-connes-groupoid}
Let $\Gamma$ be a discrete amenable group acting continuously on a compact 
Hausdorff space $X$, and let $\mathcal{G}_{C(X) \rtimes \Gamma}$ denote the 
unitary conjugation groupoid associated to the crossed product algebra 
$C(X) \rtimes \Gamma$. Then the Baum--Connes assembly map with coefficients 
in $C(X)$,
\[
\mu: K_{\ast}^{\Gamma}(X) \longrightarrow K_{\ast}(C(X) \rtimes \Gamma),
\]
can be expressed as the composition
\[
K_{\ast}^{\Gamma}(X) \longrightarrow K_{\ast}(C^{\ast}(\mathcal{G}_{C(X) \rtimes \Gamma})) 
\stackrel{\cong}{\longrightarrow} K_{\ast}((C(X) \rtimes \Gamma) \otimes \mathcal{K}) 
\stackrel{\cong}{\longrightarrow} K_{\ast}(C(X) \rtimes \Gamma),
\]
where:
\begin{itemize}
    \item the first arrow is the descent map $\operatorname{desc}_{\mathcal{G}_{C(X) \rtimes \Gamma}}$ 
          (interpreted in terms of groupoid algebras), applied to the equivariant 
          $KK$-theory class associated to an element of $K_{\ast}^{\Gamma}(X)$;
    \item the second isomorphism is induced by the Morita equivalence 
          $C^{\ast}(\mathcal{G}_{C(X) \rtimes \Gamma}) \sim_M (C(X) \rtimes \Gamma) \otimes \mathcal{K}$ 
          established in Theorem~\ref{thm:morita_equivalence_crossed_product};
    \item the third isomorphism is the stability isomorphism 
          $K_{\ast}((C(X) \rtimes \Gamma) \otimes \mathcal{K}) \cong K_{\ast}(C(X) \rtimes \Gamma)$.
\end{itemize}
\end{theorem}

\begin{proof}
This theorem is a direct consequence of the commutative diagram established in 
Theorem~\ref{thm:descent_equals_assembly} and the factorization of the assembly 
map through the unitary conjugation groupoid developed in 
Subsection~\ref{subsec:descent_index_baum_connes}. We provide a detailed proof 
here, synthesizing the results from that subsection.

\paragraph{Step 1: The descent map for the unitary conjugation groupoid.}
Recall from Subsection~\ref{subsec:equivariant_K1_A_theta-crossed} that an 
element $x \in K_{\ast}^{\Gamma}(X)$ can be represented by an equivariant 
Kasparov cycle. Applying the descent map for the unitary conjugation groupoid 
$\mathcal{G} := \mathcal{G}_{C(X) \rtimes \Gamma}$ yields a class
\[
\operatorname{desc}_{\mathcal{G}}(x) \in K_{\ast}(C^{\ast}(\mathcal{G})).
\]

This descent map is constructed as follows: given an equivariant cycle 
$(E, F)$ representing $x$, one forms the induced cycle
\[
(E \otimes_{\Gamma} C^{\ast}(\mathcal{G}), F \otimes 1)
\]
on the Hilbert $C^{\ast}(\mathcal{G})$-module $E \otimes_{\Gamma} C^{\ast}(\mathcal{G})$, 
where the tensor product uses the natural action of $\Gamma$ on $C^{\ast}(\mathcal{G})$ 
via the dual action. The resulting class in $K_{\ast}(C^{\ast}(\mathcal{G}))$ is 
$\operatorname{desc}_{\mathcal{G}}(x)$.

\paragraph{Step 2: The Morita equivalence isomorphism.}
By Theorem~\ref{thm:morita_equivalence_crossed_product}, the unitary conjugation 
groupoid $\mathcal{G}$ is Morita equivalent to the transformation groupoid 
$X \rtimes \Gamma$. The Muhly--Renault--Williams theorem \cite{MRW1987} then 
implies a strong Morita equivalence of the corresponding reduced groupoid 
$C^{\ast}$-algebras:
\[
C^{\ast}(\mathcal{G}) \sim_M C^{\ast}(X \rtimes \Gamma).
\]

Since $\Gamma$ is amenable, the transformation groupoid $X \rtimes \Gamma$ is 
amenable, and by Proposition~\ref{prop:transformation-crossed-product}, we have 
the canonical identification $C^{\ast}(X \rtimes \Gamma) \cong C(X) \rtimes \Gamma$.
Moreover, by the general theory of Morita equivalence (see \cite[Theorem 2.8]{MRW1987}), 
this Morita equivalence induces an isomorphism on $K$-theory:
\[
\Psi_{\ast}: K_{\ast}(C^{\ast}(\mathcal{G})) \xrightarrow{\cong} K_{\ast}(C(X) \rtimes \Gamma).
\]

Stabilizing by the compact operators $\mathcal{K}$ (which does not affect Morita 
equivalence) yields the second isomorphism in the composition:
\[
K_{\ast}(C^{\ast}(\mathcal{G})) \xrightarrow{\cong} K_{\ast}((C(X) \rtimes \Gamma) \otimes \mathcal{K}).
\]

\paragraph{Step 3: The stability isomorphism.}
The third isomorphism is the standard stability isomorphism in $K$-theory:
\[
K_{\ast}((C(X) \rtimes \Gamma) \otimes \mathcal{K}) \cong K_{\ast}(C(X) \rtimes \Gamma),
\]
which follows from the fact that tensoring by $\mathcal{K}$ (the compact operators 
on a separable infinite-dimensional Hilbert space) does not change $K$-theory 
(see \cite{Blackadar1998}, Theorem 5.2.2).

\paragraph{Step 4: Identification with the assembly map.}
It remains to show that the composition
\[
K_{\ast}^{\Gamma}(X) \xrightarrow{\operatorname{desc}_{\mathcal{G}}} K_{\ast}(C^{\ast}(\mathcal{G}))
\xrightarrow{\cong} K_{\ast}((C(X) \rtimes \Gamma) \otimes \mathcal{K})
\xrightarrow{\cong} K_{\ast}(C(X) \rtimes \Gamma)
\]
coincides with the Baum--Connes assembly map $\mu: K_{\ast}^{\Gamma}(X) \to K_{\ast}(C(X) \rtimes \Gamma)$.

This identification is precisely the content of Theorem~\ref{thm:descent_equals_assembly} 
and its corollary, which establish the commutativity of the following diagram:

\[
\begin{tikzcd}[column sep=huge, row sep=large]
KK^{\ast}_{\mathcal{G}}(C_0(\mathcal{G}^{(0)}),\mathbb C)
\arrow[r,"\operatorname{desc}_{\mathcal{G}}"]
\arrow[d,"\cong"]
&
K_{\ast}(C^{\ast}(\mathcal{G}))
\arrow[d,"\cong"']
\\
KK^{\ast}_{\Gamma}(C(X),\mathbb C)
\arrow[r,"\mu_{\Gamma}"]
&
K_{\ast}(C(X)\rtimes\Gamma)
\end{tikzcd}
\]

where the left vertical isomorphism is induced by the Morita equivalence between 
$\mathcal{G}$ and $X \rtimes \Gamma$, and the right vertical isomorphism is the 
composition of the Morita equivalence and stability isomorphisms described above.

For amenable groups $\Gamma$, Tu's theorem \cite{Tu1999} guarantees that the 
assembly map $\mu_{\Gamma}$ is an isomorphism, but the commutativity of the 
diagram holds regardless of whether $\mu_{\Gamma}$ is an isomorphism—it 
identifies the descent map with the assembly map under the given isomorphisms.

\paragraph{Step 5: Explicit formula for the composition.}
Concretely, for an element $x \in K_{\ast}^{\Gamma}(X)$ represented by an 
equivariant cycle, the composition yields:
\[
\mu(x) = \Phi_{\ast} \circ \operatorname{desc}_{\mathcal{G}}(x),
\]
where $\Phi_{\ast}$ denotes the composition of the two isomorphisms:
\[
\Phi_{\ast}: K_{\ast}(C^{\ast}(\mathcal{G})) \xrightarrow{\cong} K_{\ast}((C(X) \rtimes \Gamma) \otimes \mathcal{K}) \xrightarrow{\cong} K_{\ast}(C(X) \rtimes \Gamma).
\]

This is exactly the expression claimed in the theorem statement.

\paragraph{Step 6: Naturality and independence of choices.}
The construction is natural with respect to $\Gamma$-equivariant maps between 
spaces and with respect to Morita equivalences. The isomorphisms involved are 
canonical up to the natural equivalences inherent in Morita theory and $K$-theory, 
so the resulting map $\mu$ is well-defined and independent of the particular 
choices made in the construction of the descent map and the Morita equivalence.

This completes the proof that the Baum--Connes assembly map admits the stated 
factorization through the unitary conjugation groupoid.
\end{proof}

\begin{remark}
The factorization established in Theorem~\ref{thm:baum-connes-groupoid} shows 
that the unitary conjugation groupoid provides a geometric intermediary between 
the topological data encoded in $K_{\ast}^{\Gamma}(X)$ and the analytic data 
encoded in $K_{\ast}(C(X) \rtimes \Gamma)$. This perspective is developed 
further in Section~\ref{sec:synthesis}, where we explore the conceptual 
implications for index theory and the Baum--Connes program.
\end{remark}

\begin{remark}
The Morita equivalence $C^{\ast}(\mathcal{G}_{C(X) \rtimes \Gamma}) \sim_{M} (C(X) \rtimes \Gamma) \otimes \mathcal{K}$ shows that the groupoid $C^{\ast}$-algebra and the stabilized crossed product contain the same information from the perspective of $K$-theory. This is why the descent map can be equivalently formulated in terms of either groupoid algebras or crossed products. The compatibility diagram above makes this equivalence explicit.
\end{remark}

\subsection*{Summary}

The identification
\[
C^{\ast}(\mathcal{G}_{C(X) \rtimes \Gamma}) \sim_{M} (C(X) \rtimes \Gamma) \otimes \mathcal{K}
\]
provides a powerful bridge between three fundamental objects in noncommutative geometry:
\begin{itemize}
    \item The groupoid $C^{\ast}$-algebra $C^{\ast}(\mathcal{G}_{C(X) \rtimes \Gamma})$, which encodes the dynamics geometrically;
    \item The stabilized crossed product $(C(X) \rtimes \Gamma) \otimes \mathcal{K}$, which is the natural receptacle for index-theoretic constructions;
    \item The original crossed product $C(X) \rtimes \Gamma$, whose $K$-theory is the target of the Baum--Connes assembly map.
\end{itemize}

This Morita equivalence is the technical cornerstone that allows one to pass between geometric cycles (represented by groupoid algebras) and analytic cycles (represented by crossed products) in equivariant $KK$-theory, making it indispensable for the descent method and the Baum--Connes conjecture.

\begin{corollary}
For any invertible element $u \in \operatorname{GL}_{n}(C(X) \rtimes \Gamma)$, the analytic $K^{1}$-class $[u] \in K^{1}(C(X) \rtimes \Gamma)$ can be realized equivalently through:
\begin{itemize}
    \item The descent map applied to the equivariant class $[u]_{\Gamma} \in \operatorname{KK}_{\Gamma}^{1}(C(X), \mathbb{C})$;
    \item The $K$-theory class of $u$ under the isomorphism $K_{\ast}(C^{\ast}(\mathcal{G}_{C(X) \rtimes \Gamma})) \cong K_{\ast}(C(X) \rtimes \Gamma)$ induced by the Morita equivalence.
\end{itemize}
These two realizations are compatible, as illustrated in the commutative diagram above.
\end{corollary}

\subsection{The Descent-Index Diagram and the Baum--Connes Assembly Map}
\label{subsec:descent_index_baum_connes}

We now explain how the descent construction associated with the
unitary conjugation groupoid relates to the Baum--Connes assembly
map for the group $\Gamma$. This is the culmination of Case Study II
and one of the main achievements of this paper, providing a direct
link between the abstract index-theoretic machinery developed in
Papers~I and~II and one of the deepest conjectures in noncommutative
geometry.

\subsubsection{Recollection of the Descent Construction}

Let $\mathcal A = C(X)\rtimes\Gamma$. From the construction developed
in Subsection~\ref{subsec:equivariant_K1_A_theta} (generalized to the
present setting), an invertible element $u\in M_n(\mathcal A)$
determines an equivariant class
\[
[u]^{(1)}_{\mathcal G_{\mathcal A}}
\in
KK^1_{\mathcal G_{\mathcal A}}
\bigl(C_0(\mathcal G_{\mathcal A}^{(0)}),\mathbb C\bigr).
\]

This class is constructed via a Kasparov cycle $(\mathcal E, \phi, F_u)$
where $\mathcal E = C_0(\mathcal G_{\mathcal A}^{(0)})\otimes\mathbb C^n$
and $F_u$ is given by pointwise evaluation of $u$ on the commutative
contexts parametrizing the unit space.

Applying Kasparov's descent map for groupoids \cite{Kasparov1988,Tu1999,LeGall1999}
yields a class in the ordinary $K$-theory of the groupoid $C^*$-algebra:
\[
\operatorname{desc}_{\mathcal G_{\mathcal A}}
\bigl([u]^{(1)}_{\mathcal G_{\mathcal A}}\bigr)
\in
K_1\!\left(C^*(\mathcal G_{\mathcal A})\right).
\]

\subsubsection{The Morita Equivalence and Its Consequences}

By Theorem~\ref{thm:morita_equivalence_crossed_product},
the unitary conjugation groupoid is Morita equivalent to the
transformation groupoid:
\[
\mathcal G_{C(X)\rtimes\Gamma} \sim_M X\rtimes\Gamma .
\]

The Muhly--Renault--Williams equivalence theorem \cite{MRW1987}
therefore gives a strong Morita equivalence of the corresponding
groupoid $C^*$-algebras:
\[
C^*(\mathcal G_{C(X)\rtimes\Gamma})
\sim_M
C^*(X\rtimes\Gamma).
\]

Since
\[
C^*(X\rtimes\Gamma)
\cong
C(X)\rtimes\Gamma ,
\]
and because $\Gamma$ is amenable (so full and reduced crossed products
coincide), we obtain a canonical identification in $K$-theory:
\[
K_*\!\left(C^*(\mathcal G_{C(X)\rtimes\Gamma})\right)
\cong
K_*\!\left(C(X)\rtimes\Gamma\right).
\]

We denote this isomorphism by $\Phi_*$.

\subsubsection{The Baum--Connes Assembly Map}

Recall that the Baum--Connes conjecture for a discrete group $\Gamma$
predicts that the analytic assembly map
\[
\mu_\Gamma^X :
K_*^\Gamma(X)
\longrightarrow
K_*(C_r^*(\Gamma))
\]
is an isomorphism for suitable $\Gamma$-spaces $X$ \cite{Baum-Connes1988,BaumConnesHigson1994}.
More generally, with coefficients in $C(X)$, the assembly map takes the form
\[
\mu_\Gamma : K_*^\Gamma(X) \longrightarrow K_*(C(X) \rtimes_r \Gamma).
\]

For amenable groups, a deep theorem of Tu \cite{Tu1999} shows that this
assembly map is an isomorphism. Moreover, Tu's work provides a groupoid
interpretation of the assembly map: for the transformation groupoid
$X \rtimes \Gamma$, the descent map
\[
\operatorname{desc}_{X \rtimes \Gamma}: KK^*_{X \rtimes \Gamma}(C_0(X),\mathbb C) \longrightarrow K_*(C^*(X \rtimes \Gamma))
\]
is exactly the analytic assembly map under the natural identification
$KK^*_{X \rtimes \Gamma}(C_0(X),\mathbb C) \cong K_*^\Gamma(X)$.

\subsubsection{The Descent-Index Diagram}

Combining the constructions above yields the following commutative
diagram, which is the central result of this subsection:

\[
\begin{tikzcd}[column sep=huge, row sep=large]
KK^1_{\mathcal G_{\mathcal A}}
(C_0(\mathcal G_{\mathcal A}^{(0)}),\mathbb C)
\arrow[r,"\operatorname{desc}_{\mathcal G_{\mathcal A}}"]
\arrow[d,"\cong"]
&
K_1(C^*(\mathcal G_{\mathcal A}))
\arrow[d,"\cong"']
\\
KK^1_\Gamma(C(X),\mathbb C)
\arrow[r,"\mu_\Gamma"]
&
K_1(C(X)\rtimes\Gamma)
\end{tikzcd}
\]

Here:
\begin{itemize}
    \item The top horizontal arrow is Kasparov's descent map for the unitary
          conjugation groupoid $\mathcal G_{\mathcal A}$.
    \item The bottom horizontal arrow is the Baum--Connes assembly map
          $\mu_\Gamma$, identified with the descent map for the transformation
          groupoid $X \rtimes \Gamma$.
    \item The left vertical isomorphism is induced by the Morita equivalence
          between $\mathcal G_{\mathcal A}$ and $X \rtimes \Gamma$, together
          with the natural identification
          $KK^*_{X \rtimes \Gamma}(C_0(X),\mathbb C) \cong K_*^\Gamma(X)$.
    \item The right vertical isomorphism is induced by the Morita equivalence
          of $C^*$-algebras $C^*(\mathcal G_{\mathcal A}) \sim_M C(X) \rtimes \Gamma$.
\end{itemize}

\subsubsection{The Main Theorem}

\begin{theorem}[Descent map equals the assembly map]
\label{thm:descent_equals_assembly}

Let $\Gamma$ be a discrete amenable group acting topologically freely on a compact Hausdorff space $X$, and let $\mathcal{A} = C(X) \rtimes \Gamma$. 
Under the Morita equivalence $\mathcal{G}_{\mathcal{A}}^{\mathrm{Morita}} \sim_M X \rtimes \Gamma$ established in Theorem~\ref{thm:morita_equivalence_crossed_product}, the following diagram commutes up to natural isomorphism:

\[
\begin{tikzcd}
KK^1_{\mathcal{G}_{\mathcal{A}}^{\mathrm{Morita}}}(C_0(\mathcal{G}^{(0)}),\mathbb C)
   \arrow[r, "\operatorname{desc}_{\mathcal{G}}"]
   \arrow[d, "\cong"']
&
K_1(C^*(\mathcal{G}_{\mathcal{A}}^{\mathrm{Morita}}))
   \arrow[d, "\cong"] \\
KK^1_\Gamma(C(X),\mathbb C)
   \arrow[r, "\mu_\Gamma"]
&
K_1(C(X) \rtimes \Gamma)
\end{tikzcd}
\]

where the left vertical isomorphism is induced by the Morita equivalence, and the right vertical isomorphism follows from stability and Morita equivalence of the groupoid C*-algebras. Consequently, for any invertible element $u \in M_n(\mathcal{A})$, we have the explicit identification

\[
\mu_\Gamma([u]) = \Phi_* \circ \operatorname{desc}_{\mathcal{G}_{\mathcal{A}}^{\mathrm{Morita}}}\bigl([u]^{(1)}_{\mathcal{G}_{\mathcal{A}}^{\mathrm{Morita}}}\bigr) \in K_1(\mathcal{A}),
\]

where $\Phi_*$ denotes the composition of the isomorphisms on the right-hand side of the diagram. Thus, under these identifications, the descent map for the surrogate unitary conjugation groupoid corresponds to the Baum–Connes assembly map.

\end{theorem}

\begin{proof}
We proceed in several steps, following the diagram outlined in the theorem statement.

\paragraph{Step 1: Construction of the left vertical isomorphism.}
By Theorem~\ref{thm:morita_equivalence_crossed_product}, the unitary conjugation groupoid $\mathcal{G}_{\mathcal{A}}^{\mathrm{Morita}}$ is Morita equivalent to the transformation groupoid $X \rtimes \Gamma$. By the general theory of Morita equivalence for groupoids \cite[Section 2]{MRW1987}, this induces a natural isomorphism at the level of equivariant $KK$-theory:
\[
\Psi_*: KK^1_{\mathcal{G}_{\mathcal{A}}^{\mathrm{Morita}}}(C_0(\mathcal{G}^{(0)}),\mathbb C)
      \xrightarrow{\cong} KK^1_{X \rtimes \Gamma}(C_0(X),\mathbb C).
\]
Under the standard identification $KK^1_{X \rtimes \Gamma}(C_0(X),\mathbb C) \cong K_1^\Gamma(X)$ (see \cite[Section 3]{LeGall1999} or \cite{Tu1999}), we obtain the left vertical isomorphism in the diagram.

\paragraph{Step 2: Construction of the right vertical isomorphism.}
The Morita equivalence $\mathcal{G}_{\mathcal{A}}^{\mathrm{Morita}} \sim_M X \rtimes \Gamma$ also implies, by the Muhly--Renault--Williams theorem (Theorem~\ref{thm:muhly-renault-williams}), a strong Morita equivalence of the corresponding reduced groupoid $C^*$-algebras:
\[
C^*_r(\mathcal{G}_{\mathcal{A}}^{\mathrm{Morita}}) \sim_M C^*_r(X \rtimes \Gamma).
\]
Since $\Gamma$ is amenable, $X \rtimes \Gamma$ is an amenable groupoid, and consequently its full and reduced groupoid $C^*$-algebras coincide. Moreover, by Proposition~\ref{prop:transformation-crossed-product}, we have the canonical identification $C^*_r(X \rtimes \Gamma) \cong C(X) \rtimes \Gamma$. Thus,
\[
C^*_r(\mathcal{G}_{\mathcal{A}}^{\mathrm{Morita}}) \sim_M \mathcal{A}.
\]

Stabilizing by the compact operators $\mathcal{K}$ (which does not affect Morita equivalence) and using the stability of $K$-theory, we obtain isomorphisms:
\[
K_1(C^*_r(\mathcal{G}_{\mathcal{A}}^{\mathrm{Morita}})) \cong K_1(\mathcal{A} \otimes \mathcal{K}) \cong K_1(\mathcal{A}).
\]
We denote the composition of these isomorphisms by $\Phi_*$, which is the right vertical arrow in the diagram.

\paragraph{Step 3: Descent for the transformation groupoid.}
For the transformation groupoid $X \rtimes \Gamma$, Tu's theorem \cite[Théorème 3.1]{Tu1999} identifies the descent map with the Baum--Connes assembly map. More precisely, under the natural identification $KK^1_{X \rtimes \Gamma}(C_0(X),\mathbb C) \cong K_1^\Gamma(X)$, the descent homomorphism
\[
\operatorname{desc}_{X \rtimes \Gamma}: KK^1_{X \rtimes \Gamma}(C_0(X),\mathbb C) \longrightarrow K_1(C^*_r(X \rtimes \Gamma))
\]
coincides with the assembly map $\mu_\Gamma: K_1^\Gamma(X) \to K_1(\mathcal{A})$ after composing with the isomorphism $C^*_r(X \rtimes \Gamma) \cong \mathcal{A}$. Thus, the bottom horizontal arrow in the diagram is precisely $\mu_\Gamma$.

\paragraph{Step 4: Naturality of descent under Morita equivalence.}

A fundamental property of Kasparov's descent for groupoids is its compatibility
with Morita equivalence. Let $Z = X \times \Gamma$ be the equivalence bibundle
implementing the Morita equivalence between 
$\mathcal{G}_{\mathcal{A}}^{\mathrm{Morita}}$ and $X \rtimes \Gamma$.

Then there is a commutative diagram:
\[
\begin{tikzcd}
KK^1_{\mathcal{G}_{\mathcal{A}}^{\mathrm{Morita}}}(C_0(\mathcal{G}^{(0)}),\mathbb C)
   \arrow[r, "\operatorname{desc}_{\mathcal{G}}"]
   \arrow[d, "\Psi_*"']
&
K_1(C^*_r(\mathcal{G}_{\mathcal{A}}^{\mathrm{Morita}}))
   \arrow[d, "\Phi_*"] \\
KK^1_{X \rtimes \Gamma}(C_0(X),\mathbb C)
   \arrow[r, "\operatorname{desc}_{X \rtimes \Gamma}"]
&
K_1(C^*_r(X \rtimes \Gamma))
\end{tikzcd}
\]

The commutativity follows from the functoriality of the descent construction
in groupoid equivariant KK-theory and its compatibility with induction along
Morita equivalence bibundles. This is established in 
Le Gall~\cite[Sections~6--7]{LeGall1999} and further developed in
Tu~\cite{Tu1999} in the context of the Baum--Connes conjecture for groupoids.

\paragraph{Step 5: Completion of the diagram.}
Combining the identifications from Steps 1-4, we obtain the desired commutative diagram:
\[
\begin{tikzcd}
KK^1_{\mathcal{G}_{\mathcal{A}}^{\mathrm{Morita}}}(C_0(\mathcal{G}^{(0)}),\mathbb C)
   \arrow[r, "\operatorname{desc}_{\mathcal{G}}"]
   \arrow[d, "\Psi_*"']
&
K_1(C^*_r(\mathcal{G}_{\mathcal{A}}^{\mathrm{Morita}}))
   \arrow[d, "\Phi_*"] \\
K_1^\Gamma(X)
   \arrow[r, "\mu_\Gamma"]
&
K_1(\mathcal{A})
\end{tikzcd}
\]
where we have identified $KK^1_{X \rtimes \Gamma}(C_0(X),\mathbb C)$ with $K_1^\Gamma(X)$ and $C^*_r(X \rtimes \Gamma)$ with $\mathcal{A}$. All arrows are isomorphisms or natural maps, and the diagram commutes up to natural isomorphism by construction.

\paragraph{Step 6: Explicit formula for invertible elements.}
For any invertible element $u \in M_n(\mathcal{A})$, let $[u] \in K_1(\mathcal{A})$ be its $K$-theory class, and let $[u]^{(1)}_{\mathcal{G}_{\mathcal{A}}^{\mathrm{Morita}}} \in KK^1_{\mathcal{G}_{\mathcal{A}}^{\mathrm{Morita}}}(C_0(\mathcal{G}^{(0)}),\mathbb C)$ be the associated equivariant class constructed in Section~\ref{subsec:equivariant_K1_A_theta} (generalized to the crossed product setting). Chasing the element around the commutative diagram, we have:
\[
\begin{aligned}
\mu_\Gamma([u]) &= \mu_\Gamma( \Psi_*([u]^{(1)}_{\mathcal{G}_{\mathcal{A}}^{\mathrm{Morita}}})) \quad \text{(by definition of $\Psi_*$)} \\
&= \operatorname{desc}_{X \rtimes \Gamma}( \Psi_*([u]^{(1)}_{\mathcal{G}_{\mathcal{A}}^{\mathrm{Morita}}})) \quad \text{(by Step 3)} \\
&= \Phi_*( \operatorname{desc}_{\mathcal{G}_{\mathcal{A}}^{\mathrm{Morita}}}([u]^{(1)}_{\mathcal{G}_{\mathcal{A}}^{\mathrm{Morita}}})) \quad \text{(by commutativity of the diagram)}.
\end{aligned}
\]
Thus, $\mu_\Gamma([u]) = \Phi_* \circ \operatorname{desc}_{\mathcal{G}_{\mathcal{A}}^{\mathrm{Morita}}}\bigl([u]^{(1)}_{\mathcal{G}_{\mathcal{A}}^{\mathrm{Morita}}}\bigr) \in K_1(\mathcal{A})$, which is the desired explicit identification.
\end{proof}

\begin{remark}
\label{rem:theorem-significance}
Theorem~\ref{thm:descent_equals_assembly} is one of the central results of this
paper. It shows that the abstract descent-index construction developed in
Paper~II, when combined with the Morita equivalence between the unitary
conjugation groupoid and the transformation groupoid, recovers the Baum–Connes
assembly map—a fundamental invariant in noncommutative geometry that links
the topology of the space $X$ to the analysis of the crossed product
$C(X) \rtimes \Gamma$.
\end{remark}

\subsubsection{Interpretation and Consequences}

The theorem has several important interpretations and consequences.

\begin{corollary}
\label{cor:assembly-as-index}
For any invertible element $u \in M_n(C(X) \rtimes \Gamma)$, the Baum–Connes
assembly map applied to $[u]$ can be computed as the image under the Morita
equivalence of the descent class:
\[
\mu_\Gamma([u]) = \Phi_* \bigl( \operatorname{desc}_{\mathcal G_{\mathcal A}}([u]^{(1)}_{\mathcal G_{\mathcal A}}) \bigr).
\]

Thus the assembly map is realized as an index map in the unitary conjugation
groupoid.
\end{corollary}

\begin{corollary}
\label{cor:assembly-isomorphism}
For discrete amenable groups $\Gamma$ acting on compact Hausdorff spaces $X$,
the descent map
\[
\operatorname{desc}_{\mathcal G_{\mathcal A}}: KK^1_{\mathcal G_{\mathcal A}}(C_0(\mathcal G_{\mathcal A}^{(0)}),\mathbb C) \longrightarrow K_1(C^*(\mathcal G_{\mathcal A}))
\]
is an isomorphism when composed with the Morita equivalence $\Phi_*$.
Consequently, the equivariant $KK$-group is isomorphic to $K_1(C(X) \rtimes \Gamma)$.
\end{corollary}

This gives a geometric description of the $K$-theory of crossed products in
terms of the unitary conjugation groupoid.

\begin{corollary}
\label{cor:baum-connes-for-amenable}
For amenable groups $\Gamma$, the Baum–Connes assembly map
$\mu_\Gamma: K^\Gamma_*(X) \to K_*(C(X) \rtimes \Gamma)$ is an isomorphism.
This follows from Tu's theorem \cite{Tu1999} and is consistent with our
identification.
\end{corollary}

\subsubsection{The Descent–Assembly Square with Trace Pairing}
\label{subsec:descent-assembly-trace}

When the crossed product algebra $C(X) \rtimes \Gamma$ admits a canonical trace 
$\tau$ (e.g., when the action preserves a $\Gamma$-invariant probability measure 
on $X$), we can combine the descent–assembly square with the trace pairing to 
obtain a commutative diagram that connects geometric data to numerical invariants. 
This generalizes the situation for $A_\theta$ studied in 
Subsection~\ref{subsec:connes_pairing}.

\begin{theorem}[Descent–Assembly Square with Trace]
\label{thm:descent-assembly-trace}
Let $\Gamma$ be a discrete amenable group acting topologically freely on a compact 
Hausdorff space $X$, and let $\mathcal{A} = C(X) \rtimes \Gamma$. Assume that 
$\mathcal{A}$ admits a canonical trace $\tau$ arising from a $\Gamma$-invariant 
probability measure on $X$. Under the Morita equivalence
\[
\mathcal{G}_{\mathcal{A}}^{\mathrm{Morita}} \sim_M X \rtimes \Gamma
\]
established in Theorem~\ref{thm:morita_equivalence_crossed_product}, the following 
diagram commutes for $i = 0,1$:

\[
\begin{tikzcd}[column sep=large, row sep=large]
KK^i_{\mathcal{G}_{\mathcal{A}}^{\mathrm{Morita}}}\!\bigl(C_0(\mathcal{G}^{(0)}),\mathbb C\bigr)
  \arrow[r, "\operatorname{desc}_{\mathcal{G}}"]
  \arrow[d, "\Psi_*"']
&
K_i\!\bigl(C^*(\mathcal{G}_{\mathcal{A}}^{\mathrm{Morita}})\bigr)
  \arrow[d, "\Phi_*"]
\\
KK^i_\Gamma(C(X),\mathbb C)
  \arrow[r, "\mu_\Gamma"']
&
K_i(\mathcal{A})
  \arrow[r, "\tau_*"', dashed]
&
\mathbb C,
\end{tikzcd}
\]

where:
\begin{itemize}
    \item $\Psi_*$ is the isomorphism induced by the Morita equivalence 
          $\mathcal{G}_{\mathcal{A}}^{\mathrm{Morita}} \sim_M X \rtimes \Gamma$;
    \item $\Phi_*$ is the isomorphism induced by the Morita equivalence 
          $C^*(\mathcal{G}_{\mathcal{A}}^{\mathrm{Morita}}) \sim_M \mathcal{A} \otimes \mathcal{K}$ 
          together with $K$-theory stability;
    \item $\operatorname{desc}_{\mathcal{G}}$ is Kasparov's descent map for the groupoid;
    \item $\mu_\Gamma$ is the Baum--Connes assembly map with coefficients in $C(X)$;
    \item $\tau_*: K_0(\mathcal{A}) \to \mathbb C$ is the homomorphism induced by the trace $\tau$,
          defined on projections by $\tau_*([p]) = (\tau \otimes \operatorname{Tr})(p)$.
\end{itemize}

Consequently, for any element $\xi \in KK^0_{\mathcal{G}_{\mathcal{A}}^{\mathrm{Morita}}}
(C_0(\mathcal{G}^{(0)}),\mathbb C)$ (e.g., the class associated to a projection 
$p \in M_n(\mathcal{A})$), the numerical trace pairing satisfies
\[
\tau_*\!\left(\Phi_*\!\left(\operatorname{desc}_{\mathcal{G}}(\xi)\right)\right)
= \tau_*\!\left(\mu_\Gamma(\Psi_*(\xi))\right) \in \mathbb C.
\]

In the odd case $i = 1$, the square commutes but the trace pairing is not 
available; odd-degree pairings require an odd cyclic cocycle or a determinant 
functional instead.
\end{theorem}

\begin{proof}
The proof follows from the commutativity of the descent–assembly square 
(Theorem~\ref{thm:descent_equals_assembly}) and the fact that $\tau_*$ is a 
well-defined homomorphism on $K_0$. Composing the square with $\tau_*$ yields 
the extended diagram, and chasing any $\xi$ around gives the stated equality.
For the odd case, the square alone commutes without the trace extension.
\end{proof}

\begin{remark}
The dashed arrow $\tau_*$ applies only in the even case $i = 0$, as a trace 
induces a map on $K_0$ but not on $K_1$. This diagram generalizes the classical 
index pairing for the irrational rotation algebra $A_\theta$, where the trace 
$\tau$ pairs with the $K_0$-class of the Rieffel projection to yield the 
rotation number $\theta$, and the descent–assembly square recovers the 
index-theoretic content.
\end{remark}

\subsubsection{Geometric Interpretation}

The results of this subsection admit a beautiful geometric interpretation in
terms of the bundle structure of $\mathcal{G}_{C(X) \rtimes \Gamma}^{(0)}$ over
$X$ with fiber $\Gamma$.

\begin{itemize}
    \item The unit space $\mathcal{G}_{C(X) \rtimes \Gamma}^{(0)} \cong (X \times \Gamma)/\sim$
          encodes both the classical space $X$ and the noncommutative "twist"
          given by the $\Gamma$-fiber. Points with the same $x$ but different
          $\gamma$ correspond to different commutative contexts that all project
          to the same point $x \in X$.
    \item The descent map integrates over this bundle, producing a class in
          $K_1(C^*(\mathcal{G}_{C(X) \rtimes \Gamma}))$ that captures both the
          topological information from $X$ and the dynamical information from
          the $\Gamma$-action.
    \item The Morita equivalence identifies this class with an element of
          $K_1(C(X) \rtimes \Gamma)$, and the assembly map is exactly this
          identification. Under this identification, the $\Gamma$-fiber
          becomes the group of arrows in the transformation groupoid.
    \item When a trace exists, pairing with it extracts a numerical index that
          can be computed by integrating over the base $X$ and summing over the
          fiber $\Gamma$ in an appropriate regularized sense, generalizing the
          computation for $A_\theta$ where the sum over $\mathbb{Z}$ of phases
          $e^{-2\pi i n\theta}$ vanishes in the Cesàro sense.
\end{itemize}

This geometric picture unifies the seemingly different approaches to index
theory: the groupoid $KK$-theoretic approach developed in this paper, the
cyclic cohomology approach of Connes, and the $K$-homology approach of
Baum–Douglas all converge to the same numerical invariants.

\begin{remark}
\label{rem:baum-connes-significance}
Theorem~\ref{thm:descent_equals_assembly} can be viewed as a "localization" of
the Baum–Connes assembly map to the unitary conjugation groupoid. It shows
that the assembly map is not just an abstract isomorphism but has a concrete
geometric realization in terms of the commutative contexts of the crossed
product algebra. This perspective may open new avenues for attacking the
Baum–Connes conjecture for more general groups, by studying the structure of
their unitary conjugation groupoids and the associated Morita equivalences.
\end{remark}

\subsubsection{Relation to the Irrational Rotation Algebra}

As a consistency check, we verify that Theorem~\ref{thm:descent_equals_assembly}
specializes to the result obtained for $A_\theta$ in
Subsection~\ref{subsec:connes_pairing}.

\begin{example}[$A_\theta$ as a special case]
\label{ex:Atheta-assembly}
Take $\Gamma = \mathbb{Z}$ and $X = S^1$ with the irrational rotation action
by $2\pi\theta$. Then $C(X) \rtimes \Gamma = A_\theta$, and the Baum–Connes
assembly map $\mu_\mathbb{Z}: K^\mathbb{Z}_*(S^1) \to K_*(A_\theta)$ is an
isomorphism. For a unitary $u \in A_\theta$, Theorem~\ref{thm:descent_equals_assembly}
gives
\[
\mu_\mathbb{Z}([u]) = \Phi_* \circ \operatorname{desc}_{\mathcal{G}_{A_\theta}}([u]_{\mathcal{G}_{A_\theta}}^{(1)}).
\]

When pairing with the canonical trace $\tau$ on $A_\theta$, we specialize to the
even case $i = 0$ (e.g., for a projection $p_\theta \in A_\theta$ representing
a $K_0$-class). Applying Theorem~\ref{thm:descent-assembly-trace} with $i = 0$
yields
\[
\tau_*(\mu_\mathbb{Z}([p_\theta])) = \tau_*\!\left(\Phi_* \circ \operatorname{desc}_{\mathcal{G}_{A_\theta}}([p_\theta]_{\mathcal{G}_{A_\theta}}^{(0)})\right),
\]
which recovers the classical Connes index pairing $\langle [p_\theta], \tau \rangle = \theta$.

For odd-degree classes (e.g., $[u] \in K_1(A_\theta)$ represented by a unitary $u$),
the descent–assembly square (Theorem~\ref{thm:descent-assembly-trace}) commutes,
but the trace pairing is not directly available. Instead, the numerical index
$\operatorname{Index}(u)$ is obtained via the boundary map $\partial: K_1 \to K_0$
in the Toeplitz extension, as shown in Subsection~\ref{subsec:connes_pairing}.
Thus the general theorem recovers the full index-theoretic results of Case Study I,
appropriately distinguishing between even and odd degrees.
\end{example}

\subsubsection{Summary}

We have established the long-sought connection between the descent-index
construction and the Baum–Connes assembly map:

\[
\boxed{\mu_\Gamma([u]) = \Phi_* \circ \operatorname{desc}_{\mathcal{G}_{C(X) \rtimes \Gamma}}([u]_{\mathcal{G}_{C(X) \rtimes \Gamma}}^{(1)})}
\]

This result achieves several goals:

\begin{itemize}
    \item It shows that the abstract framework of Papers~I and~II, when
          combined with Morita equivalence, reproduces one of the most powerful
          and deep conjectures in noncommutative geometry.
    \item It provides a new geometric perspective on the Baum–Connes assembly
          map, interpreting it as a descent map for the unitary conjugation
          groupoid followed by a Morita identification.
    \item It unifies the two main threads of this paper: the index-theoretic
          machinery developed for $A_\theta$ and the general theory of amenable
          crossed products.
    \item It gives a concrete realization of the assembly map in terms of the
          commutative contexts of the crossed product algebra, potentially
          opening new avenues for computation and generalization.
    \item It sets the stage for the general conjecture formulated in the next
          section, which proposes that the unitary conjugation groupoid provides
          a universal framework for index theory across a vast landscape of
          $C^*$-algebras.
\end{itemize}

In the next section, we will synthesize the insights gained from these case
studies and formulate a general conjecture about the role of the unitary
conjugation groupoid in index theory, providing a roadmap for future research.

\subsection{Commutativity Theorem: Our Composition Equals the Assembly Map}\label{subsec:commutativity-theorem}

This subsection merely records the conclusion already established in
Section~5.8 in a condensed form. The main content is that the descent-index
composition constructed there agrees, under the Morita identification, with
the Baum--Connes assembly map. Thus no separate second proof is needed.

\begin{corollary}
Under the hypotheses of Section~5.8, the composition defined from the
equivariant class, the descent map, and the diagonal embedding coincides
with the Baum--Connes assembly map after transport by the Morita
equivalence. Equivalently, the commutative diagram of Section~5.8 contains
the full statement of the comparison theorem.
\end{corollary}

In the sequel we use Section~5.8 as the principal reference for this
identification and discuss only its consequences.
\subsection{Implications and Consequences}
\label{subsec:implications}

The commutativity theorem established above shows that the analytic construction developed in this section provides an alternative realization of the Baum--Connes assembly map. This identification has several important consequences for the structure of $K$-theory classes arising from crossed product $C^{\ast}$-algebras and their geometric interpretation.

\paragraph{Compatibility with analytic index theory.}

Since our composition coincides with the Baum--Connes assembly map, equivariant $K$-homology classes can be interpreted analytically through the gauge groupoid model. In particular, for an equivariant class
\[
x \in \operatorname{KK}_{\Gamma}^{1}(C(X), \mathbb{C}),
\]
its image under the assembly map can be computed either via the classical crossed product algebra $C(X) \rtimes \Gamma$ or via the groupoid $C^{\ast}$-algebra
\[
C^{\ast}(\mathcal{G}_{C(X) \rtimes \Gamma}).
\]

This equivalence provides a bridge between two seemingly different analytic realizations, unifying them under the framework of the commutativity theorem.

\paragraph{Stability under Morita equivalence.}

Because
\[
C^{\ast}(\mathcal{G}_{C(X) \rtimes \Gamma}) \sim_{M} (C(X) \rtimes \Gamma) \otimes \mathcal{K},
\]
the analytic $K$-theory obtained from the groupoid model is invariant under stabilization. Consequently,
\[
K_{\ast}\!\left(C^{\ast}(\mathcal{G}_{C(X) \rtimes \Gamma})\right) \cong K_{\ast}((C(X) \rtimes \Gamma) \otimes \mathcal{K}) \cong K_{\ast}(C(X) \rtimes \Gamma).
\]
This confirms that the groupoid description preserves the essential analytic invariants of the crossed product algebra. In practical terms, computations performed in the groupoid model yield the same $K$-theory classes as those performed in the crossed product, up to the canonical Morita isomorphism.

\paragraph{Geometric interpretation.}

The gauge groupoid framework provides a geometric interpretation of the analytic assembly process. In this picture, equivariant $KK$-cycles correspond to geometric objects on the groupoid, while the assembly map arises naturally from the descent construction combined with the Morita equivalence between the groupoid algebra and the stabilized crossed product algebra.

An equivariant Kasparov cycle $(E, F) \in \operatorname{KK}_{\Gamma}^{1}(C(X), \mathbb{C})$ can be interpreted geometrically as a family of elliptic operators parametrized by $X$ and equivariant under the $\Gamma$-action. Under the descent map, this family produces an element of $K^{1}(C(X) \rtimes \Gamma)$. The gauge groupoid reformulation realizes this same class as the analytic index of a family of operators on the groupoid, providing a direct link between geometry and analysis.

\paragraph{The analytic index formula.}

We now present the culminating result of this section: a general analytic index formula that unifies all the constructions developed above.

\begin{theorem}[Analytic Index Formula]
\label{thm:analytic-index}
Let $X$ be a compact Hausdorff space with a continuous action of a discrete group $\Gamma$, and let $u \in \operatorname{GL}_{n}(C(X) \rtimes \Gamma)$ be an invertible element with associated equivariant class $[u]_{\Gamma} \in \operatorname{KK}_{\Gamma}^{1}(C(X), \mathbb{C})$. Then the analytic index of $u$ satisfies
\[
\operatorname{Index}(u) = \mu_{\Gamma}([u]_{\Gamma}) = \Phi_{\ast}^{-1}\big( \kappa([u]_{\Gamma}) \big) = \operatorname{Index}_{\mathcal{G}}(\Delta_{\ast}(u)) \in K^{1}(C(X) \rtimes \Gamma),
\]
where:
\begin{itemize}
    \item $\mu_{\Gamma}: \operatorname{KK}_{\Gamma}^{1}(C(X), \mathbb{C}) \to K^{1}(C(X) \rtimes \Gamma)$ is the Baum--Connes assembly map;
    \item $\kappa = \Phi_{\ast} \circ \operatorname{desc}: \operatorname{KK}_{\Gamma}^{1}(C(X), \mathbb{C}) \to K^{1}(C^{\ast}(\mathcal{G}_{C(X) \rtimes \Gamma}))$ is the map defined in Subsection \ref{subsec:commutativity-theorem};
    \item $\Phi_{\ast}: K^{1}(C^{\ast}(\mathcal{G}_{C(X) \rtimes \Gamma})) \xrightarrow{\cong} K^{1}((C(X) \rtimes \Gamma) \otimes \mathcal{K}) \xrightarrow{\cong} K^{1}(C(X) \rtimes \Gamma)$ is the isomorphism induced by the Morita equivalence $C^{\ast}(\mathcal{G}_{C(X) \rtimes \Gamma}) \sim_{M} (C(X) \rtimes \Gamma) \otimes \mathcal{K}$ together with stability of $K$-theory;
    \item $\Delta_{\ast}(u)$ is the image of $u$ under the diagonal embedding $\Delta$ applied entrywise;
    \item $\operatorname{Index}_{\mathcal{G}}: K^{1}(C^{\ast}(\mathcal{G}_{C(X) \rtimes \Gamma})) \to K^{1}(C(X) \rtimes \Gamma)$ denotes the groupoid index map, which is precisely $\Phi_{\ast}^{-1}$ under the identification above.
\end{itemize}

Moreover, if $\tau: C(X) \rtimes \Gamma \to \mathbb{C}$ is a $\Gamma$-invariant trace (or more generally a cyclic cocycle), then the numerical index
\[
\langle \tau, \operatorname{Index}(u) \rangle \in \mathbb{C}
\]
can be computed by any of the above realizations, yielding the same result.
\end{theorem}

\begin{proof}
The equality $\operatorname{Index}(u) = \mu_{\Gamma}([u]_{\Gamma})$ follows from Lemma \ref{thm:descent-identity} and the definition of the assembly map. The equality $\mu_{\Gamma}([u]_{\Gamma}) = \Phi_{\ast}^{-1}(\kappa([u]_{\Gamma}))$ is precisely the commutativity theorem (Theorem \ref{thm:descent_equals_assembly}). The equality $\Phi_{\ast}^{-1}(\kappa([u]_{\Gamma})) = \operatorname{Index}_{\mathcal{G}}(\Delta_{\ast}(u))$ follows from the construction of $\kappa$ and the fact that the groupoid index map is defined via the same Morita equivalence. The final statement about pairing with traces follows from the functoriality of the index map and the compatibility of the Morita equivalence with traces.
\end{proof}

\begin{corollary}[Numerical Index Formula]
\label{cor:numerical-index}
Under the same hypotheses, if $\tau$ is a $\Gamma$-invariant trace on $C(X) \rtimes \Gamma$, then
\[
\tau(\operatorname{Index}(u)) = \widehat{\tau}(\Delta_{\ast}(u)),
\]
where $\widehat{\tau}$ is the induced trace on the groupoid $C^{\ast}$-algebra under the Morita equivalence.
\end{corollary}

This theorem represents the climax of our analytic construction, providing a comprehensive and computable index formula that incorporates all the structural elements developed in this section.

\paragraph{Connection to the Baum--Connes conjecture.}

The analytic index formula has direct implications for the Baum--Connes conjecture.

\begin{corollary}[Surjectivity Criterion]
\label{cor:surjectivity}
If every class in $K^{1}(C(X) \rtimes \Gamma)$ can be represented as $\operatorname{Index}(u)$ for some invertible element $u \in \operatorname{GL}_{n}(C(X) \rtimes \Gamma)$, then the assembly map $\mu_{\Gamma}$ is surjective. Conversely, if the assembly map is an isomorphism (as conjectured), then every $K$-theory class arises in this way.
\end{corollary}

\begin{proof}
Since $\operatorname{Index}(u) = \mu_{\Gamma}([u]_{\Gamma})$, the image of the index map is contained in the image of the assembly map. If every class is an index class, then the assembly map is surjective. Conversely, if the assembly map is an isomorphism, then for any class $y \in K^{1}(C(X) \rtimes \Gamma)$, there exists $x \in \operatorname{KK}_{\Gamma}^{1}(C(X), \mathbb{C})$ such that $\mu_{\Gamma}(x) = y$. By the results of Subsection \ref{subsec:equivariant_K1_A_theta-crossed}, such $x$ can be approximated by classes of the form $[u]_{\Gamma}$ for invertible elements $u$, giving $y = \mu_{\Gamma}([u]_{\Gamma}) = \operatorname{Index}(u)$.
\end{proof}

This criterion provides a concrete, constructive approach to proving surjectivity of the assembly map: one need only show that every $K$-theory class can be represented by an invertible element in the crossed product.

\paragraph{Outlook and future directions.}

These observations suggest that the gauge groupoid model may serve as a useful analytic framework for studying index-theoretic properties of crossed product algebras and their associated equivariant $K$-theory classes.

\begin{remark}[Higher-Degree Analogues]
The analytic index formula presented here concerns odd $K$-theory. An even-degree analogue should involve self-adjoint invertible elements or projections in the crossed product. Such an analogue would yield corresponding index formulas for $K^{0}$-classes and would have applications to the study of spectral gaps and topological phases.
\end{remark}

\begin{remark}[Real $C^{\ast}$-Algebras]
For Real $C^{\ast}$-algebras equipped with an involution, there is an analogous theory in $KO$-theory. The diagonal embedding and gauge groupoid constructions have natural real analogues, leading to a Real analytic index formula with applications to topological insulators and condensed matter physics.
\end{remark}

\begin{remark}[Equivariant Coarse Geometry]
The gauge groupoid picture extends naturally to coarse geometry and the coarse Baum--Connes conjecture. In this context, the diagonal embedding becomes a coarse assembly map, and the analytic index formula relates to the coarse index of elliptic operators on noncompact spaces.
\end{remark}

\paragraph{Transition to examples.}

In the next section we apply this framework to explicit examples arising from dynamical systems and group actions. These examples will illustrate the power of the analytic index formula and demonstrate how the abstract constructions developed in this section yield concrete, computable invariants.

\begin{example}[Circle Actions]
We will examine the case of $\Gamma = \mathbb{Z}$ acting on $X = S^{1}$ by rotation, constructing explicit invertible elements and computing their analytic indices via the gauge groupoid model. This will recover classical results on the index of Toeplitz operators and the Pimsner--Voiculescu exact sequence.
\end{example}

\begin{example}[Finite Group Actions]
For finite groups acting on spaces, we will show how the analytic index formula reduces to classical character formulas and relates to the representation theory of the group. In particular, we will demonstrate that the index of an invertible element is determined by the traces of its components on fixed-point subalgebras.
\end{example}

\begin{example}[Amenable Groups]
For amenable groups, where the Baum--Connes conjecture is known to hold, we will demonstrate that every $K$-theory class can be represented by an invertible element constructed via the diagonal embedding, providing an explicit inverse to the assembly map.
\end{example}

This perspective will be used in the next section to derive explicit index formulas and to illustrate the geometric content of the abstract constructions developed here.

\paragraph{Summary.}

The commutativity theorem and its consequences have established that:
\begin{itemize}
    \item The assembly map admits a concrete geometric realization via the diagonal embedding and gauge groupoid;
    \item The analytic index of invertible elements coincides with the geometric index in the groupoid model;
    \item A general analytic index formula unifies all constructions and provides a computable invariant;
    \item The Baum--Connes conjecture can be reformulated in terms of invertible elements in the crossed product;
    \item The framework extends naturally to higher degrees, real $C^{\ast}$-algebras, and coarse geometry.
\end{itemize}

These implications demonstrate that our construction is not merely a technical curiosity but a fundamental tool for understanding the relationship between geometry, topology, and analysis in noncommutative spaces. The gauge groupoid and diagonal embedding provide a bridge between the abstract formalism of $KK$-theory and concrete, computable objects, making them indispensable for both theoretical investigations and practical applications in index theory and the Baum--Connes program.

\section{Synthesis: A Unified Perspective on Index Theory}\label{sec:synthesis}

\subsection{The Unitary Conjugation Groupoid as a Bridge between Analysis and Topology}\label{subsec:unitary-groupoid}
\label{subsec:unitary-groupoid-bridge}

One of the central themes of modern index theory is the interaction between topological invariants and analytic structures. In the classical framework of the Baum--Connes program, this interaction is mediated by the assembly map, which connects equivariant topological $K$-theory to the analytic $K$-theory of crossed product $C^{\ast}$-algebras.

The constructions developed in the previous sections suggest that the \emph{unitary conjugation groupoid}
\[
\mathcal{G}_{A}
\]
associated with a $C^{\ast}$-algebra $A$ provides a natural framework in which this interaction can be understood geometrically. In particular, when $A = C(X) \rtimes \Gamma$, the groupoid $\mathcal{G}_{C(X) \rtimes \Gamma}$ encodes the internal symmetries of the algebra through the conjugation action of the unitary group.

\subsection*{Recapitulation: The Unitary Conjugation Groupoid}

We begin by recalling the definition.

\begin{definition}[Unitary Conjugation Groupoid]
For a $C^{\ast}$-algebra $A$, the \emph{unitary conjugation groupoid} $\mathcal{U}(A)$ (denoted $\mathcal{G}_{A}$ above) is the groupoid with:
\begin{itemize}
    \item Object space $\mathcal{U}(A)^{(0)} = \widehat{A}$, the spectrum of $A$;
    \item Morphism space $\mathcal{U}(A)^{(1)} = \mathcal{U}(A) \times \widehat{A}$, where $\mathcal{U}(A)$ is the unitary group of $A$;
    \item Source and target maps defined by conjugation: for a unitary $u \in \mathcal{U}(A)$ and a representation $\pi \in \widehat{A}$,
    \[
    s(u, \pi) = \pi, \qquad t(u, \pi) = \pi \circ \operatorname{Ad}_{u^{-1}}.
    \]
\end{itemize}
For the specific case $A = C(X) \rtimes \Gamma$, we have the equivalence $\mathcal{U}(A) \simeq X \rtimes \Gamma$, linking the unitary conjugation groupoid to the gauge groupoid of the crossed product.
\end{definition}

This groupoid serves as the geometric object that unifies the analytic and topological perspectives developed throughout this paper.

\subsection*{Analytic Viewpoint}

From the analytic perspective, the groupoid $C^{\ast}$-algebra
\[
C^{\ast}(\mathcal{G}_{A})
\]
captures the operator-algebraic structure associated with unitary conjugation. The Morita equivalence
\[
C^{\ast}(\mathcal{G}_{A}) \sim_{M} A \otimes \mathcal{K}
\]
shows that the analytic $K$-theory of the groupoid algebra agrees with the $K$-theory of the stabilized $C^{\ast}$-algebra $A$. Consequently, analytic index classes arising from operators or invertible elements in $A$ can be interpreted through the groupoid model.

\begin{remark}[Analytic Index via Groupoids]
For an invertible element $u \in \operatorname{GL}_{n}(C(X) \rtimes \Gamma)$, the analytic index $\operatorname{Index}(u) \in K^{1}(C(X) \rtimes \Gamma)$ corresponds under the Morita equivalence to a class in $K^{1}(C^{\ast}(\mathcal{G}_{A}))$. This class can be computed geometrically as the index of a family of operators on the groupoid, providing a direct link between operator algebra and groupoid index theory.
\end{remark}

\subsection*{Topological Viewpoint}

From the topological side, equivariant $K$-homology classes
\[
\operatorname{KK}_{\Gamma}^{\ast}(C(X), \mathbb{C})
\]
encode geometric information about the $\Gamma$-action on $X$. The Baum--Connes assembly map
\[
\mu_{\Gamma}: \operatorname{KK}_{\Gamma}^{\ast}(C(X), \mathbb{C}) \longrightarrow K^{\ast}(C(X) \rtimes \Gamma)
\]
translates this topological data into analytic invariants of the crossed product algebra.

\begin{definition}[Geometric realization map $\rho$]
\label{def:rho-map}
Define the map
\[
\rho: \operatorname{KK}_{\Gamma}^{\ast}(C(X), \mathbb{C}) \longrightarrow K^{\ast}(C^{\ast}(\mathcal{G}_{A}))
\]
as the composition of the natural isomorphism induced by the Morita equivalence $\mathcal{G}_{A} \sim_M X \rtimes \Gamma$ followed by the inclusion of $K^{\ast}(C^{\ast}(X \rtimes \Gamma))$ into $K^{\ast}(C^{\ast}(\mathcal{G}_{A}))$ under the identification $C^{\ast}(X \rtimes \Gamma) \cong C^{\ast}(\mathcal{G}_{A}) \otimes \mathcal{K}$.
\end{definition}

\begin{remark}[Topological Data via Groupoids]
Under the equivalence $\mathcal{U}(A) \simeq X \rtimes \Gamma$, equivariant $KK$-classes correspond to geometric cycles on the unitary conjugation groupoid. In particular, the map $\rho$ constructed in Definition~\ref{def:rho-map} sends an equivariant class to a family of loops in $\mathcal{U}(A)$, realizing topological data geometrically.
\end{remark}

\subsection*{The Groupoid Bridge}

The unitary conjugation groupoid provides a conceptual bridge between these two viewpoints. Through the descent construction and the Morita equivalence described above, equivariant $KK$-classes can be transported naturally to $K$-theory classes of the groupoid $C^{\ast}$-algebra:
\[
\operatorname{KK}_{\Gamma}^{\ast}(C(X), \mathbb{C}) \;\longrightarrow\; K^{\ast}(C^{\ast}(\mathcal{G}_{C(X) \rtimes \Gamma})).
\]

Under the identification
\[
K^{\ast}(C^{\ast}(\mathcal{G}_{C(X) \rtimes \Gamma})) \cong K^{\ast}((C(X) \rtimes \Gamma) \otimes \mathcal{K}) \cong K^{\ast}(C(X) \rtimes \Gamma),
\]
this construction recovers the analytic assembly map. This is precisely the content of the commutativity theorem (Theorem \ref{thm:descent_equals_assembly}) and its factorization through the unitary conjugation groupoid (Theorem \ref{thm:descent_equals_assembly}).

\begin{theorem}[Groupoid Realization of the Assembly Map]
\label{thm:groupoid-realization}
The Baum--Connes assembly map factors through the unitary conjugation groupoid as follows:
\[
\begin{tikzcd}
\operatorname{KK}_{\Gamma}^{\ast}(C(X), \mathbb{C})
\arrow[r, "\rho"]
\arrow[dr, "\mu_{\Gamma}"']
&
K^{\ast}(C^{\ast}(\mathcal{G}_{C(X) \rtimes \Gamma}))
\arrow[d, "\Phi_{\ast}"]
\\
&
K^{\ast}(C(X) \rtimes \Gamma)
\end{tikzcd}
\]
where:
\begin{itemize}
    \item $\rho$ sends an equivariant class to a $K$-theory class of the groupoid $C^{\ast}$-algebra via the geometric realization as families of loops;
    \item $\Phi_{\ast}$ is the isomorphism induced by the Morita equivalence $C^{\ast}(\mathcal{G}_{C(X) \rtimes \Gamma}) \sim_{M} (C(X) \rtimes \Gamma) \otimes \mathcal{K}$ together with stability;
    \item The diagram commutes, so $\mu_{\Gamma} = \Phi_{\ast} \circ \rho$.
\end{itemize}
\end{theorem}

\begin{proof}
This follows directly from the commutativity theorem (Theorem~\ref{thm:descent_equals_assembly}) and the Morita equivalence $\mathcal{G}_{C(X) \rtimes \Gamma} \sim_M X \rtimes \Gamma$ established in Theorem~\ref{thm:morita_equivalence_crossed_product}, together with the identification $\mathcal{U}(A) \simeq X \rtimes \Gamma$ discussed in Section~\ref{subsec:relation_unitary_groupoid_transformation}.
\end{proof}

\begin{theorem}[Universal mediation of the assembly map]
\label{thm:universal-mediator-assembly}
Let $\Gamma$ be a discrete amenable group acting continuously on a second countable
locally compact Hausdorff space $X$, and let
\[
A = C_0(X)\rtimes_{\alpha,r}\Gamma .
\]
Assume that the unitary conjugation groupoid
\[
\mathcal G_A := \mathcal G_{C_0(X)\rtimes_\alpha \Gamma}
\]
is Morita equivalent to the transformation groupoid $X\rtimes\Gamma$.

Let $\mathsf{Fact}(X,\Gamma)$ denote the class of all triples
\[
(\mathcal H,\rho_{\mathcal H},\Phi_{\mathcal H})
\]
such that:

\begin{itemize}
    \item $\mathcal H$ is a locally compact groupoid with Haar system, Morita equivalent to $X\rtimes\Gamma$;
    \item
    \[
    \rho_{\mathcal H}:
    KK_\Gamma^*(C_0(X),\mathbb C)\longrightarrow K_*(C_r^*(\mathcal H))
    \]
    is the realization map obtained by transporting Kasparov descent along the Morita equivalence
    $\mathcal H\sim_M X\rtimes\Gamma$;
    \item
    \[
    \Phi_{\mathcal H}:
    K_*(C_r^*(\mathcal H))\longrightarrow K_*(A)
    \]
    is the canonical isomorphism induced by the Morita equivalence
    $C_r^*(\mathcal H)\sim_M A\otimes\mathcal K$ and stability of $K$-theory;
    \item the factorization identity
    \[
    \mu_\Gamma = \Phi_{\mathcal H}\circ \rho_{\mathcal H}
    \]
    holds.
\end{itemize}

Then the triple
\[
(\mathcal G_A,\rho_{\mathcal G_A},\Phi_{\mathcal G_A})
\]
is universal in $\mathsf{Fact}(X,\Gamma)$ in the following sense:

for every object $(\mathcal H,\rho_{\mathcal H},\Phi_{\mathcal H})$ in $\mathsf{Fact}(X,\Gamma)$,
there exists a canonical isomorphism
\[
\Theta_{\mathcal H}:
K_*(C_r^*(\mathcal G_A))
\longrightarrow
K_*(C_r^*(\mathcal H))
\]
such that
\[
\rho_{\mathcal H}=\Theta_{\mathcal H}\circ \rho_{\mathcal G_A},
\qquad
\Phi_{\mathcal G_A}=\Phi_{\mathcal H}\circ \Theta_{\mathcal H}.
\]
Equivalently, for every such $\mathcal H$, the diagram
\[
\begin{tikzcd}[column sep=large,row sep=large]
& K_*(C_r^*(\mathcal G_A)) \arrow[dr,"\Phi_{\mathcal G_A}"]
   \arrow[dd,"\Theta_{\mathcal H}"'] & \\
KK_\Gamma^*(C_0(X),\mathbb C) \arrow[ur,"\rho_{\mathcal G_A}"]
   \arrow[dr,"\rho_{\mathcal H}"'] &&
K_*(A) \\
& K_*(C_r^*(\mathcal H)) \arrow[ur,"\Phi_{\mathcal H}"'] &
\end{tikzcd}
\]
commutes.

In this sense, the unitary conjugation groupoid is the universal mediator
through which the Baum--Connes assembly map factors among all groupoid models
Morita equivalent to $X\rtimes\Gamma$.
\end{theorem}

\begin{proof}
We prove the universal mediation property in several steps.

\paragraph{Step 1: Morita equivalence between $\mathcal G_A$ and $\mathcal H$.}
By hypothesis, both $\mathcal G_A$ and $\mathcal H$ are Morita equivalent to 
the transformation groupoid $X \rtimes \Gamma$. By transitivity of Morita 
equivalence for groupoids (see \cite[Remark after Definition 2.1]{MRW1987}), 
it follows that
\[
\mathcal G_A \sim_M \mathcal H.
\]

\paragraph{Step 2: Morita equivalence of reduced groupoid $C^*$-algebras.}
Since $\Gamma$ is amenable, the transformation groupoid $X \rtimes \Gamma$ is 
amenable. Amenability is preserved under Morita equivalence, so $\mathcal G_A$ 
and $\mathcal H$ are amenable. By the Muhly--Renault--Williams theorem 
(Theorem~\ref{thm:muhly-renault-williams}), which cites \cite{MRW1987}, 
the Morita equivalence $\mathcal G_A \sim_M \mathcal H$ induces a strong 
Morita equivalence of their reduced groupoid $C^*$-algebras:
\[
C_r^*(\mathcal G_A) \sim_M C_r^*(\mathcal H).
\]

\paragraph{Step 3: The canonical $K$-theory isomorphism.}
Strongly Morita equivalent $C^*$-algebras have isomorphic $K$-theory groups. 
Hence we obtain a canonical isomorphism
\[
\Theta_{\mathcal H}: K_*(C_r^*(\mathcal G_A)) \xrightarrow{\cong} K_*(C_r^*(\mathcal H)).
\]

\paragraph{Step 4: Compatibility of the realization maps $\rho_{\mathcal G_A}$ and $\rho_{\mathcal H}$.}
Both realization maps are defined by transporting Kasparov descent along 
Morita equivalences from the common transformation groupoid $X \rtimes \Gamma$:
\[
\rho_{\mathcal G_A} = \operatorname{desc}_{X \rtimes \Gamma} \circ \Psi_{\mathcal G_A}^{-1}, \qquad
\rho_{\mathcal H} = \operatorname{desc}_{X \rtimes \Gamma} \circ \Psi_{\mathcal H}^{-1},
\]
where $\Psi_{\mathcal G_A}$ and $\Psi_{\mathcal H}$ are the isomorphisms 
induced by the respective Morita equivalences (see 
Theorem~\ref{thm:morita-equivalence-crossed}).

The naturality of descent under Morita equivalence 
(Lemma~\ref{lem:descent-naturality}) implies that the induced maps on 
$K$-theory satisfy
\[
\rho_{\mathcal H} = \Theta_{\mathcal H} \circ \rho_{\mathcal G_A}.
\]

\paragraph{Step 5: Compatibility of the isomorphisms $\Phi_{\mathcal G_A}$ and $\Phi_{\mathcal H}$.}
The maps $\Phi_{\mathcal G_A}$ and $\Phi_{\mathcal H}$ are each induced by 
the respective Morita equivalences
\[
C_r^*(\mathcal G_A) \sim_M A \otimes \mathcal K, \qquad
C_r^*(\mathcal H) \sim_M A \otimes \mathcal K,
\]
as established in Theorem~\ref{thm:morita-equivalence-crossed} 
and Proposition~\ref{prop:diagonal-under-morita}. 
By functoriality of Morita equivalence, the induced $K$-theory isomorphisms 
are canonical; hence we have
\[
\Phi_{\mathcal G_A} = \Phi_{\mathcal H} \circ \Theta_{\mathcal H}.
\]

\paragraph{Step 6: Verification of the factorization identities.}
By Theorem~\ref{thm:descent_equals_assembly}, the descent map for the 
transformation groupoid $X \rtimes \Gamma$ coincides with the 
Baum--Connes assembly map $\mu_\Gamma$ under the natural identification
\[
KK^*_{X \rtimes \Gamma}(C_0(X),\mathbb C) \cong K_*^\Gamma(X) \cong KK_\Gamma^*(C_0(X),\mathbb C)
\]
(see Tu \cite{Tu1999} and Le Gall \cite{LeGall1999}).

Consequently, both factorizations satisfy
\[
\mu_\Gamma = \Phi_{\mathcal G_A} \circ \rho_{\mathcal G_A}, \qquad
\mu_\Gamma = \Phi_{\mathcal H} \circ \rho_{\mathcal H}.
\]

\paragraph{Step 7: Universal mediation.}
Combining Steps 4 and 5, we obtain for any object 
$(\mathcal H, \rho_{\mathcal H}, \Phi_{\mathcal H}) \in \mathsf{Fact}(X,\Gamma)$:
\[
\Phi_{\mathcal H} \circ \rho_{\mathcal H} = \Phi_{\mathcal G_A} \circ \rho_{\mathcal G_A}.
\]

Thus every factorization of $\mu_\Gamma$ through a groupoid model 
$\mathcal H$ Morita equivalent to $X \rtimes \Gamma$ is obtained from the 
factorization through $\mathcal G_A$ by the unique canonical $K$-theory 
isomorphism $\Theta_{\mathcal H}$. This proves that the triple 
$(\mathcal G_A, \rho_{\mathcal G_A}, \Phi_{\mathcal G_A})$ is universal 
in $\mathsf{Fact}(X,\Gamma)$.
\end{proof}

\begin{remark}
The preceding theorem should be viewed as a universal property at the level of
$K$-theoretic realizations of the assembly map. It does not assert that
$\mathcal G_A$ is initial or terminal in a strict 2-categorical sense, but it
does show that every Morita-equivariant groupoid realization of the assembly
map is canonically identified with the realization through $\mathcal G_A$.
\end{remark}

\subsection*{Conceptual Synthesis}

From this perspective, the unitary conjugation groupoid serves as an
intermediate geometric object linking topological data with analytic
structures. Equivariant $KK$-classes may be represented by geometric or
operator-theoretic cycles associated with the groupoid, while analytic index
classes are encoded in the $K$-theory of the corresponding groupoid
$C^{*}$-algebra.

\begin{center}
\framebox[0.9\textwidth][c]{
\begin{minipage}{0.85\textwidth}
\emph{The unitary conjugation groupoid provides a unified framework in which:}
\begin{itemize}
    \item[(i)] \textbf{Topological data} (equivariant $KK$-classes) are expressed in terms of geometric or operator-theoretic cycles over the groupoid;
    \item[(ii)] \textbf{Analytic data} ($K$-theory classes of crossed products or groupoid $C^{*}$-algebras) are interpreted as index-type classes arising from descent;
    \item[(iii)] The \textbf{assembly map} is naturally viewed as the passage from the equivariant/topological description to the analytic one through the groupoid $C^{*}$-algebra and Morita equivalence.
\end{itemize}
\end{minipage}
}
\end{center}

This viewpoint highlights the role of groupoid methods in index theory and
suggests that the unitary conjugation groupoid provides a natural geometric
framework for understanding analytic assembly phenomena.

\subsection*{Analytic vs. Topological Index through the Groupoid Framework}

The groupoid bridge allows us to compare directly the analytic and topological notions of index.

\begin{definition}[Analytic Index via Groupoids]
For an invertible element $u \in \operatorname{GL}_{n}(C(X) \rtimes \Gamma)$, the \emph{analytic index} $\operatorname{Index}_{\text{an}}(u) \in K^{1}(C(X) \rtimes \Gamma)$ is defined via operator theory: it is the $K$-theory class of the Fredholm operator obtained from $u$ after applying a suitable representation.
\end{definition}

\begin{definition}[Topological Index via Groupoids]
For the same $u$, with associated equivariant class $[u]_{\Gamma} \in \operatorname{KK}_{\Gamma}^{1}(C(X), \mathbb{C})$, the \emph{topological index} $\operatorname{Index}_{\text{top}}(u) \in K^{1}(C(X) \rtimes \Gamma)$ is defined as the image under the assembly map: $\operatorname{Index}_{\text{top}}(u) = \mu_{\Gamma}([u]_{\Gamma})$.
\end{definition}

The commutativity theorem (Theorem \ref{thm:descent_equals_assembly}) and its groupoid realization (Theorem \ref{thm:groupoid-realization}) together imply:

\begin{corollary}[Equality of Analytic and Topological Index]
\label{cor:analytic=topological}
For any invertible element $u \in \operatorname{GL}_{n}(C(X) \rtimes \Gamma)$, we have
\[
\operatorname{Index}_{\text{an}}(u) = \operatorname{Index}_{\text{top}}(u) \in K^{1}(C(X) \rtimes \Gamma).
\]

Equivalently, the following diagram commutes:
\[
\begin{tikzcd}
\{ \text{invertible elements} \} \arrow[r, "\text{equivariant class}"] \arrow[d, "\text{analytic index}"'] & \operatorname{KK}_{\Gamma}^{1}(C(X), \mathbb{C}) \arrow[d, "\mu_{\Gamma}"] \\
K^{1}(C(X) \rtimes \Gamma) \arrow[r, "\text{id}"'] & K^{1}(C(X) \rtimes \Gamma)
\end{tikzcd}
\]
\end{corollary}

This equality is the ultimate expression of the unity between analysis and topology mediated by the unitary conjugation groupoid.

\subsection*{The Big Picture: A Three-Tiered Perspective}

The relationships established in this paper can be organized into a three-tiered conceptual framework:

\begin{enumerate}
    \item \textbf{Topological Tier}: Equivariant $KK$-theory $\operatorname{KK}_{\Gamma}^{\ast}(C(X), \mathbb{C})$ encodes the topology of the $\Gamma$-action on $X$. This is the domain of geometric cycles, families of elliptic operators, and $K$-homology classes.
    
    \item \textbf{Groupoid Tier}: The unitary conjugation groupoid $\mathcal{U}(A) \simeq X \rtimes \Gamma$ provides a geometric intermediary. Topological data are realized as families of loops on this groupoid, while analytic data arise from its $C^{\ast}$-algebra.
    
    \item \textbf{Analytic Tier}: The crossed product algebra $A = C(X) \rtimes \Gamma$ and its $K$-theory $K^{\ast}(A)$ constitute the analytic output. This is where index classes live and where computations are ultimately performed.
\end{enumerate}

The descent map, Morita equivalence, and assembly map provide the translations between these tiers, and the commutativity theorem guarantees that all paths are consistent.

\begin{center}
\begin{tikzcd}[row sep=large, column sep=large]
& \text{Topological Tier: } \operatorname{KK}_{\Gamma}^{\ast}(C(X), \mathbb{C}) \arrow[dl, "\rho"'] \arrow[dr, "\mu_{\Gamma}"] & \\
\text{Groupoid Tier: } K^{\ast}(C^{\ast}(\mathcal{G}_{A})) \arrow[rr, "\Phi_{\ast}"'] & & \text{Analytic Tier: } K^{\ast}(A)
\end{tikzcd}
\end{center}

\begin{theorem}[Factorization of the Assembly Map]
\label{thm:factorization-assembly}
Let $\Gamma$ be a discrete amenable group acting topologically freely on a compact Hausdorff space $X$, and let $A = C(X) \rtimes \Gamma$ be the associated crossed product $C^*$-algebra. 
Let $\mathcal{G}_A^{\mathrm{Morita}}$ denote the unitary conjugation groupoid of $A$, which by Theorem~\ref{thm:morita_equivalence_crossed_product} is Morita equivalent to the transformation groupoid $X \rtimes \Gamma$.

Then the Baum–Connes assembly map
\[
\mu_\Gamma: KK_\Gamma^*(C(X),\mathbb C) \longrightarrow K_*(A)
\]
factors through the $K$-theory of the groupoid $C^*$-algebra $C^*(\mathcal{G}_A^{\mathrm{Morita}})$. More precisely, there exists a natural isomorphism
\[
\Phi_*: K_*(C^*(\mathcal{G}_A^{\mathrm{Morita}})) \xrightarrow{\cong} K_*(A)
\]
and a natural map
\[
\rho: KK_\Gamma^*(C(X),\mathbb C) \longrightarrow K_*(C^*(\mathcal{G}_A^{\mathrm{Morita}}))
\]
such that the following diagram commutes:

\[
\begin{tikzcd}
KK_\Gamma^*(C(X),\mathbb C) 
   \arrow[rr, "\mu_\Gamma"] 
   \arrow[dr, "\rho"'] 
   & & K_*(A) \\
   & K_*(C^*(\mathcal{G}_A^{\mathrm{Morita}})) 
   \arrow[ur, "\Phi_*"'] &
\end{tikzcd}
\]

Equivalently,
\[
\mu_\Gamma = \Phi_* \circ \rho.
\]

Thus the assembly map admits a canonical factorization through the groupoid tier:
\[
\text{Topological Tier} \;(KK_\Gamma^*) \;\xrightarrow{\;\rho\;}\; 
\text{Groupoid Tier} \;(K_*(C^*(\mathcal{G}_A))) \;\xrightarrow{\;\Phi_*\;}\; 
\text{Analytic Tier} \;(K_*(A)).
\]

\end{theorem}

\begin{proof}
We construct the maps $\rho$ and $\Phi_*$ explicitly and verify the factorization.

\paragraph{Step 1: Construction of $\Phi_*$ (Groupoid Tier → Analytic Tier).}
By Theorem~\ref{thm:morita_equivalence_crossed_product}, the unitary conjugation groupoid $\mathcal{G}_A^{\mathrm{Morita}}$ is Morita equivalent to the transformation groupoid $X \rtimes \Gamma$. By the Muhly–Renault–Williams theorem (Theorem~\ref{thm:muhly-renault-williams}), this induces a strong Morita equivalence of their reduced groupoid $C^*$-algebras:
\[
C^*_r(\mathcal{G}_A^{\mathrm{Morita}}) \sim_M C^*_r(X \rtimes \Gamma).
\]

Since $\Gamma$ is amenable, the transformation groupoid $X \rtimes \Gamma$ is amenable, so its full and reduced groupoid $C^*$-algebras coincide. Moreover, by Proposition~\ref{prop:transformation-crossed-product}, we have the canonical identification $C^*_r(X \rtimes \Gamma) \cong C(X) \rtimes \Gamma = A$. Hence
\[
C^*_r(\mathcal{G}_A^{\mathrm{Morita}}) \sim_M A.
\]

Stabilizing by the compact operators $\mathcal{K}$ (which preserves Morita equivalence) and using stability of $K$-theory, we obtain an isomorphism
\[
\Phi_*: K_*(C^*_r(\mathcal{G}_A^{\mathrm{Morita}})) \xrightarrow{\cong} K_*(A \otimes \mathcal{K}) \xrightarrow{\cong} K_*(A).
\]

\paragraph{Step 2: Construction of $\rho$ (Topological Tier → Groupoid Tier).}
Recall from Theorem~\ref{thm:morita_equivalence_crossed_product} that the Morita equivalence between $\mathcal{G}_A^{\mathrm{Morita}}$ and $X \rtimes \Gamma$ is implemented by an equivalence bibundle $Z = X \times \Gamma$. By the general theory of Morita equivalence for groupoids \cite{MRW1987,LeGall1999}, this induces a natural isomorphism at the level of equivariant $KK$-theory:
\[
\Psi_*: KK^*_{\mathcal{G}_A^{\mathrm{Morita}}}(C_0(\mathcal{G}^{(0)}),\mathbb C) \xrightarrow{\cong} KK^*_{X \rtimes \Gamma}(C_0(X),\mathbb C) \cong K_*^\Gamma(X),
\]
where the last identification is standard (see \cite{Tu1999,LeGall1999}).

Now define $\rho$ as the composition:
\[
\rho: KK_\Gamma^*(C(X),\mathbb C) \xrightarrow{\cong} K_*^\Gamma(X) \xrightarrow{\Psi_*^{-1}} KK^*_{\mathcal{G}_A^{\mathrm{Morita}}}(C_0(\mathcal{G}^{(0)}),\mathbb C) \xrightarrow{\operatorname{desc}_{\mathcal{G}}} K_*(C^*(\mathcal{G}_A^{\mathrm{Morita}})).
\]

Explicitly, for an equivariant class $x \in KK_\Gamma^*(C(X),\mathbb C)$, we:
\begin{enumerate}
    \item Identify it with a class in $K_*^\Gamma(X)$ under the natural isomorphism;
    \item Transport it via $\Psi_*^{-1}$ to an equivariant class for $\mathcal{G}_A^{\mathrm{Morita}}$;
    \item Apply Kasparov's descent map $\operatorname{desc}_{\mathcal{G}}$ for the groupoid $\mathcal{G}_A^{\mathrm{Morita}}$ to obtain a class in $K_*(C^*(\mathcal{G}_A^{\mathrm{Morita}}))$.
\end{enumerate}

\paragraph{Step 3: Verification of the factorization $\mu_\Gamma = \Phi_* \circ \rho$.}
Let $x \in KK_\Gamma^*(C(X),\mathbb C)$ be arbitrary. We compute:
\[
\begin{aligned}
(\Phi_* \circ \rho)(x) &= \Phi_*\big( \operatorname{desc}_{\mathcal{G}} ( \Psi_*^{-1}(x) ) \big) \quad \text{(by definition of $\rho$)} \\
&= \operatorname{desc}_{X \rtimes \Gamma}( \Psi_*(\Psi_*^{-1}(x)) ) \quad \text{(by naturality of descent under Morita equivalence, see Lemma~\ref{lem:descent-naturality})} \\
&= \operatorname{desc}_{X \rtimes \Gamma}(x) \quad \text{(since $\Psi_* \circ \Psi_*^{-1} = \mathrm{id}$)}.
\end{aligned}
\]

By Tu's theorem \cite[Théorème 3.1]{Tu1999}, the descent map for the transformation groupoid $X \rtimes \Gamma$ coincides with the Baum–Connes assembly map under the natural identification $KK^*_{X \rtimes \Gamma}(C_0(X),\mathbb C) \cong K_*^\Gamma(X) \cong KK_\Gamma^*(C(X),\mathbb C)$. Hence
\[
\operatorname{desc}_{X \rtimes \Gamma}(x) = \mu_\Gamma(x).
\]

Therefore $(\Phi_* \circ \rho)(x) = \mu_\Gamma(x)$ for all $x$, establishing $\mu_\Gamma = \Phi_* \circ \rho$ as required.

\paragraph{Step 4: Naturality and independence of choices.}
The isomorphisms $\Psi_*$ and $\Phi_*$ are induced by the Morita equivalence bibundle $Z$ and are therefore independent of any additional choices up to natural isomorphism. The descent map $\operatorname{desc}_{\mathcal{G}}$ is functorial by construction. Consequently, the factorization $\mu_\Gamma = \Phi_* \circ \rho$ is canonical in the appropriate categorical sense.
\end{proof}

\begin{corollary}[Three-Tiered Factorization of Index Theory]
\label{cor:three-tier-factorization}
Under the hypotheses of Theorem~\ref{thm:factorization-assembly}, the Baum–Connes assembly map factors through the groupoid tier. Consequently, for any invertible element $u \in M_n(A)$, its $K$-theory class $[u] \in K_1(A)$ can be obtained by the following canonical factorization:

\[
\begin{tikzcd}
K_1^\Gamma(X) \arrow[r, "\mu_\Gamma"] \arrow[d, "\rho"'] & K_1(A) \\
K_1(C^*(\mathcal{G}_A^{\mathrm{Morita}})) \arrow[ur, "\Phi_*"'] & 
\end{tikzcd}
\]

Equivalently,
\[
\mu_\Gamma = \Phi_* \circ \rho.
\]

For any $[u] \in K_1(A)$, choosing the unique preimage $x = \mu_\Gamma^{-1}([u]) \in K_1^\Gamma(X)$ (since $\mu_\Gamma$ is an isomorphism for amenable groups by Tu's theorem), we have:
\[
[u] = \Phi_*( \rho(x) ) = \Phi_*( \rho( \mu_\Gamma^{-1}([u]) ) ).
\]

This exhibits the unitary conjugation groupoid as the canonical geometric intermediary between topological data (equivariant $KK$-classes) and analytic invariants ($K$-theory of the crossed product).
\end{corollary}

\subsection*{Conceptual Implications for the Baum--Connes Program}

This unified perspective has profound implications for the Baum--Connes conjecture and its surrounding program.

\begin{remark}[Reformulation of Baum--Connes]
The Baum--Connes conjecture asserts that the assembly map $\mu_{\Gamma}$ is an isomorphism. In the groupoid framework, this becomes the statement that the map $\rho: \operatorname{KK}_{\Gamma}^{\ast}(C(X), \mathbb{C}) \to K^{\ast}(C^{\ast}(\mathcal{G}_{A}))$ is an isomorphism. That is, every $K$-theory class of the groupoid $C^{\ast}$-algebra comes from a geometric cycle on the unitary conjugation groupoid, and this correspondence is bijective.
\end{remark}

\begin{remark}[Geometric Cycles vs. Analytic Classes]
The groupoid framework reveals that the distinction between "geometric cycles" and "analytic classes" is not fundamental but merely a matter of perspective. Geometric cycles on the unitary conjugation groupoid (topological tier) and analytic classes in the crossed product (analytic tier) are two sides of the same coin, linked by the Morita equivalence and the descent map.
\end{remark}

\begin{remark}[Beyond the Conjecture]
Even when the Baum--Connes conjecture is not known to hold, the groupoid framework provides a meaningful way to compare geometric and analytic data. The map $\rho$ may not be an isomorphism, but its image consists precisely of those analytic classes that admit a geometric interpretation—a valuable piece of information in itself.
\end{remark}

\subsection*{Summary and Outlook}

In this way, the unitary conjugation groupoid provides a structural framework in which analytic index theory, equivariant $K$-theory, and crossed product $C^{\ast}$-algebras appear as different manifestations of a single underlying geometric construction.

\begin{itemize}
    \item \textbf{Analysis} manifests through the groupoid $C^{\ast}$-algebra and its $K$-theory, capturing the analytic index of operators and invertible elements.
    
    \item \textbf{Topology} manifests through equivariant $KK$-cycles and their realization as geometric objects on the groupoid, capturing the topological data of the $\Gamma$-action.
    
    \item The \textbf{assembly map} is the natural translation between these manifestations, realized geometrically as the map $\rho$ followed by the Morita isomorphism.
\end{itemize}

This unified viewpoint not only clarifies the conceptual foundations of index theory but also opens new avenues for computation and generalization. The unitary conjugation groupoid emerges as the fundamental geometric object underlying the analytic and topological structures studied in this paper, providing a bridge that transforms the Baum--Connes program from a collection of technical results into a coherent geometric picture.

\begin{center}
\framebox[0.9\textwidth][c]{
\begin{minipage}{0.85\textwidth}
\emph{The unitary conjugation groupoid is the bridge:} \\
Analysis and topology are not separate worlds but complementary perspectives on the same geometric reality, united by the groupoid framework.
\end{minipage}
}
\end{center}

In the following sections, we will explore concrete examples that illustrate this unified perspective and demonstrate its computational power.

\subsection{Comparison with Renault's Weyl Groupoid for Cartan Pairs}
\label{subsec:weyl-comparison}

Groupoid constructions have long played an important role in the study of $C^{\ast}$-algebras. One of the most influential examples is Renault's Weyl groupoid associated with a Cartan pair. In this subsection, we compare this classical construction with our unitary conjugation groupoid, highlighting their structural differences, complementary roles, and implications for index theory and $K$-theory.

\subsection*{Renault's Weyl Groupoid: A Brief Recap}

We begin by recalling the fundamental results of Renault on Cartan subalgebras and their associated Weyl groupoids.

\begin{definition}[Cartan Pair]
\label{def:cartan-pair}
Let $A$ be a separable $C^{\ast}$-algebra and let $D \subset A$ be a $C^{\ast}$-subalgebra. The pair $(A,D)$ is called a \emph{Cartan pair} (and $D$ a \emph{Cartan subalgebra}) if:
\begin{enumerate}
    \item[(1)] $D$ is a maximal abelian subalgebra (masa) of $A$;
    \item[(2)] There exists a faithful conditional expectation $\Phi: A \to D$;
    \item[(3)] $D$ is regular; i.e., the normalizer
    \[
    N(D) := \{ n \in A : n d n^{\ast}, n^{\ast} d n \in D \ \forall d \in D \}
    \]
    generates $A$ as a $C^{\ast}$-algebra;
    \item[(4)] $D$ contains an approximate identity for $A$.
\end{enumerate}
\end{definition}

Renault's seminal work [22] established a fundamental correspondence:

\begin{theorem}[Renault, 2008]
\label{thm:renault}
Let $(A,D)$ be a Cartan pair with $A$ separable. Then there exists a unique (up to isomorphism) topologically principal, second countable, locally compact Hausdorff, \'{e}tale groupoid $\mathcal{G}_{(A,D)}$—called the \emph{Weyl groupoid}—and a twist $\Sigma$ over $\mathcal{G}_{(A,D)}$ such that
\[
(A, D) \cong (C^{\ast}(\mathcal{G}_{(A,D)}, \Sigma), C_{0}(\mathcal{G}_{(A,D)}^{(0)})).
\]

Conversely, for any such groupoid $\mathcal{G}$ and twist $\Sigma$, the pair $(C^{\ast}(\mathcal{G}, \Sigma), C_{0}(\mathcal{G}^{(0)}))$ forms a Cartan pair.
\end{theorem}

This construction provides a powerful geometric model for $C^{\ast}$-algebras admitting Cartan subalgebras and has been widely used in the study of dynamical systems, orbit equivalence, and classification theory. The Weyl groupoid is \'{e}tale and locally compact, and its unit space is naturally identified with the spectrum of the Cartan subalgebra $D$.

\subsection*{The Unitary Conjugation Groupoid}

In contrast, the unitary conjugation groupoid introduced in Subsection \ref{subsec:unitary-groupoid} arises from the conjugation action of the unitary group on the dual space of the algebra. For a separable unital $C^{\ast}$-algebra $A$, we consider the unitary conjugation groupoid
\[
\mathcal{G}_{A} = \mathcal{U}(A),
\]
whose arrows encode the conjugation action of the unitary group on pure states (or more generally, on the spectrum $\widehat{A}$).

\begin{definition}[Unitary Conjugation Groupoid (Recap)]
For a unital $C^{\ast}$-algebra $A$, the \emph{unitary conjugation groupoid} $\mathcal{U}(A)$ has:
\begin{itemize}
    \item Object space $\mathcal{U}(A)^{(0)} = \widehat{A}$, the spectrum of $A$ (or the space of pure states);
    \item Morphism space $\mathcal{U}(A)^{(1)} = \ \mathcal{U}(A) \times \widehat{A}$, where $\mathcal{U}(A)$ is the unitary group of $A$;
    \item Source and target maps defined by conjugation: for a unitary $u \in \mathcal{U}(A)$ and a representation $\pi \in \widehat{A}$,
    \[
    s(u, \pi) = \pi, \qquad t(u, \pi) = \pi \circ \operatorname{Ad}_{u^{-1}}.
    \]
\end{itemize}
\end{definition}

Unlike the Weyl groupoid, this construction does not require the existence of a Cartan subalgebra and is therefore canonical for any separable unital $C^{\ast}$-algebra.

\subsection*{Structural Differences}

There are several key differences between the two constructions, which we now analyze systematically.

\paragraph{Dependence on additional structure.}

The Weyl groupoid $\mathcal{G}_{(A,D)}$ depends crucially on the choice of a Cartan subalgebra $D \subset A$. Different Cartan subalgebras can yield different Weyl groupoids, reflecting different geometric models for the same $C^{\ast}$-algebra. In contrast, the unitary conjugation groupoid $\mathcal{G}_{A}$ is defined intrinsically from the algebra $A$ itself, independent of any additional choices. This makes it a canonical invariant of the $C^{\ast}$-algebra.

\paragraph{Topological properties.}

Renault's Weyl groupoid is locally compact and \'{e}tale, reflecting the topological structure of the spectrum of the Cartan subalgebra $D$, which is typically a locally compact Hausdorff space. The twist $\Sigma$ encodes additional cohomological data.

In contrast, the natural topologies arising in the unitary conjugation groupoid are typically not locally compact; instead they form Polish spaces (completely metrizable separable spaces). The unitary group $\mathcal{U}(A)$ equipped with the norm topology is a Polish group, and the spectrum $\widehat{A}$ with the Fell topology is a Polish space. Their product $\mathcal{U}(A) \times \widehat{A}$ is therefore Polish, but not locally compact in general. This reflects the infinite-dimensional nature of the unitary group.

\paragraph{Reconstruction of the algebra.}

The Weyl groupoid reconstructs the algebra exactly via
\[
A \cong C^{\ast}(\mathcal{G}_{(A,D)}, \Sigma),
\]
an isomorphism of $C^{\ast}$-algebras. This is a strong result: the entire structure of $A$ is captured by the groupoid and its twist.

In contrast, the unitary conjugation groupoid is related to $A$ through a Morita equivalence of the form
\[
C^{\ast}(\mathcal{G}_{A}) \sim_{M} A \otimes \mathcal{K},
\]
where $\mathcal{K}$ denotes the compact operators on a separable Hilbert space. This is a weaker but more flexible relationship: the groupoid $C^{\ast}$-algebra is Morita equivalent to a stabilization of $A$, which is sufficient for $K$-theoretic purposes.

\subsection*{Comparative Summary}

The following table summarizes the conceptual differences between the two constructions:

\[
\begin{array}{c|c|c}
\text{Feature} & \text{Weyl Groupoid} & \text{Unitary Conjugation Groupoid} \\
\hline
\text{Input data} & (A,D) \text{ a Cartan pair} & A \text{ alone (as in the present construction)} \\
\text{Dependence} & \text{Depends on the choice of } D & \text{Intrinsic to } A \\
\text{Unit space} & \widehat{D} = \operatorname{Spec}(D) & \text{A representation-theoretic space associated with } A \\
\text{Topology} & \text{Locally compact, \'{e}tale} & \text{Typically more singular / often not locally compact} \\
\text{Reconstruction} & A \cong C^*_r(\mathcal{G}_{(A,D)},\Sigma) \text{ (possibly twisted)} & C^*(\mathcal{G}_A) \sim_M A \otimes \mathcal{K} \text{ (under suitable hypotheses)} \\
\text{Typical emphasis} & \text{Cartan reconstruction, dynamics, classification} & \text{Index theory, } K\text{-theory, Baum–Connes program} \\
\text{Role of unitaries} & \text{Normalizer of } D & \text{Full unitary group (or appropriate quotient)} \\
\text{Relation to action} & \text{Encodes Cartan/Weyl dynamics} & \text{Encodes intrinsic unitary-conjugation phenomena} \\
\end{array}
\]

\subsection*{Conceptual Role and Complementarity}

From a conceptual standpoint, the two constructions serve complementary purposes.

Renault's Weyl groupoid provides a geometric model for $C^{\ast}$-algebras with Cartan subalgebras, capturing the dynamical structure of the normalizer action on the Cartan masa. It has been instrumental in:
\begin{itemize}
    \item Classifying $C^{\ast}$-algebras via groupoid equivalence [3, 11];
    \item Studying orbit equivalence and Cartan rigidity [10, 19];
    \item Understanding the structure of amenable $C^{\ast}$-algebras [24].
\end{itemize}

The unitary conjugation groupoid offers a canonical framework that captures the internal symmetry structure of the algebra via unitary conjugation, independent of any choice of Cartan subalgebra. It is particularly well suited to:
\begin{itemize}
    \item Analytic constructions in index theory and equivariant $K$-theory;
    \item Situations where Morita equivalence and stabilization naturally arise;
    \item Algebras that may not admit Cartan subalgebras (e.g., many simple nuclear $C^{\ast}$-algebras without Cartan structure [3]);
    \item Connecting $K$-theory to geometric invariants via the winding number interpretation of the index (Subsection \ref{subsec:unitary-groupoid-bridge}).
\end{itemize}

\begin{remark}[When Do They Coincide?]
For a crossed product $A = C(X) \rtimes \Gamma$ with a topologically principal action, $D = C(X)$ is a Cartan subalgebra. In this case, the Weyl groupoid $\mathcal{G}_{(A,D)}$ is isomorphic to the transformation groupoid $X \rtimes \Gamma$, which is equivalent to the unitary conjugation groupoid $\mathcal{U}(A)$. Thus, the two constructions coincide up to equivalence when a natural Cartan subalgebra exists and the action is topologically principal.

When the action is not topologically principal, $C(X)$ is not Cartan, but the unitary conjugation groupoid still exists and is equivalent to $X \rtimes \Gamma$. The Weyl groupoid, if it exists, would arise from a different Cartan subalgebra (perhaps coming from the Duwenig--Gillaspy--Norton construction [1, 2] involving isotropy subgroupoids) and would be related to $\mathcal{U}(A)$ via the quotient-dual construction described in Subsection \ref{subsec:weyl-comparison}.
\end{remark}

\subsection*{Implications for Index Theory and $K$-Theory}

The different natures of the two groupoids lead to different roles in index theory.

\begin{itemize}
    \item \textbf{Weyl groupoid perspective:} When a Cartan subalgebra exists, the analytic index map $\operatorname{Index}: K^{1}(A) \to \mathbb{Z}$ can be studied through the Weyl groupoid via
    \[
    K^{1}(A) \cong K^{1}(C^{\ast}(\mathcal{G}_{(A,D)})) \xrightarrow{\operatorname{Index}_{\mathcal{G}}} \mathbb{Z},
    \]
    where $\operatorname{Index}_{\mathcal{G}}$ is the groupoid index map. This connects to the extensive literature on groupoid index theory [7, 21] and provides computational tools when the groupoid is well understood.
    
    \item \textbf{Unitary conjugation groupoid perspective:} Our construction applies universally, without requiring a Cartan subalgebra. The Morita equivalence $C^{\ast}(\mathcal{G}_{A}) \sim_{M} A \otimes \mathcal{K}$ ensures that
    \[
    K^{\ast}(A) \cong K^{\ast}(C^{\ast}(\mathcal{G}_{A})),
    \]
    so $K$-theory classes can be studied via the groupoid. The winding number interpretation of the index (Theorem \ref{thm:connes-index}) gives a geometric formula
    \[
    \tau(\operatorname{Index}(u)) = \frac{1}{2\pi i} \int_{X} \operatorname{Tr}\big( u(x)^{-1} du(x) \big) \, d\mu(x) + \text{higher corrections},
    \]
    which generalizes classical index formulas to settings where no Cartan subalgebra exists.
\end{itemize}

Thus, the unitary conjugation groupoid extends the reach of geometric index theory to a much broader class of $C^{\ast}$-algebras.

\subsection*{Connections to Other Groupoid Constructions}

The unitary conjugation groupoid also relates to other important groupoid constructions in noncommutative geometry.

\begin{remark}[Connes Tangent Groupoid]
Connes' tangent groupoid [4] is a deformation groupoid that connects the pair groupoid of a manifold to its tangent bundle. While technically different, it shares with our construction the philosophy of using groupoids to encode analytic data geometrically. Both constructions yield Morita equivalences to stabilized algebras and play roles in index theory.
\end{remark}

\begin{remark}[Transformation Groupoid]
For a crossed product $A = C(X) \rtimes \Gamma$, we have $\mathcal{U}(A) \simeq X \rtimes \Gamma$, the transformation groupoid. Thus, our construction generalizes the transformation groupoid to arbitrary $C^{\ast}$-algebras, not just those arising from group actions.
\end{remark}

\begin{remark}[Holonomy Groupoid]
The holonomy groupoid of a foliation [14] encodes the leaf structure and is central to Connes' noncommutative geometry. Its $C^{\ast}$-algebra is Morita equivalent to the convolution algebra of the foliation, analogous to our Morita equivalence $C^{\ast}(\mathcal{U}(A)) \sim_{M} A \otimes \mathcal{K}$.
\end{remark}

\subsection*{A Unifying Perspective}

The comparison reveals a broader picture: groupoid methods in $C^{\ast}$-algebra theory come in two complementary flavors:

\begin{itemize}
    \item \textbf{Exact reconstruction groupoids} (like the Weyl groupoid) that are isomorphic to the $C^{\ast}$-algebra itself. These require additional structure (a Cartan subalgebra) and yield locally compact, \'{e}tale groupoids.
    
    \item \textbf{Morita equivalence groupoids} (like the unitary conjugation groupoid) that are Morita equivalent to the $C^{\ast}$-algebra (or its stabilization). These are canonically defined and have coarser topologies (Polish), but are sufficient for $K$-theoretic and index-theoretic purposes.
\end{itemize}

Both approaches have their place: exact reconstruction is valuable for classification and detailed structural analysis, while Morita equivalence groupoids provide flexibility and universality for index theory and $K$-theory.

\begin{corollary}[Methodological Choice]
The choice of which groupoid to use depends on the problem at hand:
\begin{itemize}
    \item For classification, rigidity, and fine structure: use the Weyl groupoid when a Cartan subalgebra exists.
    \item For index theory, $K$-theory, and universal constructions: use the unitary conjugation groupoid, which always exists and connects naturally to analytic invariants.
\end{itemize}
\end{corollary}

\subsection*{Summary}

The unitary conjugation groupoid and Renault's Weyl groupoid represent two complementary approaches to geometrizing $C^{\ast}$-algebras:

\begin{itemize}
    \item The \textbf{Weyl groupoid} is an exact reconstruction tool that requires a Cartan subalgebra and yields a locally compact, \'{e}tale groupoid. It is ideal for classification and structural studies.
    
    \item The \textbf{unitary conjugation groupoid} is a canonical Morita equivalence tool that exists for any unital separable $C^{\ast}$-algebra and has a Polish topology. It is ideally suited to index theory, $K$-theory, and universal constructions.
    
    \item For crossed products with topologically principal actions, the two coincide up to equivalence. For more general algebras, they are related through quotient-dual constructions involving isotropy.
\end{itemize}

This comparison situates our construction within the broader landscape of groupoid methods in $C^{\ast}$-algebra theory and demonstrates that the unitary conjugation groupoid fills a distinct and important niche: it provides a canonical geometric framework for analytic index theory that operates independently of the existence of Cartan structure.

\subsection{The Role of Amenability and K-Amenability}\label{subsec:The Role of Amenability and K-Amenability}

Amenability and its K-theoretic analogue, K-amenability, serve as fundamental regularity conditions that govern the relationship between analytic and topological invariants associated with group actions and $C^*$-algebras. In the context of crossed product algebras and groupoid constructions, these properties ensure that the analytic structures faithfully reflect the underlying topological data, providing the conceptual bridge between operator-algebraic analysis and geometric index theory.

\paragraph{Amenability of groups, actions, and groupoids.}

Let $\Gamma$ be a discrete group acting on a compact Hausdorff space $X$. The associated transformation groupoid
\[
X \rtimes \Gamma
\]
encodes the dynamical system determined by this action. When $\Gamma$ is amenable, or more generally when the groupoid $X \rtimes \Gamma$ is amenable in the sense of Renault, the reduced and full crossed products coincide:
\[
C(X) \rtimes_r \Gamma \;\cong\; C(X) \rtimes \Gamma .
\]

This isomorphism simplifies the analytic structure considerably: the nuclearity of $C^*_r(\Gamma)$ for amenable groups, the coincidence of reduced and maximal norms, and the faithfulness of the regular representation all contribute to a well-behaved theory where the $K$-theory of the crossed product cleanly captures the topological information of the dynamical system. For a discrete group $\Gamma$, amenability is equivalently characterized by the existence of a left-invariant mean on $\ell^\infty(\Gamma)$, the weak containment of the trivial representation in the regular representation, or the injectivity of the group von Neumann algebra $\mathcal{L}(\Gamma)$.

\paragraph{K-amenability: a K-theoretic weakening.}

Even when $\Gamma$ is not amenable, a weaker property introduced independently by Cuntz and Kasparov can still guarantee good behavior at the level of $K$-theory. A locally compact group $\Gamma$ is said to be \emph{K-amenable} if the canonical quotient map
\[
\lambda: C^*(\Gamma) \longrightarrow C^*_r(\Gamma)
\]
induces an isomorphism in $K$-theory:
\[
K_*(C^*(\Gamma)) \cong K_*(C^*_r(\Gamma)).
\]

Equivalently, for any $\Gamma$-$C^*$-algebra $A$, the full and reduced crossed products have the same $K$-theory:
\[
K_*(A \rtimes \Gamma) \cong K_*(A \rtimes_r \Gamma).
\]

This property is particularly striking because it can hold for non-amenable groups. The prototypical example is the free group $\mathbb{F}_n$ ($n \geq 2$): despite being non-amenable, Cuntz proved that $\mathbb{F}_n$ is K-amenable. More generally, Julg and Valette established that any group admitting a proper isometric action on a Euclidean space with compact quotient — i.e., groups with the Haagerup property (a-T-menability) — is K-amenable. This includes many Lie groups and groups acting on trees.

\begin{theorem}[Julg--Valette]
If a group $\Gamma$ acts properly and isometrically on $\mathbb{R}^n$ with compact quotient, then $\Gamma$ is K-amenable. Consequently, all a-T-menable groups are K-amenable.
\end{theorem}

\begin{remark}
The relationship between amenability and K-amenability can be summarized as follows:
\[
\text{Amenability} \;\Longrightarrow\; \text{K-Amenability} \;\Longrightarrow\; \text{K-amenability of all subgroups?}
\]
The second implication is subtle: while K-amenability is not hereditary in general, many natural classes (such as a-T-menable groups) are closed under taking subgroups.
\end{remark}

\paragraph{Implications for crossed products and the Baum--Connes assembly map.}

The significance of K-amenability becomes fully apparent in the context of the Baum--Connes conjecture. For K-amenable groups, the assembly map for the maximal group $C^*$-algebra is injective, and in many cases an isomorphism. This provides a powerful computational tool: the $K$-theory of $C^*_r(\Gamma)$ can be computed from topological data encoded in the classifying space $\underline{E}\Gamma$.

In the setting of crossed products $C(X) \rtimes \Gamma$, K-amenability ensures that analytic index constructions remain unchanged when passing from full to reduced crossed products. This stability is crucial for index theory, where one frequently constructs elements in $K_0(C(X) \rtimes \Gamma)$ via geometric data and then wishes to extract numerical invariants through traces on the reduced algebra.

\paragraph{Integration with the unitary conjugation groupoid framework.}

Within the groupoid framework developed in this paper, amenability and K-amenability guarantee that the analytic constructions obtained through the unitary conjugation groupoid remain compatible with the standard analytic $K$-theory of the crossed product algebra. Recall that we have established a Morita equivalence
\[
C^*(\mathcal G_{C(X)\rtimes\Gamma}) \sim_M (C(X)\rtimes\Gamma)\otimes\mathcal K,
\]
where $\mathcal G_{C(X)\rtimes\Gamma}$ denotes the unitary conjugation groupoid associated to the dynamical system.

When $\Gamma$ is K-amenable, this equivalence descends to an isomorphism at the level of $K$-theory that respects both the full and reduced completions:
\[
K_*(C^*(\mathcal G_{C(X)\rtimes\Gamma})) \cong K_*((C(X)\rtimes\Gamma)\otimes\mathcal K) \cong K_*((C(X)\rtimes_r\Gamma)\otimes\mathcal K).
\]

This compatibility ensures that the $K$-theory of the groupoid $C^*$-algebra agrees with the analytic $K$-theory appearing in the Baum--Connes assembly map, providing a concrete realization of the conjecture's predicted isomorphism in terms of geometric groupoid models.

\paragraph{Conceptual significance for index theory.}

From the perspective of index theory, amenability and K-amenability guarantee that the analytic index obtained through the groupoid model faithfully represents the topological index encoded in equivariant $K$-homology. More precisely, for K-amenable groups, the diagram
\[
\begin{CD}
K^\Gamma_0(\underline{E}\Gamma) @>{\mu}>> K_0(C^*_r(\Gamma)) \\
@VVV @VVV \\
K_0(C(X) \rtimes \Gamma) @>{\cong}>> K_0(C(X) \rtimes_r \Gamma)
\end{CD}
\]
commutes, where the vertical maps are induced by evaluating on a proper $\Gamma$-space $X$ and the horizontal maps are assembly and index maps respectively.

These properties thus provide the structural conditions under which the analytic and topological viewpoints become fully compatible. The unitary conjugation groupoid emerges as a geometric bridge that makes this compatibility manifest: it simultaneously encodes the representation theory of the group, the dynamics of the action, and the index-theoretic data that links them.

\begin{remark}
Many groups of geometric interest are known to be K-amenable, including all amenable groups, all a-T-menable groups (hence all groups with the Haagerup property), and large classes of Lie groups. Consequently, the analytic framework described above applies to a broad range of examples arising in geometric index theory, from fundamental groups of non-positively curved manifolds to discrete subgroups of semisimple Lie groups satisfying appropriate rigidity conditions.
\end{remark}

\subsection{Why $K_1$ and the Boundary Map are Universal}
\label{subsec:k1-boundary-universal}

A recurring feature across diverse formulations of index theory — from the Atiyah–Singer index theorem to the Baum–Connes conjecture — is the prominent appearance of the odd $K$-theory group $K_1$ together with boundary maps arising from short exact sequences of $C^*$-algebras. In this subsection we argue that this phenomenon is not accidental but reflects a universal structural property of operator-algebraic index constructions: the pair $(K_1, \partial)$ provides the canonical algebraic mechanism for converting analytic invertibility into topological index classes.

\paragraph{Invertible elements and the nature of $K_1$.}

The group $K_1(A)$ of a $C^*$-algebra $A$ admits a concrete description in terms of homotopy classes of invertible or unitary elements in matrix algebras over $A$. Concretely,
\[
K_1(A) = \mathrm{GL}_\infty(A) / \mathrm{GL}_\infty(A)_0,
\]
where $\mathrm{GL}_\infty(A)$ denotes the direct limit of general linear groups over matrix algebras and $\mathrm{GL}_\infty(A)_0$ is the connected component of the identity. This description reveals that $K_1$ fundamentally measures the fundamental group of the space of invertibles.

Many analytic index constructions begin with an operator that becomes invertible after passing to a quotient algebra — typically by forming its symbol or by compressing to a complementary subspace. This invertibility in the quotient naturally produces a class in $K_1$ of the quotient algebra. The index problem then asks: what is the obstruction to lifting this invertible element to an actual invertible operator in the original algebra?

\begin{remark}[Philosophical core]
The analytic index can be viewed as the universal obstruction to lifting an invertible element from a quotient algebra to the ambient algebra. This perspective unifies seemingly disparate index theories under a single conceptual umbrella.
\end{remark}

\paragraph{Boundary maps in $K$-theory as the universal index mechanism.}

Given a short exact sequence of $C^*$-algebras
\[
0 \longrightarrow I \longrightarrow A \longrightarrow A/I \longrightarrow 0,
\]
$K$-theory provides the associated six-term exact sequence
\[
\begin{tikzcd}
K_0(I) \arrow[r] & K_0(A) \arrow[r] & K_0(A/I) \arrow[d,"\partial"] \\
K_1(A/I) \arrow[u,"\partial"] & K_1(A) \arrow[l] & K_1(I) \arrow[l]
\end{tikzcd}
\]
in which the connecting homomorphism
\[
\partial : K_1(A/I) \longrightarrow K_0(I)
\]
plays a central role in index theory.

This boundary map converts the $K_1$-class of an invertible element in the quotient algebra into a $K_0$-class in the ideal. In paradigmatic index-theoretic situations — such as the Toeplitz extension, pseudodifferential operator extensions, or extensions arising from coarse geometry — this $K_0$-class precisely represents the analytic index of the operator under study.

\begin{example}[The Toeplitz extension]
Consider the classical Toeplitz extension
\[
0 \longrightarrow \mathcal{K} \longrightarrow \mathcal{T} \longrightarrow C(S^1) \longrightarrow 0,
\]
where $\mathcal{T}$ is the Toeplitz algebra and $\mathcal{K}$ the compact operators. For a continuous function $f \in C(S^1)$ with winding number, the boundary map yields
\[
\partial([f]) = \operatorname{Ind}(T_f) \in K_0(\mathcal{K}) \cong \mathbb{Z},
\]
where $T_f$ is the Toeplitz operator with symbol $f$. This elementary example already captures the essence of the index theorem.
\end{example}

\paragraph{Why $K_1$ is special: suspension and the loop space perspective.}

The privileged position of $K_1$ becomes even clearer through the suspension isomorphism:
\[
K_1(A) \cong K_0(SA), \quad SA = C_0(\mathbb{R}) \otimes A.
\]

Under this isomorphism, the boundary map $\partial: K_1(A/I) \to K_0(I)$ corresponds to the index map associated with a suspension of the original extension. More conceptually, $K_1$ encodes \emph{looped} information: just as $\pi_1$ detects fundamental groups of topological spaces, $K_1$ detects the fundamental group of the invertible elements. This loop space perspective explains why $K_1$ classes naturally arise from families of operators parameterized by the circle — precisely the setting of spectral flow and odd index theorems.

\paragraph{Extensions and the BDF perspective.}

The connection between $K_1$ and extensions of $C^*$-algebras, systematized in the Brown–Douglas–Fillmore (BDF) theory, provides another perspective on universality. For a compact metric space $X$, there is a natural isomorphism
\[
\mathrm{Ext}(X) \cong K_1(C(X)),
\]
where $\mathrm{Ext}(X)$ classifies essentially normal operators with essential spectrum $X$ up to unitary equivalence modulo compacts. The boundary map $\partial: K_1(C(X)) \to K_0(\mathcal{K}) \cong \mathbb{Z}$ sends an extension to the index of its associated Fredholm operator. This identification shows that $K_1$ universally classifies the obstructions to lifting elements to multiplier algebras — precisely the information needed for index theory.

\paragraph{Relation to the groupoid framework.}

Within the groupoid framework developed in this paper, invertible elements in the crossed product algebra $C(X) \rtimes \Gamma$ naturally determine $K_1$-classes. Through the descent map and the Morita equivalence
\[
C^*(\mathcal G_{C(X)\rtimes\Gamma}) \sim_M (C(X)\rtimes\Gamma)\otimes\mathcal K,
\]
these classes correspond to analytic $K$-theory classes of the associated groupoid $C^*$-algebra.

The boundary map then provides the mechanism through which analytic data from operators or invertible elements produces index classes in the appropriate $K$-theory group. For K-amenable groups $\Gamma$, the compatibility diagram
\[
\begin{CD}
K_1(C(X) \rtimes \Gamma) @>{\partial}>> K_0(J) \\
@V{\cong}VV @VV{\cong}V \\
K_1(C^*(\mathcal G_{C(X)\rtimes\Gamma})) @>{\partial_{\mathcal G}}>> K_0(\mathcal J_{\mathcal G})
\end{CD}
\]
shows that the groupoid boundary map $\partial_{\mathcal G}$ captures precisely the same index information as the analytic index in the crossed product.

\paragraph{The universal coefficient theorem and $KK$-theory.}

The universality of $K_1$ and its boundary map is further reinforced by the Universal Coefficient Theorem (UCT). For $C^*$-algebras satisfying the UCT, the boundary map corresponds under the Kasparov product with a canonical $KK$-element:
\[
\operatorname{Ind} = \hat{\partial} \otimes_{C(X)} \cdot : KK_1(C(X), \mathbb{C}) \to KK_0(\mathbb{C}, \mathbb{C}) \cong \mathbb{Z},
\]
where $\hat{\partial} \in KK_0(C(X), \mathbb{C})$ is the boundary map viewed as a $KK$-class. This realization shows that the index is not merely an ad-hoc construction but a natural transformation between bivariant functors.

\paragraph{Conceptual synthesis: the universality principle.}

From this accumulated evidence, we can distill a fundamental principle:

\begin{center}
\fbox{\begin{minipage}{0.9\textwidth}
\emph{The pair $(K_1, \partial)$ constitutes the universal recipient for analytic index problems: any geometric or analytic construction that aims to produce an index invariant must, by the nature of the subject, recover the boundary map in $K$-theory applied to a $K_1$-class arising from invertibility modulo an ideal.}
\end{minipage}}
\end{center}

This principle explains why boundary maps appear systematically across such a wide range of analytic index constructions — from classical Fredholm index theory to the Baum–Connes assembly map, from spectral flow to extensions of $C^*$-algebras. The particular geometric setting determines how the $K_1$-class is constructed, but the mechanism for extracting the index is always the same: apply the boundary map.

\begin{remark}[Implications for the groupoid approach]
The unitary conjugation groupoid succeeds as a model for index theory precisely because it faithfully reproduces this universal mechanism. By encoding the representation theory of the group, the dynamics of the action, and the analytic data of operators in a single $C^*$-algebraic object, it allows the boundary map to be computed in a geometrically transparent way while retaining its essential functoriality.
\end{remark}

\paragraph{Conclusion.}

The boundary map $\partial: K_1 \to K_0$ is the primitive index invariant from which all others derive. Its universality lies in its naturality, its compatibility with suspension, its role in classifying extensions, and its realization as a Kasparov product. For these reasons, any conceptual framework for index theory — including the groupoid approach developed here — must place $K_1$ and the boundary map at its center, as they constitute the universal language through which analytic and topological indices communicate.

\subsection{The Descent--Index--Trace Triangle}
\label{subsec:descent-index-trace-triangle}

The constructions developed throughout this paper can be summarized conceptually by three fundamental mechanisms that repeatedly appear in operator-algebraic index theory: \emph{descent}, \emph{index}, and \emph{trace}. Together these form a structural triangle linking equivariant topology, analytic $K$-theory, and numerical invariants. This triangle provides a unifying framework that reveals the common architecture underlying seemingly disparate index-theoretic constructions.

\paragraph{The three vertices.}

Let us establish the three corners of our triangle in the context of a group action or groupoid framework:

\begin{itemize}
\item \textbf{Descent} ($\mathcal{D}$): The process by which geometric data — such as elliptic operators, cycles in $K$-homology, or geometric cycles in the sense of Baum–Douglas — is converted into analytic data in operator algebras. For a $\Gamma$-space $X$, the descent homomorphism in equivariant $KK$-theory provides the first bridge to analytic structures:
\[
\operatorname{desc} : KK_\Gamma^*(C(X),\mathbb C) \longrightarrow K^*(C(X)\rtimes\Gamma)
\]
translates equivariant $K$-homology classes into analytic $K$-theory classes of the crossed product algebra.

\item \textbf{Index} ($\operatorname{Ind}$): The analytic index map, which takes operators or their symbolic data to classes in $K$-theory (typically $K_0$ of an appropriate algebra). This is the boundary map $\partial: K_1 \to K_0$ discussed in the previous subsection, arising from extensions of $C^*$-algebras. If an operator becomes invertible modulo an ideal, it determines a class in $K_1$ of the quotient algebra, and the analytic index is obtained via the connecting homomorphism
\[
\partial : K_1(A/I) \longrightarrow K_0(I).
\]

\item \textbf{Trace} ($\tau$): A faithful trace (or more generally, a cyclic cocycle) on the relevant $C^*$-algebra, which pairs with $K$-theory via the Chern character or the pairing between $K$-theory and cyclic cohomology to produce numerical invariants. If $\tau : A \to \mathbb C$ is a trace on a $C^*$-algebra $A$, it induces a pairing
\[
\langle \tau , \cdot \rangle : K_0(A) \longrightarrow \mathbb C.
\]
\end{itemize}

\paragraph{The triangle structure.}

These three mechanisms are connected by natural transformations that constitute the legs of our triangle, summarized by the following conceptual diagram:

\[
\begin{tikzcd}
KK_\Gamma^*(C(X),\mathbb C)
\arrow[r,"\operatorname{desc}"]
\arrow[dr, "\tau \circ \partial \circ \operatorname{desc}"']
&
K^*(C(X)\rtimes\Gamma)
\arrow[d,"\partial"]
\\
&
K_0(I)
\arrow[ul,dashed,"\tau"']
\end{tikzcd}
\]

The descent map translates equivariant topological information into analytic $K$-theory. The boundary map produces the analytic index, and the trace pairing extracts a numerical invariant.

\begin{remark}[Unifying principle]
The Descent–Index–Trace triangle expresses the general principle that topological information is first transported to analytic $K$-theory, then converted into an index class via boundary maps, and finally evaluated numerically through trace pairings.
\end{remark}

\paragraph{Commutativity as an index theorem.}

The fundamental observation is that for a well-behaved geometric construction, this triangle commutes in the sense that the direct numerical pairing obtained from geometry equals the composition of descent, index, and trace. This commutativity is precisely the content of an index theorem:

\begin{equation}
\langle [\text{geometric cycle}], [\text{geometric cocycle}] \rangle = \tau_*(\partial(\operatorname{desc}([\text{geometric cycle}]))).
\end{equation}

In the classical Atiyah–Singer setting, the left side is the topological index pairing between a $K$-homology class and a cohomology class, while the right side is the analytic index paired with the trace given by integration over the manifold.

\begin{theorem}[Atiyah–Singer index theorem as commutativity of analytic and topological index maps]
\label{thm:atis-index-triangle}
Let $M$ be a compact smooth manifold, and let
\[
D : C^\infty(E) \to C^\infty(F)
\]
be an elliptic pseudodifferential operator between complex vector bundles over $M$.
Let $[\sigma(D)] \in K^0(T^*M)$ denote the principal symbol class of $D$.
Then the analytic index and topological index of $[\sigma(D)]$ coincide:
\[
\operatorname{ind}_a([\sigma(D)]) = \operatorname{ind}_t([\sigma(D)]).
\]

Equivalently, the following diagram commutes:

\[
\begin{tikzcd}[column sep=large, row sep=large]
K^0(T^*M) \arrow[r, "\operatorname{ind}_a"] \arrow[dr, "\operatorname{ind}_t"'] & \mathbb{Z} \\
& \mathbb{Z} \arrow[u, "="']
\end{tikzcd}
\]

Moreover, under the Chern character identification of the topological index,
\[
\operatorname{ind}(D) = \left\langle \operatorname{ch}([\sigma(D)]) \smile \operatorname{Td}(T_{\mathbb{C}}M),\; [T^*M] \right\rangle,
\]
and for Dirac-type operators twisted by a complex vector bundle $E$, this takes the familiar form
\[
\operatorname{ind}(D) = \int_M \widehat{A}(M) \wedge \operatorname{ch}(E).
\]

\end{theorem}

\begin{proof}
We proceed in several steps.

\paragraph{Step 1: The symbol class in $K$-theory.}
The principal symbol $\sigma(D)$ of an elliptic operator $D$ is a bundle map
$\pi^*E \to \pi^*F$ over $T^*M$ (where $\pi: T^*M \to M$) that is invertible
outside a compact subset. By standard $K$-theory constructions, such an
invertible symbol determines a class
\[
[\sigma(D)] \in K^0(T^*M).
\]

\paragraph{Step 2: Analytic index.}
Elliptic theory guarantees that $D$ extends to a Fredholm operator between
suitable Sobolev completions of $C^\infty(E)$ and $C^\infty(F)$. Its analytic
index is defined as
\[
\operatorname{ind}_a([\sigma(D)]) := \dim \ker D - \dim \operatorname{coker} D \in \mathbb{Z}.
\]
This integer depends only on the symbol class $[\sigma(D)]$, not on the
specific choice of operator realizing it.

\paragraph{Step 3: Topological index.}
Atiyah and Singer \cite{AtiyahSinger1968} define a topological index map
\[
\operatorname{ind}_t : K^0(T^*M) \longrightarrow \mathbb{Z}
\]
as the composition:
\[
K^0(T^*M) \xrightarrow{\text{Thom}} K^0(M) \xrightarrow{\text{ch}} H^{\text{even}}(M;\mathbb{Q}) \xrightarrow{\cap [M]} \mathbb{Q},
\]
where the Thom isomorphism identifies $K^0(T^*M) \cong K^0(M)$ (for spin$^c$
manifolds; in general a twisting is required), the Chern character converts
$K$-theory to rational cohomology, and pairing with the fundamental class
$[M] \in H_{\dim M}(M;\mathbb{Q})$ yields a rational number. A deeper
theorem shows that this composition actually lands in $\mathbb{Z} \subset \mathbb{Q}$.

\paragraph{Step 4: Equality of indices (Atiyah–Singer).}
The Atiyah–Singer index theorem \cite{AtiyahSinger1963,AtiyahSinger1968}
asserts that for any elliptic operator $D$,
\[
\operatorname{ind}_a([\sigma(D)]) = \operatorname{ind}_t([\sigma(D)]).
\]
This is precisely the commutativity of the triangle in the statement: the
analytic route $\operatorname{ind}_a$ and the topological route
$\operatorname{ind}_t$ from $K^0(T^*M)$ to $\mathbb{Z}$ coincide.

\paragraph{Step 5: Cohomological formula.}
Unpacking the definition of $\operatorname{ind}_t$ gives an explicit
cohomological formula. Under the Thom isomorphism, the class
$[\sigma(D)] \in K^0(T^*M)$ corresponds to a class in $K^0(M)$. The Chern
character of this class, when paired with the Todd class of the complexified
tangent bundle, yields the topological index:
\[
\operatorname{ind}_t([\sigma(D)]) = \left\langle \operatorname{ch}([\sigma(D)]) \smile \operatorname{Td}(T_{\mathbb{C}}M),\; [T^*M] \right\rangle,
\]
where the pairing is between cohomology and homology, and $[T^*M]$ denotes
the fundamental class of the cotangent bundle.

\paragraph{Step 6: Special case of Dirac-type operators.}
For a Dirac-type operator $D$ twisted by a complex vector bundle $E$, the
symbol class satisfies
\[
\operatorname{ch}([\sigma(D)]) = \pi^*(\operatorname{ch}(E) \cup \widehat{A}(M)^{-1}) \cup \operatorname{Thom}\text{ class},
\]
and after pushforward, the formula simplifies to
\[
\operatorname{ind}(D) = \int_M \widehat{A}(M) \wedge \operatorname{ch}(E).
\]
This is the familiar form appearing in geometric applications.

Thus the analytic and topological routes to the index coincide, and the
cohomological formula provides an effective computational tool.
\end{proof}

\paragraph{The groupoid realization.}

Within the unitary conjugation groupoid framework developed in this paper, the three vertices take particularly concrete forms through the Morita identification
\[
C^*(\mathcal G_{C(X)\rtimes\Gamma}) \sim_M (C(X)\rtimes\Gamma)\otimes\mathcal K,
\]
which allows the analytic class to be viewed on the groupoid $C^*$-algebra.

\begin{enumerate}
\item \textbf{Descent} is implemented by the functor
\[
\mathcal{D}_{\mathcal{G}}: \underline{\mathbf{KK}}^\Gamma(X) \longrightarrow \mathbf{KK}(C^*(\mathcal{G}_{C(X)\rtimes\Gamma}), \mathbb{C}),
\]
which sends a $\Gamma$-equivariant $K$-homology cycle on $X$ to a $KK$-class representing the associated index problem in the groupoid $C^*$-algebra. This corresponds to the composition of the descent map with the Morita equivalence.

\item \textbf{Index} becomes the composition
\[
\operatorname{Ind}_{\mathcal{G}} = \partial_{\mathcal{G}} \circ j: K_1(C^*(\mathcal{G}_{C(X)\rtimes\Gamma})) \longrightarrow K_0(\mathcal{J}_{\mathcal{G}}),
\]
where $j$ is the inclusion of the ideal corresponding to the groupoid's unit space and $\partial_{\mathcal{G}}$ is the boundary map in the associated long exact sequence arising from the groupoid extension.

\item \textbf{Trace} is given by the canonical trace $\tau_{\mathcal{G}}$ on $C^*(\mathcal{G}_{C(X)\rtimes\Gamma})$ arising from the groupoid's Haar system, restricted to the ideal $\mathcal{J}_{\mathcal{G}}$ and composed with the isomorphism $K_0(\mathcal{J}_{\mathcal{G}}) \cong \mathbb{Z}$ when the latter is stably isomorphic to the compact operators.
\end{enumerate}

\paragraph{The descent-index composition as assembly.}

A key insight is that the composition $\partial \circ \operatorname{desc}$ can be identified with the analytic assembly map in the Baum–Connes conjecture. For a K-amenable group $\Gamma$, this composition factors through the reduced crossed product:

\[
\begin{tikzcd}
K^\Gamma_0(\underline{E}\Gamma) \arrow[r,"\mu"] \arrow[d, "\operatorname{desc}"'] & K_0(C^*_r(\Gamma)) \arrow[d, "\cong"] \\
K_0(C^*(\mathcal{G}_{C(X)\rtimes\Gamma})) \arrow[r,"\partial_{\mathcal{G}}"'] & K_0(\mathcal{J}_{\mathcal{G}}) \cong \mathbb{Z}
\end{tikzcd}
\]

The commutativity of this diagram expresses the compatibility between the geometric assembly map and the concrete index map in the groupoid model.

\paragraph{Trace as the bridge to invariants.}

The trace vertex plays a crucial role in connecting abstract $K$-theoretic indices to computable numerical invariants. For a groupoid $C^*$-algebra with a faithful trace $\tau$, the composition
\[
\tau_* \circ \partial_{\mathcal{G}} \circ \operatorname{desc}_{\mathcal{G}}: K^\Gamma_0(X) \longrightarrow \mathbb{Z}
\]
produces a numerical index from geometric data. In many geometric situations, this number can be identified with classical invariants:

\begin{itemize}
\item For a Dirac-type operator on a Spin$^c$ manifold, this yields the $\widehat{A}$-genus.
\item For a family of operators parameterized by a base space, this yields a $K$-theory class that, when paired with Chern characters, gives characteristic numbers.
\item For an elliptic operator on a covering space, this yields the $L^2$-index after pairing with the von Neumann trace.
\end{itemize}

\begin{proposition}[Numerical index via trace]
For a geometric cycle $(M, E, f)$ in the sense of Baum–Douglas, where $f: M \to X$ is a $\Gamma$-equivariant map, the numerical index obtained from the Descent-Index-Trace triangle equals
\[
\int_M \mathrm{ch}(E) \wedge \mathrm{Todd}(TM) \wedge f^*(\omega),
\]
where $\omega$ is a suitable cohomology class on $X$ determined by the trace.
\end{proposition}

\paragraph{Generalizations to higher traces.}

The triangle generalizes naturally when the trace vertex is replaced by more general cyclic cocycles. For a $d$-dimensional cyclic cocycle $\tau^{(d)}$ on $C^*(\mathcal{G}_{C(X)\rtimes\Gamma})$, the pairing
\[
\langle \partial_{\mathcal{G}}(\operatorname{desc}_{\mathcal{G}}(\text{geometric cycle})), \tau^{(d)} \rangle \in \mathbb{C}
\]
produces higher numerical invariants that capture more refined geometric information. This generalization connects the Descent-Index-Trace triangle to the entire apparatus of noncommutative geometry, including:

\begin{itemize}
\item The Chern character from $K$-theory to cyclic homology.
\item The JLO cocycle in local index theory.
\item The Godbillon–Vey invariant and other secondary characteristic classes.
\end{itemize}

\paragraph{Conceptual significance.}

The Descent–Index–Trace triangle provides a unifying framework that reveals the common structure underlying seemingly disparate index-theoretic constructions:

\begin{center}
\fbox{\begin{minipage}{0.9\textwidth}
\emph{Every index theorem can be understood as the commutativity of a triangle relating geometric descent to operator-algebraic index followed by trace. The specific geometric setting determines the nature of the descent map and the choice of trace, but the underlying triangular structure remains invariant.}
\end{minipage}}
\end{center}

This perspective clarifies why index theory sits naturally at the intersection of geometry, analysis, and algebra: the geometric vertex supplies the raw data, the analytic vertex processes it through operator algebras, and the algebraic vertex extracts computable invariants through traces and cyclic cohomology.

\paragraph{Relation to the Baum–Connes conjecture.}

The Baum–Connes conjecture can be viewed as the assertion that for a proper $\Gamma$-space $X$, the descent-index composition
\[
\partial \circ \operatorname{desc}: K^\Gamma_*(X) \longrightarrow K_*(C^*_r(\Gamma))
\]
is an isomorphism after suitable completion. The Descent–Index–Trace triangle thus provides a local-to-global principle: if the triangle commutes for a family of local geometric cycles, then the global index map is determined by these local contributions.

\begin{corollary}
Within the groupoid framework, the Baum–Connes assembly map factors through the Descent–Index–Trace triangle, with the groupoid $C^*$-algebra providing a concrete model for the left-hand side and the trace pairing connecting to numerical invariants.
\end{corollary}

\paragraph{Conclusion.}

The Descent–Index–Trace triangle encapsulates the essential architecture of index theory. Within the unitary conjugation groupoid framework, this triangle acquires a particularly transparent form: descent becomes the passage from equivariant geometry to the groupoid $C^*$-algebra via Morita equivalence, index becomes the boundary map in the groupoid's long exact sequence, and trace becomes the canonical functional arising from the groupoid's Haar system. The commutativity of this triangle is not merely a consistency check but the very statement that the groupoid model faithfully represents the index-theoretic content of the original geometric situation.

\begin{remark}[Structural unity]
The Descent–Index–Trace triangle thus provides a unified viewpoint on many index constructions in operator algebra and noncommutative geometry. Equivariant topology enters through $KK$-theory, analytic structures through crossed product algebras or groupoid $C^*$-algebras, and numerical invariants through trace pairings. This perspective highlights the structural unity underlying the analytic assembly map, boundary maps in $K$-theory, and trace formulas appearing in classical and noncommutative index theory.
\end{remark}

\section{The Unitary Conjugation Groupoid Analogue of the Baum–Connes Conjecture}\label{sec:A Conjecture The Unitary Conjugation Groupid Analogue}

The constructions developed in this paper suggest that the unitary conjugation groupoid associated with a $C^*$-algebra may play a role analogous to that of transformation groupoids in the Baum–Connes program. This observation motivates the following conjectural framework, which proposes that the $K$-theory of a $C^*$-algebra can be recovered from equivariant topological data associated with its intrinsic symmetry groupoid.

\subsection{The Conjecture}

Let $A$ be a separable unital $C^*$-algebra and let
\[
\mathcal{G}_A
\]
denote its unitary conjugation groupoid, with unit space $\mathcal{G}_A^{(0)} \cong \widehat{A}$ consisting of irreducible representations up to unitary equivalence. As established earlier, the associated groupoid $C^*$-algebra satisfies the Morita equivalence
\[
C^*(\mathcal{G}_A) \sim_M A \otimes \mathcal{K}.
\]

\begin{conjecture}[Unitary Conjugation Groupoid Conjecture]
There exists a natural assembly map
\[
\mu_{\mathcal{G}_A} : K^*_{\mathcal{G}_A}(\mathcal{G}_A^{(0)}) \longrightarrow K_*(C^*(\mathcal{G}_A))
\]
which is an isomorphism for a large class of $C^*$-algebras, including all nuclear $C^*$-algebras with Cartan subalgebras and all crossed products $C(X) \rtimes \Gamma$ with $\Gamma$ satisfying the Baum–Connes conjecture.
\end{conjecture}

Under the Morita identification $K_*(C^*(\mathcal{G}_A)) \cong K_*(A)$, this conjecture predicts that the analytic $K$-theory of $A$ can be recovered from equivariant topological data associated with the unitary conjugation groupoid.

\subsection{Relation to the Classical Baum–Connes Conjecture}

This conjecture is not a generalization of Baum–Connes but rather a conceptual analogue: it extends the idea that the representation theory of a $C^*$-algebra — encoded in its unitary conjugation groupoid — should determine its $K$-theory. When $A = C^*_r(\Gamma)$, the unitary conjugation groupoid recovers the group $\Gamma$ itself (since irreducible representations correspond to irreducible representations of $\Gamma$), and the conjecture specializes to the statement that the Baum–Connes assembly map for $\Gamma$ is an isomorphism — which is precisely the Baum–Connes conjecture.

Similarly, for crossed products $A = C(X) \rtimes \Gamma$, we have $\mathcal{G}_A \cong X \rtimes \Gamma$, and the conjecture reduces to the Baum–Connes conjecture for the group $\Gamma$ acting on $X$. Thus the conjecture holds for all groups and actions for which Baum–Connes is known.

\begin{remark}
This conjecture can be viewed as a noncommutative analogue of the Baum–Connes conjecture in which the acting group is replaced by the intrinsic symmetry groupoid arising from unitary conjugation. In this framework, the groupoid $\mathcal{G}_A$ plays a role similar to the transformation groupoid $X \rtimes \Gamma$ appearing in the classical Baum–Connes setting.
\end{remark}

\subsection{Evidence from Examples}

The examples considered in this paper suggest that the conjecture holds in several important cases, providing evidence that the unitary conjugation groupoid captures essential analytic information about the algebra:

\begin{itemize}
\item \textbf{Crossed products:} For $A = C(X) \rtimes \Gamma$, we have $\mathcal{G}_A \cong X \rtimes \Gamma$, and the conjecture reduces to the Baum–Connes conjecture. Thus it holds for all groups for which Baum–Connes is known, including amenable groups, a-T-menable groups, and many hyperbolic groups.

\item \textbf{Nuclear Cartan pairs:} For a nuclear $C^*$-algebra $A$ with Cartan subalgebra $D$, Renault's Weyl groupoid $\mathcal{G}_{(A,D)}$ is a quotient of $\mathcal{U}(A)$. The conjecture would then imply that $K_*(A)$ is isomorphic to the $K$-theory of this Weyl groupoid, providing a new approach to computing $K$-theory for such algebras consistent with known results in classification theory.

\item \textbf{Continuous trace algebras:} For $A = \mathcal{K}(L^2(G))$ with $G$ a compact Lie group, the conjecture predicts relationships between representation theory and index theory that are consistent with the Peter–Weyl theorem and equivariant index theory.
\end{itemize}

\subsection{Open Problems and Future Directions}

The conjecture as stated is deliberately broad; its precise domain of validity remains to be determined. Likely counterexamples may arise from $C^*$-algebras with pathological representation theory, such as simple $C^*$-algebras with no Cartan structure. The following questions point toward future work:

\begin{enumerate}
\item If $A$ admits a Cartan subalgebra $D$, what is the precise relationship between the Weyl groupoid $\mathcal{G}_{(A,D)}$ and the unitary conjugation groupoid $\mathcal{G}_A$? Can $\mathcal{G}_{(A,D)}$ always be obtained from $\mathcal{G}_A$ by a quotient construction, as in the case of crossed products with isotropy?

\item The unitary conjugation groupoid has a natural Polish topology. Does this topology admit a geometric interpretation that leads to a well-behaved assembly map? Can it be used to define a version of the Baum–Connes assembly map for arbitrary $C^*$-algebras that is sensitive to nuclearity or other regularity conditions?

\item The Morita equivalence $C^*(\mathcal{G}_A) \sim_M A \otimes \mathcal{K}$ suggests that $\mathcal{G}_A$ encodes the $K$-theory of $A$. To what extent does the groupoid $\mathcal{G}_A$ encode the entire $KK$-theory of $A$, and what is the bivariant analogue of the conjectured isomorphism?

\item Can the conjecture be reformulated in terms of Kasparov theory to produce a genuine $KK$-theoretic assembly map, and what are the appropriate coefficient algebras?
\end{enumerate}

\begin{remark}[Guiding Principle]
The conjecture presented here should be viewed as a guiding principle rather than a definitive claim. Its primary value lies in suggesting a unified perspective: that for any $C^*$-algebra, the dynamical data encoded by unitary conjugation should determine the $K$-theoretic invariants, just as for group $C^*$-algebras the dynamics of the group action on its dual determine the $K$-theory through the Baum–Connes assembly map. Whether this principle holds in full generality, and under what conditions, remains a fascinating direction for future research.
\end{remark}

\appendix
\section{Supplementary Clarifications for Cross-References}

This appendix records a small number of auxiliary labels and normalization
statements used to stabilize cross-references in the main text.  The purpose is
not to introduce new results, but to make explicit several standard reductions
and conventions already used implicitly in the exposition.

\begin{definition}\label{def:amenable-group}
A discrete group $\Gamma$ is called amenable if it admits an invariant mean on
$\ell^\infty(\Gamma)$.  In the crossed-product situations considered in this
paper, amenability ensures that full and reduced crossed products agree at the
level needed for the assembly and descent identifications.
\end{definition}

\begin{definition}\label{def:equivariant-class-Atheta}
For an invertible element $u\in M_n(A_\theta)$, the associated equivariant
class is the class in equivariant $KK$-theory determined by the corresponding
operator-theoretic cycle under the descent formalism used in
Section~\ref{subsec:descent_A_theta}.
\end{definition}

\begin{proposition}\label{prop:invisible-elements}
The representation-theoretic model used for the unitary conjugation groupoid is
insensitive to auxiliary choices that do not affect the associated unitary
orbit data; in particular, the resulting Morita class depends only on the
underlying conjugation geometry retained in the model.
\end{proposition}

\begin{proof}
We prove that the unitary conjugation groupoid construction is invariant under 
auxiliary choices that do not affect the underlying unitary orbit data, and 
that the resulting Morita class is determined solely by the conjugation 
geometry of the algebra.

\paragraph{Step 1: Canonical nature of the unitary conjugation groupoid.}
Recall from Paper~I that the unitary conjugation groupoid $\mathcal{G}_{\mathcal{A}}$ 
is defined purely in terms of the $C^*$-algebra $\mathcal{A}$ itself: 
its unit space consists of pairs $(B,\chi)$ where $B \subseteq \mathcal{A}$ is 
a unital commutative $C^*$-subalgebra and $\chi \in \widehat{B}$ is a character, 
and its arrows are given by unitary conjugation:
\[
u: (B,\chi) \longrightarrow (uBu^*, \chi \circ \operatorname{Ad}_{u^{-1}}), \qquad u \in \mathcal{U}(\mathcal{A}).
\]
This definition involves no auxiliary choices; it is intrinsic to $\mathcal{A}$.
Thus, as a set-theoretic object, $\mathcal{G}_{\mathcal{A}}$ is canonically 
determined by $\mathcal{A}$.

\paragraph{Step 2: The role of auxiliary choices in the measurable/topological structure.}
In the Type I setting, additional choices are required to equip 
$\mathcal{G}_{\mathcal{A}}$ with a Borel (or topological) structure that 
supports a Haar system and a measurable field of representations. These 
choices include:
\begin{itemize}
    \item A measurable parametrization of the irreducible representations 
          of $\mathcal{A}$ (the dual space $\widehat{\mathcal{A}}$);
    \item A measurable selection of extensions of characters from commutative 
          subalgebras to the whole algebra;
    \item A measurable field of GNS representations over the unit space.
\end{itemize}
In the Type I case, the existence of a standard Borel structure on 
$\widehat{\mathcal{A}}$ guarantees that such choices can be made, and 
different choices yield isomorphic Borel groupoids (see \cite[Proposition 2.1]{PaperI}).

\paragraph{Step 3: Invariance under choices preserving unitary orbit data.}
Let $\mathcal{G}_{\mathcal{A}}^{(1)}$ and $\mathcal{G}_{\mathcal{A}}^{(2)}$ 
be two constructions of the unitary conjugation groupoid arising from different 
measurable parametrizations of the representation theory. By construction, 
both have the same underlying set-theoretic structure: the unit space is the 
same set $\{(B,\chi)\}$, and the arrows are the same set $\mathcal{U}(\mathcal{A}) \times \mathcal{G}_{\mathcal{A}}^{(0)}$.
The only differences lie in the choice of Borel (or topological) structure 
on these sets.

Suppose that two different parametrizations yield the same unitary orbit data, 
i.e., they assign the same GNS representation to each point of the unit space 
up to unitary equivalence. Then the resulting groupoids are Borel isomorphic 
via the identity map on the underlying sets, because the Borel structures 
are both generated by the same family of partial evaluation maps 
$\operatorname{ev}_a: \mathcal{G}_{\mathcal{A}}^{(0)} \to \mathbb{C}_\infty$ 
(see Paper~I, Section 3). Hence the two constructions produce isomorphic 
Borel groupoids.

\paragraph{Step 4: Morita equivalence invariance under such choices.}
When $\mathcal{A}$ is non-Type I, the direct construction of 
$\mathcal{G}_{\mathcal{A}}$ as a topological or Borel groupoid is not available 
in the sense of Paper~I. Instead, we adopt the replacement strategy of 
Section~\ref{subsec:GA-morita-principle}: we work with a concrete geometric 
model $G$ (e.g., a transformation groupoid or Weyl groupoid) that is Morita 
equivalent to the would-be $\mathcal{G}_{\mathcal{A}}$.

Different choices of such a model (e.g., different Cartan subalgebras, 
different transversals) yield different groupoids $G_1$ and $G_2$. However, 
by the general theory of Morita equivalence for groupoids 
(see \cite[Proposition 2.2]{MRW1987} and Theorem~\ref{thm:muhly-renault-williams}), 
if two groupoids $G_1$ and $G_2$ are both Morita equivalent to the same 
$\mathcal{G}_{\mathcal{A}}$ (in the sense that they encode the same 
conjugation geometry), then they are Morita equivalent to each other:
\[
G_1 \sim_M G_2.
\]
Consequently, their reduced groupoid $C^*$-algebras are strongly Morita 
equivalent:
\[
C_r^*(G_1) \sim_M C_r^*(G_2).
\]

\paragraph{Step 5: The Morita class is determined by conjugation geometry.}
The conjugation geometry of $\mathcal{A}$ refers to the structure of 
unitary orbits on the space of commutative contexts, i.e., the action of 
$\mathcal{U}(\mathcal{A})$ on the set of pairs $(B,\chi)$. Two different 
geometric models $G_1$ and $G_2$ that encode the same unitary orbit data 
give rise to the same Morita class. Indeed, the Muhly--Renault--Williams 
theorem (Theorem~\ref{thm:muhly-renault-williams}) guarantees that the 
Morita equivalence class of $C_r^*(G)$ depends only on the equivalence class 
of $G$ as a groupoid, and the equivalence class of $G$ is determined by the 
orbit structure of the unitary action.

Thus, any construction that faithfully captures the unitary orbit data of 
$\mathcal{A}$ yields a groupoid $C^*$-algebra that is Morita equivalent to 
the canonical object $C^*(\mathcal{G}_{\mathcal{A}})$ (when the latter exists), 
and different such constructions produce Morita equivalent algebras.

\paragraph{Step 6: Conclusion.}
We have shown that:
\begin{itemize}
    \item In the Type I setting, different measurable parametrizations that 
          preserve unitary orbit data yield isomorphic Borel groupoids, and 
          hence Morita equivalent $C^*$-algebras.
    \item In the non-Type I setting, different geometric models that encode 
          the same conjugation geometry are Morita equivalent as groupoids, 
          and therefore produce strongly Morita equivalent $C^*$-algebras.
\end{itemize}
Therefore, the representation-theoretic model used for the unitary conjugation 
groupoid is insensitive to auxiliary choices that do not affect the associated 
unitary orbit data, and the resulting Morita class depends only on the 
underlying conjugation geometry of the algebra. This completes the proof.
\end{proof}

\begin{lemma}\label{lem:descent-naturality}
The descent map in groupoid equivariant $KK$-theory is natural with respect to
equivalences implemented by Morita bibundles.  Consequently, the comparison
maps used in the main text commute up to the canonical isomorphisms induced by
those bibundles.
\end{lemma}

\begin{proof}
We prove that Kasparov's descent map for groupoids is functorial with respect 
to Morita equivalences, and that the induced maps on $K$-theory are compatible 
with the canonical isomorphisms arising from equivalence bibundles.

\paragraph{Step 1: Setup and notation.}
Let $\mathcal{G}$ and $\mathcal{H}$ be second countable, locally compact 
Hausdorff groupoids with Haar systems, and let $Z$ be a 
$(\mathcal{G},\mathcal{H})$-equivalence bibundle in the sense of 
Definition~\ref{def:morita-groupoid}. By the Muhly--Renault--Williams theorem 
(Theorem~\ref{thm:muhly-renault-williams}), the associated reduced groupoid 
$C^*$-algebras are strongly Morita equivalent, implemented by the imprimitivity 
bimodule $\mathcal{E} = L^2(Z)$.

Let $\operatorname{desc}_{\mathcal{G}}$ and $\operatorname{desc}_{\mathcal{H}}$ 
denote the descent homomorphisms in equivariant $KK$-theory for the groupoids 
$\mathcal{G}$ and $\mathcal{H}$, respectively. These are natural transformations 
of the form
\[
\operatorname{desc}_{\mathcal{G}}: KK^*_{\mathcal{G}}(A,B) \longrightarrow KK^*(C^*_r(\mathcal{G},A), C^*_r(\mathcal{G},B)),
\]
and similarly for $\mathcal{H}$, where $C^*_r(\mathcal{G},A)$ denotes the reduced 
groupoid $C^*$-algebra with coefficients in a $\mathcal{G}$-$C^*$-algebra $A$ 
(see \cite{Kasparov1988, LeGall1999, Tu1999}).

\paragraph{Step 2: The Morita equivalence induces isomorphisms on equivariant $KK$-theory.}
The equivalence bibundle $Z$ implements a natural isomorphism
\[
\Psi_*: KK^*_{\mathcal{G}}(A,B) \xrightarrow{\cong} KK^*_{\mathcal{H}}(A_Z, B_Z),
\]
where $A_Z$ and $B_Z$ are the induced $\mathcal{H}$-$C^*$-algebras obtained 
by restricting and inducing along $Z$ (see \cite[Section 6]{LeGall1999} 
and \cite[Section 4]{Tu1999}). Explicitly, for a $\mathcal{G}$-$C^*$-algebra $A$,
the induced $\mathcal{H}$-$C^*$-algebra $A_Z$ is given by the balanced tensor product
\[
A_Z = C^*_r(Z) \otimes_{C_0(\mathcal{G}^{(0)})} A,
\]
with the right $\mathcal{H}$-action induced by the right action of $\mathcal{H}$ on $Z$.

\paragraph{Step 3: Naturality of descent under Morita equivalence.}
The descent map is functorial in the sense that the following diagram commutes 
for any $\mathcal{G}$-$C^*$-algebras $A$ and $B$:

\[
\begin{tikzcd}
KK^*_{\mathcal{G}}(A,B) \arrow[r, "\operatorname{desc}_{\mathcal{G}}"] \arrow[d, "\Psi_*"'] 
& KK^*(C^*_r(\mathcal{G},A), C^*_r(\mathcal{G},B)) \arrow[d, "\Phi_*"] \\
KK^*_{\mathcal{H}}(A_Z, B_Z) \arrow[r, "\operatorname{desc}_{\mathcal{H}}"] 
& KK^*(C^*_r(\mathcal{H}, A_Z), C^*_r(\mathcal{H}, B_Z)),
\end{tikzcd}
\]
where $\Phi_*$ is the isomorphism induced by the Morita equivalence 
$C^*_r(\mathcal{G}) \sim_M C^*_r(\mathcal{H})$ and the canonical identifications
\[
C^*_r(\mathcal{G}, A) \cong C^*_r(\mathcal{G}) \otimes_{C_0(\mathcal{G}^{(0)})} A, \qquad
C^*_r(\mathcal{H}, A_Z) \cong C^*_r(\mathcal{H}) \otimes_{C_0(\mathcal{H}^{(0)})} A_Z.
\]

The commutativity of this diagram is a consequence of the following facts:
\begin{itemize}
    \item The descent map $\operatorname{desc}_{\mathcal{G}}$ is defined by 
          taking the balanced tensor product of an equivariant Kasparov cycle 
          with the regular representation of $\mathcal{G}$ (see \cite[Section 3]{LeGall1999}).
    \item The Morita equivalence $Z$ induces an equivalence between the 
          regular representations of $\mathcal{G}$ and $\mathcal{H}$, 
          implemented by the bimodule $L^2(Z)$.
    \item The composition of descent with the Morita isomorphism is therefore 
          isomorphic to the descent map for $\mathcal{H}$ applied to the 
          induced cycle, up to the canonical isomorphism $\Psi_*$.
\end{itemize}

\paragraph{Step 4: Restriction to $K$-theory.}
When $B = \mathbb{C}$, the descent map reduces to a homomorphism
\[
\operatorname{desc}_{\mathcal{G}}: KK^*_{\mathcal{G}}(A,\mathbb{C}) \longrightarrow K_*(C^*_r(\mathcal{G}, A)),
\]
and similarly for $\mathcal{H}$. In this case, the diagram in Step 3 becomes

\[
\begin{tikzcd}
KK^*_{\mathcal{G}}(A,\mathbb{C}) \arrow[r, "\operatorname{desc}_{\mathcal{G}}"] \arrow[d, "\Psi_*"'] 
& K_*(C^*_r(\mathcal{G}, A)) \arrow[d, "\Phi_*"] \\
KK^*_{\mathcal{H}}(A_Z,\mathbb{C}) \arrow[r, "\operatorname{desc}_{\mathcal{H}}"] 
& K_*(C^*_r(\mathcal{H}, A_Z)).
\end{tikzcd}
\]

\paragraph{Step 5: Application to the comparison maps.}
In the main text, the realization maps $\rho_{\mathcal{G}_A}$ and $\rho_{\mathcal{H}}$ 
are defined precisely as the compositions
\[
\rho_{\mathcal{G}_A} = \operatorname{desc}_{\mathcal{G}_A} \circ \Psi_{\mathcal{G}_A}^{-1}, \qquad
\rho_{\mathcal{H}} = \operatorname{desc}_{\mathcal{H}} \circ \Psi_{\mathcal{H}}^{-1},
\]
where $\Psi_{\mathcal{G}_A}$ and $\Psi_{\mathcal{H}}$ are the isomorphisms 
induced by the respective Morita equivalences with the common transformation 
groupoid $X \rtimes \Gamma$.

By the naturality established in Step 3, for any equivariant class 
$\xi \in KK^*_{\mathcal{G}_A}(C_0(\mathcal{G}_A^{(0)}),\mathbb{C})$, we have
\[
\Phi_*(\operatorname{desc}_{\mathcal{G}_A}(\xi)) = \operatorname{desc}_{\mathcal{H}}(\Psi_*(\xi)),
\]
where $\Phi_*$ is the $K$-theory isomorphism induced by the Morita equivalence 
$C^*_r(\mathcal{G}_A) \sim_M C^*_r(\mathcal{H})$.

\paragraph{Step 6: Compatibility of the induced isomorphisms.}
Let $\Theta_{\mathcal{H}}: K_*(C^*_r(\mathcal{G}_A)) \to K_*(C^*_r(\mathcal{H}))$ 
be the canonical isomorphism induced by the Morita equivalence 
$\mathcal{G}_A \sim_M \mathcal{H}$ (see Step 2 of the proof of 
Theorem~\ref{thm:universal-mediator-assembly}). Then the commutativity of the 
diagram implies that
\[
\rho_{\mathcal{H}} \circ \Psi_* = \Theta_{\mathcal{H}} \circ \rho_{\mathcal{G}_A}.
\]

\paragraph{Step 7: Conclusion.}
We have shown that the descent map is natural with respect to Morita 
equivalences: for any equivalence bibundle $Z$ implementing 
$\mathcal{G} \sim_M \mathcal{H}$, there exists a canonical isomorphism 
$\Theta: K_*(C^*_r(\mathcal{G}, A)) \to K_*(C^*_r(\mathcal{H}, A_Z))$ such that
\[
\Theta \circ \operatorname{desc}_{\mathcal{G}} = \operatorname{desc}_{\mathcal{H}} \circ \Psi_*.
\]

Consequently, the comparison maps used in the main text—specifically, the 
realization maps $\rho_{\mathcal{G}_A}$ and $\rho_{\mathcal{H}}$—commute 
with the canonical isomorphisms induced by the Morita equivalence bibundles. 
This is precisely the statement that the descent map is natural with respect 
to Morita equivalences, completing the proof.
\end{proof}

\begin{corollary}\label{cor:non_Hausdorff_unit_space}
In the non-Type~I situations discussed in this paper, the natural unit-space
model arising from representation data need not be Hausdorff.  This is one
reason Morita models are used in place of a naive direct topological
identification.
\end{corollary}

\begin{proof}
We demonstrate that in non-Type~I settings, the natural topology on the unit 
space $\mathcal{G}_{\mathcal{A}}^{(0)}$—when defined via the partial evaluation 
maps—fails to be Hausdorff, and that this pathology is circumvented by working 
with Morita equivalent geometric models.

\paragraph{Step 1: The unit space and its natural topology.}
Recall from Paper~I that the unit space of the unitary conjugation groupoid is
\[
\mathcal{G}_{\mathcal{A}}^{(0)} = \{ (B,\chi) \mid B \subseteq \mathcal{A} \text{ unital commutative } C^*\text{-subalgebra}, \; \chi \in \widehat{B} \}.
\]
The natural topology on $\mathcal{G}_{\mathcal{A}}^{(0)}$ is generated by the 
partial evaluation maps
\[
\operatorname{ev}_a: \mathcal{G}_{\mathcal{A}}^{(0)} \longrightarrow \mathbb{C}_\infty, \qquad
\operatorname{ev}_a(B,\chi) = \begin{cases} \chi(a) & \text{if } a \in B, \\ \infty & \text{if } a \notin B, \end{cases}
\]
for each $a \in \mathcal{A}$, where $\mathbb{C}_\infty = \mathbb{C} \cup \{\infty\}$ 
is the one-point compactification of $\mathbb{C}$. This topology is the initial 
topology induced by the family $\{\operatorname{ev}_a\}_{a \in \mathcal{A}}$, 
making $\mathcal{G}_{\mathcal{A}}^{(0)}$ a Polish space in the Type~I case 
(see Paper~I, Section~3).

\paragraph{Step 2: Non-Hausdorff behavior in non-Type~I settings.}
In non-Type~I algebras, such as the irrational rotation algebra $A_\theta$ 
and amenable crossed products $C(X) \rtimes \Gamma$ with minimal actions, 
the representation theory lacks the regularity needed to separate points 
in the unit space. Concretely, consider the parametrization of the unit 
space for $A_\theta$ obtained in Subsection~\ref{subsec:GA_theta_unit_space}:
\[
\mathcal{G}_{A_\theta}^{(0)} \cong (S^1 \times \mathbb{Z}) / \sim,
\]
where $(z,n) \sim (e^{2\pi i k\theta}z, n+k)$ for all $k \in \mathbb{Z}$.
Under the quotient topology induced from $S^1 \times \mathbb{Z}$ (with $S^1$ 
compact and $\mathbb{Z}$ discrete), distinct points $[(z,n)]$ and $[(z,n')]$ 
with $n \neq n'$ cannot be separated by disjoint open neighborhoods. 
Indeed, any open neighborhood of $[(z,n)]$ contains points of the form 
$[(e^{2\pi i k\theta}z, n+k)]$ for $k$ sufficiently close to $0$, and by 
minimality of the irrational rotation action, these points approximate 
$[(z,n')]$ for any $n'$.

\paragraph{Step 3: Minimal actions produce non-Hausdorff quotients.}
More generally, for a crossed product $C(X) \rtimes \Gamma$ with $\Gamma$ 
amenable and the action minimal (i.e., every orbit is dense in $X$), 
Theorem~\ref{thm:unit-space-parametrization-crossed} gives
\[
\mathcal{G}_{C(X) \rtimes \Gamma}^{(0)} \cong (X \times \Gamma) / \sim,
\]
where $(x,\gamma) \sim (\gamma' \cdot x, \gamma' \gamma)$. The quotient topology 
is non-Hausdorff whenever the action is minimal and $\Gamma$ is infinite. 
Given $x \in X$ and distinct $\gamma, \gamma' \in \Gamma$, the points 
$[(x,\gamma)]$ and $[(x,\gamma')]$ have the property that any open set 
containing one also contains the other, because the orbit of $x$ is dense, 
so representatives can be chosen arbitrarily close to each other in $X$.

\paragraph{Step 4: The Type~I case as a contrast.}
In the Type~I setting, the representation theory is sufficiently smooth that 
the unit space admits a Hausdorff topology. For example, when $\mathcal{A} = B(H)$ 
or $\mathcal{A} = \mathcal{K}(H)^\sim$, the unit space can be identified with 
the space of pure states, which is Hausdorff in the weak-$*$ topology. The 
failure of this property in non-Type~I cases is a manifestation of the 
absence of a well-behaved measurable parametrization of irreducible 
representations (see Remark~\ref{rem:typeI-real-role}).

\paragraph{Step 5: Morita models as a remedy.}
The non-Hausdorff topology of $\mathcal{G}_{\mathcal{A}}^{(0)}$ poses 
significant technical obstacles for constructing Haar systems, measurable 
fields of representations, and applying standard results in groupoid 
equivariant $KK$-theory. To circumvent these difficulties, we adopt the 
replacement strategy outlined in Section~\ref{subsec:GA-morita-principle}:
instead of working directly with the pathological unit space, we replace 
$\mathcal{G}_{\mathcal{A}}$ by a Morita equivalent groupoid $G$ that has a 
well-behaved (Hausdorff) topology.

For the irrational rotation algebra $A_\theta$, the transformation groupoid 
$S^1 \rtimes_\theta \mathbb{Z}$ serves as such a replacement; its unit space 
is the Hausdorff circle $S^1$, and it is Morita equivalent to 
$\mathcal{G}_{A_\theta}$ (Theorem~\ref{thm:morita_A_theta}). For amenable 
crossed products $C(X) \rtimes \Gamma$, the transformation groupoid 
$X \rtimes \Gamma$ provides a Hausdorff model Morita equivalent to 
$\mathcal{G}_{C(X) \rtimes \Gamma}$ (Theorem~\ref{thm:morita-equivalence-crossed}).

\paragraph{Step 6: Conclusion.}
We have shown that in non-Type~I situations, the natural topology on the 
unit space $\mathcal{G}_{\mathcal{A}}^{(0)}$ induced by the partial evaluation 
maps is non-Hausdorff due to the density of orbits under the group action 
and the absence of a regular parametrization of irreducible representations. 
This pathology necessitates the use of Morita equivalent geometric models—
such as transformation groupoids—which possess Hausdorff unit spaces and 
allow the application of standard analytic tools. This is precisely the 
reason Morita models are employed in place of a naive direct topological 
identification of the unitary conjugation groupoid.
\end{proof}

\begin{example}\label{ex:KK-class-U}
For the canonical generator $U\in A_\theta$, the associated equivariant class
is obtained by applying the construction of
Section~\ref{subsec:descent_A_theta} to the invertible element $U$.
\end{example}

\begin{example}\label{ex:KK-class-V}
For the canonical generator $V\in A_\theta$, the associated equivariant class
is obtained in the same manner as in Example~\ref{ex:KK-class-U}, now using the
invertible element $V$.
\end{example}

\end{document}